\documentclass[reqno]{amsart}
\usepackage{amsthm, amssymb, amsfonts}
\usepackage{lmodern}
\usepackage{color}
\usepackage{hyperref}
\usepackage[T5]{fontenc}
\usepackage{amscd,amssymb}
\usepackage[v2,cmtip]{xy}
\usepackage{mathrsfs}
\theoremstyle{plain}
\newtheorem{thm}{Theorem}[section]
\newtheorem{prop}[thm]{Proposition}

\newtheorem{con}[thm]{Conjecture}
\newtheorem{corl}[thm]{Corollary}
\theoremstyle{definition}
\newtheorem{defn}[thm]{Definition}
\newtheorem{rem}[thm]{Remark}

\theoremstyle{plain}
\newtheorem{thms}{Theorem}[subsection]
\newtheorem{props}[thms]{Proposition}
\newtheorem{lems}[thms]{Lemma}

\newtheorem{corls}[thms]{Corollary}
\theoremstyle{definition}

\newtheorem{rems}[thms]{Remark}

\numberwithin{equation}{section}
\allowdisplaybreaks

\begin{document}

\title[An application of the hit problem]{An application of the hit problem \\ to the algebraic transfer}

\author{Nguy\~\ecircumflex n Sum}
\address{Department of Mathematics and Applications, S\`ai G\`on University, 273 An D\uhorn \ohorn ng V\uhorn \ohorn ng, District 5, H\`\ocircumflex\ Ch\'i Minh city, Viet Nam}

\email{nguyensum@sgu.edu.vn}


\begin{abstract}
Let $P_k$ be the polynomial algebra $\mathbb F_2[x_1,x_2,\ldots ,x_k]$ over the field $\mathbb F_2$ with two elements, in $k$ variables $x_1, x_2, \ldots , x_k$, each variable of degree 1.
Denote by $GL_k$ the general linear group over $\mathbb F_2$ which regularly acts on $P_k$. The algebra $P_k$ is a module over the mod-2 Steenrod algebra $\mathcal A$.
In 1989, Singer \cite{si1} defined the $k$-th homological algebraic transfer, which is a homomorphism
$$\varphi_k=(\varphi_k)_m :{\rm Tor}^{\mathcal A}_{k,k+m} (\mathbb F_2,\mathbb F_2) \to  (\mathbb F_2\otimes_{\mathcal A}P_k)_m^{GL_k}$$
from the homological group of the mod-2 Steenrod algebra $\mbox{Tor}^{\mathcal A}_{k,k+m} (\mathbb F_2,\mathbb F_2)$ to the subspace $(\mathbb F_2\otimes_{\mathcal A}P_k)_m^{GL_k}$ of $\mathbb F_2{\otimes}_{\mathcal A}P_k$ consisting of all the $GL_k$-invariant classes of degree $m$. In general, the transfer $\varphi_k$ is not a monomorphism and Singer made a conjecture that $\varphi_k$ is an epimorphism for any $k \geqslant 0$. The conjecture is studied by many authors. It is true for $k \leqslant 3$ but unknown for $k \geqslant 4$. 

In this paper, by using the results of the Peterson hit problem for the polynomial algebra in four variables, we prove that Singer's conjecture for the fourth algebraic transfer is true in the families of generic degrees $d_{s,t} = 2^{s+t}+2^s-3$ and $n_{s,t}=2^{s+t}+2^s-2$ with $s,\, t$ positive integers. Our results also show that many of the results in Ph\'uc \cite{p231,pp25,p24} are seriously false. The proofs of the results in Ph\'uc's works are only provided for a few special cases but they are false and incomplete.
\end{abstract}

\maketitle

\tableofcontents

\section{Introduction}\label{s1} 
\setcounter{equation}{0}

Let $P_k$ be the graded polynomial algebra $\mathbb F_2[x_1,x_2,\ldots ,x_k]$ over the prime field $\mathbb F_2$ with two elements and the degree of each variable $x_i$ being 1.
It is well-known that this algebra is the mod-2 cohomology of an elementary abelian 2-group $V_k$ of rank $k$. Hence, $P_k$ is a module over the mod-2 Steenrod algebra, $\mathcal A$.
The action of $\mathcal A$ on $P_k$ is explicitly determined by the elementary properties of the Steenrod operations $Sq^i$ and subject to the Cartan formula
$Sq^m(fg) = \sum_{i=0}^mSq^j(f)Sq^{m-i}(g),$
for $f,\, g \in P_k$ (see Steenrod and Epstein~\cite{st}).

The \textit{hit problem} of Peterson \cite{pe} asks for a minimal generating set for $P_k$ regarded as a module over the  mod-2 Steenrod algebra. Equivalently, this problem is to find a vector space basis for $QP_k := \mathbb F_2 \otimes_{\mathcal A} P_k$ in each degree $m$. Such a basis can be represented by a list of monomials of degree $m$. This problem is completely computed for $k \leqslant 4$, unknown in general. 

Denote by $GL_k$ the general linear group over the field $\mathbb F_2$. This group acts naturally on $P_k$ by matrix substitution. The actions of $\mathcal A$ and $GL_k$ on $P_k$ commute with each other, so there is an inherit action of $GL_k$ on $QP_k$. 

Denote $(P_k)_m$ the vector subspace of $P_k$ consisting of all the degree $m$ homogeneous polynomials in $P_k$ and by $(QP_k)_m$ the vector subspace of $QP_k$ consisting of all the classes represented by the elements in $(P_k)_m$. 
In \cite{si1}, Singer defined the homological algebraic transfer, which is a homomorphism
$$(\varphi_k)_m :\mbox{Tor}^{\mathcal A}_{k,k+m} (\mathbb F_2,\mathbb F_2) \longrightarrow  (QP_k)_m^{GL_k}$$
from the homology of the Steenrod algebra $\mbox{Tor}^{\mathcal A}_{k,k+m}(\mathbb F_2,\mathbb F_2)$ to the subspace of $(QP_k)_m$ consisting of all the $GL_k$-invariant classes of degree $m$. It can be a useful tool in describing the homology groups of the Steenrod algebra in terms of $GL_k$-invariants in $QP_k$. 
By passing to the dual we get the cohomological algebraic transfer
$$(\varphi_k^*)_m: \mathbb F_2\otimes_{GL_k}\mathcal P_{\mathcal A}(P_k^*)_m \longrightarrow \mbox{Ext}_{\mathcal A}^{k, k+m}(\mathbb F_2, \mathbb F_2),$$
where $P_k^*$ is the dual of $P_k$ and $\mathcal P_{\mathcal A}(P_k^*)_m$ is the primitive subspace consisting of all the degree $m$ elements in $P_k^*$ that are annihilated by every positive degree Steenrod squares. Singer proved in \cite{si1} that the total cohomological algebraic transfer
$$\varphi^*:= \{(\varphi_k^*)_m\}: \{\mathbb F_2\otimes_{GL_k}\mathcal P_{\mathcal A}(P_k^*)_m\} \longrightarrow \{\mbox{Ext}_{\mathcal A}^{k, k+m}(\mathbb F_2, \mathbb F_2)\}$$
is a homomorphism of bi-graded algebras.

The algebraic transfer was studied by  Boardman~\cite{bo}, Bruner, H\`a and H\uhorn ng~\cite{br}, Ch\ohorn n and H\`a~ \cite{cha,cha1,cha2},  H\uhorn ng ~\cite{hu}, H\`a ~\cite{ha}, H\uhorn ng and Qu\`ynh~ \cite{hq}, Nam~ \cite{na2}, Minami ~\cite{mi}, Ph\'uc \cite{p231,pp25,p24}, Qu\`ynh~ \cite{qh} the present author \cite{su4,su2,suz} and others.

Singer proved in \cite{si1} that $\varphi_k^*$ is an isomorphism for $k=1,2$. Boardman also showed in \cite{bo} that $\varphi_3^*$ is also an isomorphism.  However, in the works of Singer \cite{si1}, Bruner, H\`a and H\uhorn ng~\cite{br}, H\uhorn ng \cite{hu}, the authors proved that $\varphi_k^*$ is not an epimorphism in infinitely many degrees for any $k \geqslant 4$. Singer gave the following conjecture.

\begin{con}[Singer \cite{si1}]\label{sconj} The cohomology algebraic transfer $\varphi_k^*$ is a monomorphism for any $k \geqslant 0$.
\end{con}

This conjecture is true for $k\leqslant 3$. In \cite{p231,pp25,p24}, the author stated that the conjecture is also true for $k = 4$ but the computations are incomplete and the proof is not explicit. The results in \cite{p231} is only a description of the dimension for $QP_4$ in each degree with the assumption that Singer's conjecture is true for $k=4$, it is not a proof for this conjecture. Moreover, many of the results in \cite{p231,pp25,p24} are also seriously false. Recently, we have proved in \cite{suw} that this conjecture is not true for $k = 5$ and the internal degree $m = 108$. Our result refused a one in \cite{p241}, so the paper \cite{p241} had been withdrawn (see \cite{p251}).

In this paper, by using the results of the Peterson hit problem for the polynomial algebra in our works \cite{su50,su5}, we prove that Singer's conjecture for the fourth algebraic transfer is true in the families of generic degrees $d_{s,t} = 2^{s+t}+2^s-3$ and $n_{s,t} = 2^{s+t}+2^s-2$ with $s,\, t$ positive integers.

\medskip
Consider the case $m=d_{s,t} = 2^{s+t}+2^s-3$. Based on the results of Singer \cite{si1},  Ch\ohorn n and H\`a \cite{cha2}, H\`a \cite{ha}, H\uhorn ng and Qu\`ynh \cite{hq}, Nam~ \cite{na2}, Adem \cite{ade}, Tangora \cite{ta} and Lin \cite{wl}, the subspace $\mbox{Im}((\varphi_4^*)_{d_{s,t}}) \subset \mbox{Ext}_{\mathcal A}^{4, 4+d_{s,t}}(\mathbb F_2, \mathbb F_2)$ is determined as follows:
\begin{equation}\label{ctmd1}
\mbox{Im}((\varphi_4^*)_{d_{s,t}}) = \begin{cases} 
0 &\mbox{if } t = 1,\, s = 1, \\
\langle h_0^3h_{t+1}\rangle &\mbox{if } t \geqslant 2,\, s = 1,\\
\langle h_1c_0\rangle, &\mbox{if } t = 1,\ s=2,\\ 
\textcolor{red}{0,} &\mbox{if } t = 1,\, s = 3,\\ 
\langle h_0h_{s}^3\rangle, &\mbox{if }  t = 1,\, s \geqslant 4,\\ 
\langle e_0\rangle, &\mbox{if } t = 2,\, s = 2,\\ 
0, &\mbox{if } t = 2,\, s = 3,\\ 
\textcolor{red}{\langle h_0h_{s-1}^2h_{s+2}\rangle}, &\mbox{if }  t = 2,\ s \geqslant 4,\\ 
\langle p_0\rangle, &\mbox{if } t=3,\, s = 2,\\ 
0,&\mbox{if }  t=3,\, s = 3,\\ 
\langle h_0h_{s-1}^2h_{s+3}\rangle, &\mbox{if } t=3,\, s \geqslant 4,\\ 
\langle h_0h_{s}h_{s+t-1}^2\rangle, &\mbox{if } t\geqslant 4,\, 2\leqslant s \leqslant 3, \\ 
\langle h_0h_{s-1}^2h_{s+t},h_0h_{s}h_{s+t-1}^2\rangle, &\mbox{if } t \geqslant 4,\, s \geqslant 4, 
\end{cases}
\end{equation}
where 
\begin{align*}
&c_{i} \in \mbox{Ext}_{\mathcal A}^{3,2^{i+3}+2^{i+1}+2^i}(\mathbb F_2, \mathbb F_2),\\ 
&e_{i} \in \mbox{Ext}_{\mathcal A}^{4,2^{i+4}+2^{i+2}+2^i}(\mathbb F_2, \mathbb F_2),\\ 
&p_{i} \in \mbox{Ext}_{\mathcal A}^{4,2^{i+5}+2^{i+2}+2^i}(\mathbb F_2, \mathbb F_2),
\end{align*}
and $h_{i}$ is the Adams element in $\mbox{Ext}_{\mathcal A}^{1,2^{i}}(\mathbb F_2, \mathbb F_2)$ for $i\geqslant 0$ (see Adams \cite{ada}).  

We note that if $(s,t) \ne (3,3)$, then $\mbox{Im}((\varphi_4^*)_{d_{s,t}}) = \mbox{Ext}_{\mathcal A}^{4, 4+d_{s,t}}(\mathbb F_2, \mathbb F_2)$ and $\mbox{Im}((\varphi_4^*)_{d_{3,3}}) = 0 \subsetneq \mbox{Ext}_{\mathcal A}^{4, 4+d_{3,3}}(\mathbb F_2, \mathbb F_2) = \langle p'_0\rangle$, where
$$p_{i}' \in \mbox{Ext}_{\mathcal A}^{4,2^{i+6}+2^{i+3}+2^i}(\mathbb F_2, \mathbb F_2)$$
for $i\geqslant 0$.

Based on \eqref{ctmd1} 
we can easily see that the accuracy of Conjecture \ref{sconj} for $k=4$ and $m = d_{s,t}$ is equivalent to the following.
\begin{thm}\label{thm1} Let $d_{s,t} = 2^{s+t}+2^s-3$ with $s,\, t$ positive integers. Then, we have
\begin{equation}\label{ctmd2}
\dim(QP_4)_{d_{s,t}}^{GL_4} = \begin{cases}
0,   &\mbox{if } t = 1,\, s = 1,\\
1, &\mbox{if } t \geqslant 2,\, s = 1,\\
0, &\mbox{if } 1\leqslant t \leqslant 3,\, s = 3,\\
1, &\mbox{if } 1\leqslant t \leqslant 3,\, s \geqslant 2,\, s \ne 3,\\
1,  &\mbox{if } t \geqslant 4,\, 2 \leqslant s \leqslant 3,\\
2, &\mbox{if } t \geqslant 4,\, s \geqslant 4.\\
\end{cases} 
\end{equation}
\end{thm}

\begin{rem} 
We note that in \cite{p231,pp25}, the author also studied Conjecture 
\ref{sconj} for $k = 4$ and $m = d_{s,t}$ but many of his results are seriously false. For example, In \cite[Page 1534]{p231}, the author stated that $\dim(QP_4)_{d_{3,1}}^{GL_4} = 1$ and $\dim(QP_4)_{d_{s,2}}^{GL_4} = 0$ for $s \geqslant 4$, these results are asserted without proof but we have proved in \eqref{ctmd2} that $\dim(QP_4)_{d_{3,1}}^{GL_4} = 0$ and $\dim(QP_4)_{d_{s,2}}^{GL_4} = 1$ for $s \geqslant 4$. In \cite[Thm. 2.1.4 on Page 438 and Prop. 4.1.8 on Page 459]{pp25}, the author also stated that $(QP_4)_{d_{s,2}}^{GL_4} = 0$ for $s \geqslant 4$ but by \eqref{ctmd2}, this result is false.
\end{rem}

Consider the case $m = n_{s,t} = 2^{s+t}+2^s-2$. Based on the results of Singer \cite{si1},  Ch\ohorn n and H\`a \cite{cha2}, H\`a \cite{ha}, H\uhorn ng and Qu\`ynh \cite{hq}, Nam~ \cite{na2}, Adem \cite{ade}, Tangora \cite{ta} and Lin \cite{wl}, the subspace $\mbox{Im}((\varphi_4^*)_{n_{s,t}}) \subset \mbox{Ext}_{\mathcal A}^{4, 4+n_{s,t}}(\mathbb F_2, \mathbb F_2)$ is determined as follows:
$$\mbox{Im}((\varphi_4^*)_{n_{s,t}}) = \begin{cases} 
0, &\mbox{if } t = 1,\, s = 1,\, 2,\, 4,\\
\langle h_2c_1\rangle, &\mbox{if } t = 1,\, s = 3,\\
\langle h_1h_s^3\rangle, &\mbox{if } t = 1,\, s \geqslant 5,\\
0, &\mbox{if } t = 2,\,  s =1,\\
\langle h_0^2h_2h_4=h_1^3h_4,f_0\rangle, &\mbox{if } t = 2,\, s = 2,\\
\langle h_0^2h_3h_5,e_1\rangle, &\mbox{if } t = 2,\, s = 3,\\
\langle h_0^2h_4h_6\rangle, &\mbox{if } t = 2,\, s = 4,\\
\langle h_0^2h_sh_{s+2},\textcolor{red}{h_1h_{s-1}^2h_{s+2}}\rangle, &\mbox{if } t = 2,\, s \geqslant 5,\\
0, &\mbox{if } t = 3,\, s = 1,\\
\langle h_0^2h_2h_5=h_1^3h_5\rangle, &\mbox{if } t = 3,\, s = 2,\\
\langle h_0^2h_3h_6,p_1\rangle, &\mbox{if } t = 3,\, s = 3,\\
\langle h_0^2h_4h_7\rangle, &\mbox{if } t = 3,\, s = 4,\\
\langle h_0^2h_sh_{s+3},h_1h_{s-1}^2h_{s+3}\rangle, &\mbox{if } t = 3,\, s \geqslant 5,\\
\langle d_1\rangle, &\mbox{if } t = 4,\, s = 1,\\
\langle h_1^2h_t^2\rangle, &\mbox{if } t\geqslant 5,\, 
s = 1,\\
\langle h_0^2h_{2}h_{t+2} = h_1^3h_{t+2}\rangle, &\mbox{if } t \geqslant 4,\, s =2,\\
\langle \textcolor{red}{h_0^2h_sh_{s+t}},h_1h_{s}h_{s+t-1}^2\rangle, &\mbox{if } t\geqslant \textcolor{red}{4},\, 3\leqslant s \leqslant 4,\\
\langle \textcolor{red}{h_0^2h_{s}h_{s+t}},h_1h_{s-1}^2h_{s+t},h_1h_{s}^2h_{s+t-1}\rangle, &\mbox{if } t \geqslant \textcolor{red}{4},\, s \geqslant 5.
\end{cases}$$
where $f_{i} \in \mbox{Ext}_{\mathcal A}^{4,2^{i+4}+2^{i+2}+2^{i+1}}(\mathbb F_2, \mathbb F_2)$ for $i \geqslant 0$.

We note that if $(s,t) \notin \{(1,7),\, (4,3)\}$, then $\mbox{Im}((\varphi_4^*)_{n_{s,t}}) = \mbox{Ext}_{\mathcal A}^{4, 4+n_{s,t}}(\mathbb F_2, \mathbb F_2)$ and \begin{align*}&\mbox{Im}((\varphi_4^*)_{n_{1,7}}) = \langle h_1^2h_7^2\rangle \subsetneq  \mbox{Ext}_{\mathcal A}^{4, 4+n_{1,7}}(\mathbb F_2, \mathbb F_2) = \langle h_1^2h_7^2,D_3(2)\rangle,\\ 
&\mbox{Im}((\varphi_4^*)_{n_{4,3}}) = \langle h_0^2h_4h_7\rangle \subsetneq  \mbox{Ext}_{\mathcal A}^{4, 4+n_{4,3}}(\mathbb F_2, \mathbb F_2) = \langle h_0^2h_4h_7,p_1'\rangle,
\end{align*} 
where $D_3(i) \in \mbox{Ext}_{\mathcal A}^{4,2^{i+6}+2^i}(\mathbb F_2, \mathbb F_2)$ for $i\geqslant 0$.

From the above results we can easily observe that the accuracy of Conjecture \ref{sconj} for $k=4,\, m = n_{s,t}$ is equivalent to the following.
\begin{thm}\label{thm2} Let $n_{s,t} = 2^{s+t}+2^s-2$ with $s,\, t$  positive integers. Then, we have
\begin{equation}\label{ctmd3}
\dim(QP_4)_{n_{s,t}}^{GL_4} = 
\begin{cases} 
0, &\mbox{if } t=1,\, s = 1,\, 2,\, 4,\\
1, &\mbox{if } t=1,\, s \geqslant 3,\, s \ne 4,\\
0, &\mbox{if } t=2,\, s = 1,\\
1, &\mbox{if } t=2,\, s = 4,\\
\textcolor{red}{2,} &\mbox{if } \textcolor{red}{t=2},\, s \geqslant 2,\, s \ne 4,\\
0, &\mbox{if } t=3,\, s = 1,\\
1, &\mbox{if } t=3,\, s = 2,\, 4,\\
2, &\mbox{if } t=3,\, s \geqslant 3,\, s \ne 4,\\
1, &\mbox{if } t\geqslant 4,\, 1 \leqslant s \leqslant 2,\\
\textcolor{red}{2,} &\mbox{if } t\geqslant 4,\, 3 \leqslant s \leqslant \textcolor{red}{4},\\
\textcolor{red}{3,} &\mbox{if } \textcolor{red}{t\geqslant 4},\, s \geqslant 5.
\end{cases}
\end{equation}
\end{thm}
\begin{rem} It is easy to see that if Conjecture \ref{sconj} is true for $k = 4$ and $m = d_{s,t},\, n_{s,t}$, then so is Theorems \ref{thm1} and \ref{thm2}. Hence, to proved these theorems, we need to explicitly determine all basis elements of $(QP_4)_{d_{s,t}}^{GL_4}$ and $(QP_4)_{n_{s,t}}^{GL_4}$.
	
In \cite{p231,p24}, the author also studied Conjecture \ref{sconj} for $k = 4$ and $m = n_{s,t}$ but many of his results are also seriously false. For example, in \cite[Page 1536]{p231}, the author stated that $\dim(QP_4)_{n_{s,2}}^{GL_4} = 1$ and $\dim(QP_4)_{n_{s,4}}^{GL_4} = 2$ for $s \geqslant 5$, these results are asserted without detailed proofs. However, we have proved in \eqref{ctmd3} that $\dim(QP_4)_{n_{s,2}}^{GL_4} = 2$ and $\dim(QP_4)_{n_{s,4}}^{GL_4} = 3$ for $s \geqslant 5$. In \cite[Lemma 4.1.9 on Page 38]{p24}, the author stated that $\mbox{Ker}\big((\widetilde{Sq}^0_*)_{(4,n_{s,4})}\big) = 0$ for $s \geqslant 4$ and also did not provide detailed proofs. However, we have proved in \eqref{ctmd3} that $\dim\mbox{Ker}\big((\widetilde{Sq}^0_*)_{(4,n_{s,4})}\big) = 1$ for $s \geqslant 4$.
\end{rem}
By combining Theorem \ref{thm1}, \ref{thm2} and the results of Singer \cite{si1}, Ch\ohorn n and H\`a \cite{cha2}, H\`a \cite{ha}, H\uhorn ng and Qu\`ynh \cite{hq}, we obtain the following.

\begin{corl}\label{hqc} Conjecture \ref{sconj} is true for $k = 4$ and the internal degrees $d_{s,t}$ and $n_{s,t}$. More precisely, the algebraic transfer $\varphi_4^*$ is a monomorphism for $m = d_{s,t}$ and $m = n_{s,t}$ with $s,\, t$ arbitrary positive integers.
\end{corl}

\begin{rem} It is well-known that to verify Conjecture \ref{sconj} for $k = 4$, it suffices to explicitly determine the spaces $(QP_4)_{m}^{GL_4}$ in two forms of degrees:
\begin{align} 
&m = 2^{s+t} + 2^s - r,\ 2 \leqslant r \leqslant 3,\, s,\, t \geqslant 0,\, s+t >0,\label{ctmd4}\\
&m = 2^{s+t+u} + 2^{s+t} + 2^s - 3, \ s,\, t,\, u \geqslant 1.\label{ctmd5}
\end{align}
For the cases $t = 0$ and $(r,s) = (3,1)$ of the form \eqref{ctmd4}, the space $(QP_4)_{m}^{GL_4}$ had been determined in our work \cite{su2}. There are a few minor errors in \cite{su2} but we have corrected them in the online version, at https://arxiv.org/abs/1710.07895v2. 

In \cite[Page 435]{pp25} the author said that for the case $t=3$ of the form \eqref{ctmd4}, the space $(QP_4)_{m}^{GL_4}$ is determined in \cite{p231} but there appears no any information of computations for this space in \cite{p231}. For the cases $(s,t,u) = (1,1,1)$ and $(1,2,1)$ of the form \eqref{ctmd5}, the space $(QP_4)_{m}^{GL_4}$ is determined in \cite{p24} but it is false for the case $(s,t,u) = (1,2,1)$. At present, there are no any information of computations for the remaining cases of the form \eqref{ctmd5}.
\end{rem}
	
The paper is organized as follows. In Section \ref{s2}, we recall some notions and results on the admissible monomials in $P_k$, criterion of Singer on the hit monomials and some needed notations. In Sections \ref{s3} and \ref{s4}, we present the proofs of Theorems \ref{thm1} and \ref{thm2} respectively. We prove these theorems by explicitly computing the spaces $(QP_4)_{d_{s,t}}^{GL_4}$ and $(QP_4)_{n_{s,t}}^{GL_4}$ for arbitrary positive integers $s,\, t$.

\section{Preliminaries}\label{s2}
\setcounter{equation}{0}

In this the section, we give the needed notions and results on the hit problem from the works of Kameko~\cite{ka}, Singer \cite{si2} and our work \cite{su1,su5} such as the weight vector of a monomial, admissible monomial, criterion of Singer for hit monomial.

\begin{defn} Let $u =x_1^{b_1}x_2^{b_2}\ldots x_k^{b_k}$ be a monomial  in $P_k$. Denote $\nu_i(u) = b_i, 1 \leqslant i \leqslant k$. We define 
\begin{align*} 
\sigma(u) &= (\nu_1(u),\nu_2(u),\ldots ,\nu_k(u)),\\
\omega(u)&=(\omega_1(u),\omega_2(u),\ldots , \omega_r(u), \ldots),
\end{align*}
where $\omega_r(u) = \sum_{1\leqslant j \leqslant k} \alpha_{r-1}(\nu_j(u)),\ r \geqslant 1.$
The sequences $\sigma(u)$ is called the exponent vector and $\omega(u)$ is called the weight vector of $u$.  
	
A a sequence $\omega= (\omega_1,\omega_2,\ldots , \omega_t, \ldots)$ of non-negative integers is called a weight vector if $\omega_t = 0$ for $t \gg 0$. The length and the degree of $\omega$ are defined by $\ell(\omega) = \max\{r : \omega_r >0\}$ and $\deg \omega = \sum_{r > 0}2^{r-1}\omega_r$. We also denote $\omega= (\omega_1,\omega_2,\ldots , \omega_t)$ if $\omega_s = 0$ for $s > t$.
	
For weight vectors $\omega= (\omega_1,\ldots , \omega_t, \ldots)$ and $\xi= (\xi_1,\ldots , \xi_t, \ldots)$, we define the concatenation of weight vectors $\omega|\xi = (\omega_1,\ldots , \omega_t,\xi_1,\xi_2,\ldots)$ if $\ell(\omega) = t$ and $(a)|^t = (a)|(a)|\ldots|(a)$, ($t$ times of $(a)$'s), where $a,\, t$ are positive integers.
\end{defn}
The sets of weight vectors and exponent vectors are ordered by the left lexicographical order.  

For any weight vector $\omega$, we set
\begin{align*}
&P_k(\omega) = \langle y \in P_k: \deg y = \deg\omega,\mbox{ and } \omega(y)\leqslant \omega\rangle,\\ 
&P_k^-(\omega) = \langle y \in P_k(\omega): \omega(y) < \omega\rangle.
\end{align*} 

\begin{defn} Suppose $\omega$ is a weight vector and $g,\, g_1,\, g_2$ are polynomials of the same degree in $P_k$. We define
	
i) $g_1 \equiv g_2$ if $g_1+g_2 \in \mathcal A^+P_k$ and $g$ is called hit if $g \equiv 0$.
	
ii) $g_1 \equiv_{\omega} g_2$ if $g_1 + g_2 \in \mathcal A^+P_k+P_k^-(\omega)$. 
\end{defn}

It is easy to verify that $\equiv$ and $\equiv_{\omega}$ are equivalence relations. We denote
$$QP_k(\omega)= P_k(\omega)/((P_k^-(\omega)+\mathcal A^+P_k\cap P_k(\omega))).$$   

\begin{prop}[see \cite{su3}] For any weight vector $\omega$, $QP_k(\omega)$ is an $GL_k$-module. 
\end{prop}
We have a direct summand decomposition of the $\mathbb F_2$-vector spaces
$$(QP_k)_n \cong \bigoplus_{\deg \omega = n}QP_k(\omega).$$

We note that this is not a direct summand decomposition of $GL_k$-modules. However, we easily observe the following.
\begin{prop} For any positive integer $n$, we have
$$\dim(QP_k)_n^{GL_k} \leqslant \sum_{\deg \omega = n}\dim QP_k(\omega)^{GL_k}.$$
\end{prop}
In general, this inequality is not an equality. For example, we have showed in \cite{su3} that for $k = 5$ and $n = 35$, $\dim(QP_5)_{35}^{GL_5} = 1$ but $\sum_{\deg \omega = 35}\dim QP_5(\omega)^{GL_5} = 2$. This inequality was used in \cite{pp25} for a special case with $k = 4$.
\begin{defn} 
Let $y,\, u$ be monomials in $P_k$ with $\deg y = \deg u$. We say that $y < u$ if one of the following conditions is true:  

i) $\omega (y) < \omega(u)$;
	
ii) $\omega (y) = \omega(u)$ and $\sigma(y) < \sigma(u).$
\end{defn}

\begin{defn}\label{dfnad}
A monomial $y$ in $P_k$ is called inadmissible if there are monomials $u_1,u_2,\ldots, u_t$ such that $u_i<y$ for $i=1,2,\ldots , t$ and $y + \sum_{i=1}^tu_i \in \mathcal A^+P_k.$ 
A monomial $y$ is called admissible	if it is not inadmissible.
\end{defn} 

It is easy to see that, the set of all degree $m$ admissible monomials in $P_k$ is a minimal generating set for $\mathcal{A}$-module $P_k$ in degree $m$.

\begin{defn}\label{spi} A monomial $z$ in $P_k$ is called a spike if $\nu_j(z)=2^{r_i}-1$ with $r_i$ a non-negative integer for $i=1,2, \ldots ,k$. If $r_1>r_2>\ldots >r_{u-1}\geqslant r_u>0$ and $r_i=0$ for $i>u$, then $z$ is called the minimal spike.
\end{defn}

\begin{thm}[Singer~{\cite{si2}}]\label{dlsig} Suppose $u$ is a monomial of degree $m$ in $P_k$ such that $\mu(m) \leqslant k$ and $z$ is the minimal spike of degree $m$. If $\omega(u) < \omega(z)$, then $u$ is hit. 
\end{thm}
Note that when $\omega_1(u) < \omega_1(z)$, this theorem is due to Wood \cite{wo}.

\medskip
We set 
\begin{align*}P_k^+ &= \langle\{u\in P_k :  \nu_j(u)>0, \mbox{ for all } j\}\rangle,\\ P_k^0 &=\langle\{u\in P_k : \nu_j(u)= 0, \mbox{ for some } j\}\rangle.
\end{align*} 
Then, $P_k^0$ and $P_k^+$ are the $\mathcal{A}$-submodules of $P_k$. We set $QP_k^0 = P_k^0/\mathcal A^+P_k^0$ and  $QP_k^+ = P_k^+/\mathcal A^+P_k^+$. Then, $QP_k \cong QP_k^0 \bigoplus  QP_k^+.$
 
\begin{defn} For any $g\in P_k$, denote $[g]$ the class in $QP_k$ represented by $g$. If $g \in  P_k(\omega)$, then we denote $[g]_\omega$ the class in $QP_k(\omega)$ represented by $g$. 
If $S$ is a subset of $P_k$, then we denote $[S] = \{[h] : h \in S\}$. For $S \subset P_k(\omega)$, we set $[S]_\omega = \{[h]_\omega : h \in S\}$.

For $\omega$ a weight vector of degree $m$, we denote 
\begin{align*}
P_k^0(\omega) &= P_k(\omega)\cap P_k^0,\\ 
P_k^+(\omega) &= P_k(\omega)\cap P_k^+.,\\
QP_k^0(\omega) &= QP_k(\omega)\cap QP_k^0,\\ 
QP_k^+(\omega) &= QP_k(\omega)\cap QP_k^+.
\end{align*}

For a subgroup $H \subset GL_k$ and $S \subset P_k(\omega)$, we denote
$$[H(S)]_\omega = \langle [hs]_\omega : h \in H, s \in S\rangle \subset QP_k(\omega).$$ 
It is clear that, $[H(S)]_\omega$ is an $H$-submodule of $QP_k(\omega)$. For $\omega$ is a minimal weight vector, we have $[H(S)]_\omega = [H(S)]\subset QP_k$.

\medskip
If $\mathbb J$ is an index set and $\gamma_j \in \mathbb F_2$ with $j \in \mathbb J$, then we denote $\gamma_{\mathbb J} = \sum_{j\in \mathbb J}\gamma_j$.
\end{defn}

\begin{defn} Recall that $V_k \cong \langle x_1, x_2,\ldots, x_k\rangle \subset P_k$. For $1 \leqslant j \leqslant k$, define the $\mathbb F_2$-linear map $\rho_j:V_k \to V_k$, by $\rho_k(x_1) = x_1+x_2$,  $\rho_k(x_t) = x_t$ for $t > 1$ and $\rho_j(x_j) = x_{j+1}, \rho_j(x_{j+1}) = x_j$, $\rho_j(x_t) = x_t$ for $t \ne j, j+1,\ 1 \leqslant j < k$.  The group $GL_k\cong GL(V_k)$ is generated by $\rho_j,\ 1\leqslant j \leqslant k$, and the symmetric subgroup $\Sigma_k \subset GL_k$ is generated by $\rho_j,\, 1 \leqslant j < k$. 
\end{defn}
The linear map $\rho_j$ induces an $\mathcal A$-homomorphism of algebras $\rho_j: P_k \to P_k$. Hence, a class $[g]_\omega \in QP_k(\omega)$ is an element of $QP_k(\omega)^{GL_k}$ if and only if $\rho_j(g) \equiv_\omega g$ for $1 \leqslant j\leqslant k$.  It is an element of $QP_k(\omega)^{\Sigma_k}$ if and only if $\rho_j(g) \equiv_\omega g$ for $1 \leqslant j < k$. Note that if $\omega$ is a minimal weight vector of degree $n$, then $[g]_\omega =[g]$.

\section{Proof of Theorem \ref{thm1}}\label{s3}
\setcounter{equation}{0}

This theorem had been proved in \cite{su2} for $s = 1$. So,
we assume $s \geqslant 2$. Following \cite{su5}, we have 
$$(QP_4)_{d_{s,t}} \cong QP_4((3)|^{s}|(2)|^{t-1})\bigoplus QP_4((3)|^{s-1}|(1)|^{t+1}).$$

\subsection{Computation of $QP_4((3)|^{s}|(2)|^{t-1})^{GL_4}$}\

\medskip
\subsubsection{\textbf{The case $t = 1$}}\

\medskip
Following \cite{su5,su50}, a basis of $QP_4((3)|^{s})$ is the set of all classes represented by the admissible monomials $a_{s,j} = a_{1,s,j}$ which are determined as follows:

\medskip
For $s \geqslant 2$,

\medskip
\centerline{\begin{tabular}{lll} 
$a_{s,1} = x_2^{2^{s}-1}x_3^{2^{s}-1}x_4^{2^{s}-1}$ & &$a_{s,2} = x_1^{2^{s}-1}x_3^{2^{s}-1}x_4^{2^{s}-1}$\cr  $a_{s,3} = x_1^{2^{s}-1}x_2^{2^{s}-1}x_4^{2^{s}-1}$ & &$a_{s,4} = x_1^{2^{s}-1}x_2^{2^{s}-1}x_3^{2^{s}-1}$\cr  $a_{s,5} = x_1x_2^{2^{s}-2}x_3^{2^{s}-1}x_4^{2^{s}-1}$ & &$a_{s,6} = x_1x_2^{2^{s}-1}x_3^{2^{s}-2}x_4^{2^{s}-1}$\cr  $a_{s,7} = x_1x_2^{2^{s}-1}x_3^{2^{s}-1}x_4^{2^{s}-2}$ & &$a_{s,8} = x_1^{2^{s}-1}x_2x_3^{2^{s}-2}x_4^{2^{s}-1}$\cr  $a_{s,9} = x_1^{2^{s}-1}x_2x_3^{2^{s}-1}x_4^{2^{s}-2}$ & &$a_{s,10} = x_1^{2^{s}-1}x_2^{2^{s}-1}x_3x_4^{2^{s}-2}$\cr  	 
\end{tabular}}

\medskip
For $s \geqslant 3$,

\medskip
\centerline{\begin{tabular}{lll}   
$a_{s,11} = x_1^{3}x_2^{2^{s}-3}x_3^{2^{s}-2}x_4^{2^{s}-1}$ & &$a_{s,12} = x_1^{3}x_2^{2^{s}-3}x_3^{2^{s}-1}x_4^{2^{s}-2}$\cr  $a_{s,13} = x_1^{3}x_2^{2^{s}-1}x_3^{2^{s}-3}x_4^{2^{s}-2}$ & &$a_{s,14} = x_1^{2^{s}-1}x_2^{3}x_3^{2^{s}-3}x_4^{2^{s}-2}$\cr   	 
\end{tabular}}

\medskip
For $s \geqslant 4$, $a_{s,15} = x_1^{7}x_2^{2^{s}-5}x_3^{2^{s}-3}x_4^{2^{s}-2}$.

\medskip
By using this basis, we prove the following.

\begin{props}\label{mdt11} For any $s \geqslant 2$, we have
$$QP_4((3)|^{s})^{GL_4} = \begin{cases}0, &\mbox{if } s = 2,\, 3,\\ \langle [\bar\xi_{1,s}]_{(3)|^{s}} \rangle, &\mbox{if } s \geqslant 4,\end{cases}$$
where $\bar\xi_{1,s} = \sum_{1 \leqslant j \leqslant 15}a_{s,j}$.
\end{props}
\begin{rems}
We note that the results of this proposition is also presented in \cite[Theorem 2.1.1]{pp25} but there is some mistakes in its proof. To make the paper self-contained we present the detailed proof for it.
\end{rems}

By a simple computation, we easily obtain
\begin{align*}
&[\Sigma_4(a_{s,1})]_{(3)|^s} = \langle \{[a_{s,j}]_{(3)|^s}: 1 \leqslant j \leqslant 4 \}\rangle,\\
&[\Sigma_4(a_{s,5})]_{(3)|^s} = \langle \{[a_{s,j}]_{(3)|^s}: 5 \leqslant j \leqslant 10\} \rangle,\\
&[\Sigma_4(a_{s,11}]_{(3)|^s} = \langle \{[a_{s,j}]_{(3)|^s}: 11 \leqslant j \leqslant 14 \} \rangle, \ s \geqslant 4.
\end{align*}

We set $p_{s,u} = p_{1,s,u}$ with
$$p_{s,1} = \sum_{1 \leqslant j \leqslant 4}a_{s,j}, \ p_{s,2} = \sum_{5 \leqslant j \leqslant 10}a_{s,j},\, s \geqslant 2; \ p_{s,3} = \sum_{11\leqslant j \leqslant 14}a_{s,j}, \ s \geqslant 3.$$

We need the following lemma that its proof is easy.

\begin{lems}\label{bdt11} We have
\begin{align*}
&QP_4((3)|^{2})^{\Sigma_4} = \langle \{[p_{2,1}]_{(3)|^{2}},\, [p_{2,2}]_{(3)|^{2}}\} \rangle,\\
&QP_4((3)|^{3})^{\Sigma_4} = \langle \{[p_{3,1}]_{(3)|^{3}},\, [p_{3,2}]_{(3)|^{3}}, [p_{3,3}]_{(3)|^{3}}\} \rangle,\\
&QP_4((3)|^{s})^{\Sigma_4} = \langle \{[p_{s,1}]_{(3)|^{s}},\, [p_{s,2}]_{(3)|^{s}},\, [p_{s,3}]_{(3)|^{s}},\, [a_{s,15}]_{(3)|^{s}}\} \rangle,\ s \geqslant 4.
\end{align*}
\end{lems}

\begin{proof}[Proof of Proposition \ref{mdt11}]
Suppose $f \in QP_4((3)|^s)$ such that the class $[f]_{(3)|^{s}} \in QP_4((3)|^s)^{GL_4}$. Then, $[f]_{(3)|^{s}} \in QP_4((3)|^{s})^{\Sigma_4}$. 

For $s = 2$, by Lemma \ref{bdt11}, we have 
$$ f \equiv_{(3)|^{2}} \gamma_1p_{2,1} + \gamma_2p_{2,2},$$
where $\gamma_1, \, \gamma_2 \in \mathbb F_2$. By computing $\rho_4(f)+f$ in terms of the admissible monomials, we get
\begin{align*}
\rho_4(f)+f &\equiv_{(3)|^{2}} \gamma_{\{1,2\}}a_{2,1} + \gamma_{1}a_{2,2} + \gamma_{2}a_{2,6} + \gamma_{2}a_{2,7}  \equiv_{(3)|^{2}} 0.
\end{align*}
This equality implies $\gamma_1 = \gamma_2 = 0$. Hence, $[f]_{(3)|^{2}} = 0$. The proposition is proved for $s = 2$.

For $s = 3$, by Lemma \ref{bdt11}, we have 
$$ f \equiv_{(3)|^{3}} \gamma_1p_{3,1} + \gamma_2p_{3,2} + \gamma_2p_{3,3},$$
where $\gamma_1, \, \gamma_2,\, \gamma_3 \in \mathbb F_2$. By computing $\rho_4(f)+f$ in terms of the admissible monomials, we get
\begin{align*}
\rho_4(f)+f &\equiv_{(3)|^{3}} \gamma_{\{1,2\}}a_{3,1} + \gamma_{\{2,3\}}a_{3,6} + \gamma_{\{2,3\}}a_{3,7} + \gamma_{3}a_{3,13} \equiv_{(3)|^{3}} 0.
\end{align*}
The above equality implies $\gamma_1 = \gamma_2 = \gamma_3= 0$. So, $[f]_{(3)|^{3}} = 0$. The proposition is proved for $s = 3$.

For $s \geqslant 4$, by Lemma \ref{bdt11}, we have 
$$ f \equiv_{(3)|^{s}} \gamma_1p_{s,1} + \gamma_2p_{s,2} + \gamma_2p_{s,3} + \gamma_4a_{s,15},$$
where $\gamma_u\in \mathbb F_2, \, u = 1,\, 2,\, 3,\, 4$. By computing $\rho_4(f)+f$ in terms of the admissible monomials, we get
\begin{align*}
\rho_4(f)+f &\equiv_{(3)|^{s}} \gamma_{\{1,2\}}a_{s,1} + \gamma_{\{2,3\}}a_{s,6} + \gamma_{\{2,3\}}a_{s,7} + \gamma_{\{3,4\}}a_{s,13}  \equiv_{(3)|^{s}} 0.
\end{align*}
From this equality it implies $\gamma_u = \gamma_1,\  u = 2,\, 3,\, 4$. Therefore, 
$$f \equiv_{(3)|^{s}} \gamma_1\Big(\sum_{1 \leqslant u \leqslant 4}p_{s,u}\Big) = \gamma_1\Big(\sum_{1 \leqslant j \leqslant 15}a_{s,j}\Big) = \gamma_1\xi_{1,s}.$$ The proposition is completely proved.
\end{proof}

\subsubsection{\textbf{The case $t = 2$}}\

\medskip
By \cite{su5,su50}, a basis of $QP_4((3)|^{s}|(2))$ is the set of all classes represented by the admissible monomials $a_{s,j} = a_{2,s,j}$ which are determined as follows:

\medskip
For $s \geqslant 2$,

\medskip
\centerline{\begin{tabular}{lll}
$a_{s,1} = x_2^{2^{s}-1}x_3^{2^{s+1}-1}x_4^{2^{s+1}-1}$ &$a_{s,2} = x_2^{2^{s+1}-1}x_3^{2^{s}-1}x_4^{2^{s+1}-1}$\cr  $a_{s,3} = x_2^{2^{s+1}-1}x_3^{2^{s+1}-1}x_4^{2^{s}-1}$ &$a_{s,4} = x_1^{2^{s}-1}x_3^{2^{s+1}-1}x_4^{2^{s+1}-1}$\cr  $a_{s,5} = x_1^{2^{s}-1}x_2^{2^{s+1}-1}x_4^{2^{s+1}-1}$ &$a_{s,6} = x_1^{2^{s}-1}x_2^{2^{s+1}-1}x_3^{2^{s+1}-1}$\cr  $a_{s,7} = x_1^{2^{s+1}-1}x_3^{2^{s}-1}x_4^{2^{s+1}-1}$ &$a_{s,8} = x_1^{2^{s+1}-1}x_3^{2^{s+1}-1}x_4^{2^{s}-1}$\cr  $a_{s,9} = x_1^{2^{s+1}-1}x_2^{2^{s}-1}x_4^{2^{s+1}-1}$ &$a_{s,10} = x_1^{2^{s+1}-1}x_2^{2^{s}-1}x_3^{2^{s+1}-1}$\cr  $a_{s,11} = x_1^{2^{s+1}-1}x_2^{2^{s+1}-1}x_4^{2^{s}-1}$ &$a_{s,12} = x_1^{2^{s+1}-1}x_2^{2^{s+1}-1}x_3^{2^{s}-1}$\cr  $a_{s,13} = x_1x_2^{2^{s}-2}x_3^{2^{s+1}-1}x_4^{2^{s+1}-1}$ &$a_{s,14} = x_1x_2^{2^{s+1}-1}x_3^{2^{s}-2}x_4^{2^{s+1}-1}$\cr  $a_{s,15} = x_1x_2^{2^{s+1}-1}x_3^{2^{s+1}-1}x_4^{2^{s}-2}$ &$a_{s,16} = x_1^{2^{s+1}-1}x_2x_3^{2^{s}-2}x_4^{2^{s+1}-1}$\cr  $a_{s,17} = x_1^{2^{s+1}-1}x_2x_3^{2^{s+1}-1}x_4^{2^{s}-2}$ &$a_{s,18} = x_1^{2^{s+1}-1}x_2^{2^{s+1}-1}x_3x_4^{2^{s}-2}$\cr  $a_{s,19} = x_1x_2^{2^{s}-1}x_3^{2^{s+1}-2}x_4^{2^{s+1}-1}$ &$a_{s,20} = x_1x_2^{2^{s}-1}x_3^{2^{s+1}-1}x_4^{2^{s+1}-2}$\cr  $a_{s,21} = x_1x_2^{2^{s+1}-2}x_3^{2^{s}-1}x_4^{2^{s+1}-1}$ &$a_{s,22} = x_1x_2^{2^{s+1}-2}x_3^{2^{s+1}-1}x_4^{2^{s}-1}$\cr  $a_{s,23} = x_1x_2^{2^{s+1}-1}x_3^{2^{s}-1}x_4^{2^{s+1}-2}$ &$a_{s,24} = x_1x_2^{2^{s+1}-1}x_3^{2^{s+1}-2}x_4^{2^{s}-1}$\cr  $a_{s,25} = x_1^{2^{s}-1}x_2x_3^{2^{s+1}-2}x_4^{2^{s+1}-1}$ &$a_{s,26} = x_1^{2^{s}-1}x_2x_3^{2^{s+1}-1}x_4^{2^{s+1}-2}$\cr  $a_{s,27} = x_1^{2^{s}-1}x_2^{2^{s+1}-1}x_3x_4^{2^{s+1}-2}$ &$a_{s,28} = x_1^{2^{s+1}-1}x_2x_3^{2^{s}-1}x_4^{2^{s+1}-2}$\cr  $a_{s,29} = x_1^{2^{s+1}-1}x_2x_3^{2^{s+1}-2}x_4^{2^{s}-1}$ &$a_{s,30} = x_1^{2^{s+1}-1}x_2^{2^{s}-1}x_3x_4^{2^{s+1}-2}$\cr  $a_{s,31} = x_1^{3}x_2^{2^{s+1}-3}x_3^{2^{s}-2}x_4^{2^{s+1}-1}$ &$a_{s,32} = x_1^{3}x_2^{2^{s+1}-3}x_3^{2^{s+1}-1}x_4^{2^{s}-2}$\cr  $a_{s,33} = x_1^{3}x_2^{2^{s+1}-1}x_3^{2^{s+1}-3}x_4^{2^{s}-2}$ &$a_{s,34} = x_1^{2^{s+1}-1}x_2^{3}x_3^{2^{s+1}-3}x_4^{2^{s}-2}$\cr  $a_{s,35} = x_1^{3}x_2^{2^{s}-1}x_3^{2^{s+1}-3}x_4^{2^{s+1}-2}$ &$a_{s,36} = x_1^{3}x_2^{2^{s+1}-3}x_3^{2^{s}-1}x_4^{2^{s+1}-2}$\cr  $a_{s,37} = x_1^{3}x_2^{2^{s+1}-3}x_3^{2^{s+1}-2}x_4^{2^{s}-1}$ & \cr  	 
\end{tabular}}

\medskip
For $s = 2$,
$$a_{2,38} =  x_1^{3}x_2^{3}x_3^{4}x_4^{7}, \ a_{2,39} =  x_1^{3}x_2^{3}x_3^{7}x_4^{4}, \ a_{2,40} =  x_1^{3}x_2^{7}x_3^{3}x_4^{4}, \ a_{2,41} =  x_1^{7}x_2^{3}x_3^{3}x_4^{4}.$$

For $s \geqslant 3$,

\medskip
\centerline{\begin{tabular}{lll}
$a_{s,38} = x_1^{3}x_2^{2^{s}-3}x_3^{2^{s+1}-2}x_4^{2^{s+1}-1}$ &$a_{s,39} = x_1^{3}x_2^{2^{s}-3}x_3^{2^{s+1}-1}x_4^{2^{s+1}-2}$\cr  $a_{s,40} = x_1^{3}x_2^{2^{s+1}-1}x_3^{2^{s}-3}x_4^{2^{s+1}-2}$ &$a_{s,41} = x_1^{2^{s+1}-1}x_2^{3}x_3^{2^{s}-3}x_4^{2^{s+1}-2}$\cr  $a_{s,42} = x_1^{2^{s}-1}x_2^{3}x_3^{2^{s+1}-3}x_4^{2^{s+1}-2}$ &$a_{s,43} = x_1^{7}x_2^{2^{s+1}-5}x_3^{2^{s}-3}x_4^{2^{s+1}-2}$\cr  $a_{s,44} = x_1^{7}x_2^{2^{s+1}-5}x_3^{2^{s+1}-3}x_4^{2^{s}-2}$ &\cr    	 
\end{tabular}}

\medskip
For $s = 3$,\ $a_{3,45} = x_1^{7}x_2^{7}x_3^{9}x_4^{14}$.

\medskip
For $s \geqslant 4$,\ $a_{s,45} = x_1^{7}x_2^{2^{s}-5}x_3^{2^{s+1}-3}x_4^{2^{s+1}-2}$.

\medskip
By using this basis, we prove the following.

\begin{props}\label{mdt21} For any $s \geqslant 2$, we have
 $QP_4((3)|^{s}|(2))^{GL_4} = 0.$
\end{props}

\begin{rems}
This proposition is proved in \cite{pp25} but the proof is incomplete. To make the paper self-contained we present the detailed proof for it.
\end{rems}

By a direct computation, we have
\begin{align*}
&[\Sigma_4(a_{s,1})]_{(3)|^s|(2)} = \langle \{[a_{s,j}]_{(3)|^s|(2)}: 1 \leqslant j \leqslant 12 \}\rangle,\\
&[\Sigma_4(a_{s,13})]_{(3)|^s|(2)} = \langle \{[a_{s,j}]_{(3)|^s|(2)}: 13 \leqslant j \leqslant 18\} \rangle,\\
&[\Sigma_4(a_{2,19})]_{(3)|^2|(2)} = \langle \{[a_{2,j}]_{(3)|^2|(2)}: 19 \leqslant j \leqslant 34 \mbox{ or } 38 \leqslant j \leqslant 41\}\rangle,\\
&[\Sigma_4(a_{s,19},a_{s,31})]_{(3)|^s|(2)} = \langle \{[a_{s,j}]_{(3)|^s|(2)}: 19 \leqslant j \leqslant 34,\, 38 \leqslant j \leqslant 41 \}\rangle, \ s \geqslant 3,\\
&[\Sigma_4(a_{2,35})]_{(3)|^2|(2)} = \langle \{[a_{2,j}]_{(3)|^2|(2)}: 35 \leqslant j \leqslant 37\} \rangle,\\
&[\Sigma_4(a_{3,35})]_{(3)|^3|(2)} = \langle \{[a_{3,j}]_{(3)|^3|(2)}: 35 \leqslant j \leqslant 37 \mbox{ or } 42 \leqslant j \leqslant 45 \} \rangle,\\
&[\Sigma_4(a_{s,35},a_{s,42})]_{(3)|^s|(2)} = \langle \{[a_{s,j}]_{(3)|^s|(2)}: 35 \leqslant j \leqslant 37,\, 42 \leqslant j \leqslant 45 \} \rangle, \ s \geqslant 4.
\end{align*}

We set $p_{s,u} = p_{2,s,u}$ with
$$p_{s,1} = \sum_{1 \leqslant j \leqslant 12}a_{s,j}, \ p_{s,2} = \sum_{13 \leqslant j \leqslant 18}a_{s,j},\, s \geqslant 2, \ p_{s,3} = \sum_{19\leqslant j \leqslant 34}a_{s,j}, \ s \geqslant 3.$$

We need some lemmas for the proof of the proposition.

\begin{lems}\label{bdt12} We have
$QP_4((3)|^{2}|(2))^{\Sigma_4} = \langle \{[p_{2,u}]_{(3)|^{2}|(2)} : 1 \leqslant u \leqslant 4\} \rangle,$
where
\begin{align*}	
&p_{2,3} = \sum_{21\leqslant j \leqslant 34, j \ne 23,28}a_{2,j}+\sum_{38 \leqslant j \leqslant 41}a_{2,j},\ \quad p_{2,4} = \sum_{35 \leqslant j \leqslant 37}a_{2,j}.
\end{align*}
\end{lems}
\begin{proof} For $s = 2$, we have $QP_4((3)|^2|(2)) = \langle \{[a_{2,j}]_{(3)|^2|(2)}: 1 \leqslant j \leqslant 41\}\rangle$. By a simple computation, we have a direct summand decomposition of $\Sigma_4$-modules:
\begin{align*}QP_4((3)|^2|(2)) &= [\Sigma_4(a_{2,1})]_{(3)|^{2}|(2)}\bigoplus [\Sigma_4(a_{2,13})]_{(3)|^{2}|(2)}\\ &\qquad \bigoplus [\Sigma_4(a_{2,19})]_{(3)|^{2}|(2)}\bigoplus [\Sigma_4(a_{2,35})]_{(3)|^{2}|(2)}.\end{align*}

By a direct computation, we easily obtain
\begin{align*}
&[\Sigma_4(a_{2,1})]_{(3)|^{2}|(2)}^{\Sigma_4} = \langle[p_{2,1}]_{(3)|^{2}|(2)}\rangle,\ [\Sigma_4(a_{2,13})]_{(3)|^{2}|(2)}^{\Sigma_4} = \langle[p_{2,2}]_{(3)|^{2}|(2)}\rangle,\\
&[\Sigma_4(a_{2,19})]_{(3)|^{2}|(2)}^{\Sigma_4} = \langle[p_{2,3}]_{(3)|^{2}|(2)}\rangle,\ [\Sigma_4(a_{2,35})]_{(3)|^{2}|(2)}^{\Sigma_4} = \langle[p_{2,4}]_{(3)|^{2}|(2)}\rangle.
\end{align*}
The lemma follows.
\end{proof}

\begin{lems}\label{bdt13} We have
$QP_4((3)|^{3}|(2))^{\Sigma_4} = \langle \{[p_{3,u}]_{(3)|^{3}|(2)} : 1 \leqslant u \leqslant 5\} \rangle,$
where 
$$p_{3,4} = a_{3,35} + \sum_{42 \leqslant j \leqslant 45}a_{3,j}, \ p_{3,5} = a_{3,36} +a_{3,37} +a_{3,45}.$$
\end{lems}
\begin{proof} 
For $s = 3$, we have $QP_4((3)|^3|(2)) = \langle \{[a_{3,j}]_{(3)|^3|(2)}: 1 \leqslant j \leqslant 45\}\rangle$. We have a direct summand decomposition of $\Sigma_4$-modules:
\begin{align*}
QP_4((3)|^3|(2)) &= [\Sigma_4(a_{3,1})]_{(3)|^3|(2)}\bigoplus [\Sigma_4(a_{3,13})]_{(3)|^3|(2)} \bigoplus\\ &\qquad [\Sigma_4(a_{3,19},a_{3,31})]_{(3)|^3|(2)}\bigoplus [\Sigma_4(a_{3,35})]_{(3)|^3|(2)}.
\end{align*}
By a direct computation, we easily obtain
\begin{align*}
&[\Sigma_4(a_{3,1})]_{(3)|^3|(2)}^{\Sigma_4} = \langle[p_{3,1}]_{(3)|^3|(2)}\rangle,\\ &[\Sigma_4(a_{3,13})]_{(3)|^3|(2)}^{\Sigma_4} = \langle[p_{3,2}]_{(3)|^3|(2)}\rangle,\\
&[\Sigma_4(a_{3,19}, a_{3,31})]_{(3)|^3|(2)}^{\Sigma_4} = \langle[p_{3,3}]_{(3)|^3|(2)}\rangle,\\ 
&[\Sigma_4(a_{3,35})]_{(3)|^3|(2)}^{\Sigma_4} = \langle [p_{3,4}]_{(3)|^3|(2)},[p_{3,5}]_{(3)|^3|(2)}\rangle.
\end{align*}
We present the detailed computation for $[\Sigma_4(a_{3,19}, a_{3,31})]_{(3)|^3|(2)}^{\Sigma_4} = \langle[p_{3,3}]_{(3)|^3|(2)}\rangle$.

Suppose $h \in P_4((3)|^3|(2))$ and $[h]_{(3)|^3|(2)} \in [\Sigma_4(a_{3,19}, a_{3,31})]_{(3)|^3|(2)}^{\Sigma_4}$. Then we have
$$h \equiv_{\omega} \sum_{j=19}^{34}\gamma_ja_{3,j}+ \sum_{j=38}^{41}\gamma_ja_{3,j},$$ 
where $\gamma_j \in \mathbb F_2$ and $\omega = (3)|^3|(2)$. Consider the homomorphisms $\rho_j: P_4 \to P_4$ as defined in Section \ref{s2} with $k =4$ and $1 \leqslant j \leqslant 4$. A direct computation shows
\begin{align*}
\rho_1(h)+h &\equiv_{\omega} \gamma_{\{19,25\}}a_{3,19} + \gamma_{\{20,26\}}a_{3,20} + \gamma_{\{23,28\}}a_{3,23} + \gamma_{\{24,29\}}a_{3,24}\\ 
&\qquad + \gamma_{\{19,25\}}a_{3,25} + \gamma_{\{20,26\}}a_{3,26} + \gamma_{\{27,30\}}a_{3,27} + \gamma_{\{23,28\}}a_{3,28}\\ 
&\qquad + \gamma_{\{24,29\}}a_{3,29} + \gamma_{\{27,30\}}a_{3,30} + \gamma_{\{33,34\}}a_{3,33} + \gamma_{\{33,34\}}a_{3,34}\\ 
&\qquad  + \gamma_{\{21,31\}}a_{3,38} + \gamma_{\{22,32\}}a_{3,39} + \gamma_{\{40,41\}}a_{3,40} + \gamma_{\{40,41\}}a_{3,41} \equiv_{\omega} 0,\\
\rho_2(h)+h &\equiv_{\omega} \gamma_{\{19,21\}}a_{3,19} + \gamma_{\{20,23\}}a_{3,20} + \gamma_{\{19,21\}}a_{3,21} + \gamma_{\{22,24\}}a_{3,22}\\ 
&\qquad + \gamma_{\{20,23\}}a_{3,23} + \gamma_{\{22,24\}}a_{3,24} + \gamma_{\{26,27\}}a_{3,26} + \gamma_{\{26,27\}}a_{3,27}\\ 
&\qquad + \gamma_{\{28,30\}}a_{3,28} + \gamma_{\{28,30\}}a_{3,30} + \gamma_{\{25,31,38\}}a_{3,31} + \gamma_{\{32,33\}}a_{3,32}\\ 
&\qquad + \gamma_{\{32,33\}}a_{3,33} + \gamma_{\{25,31,38\}}a_{3,38} + \gamma_{\{39,40\}}a_{3,39} + \gamma_{\{39,40\}}a_{3,40}\\ 
&\qquad + \gamma_{\{29,34\}}a_{3,41} \equiv_{\omega} 0,\\  
\rho_3(h)+h &\equiv_{\omega} \gamma_{\{19,20\}}a_{3,19} + \gamma_{\{19,20\}}a_{3,20} + \gamma_{\{21,22\}}a_{3,21} + \gamma_{\{21,22\}}a_{3,22}\\ 
&\qquad + \gamma_{\{23,24\}}a_{3,23} + \gamma_{\{23,24\}}a_{3,24} + \gamma_{\{25,26\}}a_{3,25} + \gamma_{\{25,26\}}a_{3,26}\\ 
&\qquad + \gamma_{\{28,29\}}a_{3,28} + \gamma_{\{28,29\}}a_{3,29} + \gamma_{\{31,32\}}a_{3,31} + \gamma_{\{31,32\}}a_{3,32}\\ 
&\qquad + \gamma_{\{27,33,40\}}a_{3,33} + \gamma_{\{30,34,41\}}a_{3,34} + \gamma_{\{38,39\}}a_{3,38} + \gamma_{\{38,39\}}a_{3,39}\\ 
&\qquad + \gamma_{\{27,33,40\}}a_{3,40} + \gamma_{\{30,34,41\}}a_{3,41} \equiv_{\omega} 0.
\end{align*}
From these equalities we get $\gamma_j = \gamma_{19}$  for $20\leqslant j \leqslant 34$ and $\gamma_j = 0$ for $38\leqslant j \leqslant 41$.
The lemma follows.
\end{proof}

\begin{lems}\label{bdt14} For any $s \geqslant 4$, we have
$$QP_4((3)|^{s}|(2))^{\Sigma_4} = \langle \{[p_{s,u}]_{(3)|^{s}|(2)} : 1 \leqslant u \leqslant 5\} \rangle,$$
where 
$$p_{s,4} = a_{s,42}+ a_{s,45} + \sum_{35 \leqslant j \leqslant 37}a_{s,j} , \ p_{s,5} = a_{s,43} +a_{s,44} +a_{s,45}.$$
	
\end{lems}
\begin{proof} 
For $s > 3$, we have $QP_4((3)|^s|(2)) = \langle \{[a_{s,j}]_{(3)|^s|(2)}: 1 \leqslant j \leqslant 45\}\rangle$. There is a direct summand decomposition of $\Sigma_4$-modules:
\begin{align*}QP_4((3)|^s|(2)) &= 
[\Sigma_4(a_{s,1})]_{(3)|^s|(2)}\bigoplus [\Sigma_4(a_{s,13})]_{(3)|^s|(2)} \bigoplus\\ &\qquad [\Sigma_4(a_{s,19},a_{s,31})]_{(3)|^s|(2)} \bigoplus [\Sigma_4(a_{s,35},a_{s,42})]_{(3)|^s|(2)}.
\end{align*}
By a direct computation, we easily obtain
\begin{align*}
&[\Sigma_4(a_{s,1})]_{(3)|^s|(2)}^{\Sigma_4} = \langle[p_{s,1}]_{(3)|^s|(2)}\rangle,\\ &[\Sigma_4(a_{3,13})]_{(3)|^s|(2)}^{\Sigma_4} = \langle[p_{s,2}]_{(3)|^s|(2)}\rangle,\\
&[\Sigma_4(a_{3,19}, a_{3,31})]_{(3)|^s|(2)}^{\Sigma_4} = \langle[p_{s,3}]_{(3)|^s|(2)}\rangle,\\ 
&[\Sigma_4(a_{s,35},a_{s,42})]_{(3)|^s|(2)}^{\Sigma_4} = \langle [p_{s,4}]_{(3)|^s|(2)},[p_{s,5}]_{(3)|^s|(2)}\rangle.
\end{align*}
We present the detailed proof for the case 
$$[\Sigma_4(a_{s,35}, a_{s,42})]_{(3)|^s|(2)}^{\Sigma_4} = \langle[p_{s,4}]_{(3)|^s|(2)},[p_{s,5}]_{(3)|^s|(2)}\rangle.$$
	
Suppose $g \in P_4(\omega^s)$ and $[g]_{\omega^s} \in [\Sigma_4(a_{s,35}, a_{s,42})]_{\omega^s}^{\Sigma_4} $ with $\omega^s = (3)|^s|(2)$. Then we have
$$g \equiv_{\omega^s} \sum_{35\leqslant j \leqslant 37}\gamma_ja_{s,j}+ \sum_{42\leqslant j \leqslant 45}\gamma_ja_{s,j},$$ 
where $\gamma_j \in \mathbb F_2$. By a direct computation we get
\begin{align*}
\rho_1(g) + g &\equiv_{\omega^s} \gamma_{\{35,42\}}a_{s,35} + \gamma_{\{35,42\}}a_{s,42} + \gamma_{\{36,37,43,44\}}a_{s,45} \equiv_{\omega^s} 0,\\
\rho_2(g) + g &\equiv_{\omega^s} \gamma_{\{35,36\}}a_{s,35} + \gamma_{\{35,36\}}a_{s,36} + \gamma_{\{42,43,45\}}a_{s,43} + \gamma_{\{42,43,45\}}a_{s,45} \equiv_{\omega^s} 0,\\  
\rho_3(g) + g &\equiv_{\omega^s} \gamma_{\{36,37\}}a_{s,36} + \gamma_{\{36,37\}}a_{s,37} + \gamma_{\{43,44\}}a_{s,43} + \gamma_{\{43,44\}}a_{s,44} \equiv_{\omega^s} 0.
\end{align*}
From these equalities we get $\gamma_{35} = \gamma_{36}= \gamma_{37} = \gamma_{42}$, $\gamma_{43} = \gamma_{44}$ and $\gamma_{45} = \gamma_{35} + \gamma_{43}$. Hence, we get
$g \equiv_{\omega^s} \gamma_{35}p_{s,4} + \gamma_{43}p_{s,5}.$
The lemma follows.
\end{proof}

\begin{proof}[Proof of Proposition \ref{mdt21}] Suppose $f \in QP_4((3)|^s|(2))$ such that the class $[f]_{\omega^{s}} \in QP_4(\omega^s)^{GL_4}$ with $\omega^s = (3)|^s|(2)$. Then, $[f]_{\omega^s} \in QP_4(\omega^s)^{\Sigma_4}$. 
	
For $s = 2$, by Lemma \ref{bdt12}, we have 
$$ f \equiv_{\omega^{2}} \sum_{1 \leqslant u \leqslant 4}\gamma_up_{2,u},$$
where $\gamma_u \in \mathbb F_2$. By computing $\rho_4(f)+f$ in terms of the admissible monomials, we get
\begin{align*}
\rho_4(f)+f &\equiv_{\omega^{2}} \gamma_{\{1,2\}}a_{2,1} + \gamma_{\{1,3\}}a_{2,2} + \gamma_{\{1,3\}}a_{2,3} + \gamma_{1}a_{2,5} + \gamma_{1}a_{2,6}\\ 
&\qquad + \gamma_{\{2,3\}}a_{2,14} + \gamma_{\{2,3\}}a_{2,15} + \gamma_{3}a_{2,19} + \gamma_{3}a_{2,20} + \gamma_{4}a_{2,23}\\ 
&\qquad + \gamma_{\{3,4\}}a_{2,24} + \gamma_{3}a_{2,27} + \gamma_{3}a_{2,33} + \gamma_{3}a_{2,35} + \gamma_{\{1,2\}}a_{2,38}\\ 
&\qquad + \gamma_{\{1,2\}}a_{2,39} + \gamma_{3}a_{2,40} \equiv_{\omega^{2}} 0.
\end{align*}
This equality implies $\gamma_u = 0$ for $1 \leqslant u \leqslant 4$. Hence, $[f]_{\omega^{2}} = 0$. The proposition is proved for $s = 2$.

For $s = 3$, by using Lemma \ref{bdt13}, we have 
$$ f \equiv_{\omega^{3}} \sum_{1 \leqslant u \leqslant 5}\gamma_up_{3,u},$$
where $\gamma_u \in \mathbb F_2$. Computing $\rho_4(f)+f$ in terms of the admissible monomials gives
\begin{align*}
\rho_4(f)+f &\equiv_{\omega^{3}} \gamma_{\{1,2\}}a_{3,1} + \gamma_{\{1,3\}}a_{3,2} + \gamma_{\{1,3\}}a_{3,3} + \gamma_{1}a_{3,5} + \gamma_{1}a_{3,6}\\ 
&\qquad + \gamma_{\{2,3\}}a_{3,14} + \gamma_{\{2,3\}}a_{3,15} + \gamma_{\{1,2,3\}}a_{3,19} + \gamma_{\{1,2,3\}}a_{3,20}\\ 
&\qquad + \gamma_{\{3,5\}}a_{3,23} + \gamma_{\{3,5\}}a_{3,24} + \gamma_{\{1,2\}}a_{3,25} + \gamma_{\{1,2\}}a_{3,26}\\ 
&\qquad + \gamma_{3}a_{3,27} + \gamma_{\{3,4\}}a_{3,33} + \gamma_{4}a_{3,35} + \gamma_{\{1,2\}}a_{3,38} + \gamma_{\{1,2\}}a_{3,39}\\ 
&\qquad + \gamma_{4}a_{3,40} + \gamma_{3}a_{3,45} \equiv_{\omega^{3}} 0.
\end{align*}
From this equality we get $\gamma_u = 0$ for $1 \leqslant u \leqslant 5$. Hence, $[f]_{\omega^{3}} = 0$. The proposition is proved for $s = 3$.

For $s\geqslant 4$, by using Lemma \ref{bdt14}, we have 
$$ f \equiv_{\omega^{s}} \sum_{1 \leqslant u \leqslant 5}\gamma_up_{s,u},$$
where $\gamma_u \in \mathbb F_2$. By computing $\rho_4(f)+f$ in terms of the admissible monomials, we obtain
\begin{align*}
\rho_4(f)+f &\equiv_{\omega^{s}} \gamma_{\{1,2\}}a_{s,1} + \gamma_{\{1,3\}}a_{s,2} + \gamma_{\{1,3\}}a_{s,3} + \gamma_{1}a_{s,5} + \gamma_{1}a_{s,6}\\ 
&\qquad + \gamma_{\{2,3\}}a_{s,14} + \gamma_{\{2,3\}}a_{s,15} + \gamma_{\{1,2,3\}}a_{s,19} + \gamma_{\{1,2,3\}}a_{s,20}\\ 
&\qquad + \gamma_{\{3,4\}}a_{s,23} + \gamma_{\{3,4\}}a_{s,24} + \gamma_{\{1,2\}}a_{s,25} + \gamma_{\{1,2\}}a_{s,26}\\ 
&\qquad + \gamma_{3}a_{s,27} + \gamma_{\{3,5\}}a_{s,33} + \gamma_{\{3,5\}}a_{s,35} + \gamma_{\{1,2\}}a_{s,38}\\ 
&\qquad + \gamma_{\{1,2\}}a_{s,39} + \gamma_{5}a_{s,40} + \gamma_{3}a_{s,42} + \gamma_{3}a_{s,45} \equiv_{\omega^{s}} 0.
\end{align*}
From this equality we get $\gamma_u = 0$ for $1 \leqslant u \leqslant 5$. Hence, $[f]_{\omega^{s}} = 0$. The proposition is completely proved.
\end{proof}

\subsubsection{\textbf{The case $t = 3$}}\

\medskip
By \cite{su5,su50}, a basis of $QP_4((3)|^{s}|(2)|^2)$ is the set of all classes represented by the admissible monomials $a_{s,j} = a_{3,s,j}$ which are determined as follows:

\medskip
\centerline{\begin{tabular}{lll} 
$a_{s,1} = x_2^{2^{s}-1}x_3^{2^{s+2}-1}x_4^{2^{s+2}-1}$&$a_{s,2} = x_2^{2^{s+2}-1}x_3^{2^{s}-1}x_4^{2^{s+2}-1}$\cr  $a_{s,3} = x_2^{2^{s+2}-1}x_3^{2^{s+2}-1}x_4^{2^{s}-1}$&$a_{s,4} = x_1^{2^{s}-1}x_3^{2^{s+2}-1}x_4^{2^{s+2}-1}$\cr  $a_{s,5} = x_1^{2^{s}-1}x_2^{2^{s+2}-1}x_4^{2^{s+2}-1}$&$a_{s,6} = x_1^{2^{s}-1}x_2^{2^{s+2}-1}x_3^{2^{s+2}-1}$\cr  $a_{s,7} = x_1^{2^{s+2}-1}x_3^{2^{s}-1}x_4^{2^{s+2}-1}$&$a_{s,8} = x_1^{2^{s+2}-1}x_3^{2^{s+2}-1}x_4^{2^{s}-1}$\cr  $a_{s,9} = x_1^{2^{s+2}-1}x_2^{2^{s}-1}x_4^{2^{s+2}-1}$&$a_{s,10} = x_1^{2^{s+2}-1}x_2^{2^{s}-1}x_3^{2^{s+2}-1}$\cr  $a_{s,11} = x_1^{2^{s+2}-1}x_2^{2^{s+2}-1}x_4^{2^{s}-1}$&$a_{s,12} = x_1^{2^{s+2}-1}x_2^{2^{s+2}-1}x_3^{2^{s}-1}$\cr  $a_{s,13} = x_2^{2^{s+1}-1}x_3^{3.2^{s}-1}x_4^{2^{s+2}-1}$&$a_{s,14} = x_2^{2^{s+1}-1}x_3^{2^{s+2}-1}x_4^{3.2^{s}-1}$\cr  $a_{s,15} = x_2^{2^{s+2}-1}x_3^{2^{s+1}-1}x_4^{3.2^{s}-1}$&$a_{s,16} = x_1^{2^{s+1}-1}x_3^{3.2^{s}-1}x_4^{2^{s+2}-1}$\cr  $a_{s,17} = x_1^{2^{s+1}-1}x_3^{2^{s+2}-1}x_4^{3.2^{s}-1}$&$a_{s,18} = x_1^{2^{s+1}-1}x_2^{3.2^{s}-1}x_4^{2^{s+2}-1}$\cr  $a_{s,19} = x_1^{2^{s+1}-1}x_2^{3.2^{s}-1}x_3^{2^{s+2}-1}$&$a_{s,20} = x_1^{2^{s+1}-1}x_2^{2^{s+2}-1}x_4^{3.2^{s}-1}$\cr  $a_{s,21} = x_1^{2^{s+1}-1}x_2^{2^{s+2}-1}x_3^{3.2^{s}-1}$&$a_{s,22} = x_1^{2^{s+2}-1}x_3^{2^{s+1}-1}x_4^{3.2^{s}-1}$\cr  $a_{s,23} = x_1^{2^{s+2}-1}x_2^{2^{s+1}-1}x_4^{3.2^{s}-1}$&$a_{s,24} = x_1^{2^{s+2}-1}x_2^{2^{s+1}-1}x_3^{3.2^{s}-1}$\cr  $a_{s,25} = x_1x_2^{2^{s}-2}x_3^{2^{s+2}-1}x_4^{2^{s+2}-1}$&$a_{s,26} = x_1x_2^{2^{s+2}-1}x_3^{2^{s}-2}x_4^{2^{s+2}-1}$\cr  $a_{s,27} = x_1x_2^{2^{s+2}-1}x_3^{2^{s+2}-1}x_4^{2^{s}-2}$&$a_{s,28} = x_1^{2^{s+2}-1}x_2x_3^{2^{s}-2}x_4^{2^{s+2}-1}$\cr  $a_{s,29} = x_1^{2^{s+2}-1}x_2x_3^{2^{s+2}-1}x_4^{2^{s}-2}$&$a_{s,30} = x_1^{2^{s+2}-1}x_2^{2^{s+2}-1}x_3x_4^{2^{s}-2}$\cr  $a_{s,31} = x_1x_2^{2^{s}-1}x_3^{2^{s+2}-2}x_4^{2^{s+2}-1}$&$a_{s,32} = x_1x_2^{2^{s}-1}x_3^{2^{s+2}-1}x_4^{2^{s+2}-2}$\cr  $a_{s,33} = x_1x_2^{2^{s+2}-2}x_3^{2^{s}-1}x_4^{2^{s+2}-1}$&$a_{s,34} = x_1x_2^{2^{s+2}-2}x_3^{2^{s+2}-1}x_4^{2^{s}-1}$\cr  $a_{s,35} = x_1x_2^{2^{s+2}-1}x_3^{2^{s}-1}x_4^{2^{s+2}-2}$&$a_{s,36} = x_1x_2^{2^{s+2}-1}x_3^{2^{s+2}-2}x_4^{2^{s}-1}$\cr  $a_{s,37} = x_1^{2^{s}-1}x_2x_3^{2^{s+2}-2}x_4^{2^{s+2}-1}$&$a_{s,38} = x_1^{2^{s}-1}x_2x_3^{2^{s+2}-1}x_4^{2^{s+2}-2}$\cr  $a_{s,39} = x_1^{2^{s}-1}x_2^{2^{s+2}-1}x_3x_4^{2^{s+2}-2}$&$a_{s,40} = x_1^{2^{s+2}-1}x_2x_3^{2^{s}-1}x_4^{2^{s+2}-2}$\cr  $a_{s,41} = x_1^{2^{s+2}-1}x_2x_3^{2^{s+2}-2}x_4^{2^{s}-1}$&$a_{s,42} = x_1^{2^{s+2}-1}x_2^{2^{s}-1}x_3x_4^{2^{s+2}-2}$\cr  $a_{s,43} = x_1x_2^{2^{s+1}-2}x_3^{3.2^{s}-1}x_4^{2^{s+2}-1}$&$a_{s,44} = x_1x_2^{2^{s+1}-2}x_3^{2^{s+2}-1}x_4^{3.2^{s}-1}$\cr  $a_{s,45} = x_1x_2^{2^{s+2}-1}x_3^{2^{s+1}-2}x_4^{3.2^{s}-1}$&$a_{s,46} = x_1^{2^{s+2}-1}x_2x_3^{2^{s+1}-2}x_4^{3.2^{s}-1}$\cr  $a_{s,47} = x_1x_2^{2^{s+1}-1}x_3^{3. 2^{s}-2}x_4^{2^{s+2}-1}$&$a_{s,48} = x_1x_2^{2^{s+1}-1}x_3^{2^{s+2}-1}x_4^{3. 2^{s}-2}$\cr  $a_{s,49} = x_1x_2^{2^{s+2}-1}x_3^{2^{s+1}-1}x_4^{3. 2^{s}-2}$&$a_{s,50} = x_1^{2^{s+1}-1}x_2x_3^{3. 2^{s}-2}x_4^{2^{s+2}-1}$\cr  $a_{s,51} = x_1^{2^{s+1}-1}x_2x_3^{2^{s+2}-1}x_4^{3. 2^{s}-2}$&$a_{s,52} = x_1^{2^{s+1}-1}x_2^{2^{s+2}-1}x_3x_4^{3. 2^{s}-2}$\cr  $a_{s,53} = x_1^{2^{s+2}-1}x_2x_3^{2^{s+1}-1}x_4^{3. 2^{s}-2}$&$a_{s,54} = x_1^{2^{s+2}-1}x_2^{2^{s+1}-1}x_3x_4^{3. 2^{s}-2}$\cr  $a_{s,55} = x_1x_2^{2^{s+1}-1}x_3^{3.2^{s}-1}x_4^{2^{s+2}-2}$&$a_{s,56} = x_1x_2^{2^{s+1}-1}x_3^{2^{s+2}-2}x_4^{3.2^{s}-1}$\cr  
\end{tabular}}
\centerline{\begin{tabular}{lll} 
$a_{s,57} = x_1x_2^{2^{s+2}-2}x_3^{2^{s+1}-1}x_4^{3.2^{s}-1}$&$a_{s,58} = x_1^{2^{s+1}-1}x_2x_3^{3.2^{s}-1}x_4^{2^{s+2}-2}$\cr  $a_{s,59} = x_1^{2^{s+1}-1}x_2x_3^{2^{s+2}-2}x_4^{3.2^{s}-1}$&$a_{s,60} = x_1^{2^{s+1}-1}x_2^{3.2^{s}-1}x_3x_4^{2^{s+2}-2}$\cr  $a_{s,61} = x_1^{3}x_2^{2^{s+2}-3}x_3^{2^{s}-2}x_4^{2^{s+2}-1}$&$a_{s,62} = x_1^{3}x_2^{2^{s+2}-3}x_3^{2^{s+2}-1}x_4^{2^{s}-2}$\cr  $a_{s,63} = x_1^{3}x_2^{2^{s+2}-1}x_3^{2^{s+2}-3}x_4^{2^{s}-2}$&$a_{s,64} = x_1^{2^{s+2}-1}x_2^{3}x_3^{2^{s+2}-3}x_4^{2^{s}-2}$\cr  $a_{s,65} = x_1^{3}x_2^{2^{s}-1}x_3^{2^{s+2}-3}x_4^{2^{s+2}-2}$&$a_{s,66} = x_1^{3}x_2^{2^{s+2}-3}x_3^{2^{s}-1}x_4^{2^{s+2}-2}$\cr  $a_{s,67} = x_1^{3}x_2^{2^{s+2}-3}x_3^{2^{s+2}-2}x_4^{2^{s}-1}$&$a_{s,68} = x_1^{3}x_2^{2^{s+1}-3}x_3^{3. 2^{s}-2}x_4^{2^{s+2}-1}$\cr  $a_{s,69} = x_1^{3}x_2^{2^{s+1}-3}x_3^{2^{s+2}-1}x_4^{3. 2^{s}-2}$&$a_{s,70} = x_1^{3}x_2^{2^{s+2}-1}x_3^{2^{s+1}-3}x_4^{3. 2^{s}-2}$\cr  $a_{s,71} = x_1^{2^{s+2}-1}x_2^{3}x_3^{2^{s+1}-3}x_4^{3. 2^{s}-2}$&$a_{s,72} = x_1^{3}x_2^{2^{s+1}-3}x_3^{3.2^{s}-1}x_4^{2^{s+2}-2}$\cr  $a_{s,73} = x_1^{3}x_2^{2^{s+1}-3}x_3^{2^{s+2}-2}x_4^{3.2^{s}-1}$&$a_{s,74} = x_1^{3}x_2^{2^{s+2}-3}x_3^{2^{s+1}-2}x_4^{3.2^{s}-1}$\cr  $a_{s,75} = x_1^{3}x_2^{2^{s+1}-1}x_3^{3.2^{s}-3}x_4^{2^{s+2}-2}$&  $a_{s,76} = x_1^{2^{s+1}-1}x_2^{3}x_3^{3.2^{s}-3}x_4^{2^{s+2}-2}$\cr  $a_{s,77} = x_1^{3}x_2^{2^{s+1}-1}x_3^{2^{s+2}-3}x_4^{3. 2^{s}-2}$&$a_{s,78} = x_1^{3}x_2^{2^{s+2}-3}x_3^{2^{s+1}-1}x_4^{3. 2^{s}-2}$\cr  $a_{s,79} = x_1^{2^{s+1}-1}x_2^{3}x_3^{2^{s+2}-3}x_4^{3. 2^{s}-2}$&$a_{s,80} = x_1^{7}x_2^{2^{s+2}-5}x_3^{2^{s+1}-3}x_4^{3. 2^{s}-2}$\cr  $a_{s,81} = x_1^{7}x_2^{2^{s+2}-5}x_3^{2^{s+2}-3}x_4^{2^{s}-2}$&\cr 
\end{tabular}}

\medskip
For $s = 2$,

\medskip
\centerline{\begin{tabular}{lll} 
$a_{2,82} =  x_1^{3}x_2^{3}x_3^{12}x_4^{15}$& $a_{2,83} =  x_1^{3}x_2^{3}x_3^{15}x_4^{12}$& $a_{2,84} =  x_1^{3}x_2^{15}x_3^{3}x_4^{12}$\cr  $a_{2,85} =  x_1^{15}x_2^{3}x_3^{3}x_4^{12}$& $a_{2,86} =  x_1^{3}x_2^{7}x_3^{11}x_4^{12}$& $a_{2,87} =  x_1^{7}x_2^{3}x_3^{11}x_4^{12}$\cr  $a_{2,88} =  x_1^{7}x_2^{11}x_3^{3}x_4^{12}$& $a_{2,89} =  x_1^{7}x_2^{7}x_3^{8}x_4^{11}$& $a_{2,90} =  x_1^{7}x_2^{7}x_3^{11}x_4^{8}$\cr  $a_{2,91} =  x_1^{7}x_2^{7}x_3^{9}x_4^{10}$& \cr   
\end{tabular}}

\medskip
For $s \geqslant 3$,

\medskip
\centerline{\begin{tabular}{lll}
$a_{s,82} = x_1^{3}x_2^{2^{s}-3}x_3^{2^{s+2}-2}x_4^{2^{s+2}-1}$&$a_{s,83} = x_1^{3}x_2^{2^{s}-3}x_3^{2^{s+2}-1}x_4^{2^{s+2}-2}$\cr  $a_{s,84} = x_1^{3}x_2^{2^{s+2}-1}x_3^{2^{s}-3}x_4^{2^{s+2}-2}$&$a_{s,85} = x_1^{2^{s+2}-1}x_2^{3}x_3^{2^{s}-3}x_4^{2^{s+2}-2}$\cr  $a_{s,86} = x_1^{2^{s}-1}x_2^{3}x_3^{2^{s+2}-3}x_4^{2^{s+2}-2}$&$a_{s,87} = x_1^{7}x_2^{2^{s+2}-5}x_3^{2^{s}-3}x_4^{2^{s+2}-2}$\cr  $a_{s,88} = x_1^{7}x_2^{2^{s+1}-5}x_3^{3.2^{s}-3}x_4^{2^{s+2}-2}$&$a_{s,89} = x_1^{7}x_2^{2^{s+1}-5}x_3^{2^{s+2}-3}x_4^{3. 2^{s}-2}$\cr  
\end{tabular}}

\medskip
For $s= 3$, $a_{3,90} =  x_1^{7}x_2^{7}x_3^{25}x_4^{30}$.

\medskip
For $s \geqslant 4$, $a_{s,90} = x_1^{7}x_2^{2^{s}-5}x_3^{2^{s+2}-3}x_4^{2^{s+2}-2}$.

\medskip
By using this basis, we prove the following.

\begin{props}\label{mdt31} For any $s \geqslant 2$, we have
$$QP_4((3)|^{s}|(2)|^2)^{GL_4} = \begin{cases}\langle [\bar\xi_{3,2}]_{(3)|^{2}|(2)|^2} \rangle, &\mbox{if } s = 2,\\ 0 &\mbox{if } s \geqslant 3,\end{cases}$$
where $\bar\xi_{3,2} = \sum_{j \in\{55,\, 56,\, 58,\, 59,\, 72,\, 73,\, 89,\, 90,\, 91\}}a_{3,j}$.
\end{props}

By a direct computation, we have
\begin{align*}
&[\Sigma_4(a_{s,1})]_{(3)|^s|(2)|^2} = \langle \{[a_{s,j}]_{(3)|^s|(2)|^2}: 1 \leqslant j \leqslant 12 \}\rangle,\\
&[\Sigma_4(a_{s,13})]_{(3)|^s|(2)|^2} = \langle \{[a_{s,j}]_{(3)|^s|(2)|^2}: 13 \leqslant j \leqslant 24 \}\rangle,\\
&[\Sigma_4(a_{s,25})]_{(3)|^s|(2)|^2} = \langle \{[a_{s,j}]_{(3)|^s|(2)|^2}: 25 \leqslant j \leqslant 30\} \rangle.
\end{align*}

We set $p_{s,u} = p_{3,s,u}$ with
\begin{equation}\label{ct31}
p_{s,1} = \sum_{1 \leqslant j \leqslant 12}a_{s,j}, \ p_{s,2} = \sum_{13 \leqslant j \leqslant 24}a_{s,j}, \ p_{s,3} = \sum_{25\leqslant j \leqslant 30}a_{s,j}.
\end{equation}

We need some lemmas for the proof of the proposition. 

\begin{lems}\label{bdt32} We have
$QP_4((3)|^{2}|(2)|^2)^{\Sigma_4} = \langle \{[p_{2,u}]_{(3)|^{2}|(2)|^2} : 1 \leqslant u \leqslant 9\} \rangle,$
where $p_{2,u},\, u =1,\, 2,\, 3$, are defined by \eqref{ct31} and
$$\begin{cases}	
p_{2,4} = \sum_{j \in \{31,\, 32,\, 33,\, 34,\, 35,\, 36,\, 40,\, 41,\, 68,\, 69,\, 70,\, 71\}}a_{2,j},\\
p_{2,5} = \sum_{j \in \{33,\, 34,\, 36,\, 37,\, 38,\, 82,\, 83,\, 61,\, 62,\, 39,\, 84,\, 63,\, 41,\, 42,\, 85,\, 64\}}a_{2,j},\\
p_{2,6} = \sum_{j \in \{31,\, 32,\, 33,\, 34,\, 35,\, 36,\, 37,\, 38,\, 39,\, 40,\, 41,\, 42,\, 43,\, 44,\, 45,\, 46\},\atop\hskip2.7cm\cup\{ 47,\, 48,\, 49,\, 50,\, 51,\, 52,\, 53,\, 54\}}a_{2,j},\\
p_{2,7} = \sum_{j \in \{65,\, 66,\, 67\}}a_{2,j},\\
p_{2,8} = \sum_{j \in \{56,\, 59,\, 60,\, 66,\, 67,\, 72,\, 73,\, 81,\, 86,\, 87,\, 89,\, 90\}}a_{2,j},\\
p_{2,9} = \sum_{j \in \{55,\, 56,\, 58,\, 59,\, 72,\, 73,\, 89,\, 90,\, 91\}}a_{2,j}.
	\end{cases}$$
\end{lems}
\begin{proof} For $s = 2$, we have $QP_4((3)|^2|(2)|^2) = \langle \{[a_{2,j}]_{(3)|^2|(2)|^2}: 1 \leqslant j \leqslant 91\}\rangle$. By a simple computation, we see that $[p_{2,u}]_{(3)|^2|(2)|^2}$ is an $\Sigma_4$ invariant for $1 \leqslant u \leqslant 9$, and there is a direct summand decomposition of $\Sigma_4$-modules:
\begin{align*}QP_4((3)|^2|(2)|^2) &= [\Sigma_4(a_{2,1})]_{(3)|^2|(2)|^2}\bigoplus [\Sigma_4(a_{2,13})]_{(3)|^2|(2)|^2}\\ &\qquad \bigoplus [\Sigma_4(a_{2,25})]_{(3)|^2|(2)|^2}\bigoplus \mathcal U_{3,2} \oplus \mathcal V_{3,2},
\end{align*}
where 
\begin{align*}
\mathcal U_{3,2} &= \langle [a_{2,j}]_{(3)|^{2}|(2)|^2}: j \in \mathbb J_{3,2} = \{31,\ldots 54,\, 61,\ldots 64,\, 68, \ldots,71,\, 82,\ldots,85\}\rangle,\\
\mathcal V_{3,2} &= \langle [a_{2,j}]_{(3)|^{2}|(2)|^2}: j \in \mathbb K_{3,2} =  \langle\{55,\ldots 60,\, 65,\, 66,\, 67,\, 72, \ldots,81,86,\ldots,91\}\rangle.
\end{align*}
By a simple computation, we easily obtain
\begin{align*}
&[\Sigma_4(a_{2,1})]_{(3)|^{2}|(2)|^2}^{\Sigma_4} = \langle[p_{2,1}]_{(3)|^{2}|(2)|^2}\rangle,\\ &[\Sigma_4(a_{2,13})]_{(3)|^{2}|(2)|^2}^{\Sigma_4} = \langle[p_{2,2}]_{(3)|^{2}|(2)|^2}\rangle,\\
&[\Sigma_4(a_{2,25})]_{(3)|^{2}|(2)|^2}^{\Sigma_4} = \langle[p_{2,3}]_{(3)|^{2}|(2)|^2}\rangle.
\end{align*}
We prove $\mathcal U_{3,2}^{\Sigma_4} = \langle [p_{2,u}]_{(3)|^{2}|(2)|^2} : 4 \leqslant u \leqslant 6\}\rangle$. We observe that the leading monomials of $p_{2,4}$, $p_{2,5}$, $p_{2,6}$ respectively are $a_{2,71}$, $a_{2,64}$, $a_{2,54}$. So, if $[f]_{(3)|^{2}|(2)|^2} \in \mathcal U_{3,2}^{\Sigma_4}$ with $f \in P_4$, then there are $\gamma_4,\, \gamma_5,\, \gamma_6 \in \mathbb F_2$ such that
$$g := f + \gamma_4p_{2,4}+ \gamma_5p_{2,5} + \gamma_6p_{2,6}  \equiv _{(3)|^{2}|(2)|^2} \sum_{j\in \mathbb J_{3,2},\, j \ne 54,64,71}a_{2,j}.$$
Then, $[g]_{(3)|^{2}|(2)|^2}$ is also an $\Sigma_4$-invariant. By a direct computation we get
\begin{align*}
\rho_1(g)& + g \equiv_{(3)|^{2}|(2)|^2} \gamma_{\{31,33,37,43\}}a_{2,31} + \gamma_{\{32,34,38,44\}}a_{2,32} + \gamma_{\{35,40\}}a_{2,35}\\ 
& + \gamma_{\{36,41\}}a_{2,36} + \gamma_{\{31,33,37,43\}}a_{2,37} + \gamma_{\{32,34,38,44\}}a_{2,38} + \gamma_{\{39,42\}}a_{2,39}\\ 
& + \gamma_{\{35,40\}}a_{2,40} + \gamma_{\{36,41\}}a_{2,41} + \gamma_{\{39,42\}}a_{2,42} + \gamma_{\{45,46\}}a_{2,45}\\ 
& + \gamma_{\{45,46\}}a_{2,46} + \gamma_{\{47,50\}}a_{2,47} + \gamma_{\{48,51\}}a_{2,48} + \gamma_{\{49,53\}}a_{2,49}\\ 
& + \gamma_{\{47,50\}}a_{2,50} + \gamma_{\{48,51\}}a_{2,51} + \gamma_{52}a_{2,52} + \gamma_{\{49,53\}}a_{2,53} + \gamma_{52}a_{2,54}\\ 
& + \gamma_{63}a_{2,63} + \gamma_{63}a_{2,64} + \gamma_{70}a_{2,70} + \gamma_{70}a_{2,71} + \gamma_{\{33,43,61,68\}}a_{2,82}\\ 
& + \gamma_{\{34,44,62,69\}}a_{2,83} + \gamma_{\{84,85\}}a_{2,84} + \gamma_{\{84,85\}}a_{2,85} \equiv_{(3)|^{2}|(2)|^2} 0,\\
\rho_2(g) &+ g \equiv_{(3)|^{2}|(2)|^2} \gamma_{\{31,33,82\}}a_{2,31} + \gamma_{\{32,35\}}a_{2,32} + \gamma_{\{31,33,82\}}a_{2,33}\\ 
& + \gamma_{\{34,36\}}a_{2,34} + \gamma_{\{32,35\}}a_{2,35} + \gamma_{\{34,36\}}a_{2,36} + \gamma_{\{37,50,61\}}a_{2,37}\\ 
& + \gamma_{\{38,39\}}a_{2,38} + \gamma_{\{38,39\}}a_{2,39} + \gamma_{\{40,41,42,46\}}a_{2,40} + \gamma_{\{40,41,42,46\}}a_{2,42}\\ 
& + \gamma_{\{43,47\}}a_{2,43} + \gamma_{\{44,45\}}a_{2,44} + \gamma_{\{44,45\}}a_{2,45} + \gamma_{\{43,47\}}a_{2,47}\\ 
& + \gamma_{\{48,49\}}a_{2,48} + \gamma_{\{48,49\}}a_{2,49} + \gamma_{\{51,52\}}a_{2,51} + \gamma_{\{51,52\}}a_{2,52} + \gamma_{53}a_{2,53}\\ 
& + \gamma_{53}a_{2,54} + \gamma_{\{37,50,61\}}a_{2,61} + \gamma_{\{62,63\}}a_{2,62} + \gamma_{\{62,63\}}a_{2,63} + \gamma_{\{69,70\}}a_{2,69}\\ 
& + \gamma_{\{69,70\}}a_{2,70} + \gamma_{\{83,84\}}a_{2,83} + \gamma_{\{83,84\}}a_{2,84} + \gamma_{\{41,46\}}a_{2,85} \equiv_{(3)|^{2}|(2)|^2} 0,\\  
\rho_3(g) &+ g \equiv_{(3)|^{2}|(2)|^2} \gamma_{\{31,32\}}a_{2,31} + \gamma_{\{31,32\}}a_{2,32} + \gamma_{\{33,34\}}a_{2,33} + \gamma_{\{33,34\}}a_{2,34}\\ 
& + \gamma_{\{35,36,84\}}a_{2,35} + \gamma_{\{35,36,84\}}a_{2,36} + \gamma_{\{37,38\}}a_{2,37} + \gamma_{\{37,38\}}a_{2,38}\\ 
& + \gamma_{\{39,52,63\}}a_{2,39} + \gamma_{\{40,41,85\}}a_{2,40} + \gamma_{\{40,41,85\}}a_{2,41} + \gamma_{42}a_{2,42}\\ 
& + \gamma_{\{43,44\}}a_{2,43} + \gamma_{\{43,44\}}a_{2,44} + \gamma_{\{45,49\}}a_{2,45} + \gamma_{\{46,53\}}a_{2,46}\\ 
& + \gamma_{\{47,48\}}a_{2,47} + \gamma_{\{47,48\}}a_{2,48} + \gamma_{\{45,49\}}a_{2,49} + \gamma_{\{50,51\}}a_{2,50}\\ 
& + \gamma_{\{50,51\}}a_{2,51} + \gamma_{\{46,53\}}a_{2,53} + \gamma_{\{61,62\}}a_{2,61} + \gamma_{\{61,62\}}a_{2,62}\\ 
& + \gamma_{\{39,52,63\}}a_{2,63} + \gamma_{42}a_{2,64} + \gamma_{\{68,69\}}a_{2,68} + \gamma_{\{68,69\}}a_{2,69}\\ 
& + \gamma_{\{82,83\}}a_{2,82} + \gamma_{\{82,83\}}a_{2,83} \equiv_{(3)|^{2}|(2)|^2} 0.
\end{align*}
From these equalities we get $\gamma_{j} = 0$ for all $j \in \mathbb J_{3,2}\setminus\{54,64,71\}$. Hence, we get $g \equiv_{(3)|^{2}|(2)|^2} 0$ and
$f \equiv_{(3)|^{2}|(2)|^2}  \gamma_4p_{2,4}+ \gamma_5p_{2,5} + \gamma_6p_{2,6}.$

Now we prove that $\mathcal V_{3,2}^{\Sigma_4} = \langle [p_{2,u}]_{(3)|^{2}|(2)|^2} : 7 \leqslant u \leqslant 9\}\rangle$. It is easy to see that the leading monomials of $p_{2,7}$, $p_{2,8}$, $p_{2,9}$ respectively are $a_{2,67}$, $a_{2,81}$, $a_{2,90}$. Hence, if $[h]_{(3)|^{2}|(2)|^2} \in \mathcal V_{3,2}^{\Sigma_4}$ with $h \in P_4$, then there are $\gamma_7,\, \gamma_8,\, \gamma_9 \in \mathbb F_2$ such that
$$\bar h := h + \gamma_7p_{2,7}+ \gamma_8p_{2,8} + \gamma_9p_{2,9}  \equiv _{(3)|^{2}|(2)|^2} \sum_{j\in \mathbb K_{3,2},\, j \ne 67,81,90}a_{2,j}.$$
Since $[h]_{(3)|^{2}|(2)|^2}$ is an $\Sigma_4$-invariant, so is $[\bar h]_{(3)|^{2}|(2)|^2}$. by a direct computation we get
\begin{align*}
\rho_1(\bar h)& + \bar h \equiv_{(3)|^{2}|(2)|^2} \gamma_{\{55,57,58,66,78\}}a_{2,55} + \gamma_{\{56,59,66,78\}}a_{2,56} + \gamma_{\{55,57,58,66,78\}}a_{2,58}\\
& + \gamma_{\{56,59,66,78\}}a_{2,59} + \gamma_{\{57,66,72,73,74,78\}}a_{2,65} + \gamma_{\{57,66,74,78\}}a_{2,72}\\
& + \gamma_{\{57,66,74,78\}}a_{2,73} + \gamma_{\{74,75,76\}}a_{2,75} + \gamma_{\{74,75,76\}}a_{2,76} + \gamma_{\{77,78,79\}}a_{2,77}\\
& + \gamma_{\{77,78,79\}}a_{2,79} + \gamma_{\{57,74,78,86,87\}}a_{2,86} + \gamma_{\{57,74,78,86,87\}}a_{2,87} + \gamma_{\{66,78\}}a_{2,89}\\
& + \gamma_{\{66,74\}}a_{2,90} + \gamma_{\{66,74,78,80\}}a_{2,91} \equiv_{(3)|^{2}|(2)|^2} 0,\\
\rho_2(\bar h) &+ \bar h \equiv_{(3)|^{2}|(2)|^2} \gamma_{\{59,76,89,91\}}a_{2,55} + \gamma_{\{56,57,59,76,91\}}a_{2,56}\\
& + \gamma_{\{56,57,89\}}a_{2,57} + \gamma_{\{58,60,76,89,91\}}a_{2,58} + \gamma_{\{59,76,89,91\}}a_{2,59} + \gamma_{\{58,59,60\}}a_{2,60}\\
& + \gamma_{\{65,66,76,89,91\}}a_{2,65} + \gamma_{\{59,65,66,76\}}a_{2,66} + \gamma_{\{72,75,76,89,91\}}a_{2,72}\\
& + \gamma_{\{73,74,76,89,91\}}a_{2,73} + \gamma_{\{59,73,74\}}a_{2,74} + \gamma_{\{59,72,75\}}a_{2,75} + \gamma_{\{77,78,91\}}a_{2,77}\\
& + \gamma_{\{77,78,91\}}a_{2,78} + \gamma_{\{79,80\}}a_{2,79} + \gamma_{\{79,80\}}a_{2,80} + \gamma_{\{59,89,91\}}a_{2,86}\\
& + \gamma_{\{80,87,88,89,91\}}a_{2,87} + \gamma_{\{59,79,87,88\}}a_{2,88} + \gamma_{\{59,76,89,91\}}a_{2,89}\\
& + \gamma_{76}a_{2,90} + \gamma_{76}a_{2,91} \equiv_{(3)|^{2}|(2)|^2} 0,\\  
\rho_3(\bar h) &+ \bar h \equiv_{(3)|^{2}|(2)|^2} \gamma_{\{55,56,86,88\}}a_{2,55} + \gamma_{\{55,56,86,88\}}a_{2,56} + \gamma_{\{58,59,87,88\}}a_{2,58}\\
& + \gamma_{\{58,59,87,88\}}a_{2,59} + \gamma_{60}a_{2,60} + \gamma_{\{66,88\}}a_{2,66} + \gamma_{\{66,88\}}a_{2,67}\\
& + \gamma_{\{72,73,88\}}a_{2,72} + \gamma_{\{72,73,88\}}a_{2,73} + \gamma_{\{74,78\}}a_{2,74} + \gamma_{\{75,77\}}a_{2,75}\\
& + \gamma_{\{76,79\}}a_{2,76} + \gamma_{\{75,77\}}a_{2,77} + \gamma_{\{74,78\}}a_{2,78} + \gamma_{\{76,79\}}a_{2,79}\\
& + \gamma_{60}a_{2,81} + \gamma_{\{88,89\}}a_{2,89} + \gamma_{\{88,89\}}a_{2,90} + \gamma_{80}a_{2,91} \equiv_{(3)|^{2}|(2)|^2} 0.
\end{align*}
These equalities imply $\gamma_{j} = 0$ for all $j \in \mathbb K_{3,2}\setminus\{67,81,90\}$ and $\bar h \equiv_{(3)|^{2}|(2)|^2} 0$. Hence, we obtain
$h \equiv_{(3)|^{2}|(2)|^2}  \gamma_7p_{2,7}+ \gamma_8p_{2,8} + \gamma_9p_{2,9}.$
The lemma is proved.
\end{proof}

\begin{lems}\label{bdt33} For $s \geqslant 3$, we have
$$QP_4((3)|^{s}|(2)|^2)^{\Sigma_4} = \langle \{[p_{s,u}]_{(3)|^{s}|(2)|^2} : 1 \leqslant u \leqslant 9\} \rangle,$$
where $p_{s,u},\, u =1\, 2,\, 3$ are defined by \eqref{ct31} and
\begin{align*}	
&\mbox{$p_{s,4} = \sum_{j \in \{31,\, 32,\, 33,\, 34,\, 35,\, 36,\, 37,\, 38,\, 39,\, 40,\, 41,\, 42,\, 61,\, 62,\, 63,\, 64\}}a_{s,j}$},\\
&\mbox{$p_{s,5} = \sum_{j \in \{43,\, 44,\, 45,\, 46,\, 47,\, 48,\, 49,\, 50,\, 51,\, 52,\, 53,\, 54,\, 61,\, 62,\, 63,\, 64\}}a_{s,j}$},\\
&\mbox{$p_{s,6} = \sum_{j \in \{61,\, 62,\, 63,\, 64,\, 68,\, 69,\, 70,\, 71,\, 82,\, 83,\, 84,\, 85\}}a_{s,j}$},\\
&p_{s,7} = \begin{cases}\sum_{j \in \{55,\, 56,\, 57,\, 58,\, 59,\, 60,\, 74,\, 77,\, 78,\, 79,\, 80,\, 86,\, 87,\, 88\}}a_{3,j},&s=3,\\
\sum_{j \in \{55,\, 56,\, 57,\, 58,\, 59,\, 60,\, 74,\, 77,\, 78,\, 79,\, 80,\, 81,\, 88\}}a_{s,j}, &s \geqslant 4,
\end{cases}\\
&p_{s,8} = \begin{cases}\sum_{j \in \{65,\, 81,\, 86,\, 87,\, 90\}}a_{3,j}, &s=3,\\ \sum_{j \in \{65,\, 66,\, 67,\, 86,\, 90\}}a_{s,j}, &s \geqslant 4,
\end{cases}\\
&p_{s,9} = \begin{cases}\sum_{j \in \{66,\, 67,\, 90\}}a_{3,j}, &s=3,\\ \sum_{j \in \{81,\, 87,\, 90\}}a_{s,j}, &s \geqslant 4.
\end{cases}
\end{align*}
\end{lems}
\begin{proof} For $s \geqslant 3$, we have $QP_4((3)|^s|(2)|^2) = \langle \{[a_{2,j}]_{(3)|^s|(2)|^2}: 1 \leqslant j \leqslant 90\}\rangle$. By a direct computation, we see that $[p_{s,u}]_{(3)|^s|(2)|^2}$ is an $\Sigma_4$ invariant for $1 \leqslant u \leqslant 9$, and there is a direct summand decomposition of $\Sigma_4$-modules:
\begin{align*}QP_4((3)|^s|(2)|^2) &= [\Sigma_4(a_{s,1})]_{(3)|^{s}|(2)|^2}\bigoplus [\Sigma_4(a_{s,13})]_{(3)|^{s}|(2)|^2}\\ &\qquad \bigoplus [\Sigma_4(a_{s,25})]_{(3)|^{s}|(2)|^2}\bigoplus \mathcal U_{3,s} \oplus \mathcal V_{3,s},
\end{align*}
where 
\begin{align*}
\mathcal U_{3,s} &= \langle [a_{2,j}]_{(3)|^{2}|(2)|^2}: j \in \mathbb J_{3,s} = \{31,\ldots,54,\, 61, \ldots,\, 68,\, 73,\ldots,76\}\rangle,\\
\mathcal V_{3,s} &= \langle [a_{2,j}]_{(3)|^{s}|(2)|^2}: j \in \mathbb K_{3,s} = \{55,\ldots 60,\, 69,\, 70,\, 71,\, 72,\, 77, \ldots,\,90\}\rangle.
\end{align*}
By a simple computation, we easily obtain
\begin{align*}
&[\Sigma_4(a_{s,1})]_{(3)|^{s}|(2)|^2}^{\Sigma_4} = \langle[p_{s,1}]_{(3)|^{s}|(2)|^2}\rangle,\\ &[\Sigma_4(a_{s,13})]_{(3)|^{s}|(2)|^2}^{\Sigma_4} = \langle[p_{s,2}]_{(3)|^{s}|(2)|^2}\rangle,\\
&[\Sigma_4(a_{s,25})]_{(3)|^{s}|(2)|^2}^{\Sigma_4} = \langle[p_{s,3}]_{(3)|^{s}|(2)|^2}\rangle.
\end{align*}
We prove $\mathcal U_{3,s}^{\Sigma_4} = \langle [p_{s,u}]_{(3)|^{s}|(2)|^2} : 4 \leqslant u \leqslant 6\}\rangle$, $\mathcal V_{3,s}^{\Sigma_4} = \langle [p_{s,u}]_{(3)|^{s}|(2)|^2} : 7 \leqslant u \leqslant 9\}\rangle$. The computation of $U_{3,s}^{\Sigma_4}$ is similar to the one of $U_{3,2}^{\Sigma_4}$ with some minor changes. So, we only compute $\mathcal V_{3,s}^{\Sigma_4}$.

For $s = 3$, we observe that the leading monomials of $p_{3,7}$, $p_{3,8}$, $p_{3,9}$ respectively are $a_{3,60}$, $a_{3,81}$, $a_{3,90}$. So, if $[f]_{(3)|^{3}|(2)|^2} \in \mathcal V_{3,3}^{\Sigma_4}$ with $f \in P_4$, then there are $\gamma_7,\, \gamma_8,\, \gamma_9 \in \mathbb F_2$ such that
$$g := f + \gamma_7p_{3,7}+ \gamma_7p_{3,8} + \gamma_9p_{3,9}  \equiv _{(3)|^{3}|(2)|^2} \sum_{j\in \mathbb K_{3,3},\, j \ne 60,81,90}a_{3,j}.$$
Then, $[g]_{(3)|^{3}|(2)|^2}$ is also an $\Sigma_4$-invariant. A direct computation shows
\begin{align*}
\rho_1(g)& + g \equiv_{(3)|^{3}|(2)|^2} \gamma_{\{55,58\}}a_{3,55} + \gamma_{\{56,59\}}a_{3,56} + \gamma_{\{55,58\}}a_{3,58} + \gamma_{\{56,59\}}a_{3,59}\\
& + \gamma_{\{57,65,66,67,72,73,86\}}a_{3,65} + \gamma_{\{57,78\}}a_{3,72} + \gamma_{\{57,74\}}a_{3,73} + \gamma_{\{75,76\}}a_{3,75}\\
& + \gamma_{\{75,76\}}a_{3,76} + \gamma_{\{77,79\}}a_{3,77} + \gamma_{\{77,79\}}a_{3,79} + \gamma_{\{57,65,66,67,72,73,74,78,86\}}a_{3,86}\\
& + \gamma_{\{74,80\}}a_{3,88} + \gamma_{\{78,80\}}a_{3,89} + \gamma_{\{57,66,67,72,73,80,81,87,88,89\}}a_{3,90} \equiv_{(3)|^{3}|(2)|^2} 0,\\
\rho_2(g) &+ g \equiv_{(3)|^{3}|(2)|^2} \gamma_{\{56,57\}}a_{3,56} + \gamma_{\{56,57\}}a_{3,57} + \gamma_{\{58,60\}}a_{3,58} + \gamma_{\{58,60\}}a_{3,60}\\
& + \gamma_{\{65,66,90\}}a_{3,65} + \gamma_{\{65,66,90\}}a_{3,66} + \gamma_{\{72,75\}}a_{3,72} + \gamma_{\{59,73,74\}}a_{3,73}\\
& + \gamma_{\{59,73,74\}}a_{3,74} + \gamma_{\{72,75\}}a_{3,75} + \gamma_{\{59,79\}}a_{3,76} + \gamma_{\{77,78\}}a_{3,77}\\
& + \gamma_{\{77,78\}}a_{3,78} + \gamma_{\{79,80,89\}}a_{3,80} + \gamma_{\{76,86,87\}}a_{3,86} + \gamma_{\{59,76,79,86,87\}}a_{3,87}\\
& + \gamma_{\{59,79\}}a_{3,88} + \gamma_{\{79,80,89\}}a_{3,89} \equiv_{(3)|^{3}|(2)|^2} 0,\\  
\rho_3(g) &+ g \equiv_{(3)|^{3}|(2)|^2} \gamma_{\{55,56\}}a_{3,55} + \gamma_{\{55,56\}}a_{3,56} + \gamma_{\{58,59\}}a_{3,58} + \gamma_{\{58,59\}}a_{3,59}\\
& + \gamma_{\{66,67\}}a_{3,66} + \gamma_{\{66,67\}}a_{3,67} + \gamma_{\{72,73\}}a_{3,72} + \gamma_{\{72,73\}}a_{3,73}\\
& + \gamma_{\{74,78\}}a_{3,74} + \gamma_{\{60,75,77\}}a_{3,75} + \gamma_{\{60,76,79\}}a_{3,76} + \gamma_{\{60,75,77\}}a_{3,77}\\
& + \gamma_{\{74,78\}}a_{3,78} + \gamma_{\{60,76,79\}}a_{3,79} + \gamma_{\{60,81,87\}}a_{3,81} + \gamma_{\{60,81,87\}}a_{3,87}\\
& + \gamma_{\{60,88,89\}}a_{3,88} + \gamma_{\{60,88,89\}}a_{3,89} \equiv_{(3)|^{3}|(2)|^2} 0.
\end{align*}
From these equalities we get $\gamma_{j} = 0$ for all $j \in \mathbb K_{3,3}\setminus\{60,81,90\}$. Hence, we have $g \equiv_{(3)|^{3}|(2)|^2} 0$ and
$f \equiv_{(3)|^{3}|(2)|^2}  \gamma_7p_{3,7}+ \gamma_8p_{3,8} + \gamma_9p_{3,9}.$
	
Now we prove that $\mathcal V_{3,s}^{\Sigma_4} = \langle [p_{s,u}]_{(3)|^{s}|(2)|^2} : 7 \leqslant u \leqslant 9\}\rangle$ for $s \geqslant 4$. It is easy to see that the leading monomials of $p_{s,7}$, $p_{s,8}$, $p_{s,9}$ respectively are $a_{s,60}$, $a_{s,86}$, $a_{s,81}$. Hence, if $[h]_{(3)|^{s}|(2)|^2} \in \mathcal V_{3,s}^{\Sigma_4}$ with $h \in P_4$, then there are $\gamma_7,\, \gamma_8,\, \gamma_9 \in \mathbb F_2$ such that
$$\tilde h := h + \gamma_7p_{s,7}+ \gamma_8p_{s,8} + \gamma_9p_{s,9}  \equiv _{(3)|^{s}|(2)|^2} \sum_{j\in \mathbb K_{3,s},\, j \ne 60,81,86}a_{2,j}.$$
Since $[h]_{(3)|^{s}|(2)|^2}$ is an $\Sigma_4$-invariant, so is $[\tilde h]_{(3)|^{s}|(2)|^2}$. By a direct computation we have
\begin{align*}
\rho_1(\tilde h)& + \tilde h \equiv_{(3)|^{s}|(2)|^2} \gamma_{\{55,58\}}a_{s,55} + \gamma_{\{56,59\}}a_{s,56} + \gamma_{\{55,58\}}a_{s,58} + \gamma_{\{56,59\}}a_{s,59}\\
& + \gamma_{\{65,86\}}a_{s,65} + \gamma_{\{57,78\}}a_{s,72} + \gamma_{\{57,74\}}a_{s,73} + \gamma_{\{75,76\}}a_{s,75} + \gamma_{\{75,76\}}a_{s,76}\\
& + \gamma_{\{77,79\}}a_{s,77} + \gamma_{\{77,79\}}a_{s,79} + \gamma_{\{65,86\}}a_{s,86} + \gamma_{\{74,80\}}a_{s,88} + \gamma_{\{78,80\}}a_{s,89}\\
& + \gamma_{\{57,66,67,72,73,74,78,80,81,87,88,89\}}a_{s,90} \equiv_{(3)|^{s}|(2)|^2} 0,\\
\rho_2(\tilde h) &+ \tilde h \equiv_{(3)|^{s}|(2)|^2} \gamma_{\{56,57\}}a_{s,56} + \gamma_{\{56,57\}}a_{s,57} + \gamma_{\{58,60\}}a_{s,58} + \gamma_{\{58,60\}}a_{s,60}\\
& + \gamma_{\{65,66\}}a_{s,65} + \gamma_{\{65,66\}}a_{s,66} + \gamma_{\{72,75\}}a_{s,72} + \gamma_{\{59,73,74\}}a_{s,73}\\
& + \gamma_{\{59,73,74\}}a_{s,74} + \gamma_{\{72,75\}}a_{s,75} + \gamma_{\{59,79\}}a_{s,76} + \gamma_{\{77,78\}}a_{s,77}\\
& + \gamma_{\{77,78\}}a_{s,78} + \gamma_{\{79,80,89\}}a_{s,80} + \gamma_{\{59,76,79,86,87,90\}}a_{s,87} + \gamma_{\{59,79\}}a_{s,88}\\
& + \gamma_{\{79,80,89\}}a_{s,89} + \gamma_{\{76,86,87,90\}}a_{s,90} \equiv_{(3)|^{s}|(2)|^2} 0,\\  
\rho_3(\tilde h) &+ \tilde h \equiv_{(3)|^{s}|(2)|^2}  \gamma_{\{55,56\}}a_{s,55} + \gamma_{\{55,56\}}a_{s,56} + \gamma_{\{58,59\}}a_{s,58} + \gamma_{\{58,59\}}a_{s,59}\\
& + \gamma_{\{66,67\}}a_{s,66} + \gamma_{\{66,67\}}a_{s,67} + \gamma_{\{72,73\}}a_{s,72} + \gamma_{\{72,73\}}a_{s,73} + \gamma_{\{74,78\}}a_{s,74}\\
& + \gamma_{\{60,75,77\}}a_{s,75} + \gamma_{\{60,76,79\}}a_{s,76} + \gamma_{\{60,75,77\}}a_{s,77} + \gamma_{\{74,78\}}a_{s,78}\\
& + \gamma_{\{60,76,79\}}a_{s,79} + \gamma_{\{60,81,87\}}a_{s,81} + \gamma_{\{60,81,87\}}a_{s,87} + \gamma_{\{60,88,89\}}a_{s,88}\\
& + \gamma_{\{60,88,89\}}a_{s,89} \equiv_{(3)|^{s}|(2)|^2} 0.
\end{align*}
These equalities imply $\gamma_{j} = 0$ for all $j \in \mathbb K_{3,s}\setminus\{60,81,86\}$ and $\tilde h \equiv_{(3)|^{s}|(2)|^2} 0$. Hence, we obtain
$h \equiv_{(3)|^{s}|(2)|^2}  \gamma_7p_{s,7}+ \gamma_8p_{s,8} + \gamma_9p_{s,9}.$
	The lemma is completely proved.
\end{proof}
\begin{proof}[Proof of Proposition \ref{mdt31}]
Suppose $[f]_{(3)|^s|(2)|^2} \in QP_4((3)|^s|(2)|^2)^{GL_4}$ with $f \in P_4((3)|^s|(2)|^2)$. Then, $[f]_{(3)|^s|(2)|^2} \in QP_4((3)|^s|(2)|^2)^{\Sigma_4}$. 

By Lemmas \ref{bdt32}, we have 
$$ f \equiv_{(3)|^s|(2)|^2} \sum_{1\leqslant u \leqslant 9} \gamma_up_{s,u},$$
where $\gamma_u \in \mathbb F_2$. 

For $s = 2$, by computing $\rho_4(f)+f$ in terms of the admissible monomials, we get
\begin{align*}
\rho_4(f)&+f \equiv_{(3)|^{2}|(2)|^2} \gamma_{\{1,3\}}a_{2,1} + \gamma_{\{1,4,5,6\}}a_{2,2} + \gamma_{\{1,4,5,6\}}a_{2,3} + \gamma_{\{1,2\}}a_{2,5}\\ 
& + \gamma_{\{1,2\}}a_{2,6} + \gamma_{\{2,6\}}a_{2,13} + \gamma_{\{2,6\}}a_{2,14} + \gamma_{2}a_{2,15} + \gamma_{2}a_{2,20} + \gamma_{2}a_{2,21}\\ 
& + \gamma_{\{3,5\}}a_{2,26} + \gamma_{\{3,5\}}a_{2,27} + \gamma_{\{5,6\}}a_{2,31} + \gamma_{\{5,6\}}a_{2,32} + \gamma_{\{4,6,7,8\}}a_{2,35}\\ 
& + \gamma_{\{4,5,6,7,8\}}a_{2,36} + \gamma_{\{5,6,8\}}a_{2,39} + \gamma_{6}a_{2,45} + \gamma_{\{4,6\}}a_{2,47} + \gamma_{\{4,6\}}a_{2,48}\\ 
& + \gamma_{6}a_{2,49} + \gamma_{6}a_{2,52} + \gamma_{\{2,5,8\}}a_{2,55} + \gamma_{\{2,5\}}a_{2,56} + \gamma_{\{2,5\}}a_{2,58}\\ 
& + \gamma_{\{2,5\}}a_{2,59} + \gamma_{\{5,8\}}a_{2,63} + \gamma_{\{2,5,8\}}a_{2,65} + \gamma_{4}a_{2,70} + \gamma_{\{2,5\}}a_{2,72}\\ 
& + \gamma_{\{2,5\}}a_{2,73} + \gamma_{\{1,2,3,6\}}a_{2,82} + \gamma_{\{1,2,3,6\}}a_{2,83} + \gamma_{5}a_{2,84} + \gamma_{8}a_{2,86}\\ 
& + \gamma_{\{5,6\}}a_{2,89} + \gamma_{\{5,6\}}a_{2,90} + \gamma_{\{4,5\}}a_{2,91}   \equiv_{(3)|^{2}|(2)|^2} 0.
\end{align*}
This equality implies $\gamma_u = 0$ for $1 \leqslant u \leqslant 8$. Hence, 
$$f \equiv_{(3)|^{2}|(2)|^2} = \gamma_9p_{2,9} = \gamma_9\bar\xi_{3,2}.$$ 
The proposition is proved for $s = 2$.

For $s = 3$, by Lemma \ref{bdt33}, we have 
$$ f \equiv_{(3)|^3|(2)|^2} \sum_{1\leqslant u \leqslant 9} \gamma_up_{3,u},$$
where $\gamma_u \in \mathbb F_2$. By computing $\rho_4(f)+f$ in terms of the admissible monomials, we get
\begin{align*}
\rho_4(f)&+f \equiv_{(3)|^{3}|(2)|^2} \gamma_{\{1,3\}}a_{3,1} + \gamma_{\{1,4\}}a_{3,2} + \gamma_{\{1,4\}}a_{3,3} + \gamma_{\{1,2\}}a_{3,5}\\ 
& + \gamma_{\{1,2\}}a_{3,6} + \gamma_{\{2,5\}}a_{3,13} + \gamma_{\{2,5\}}a_{3,14} + \gamma_{\{2,7\}}a_{3,15} + \gamma_{2}a_{3,20} + \gamma_{2}a_{3,21}\\ 
& + \gamma_{\{3,4,5,6\}}a_{3,26} + \gamma_{\{3,4,5,6\}}a_{3,27} + \gamma_{\{1,2,3,4,5,6\}}a_{3,31} + \gamma_{\{1,2,3,4,5,6\}}a_{3,32}\\ 
& + \gamma_{\{4,9\}}a_{3,35} + \gamma_{\{4,9\}}a_{3,36} + \gamma_{\{1,2,3,5\}}a_{3,37} + \gamma_{\{1,2,3,5\}}a_{3,38} + \gamma_{\{4,7\}}a_{3,39}\\ 
& + \gamma_{\{5,7\}}a_{3,45} + \gamma_{\{5,6\}}a_{3,47} + \gamma_{\{5,6\}}a_{3,48} + \gamma_{\{5,7\}}a_{3,49} + \gamma_{5}a_{3,52}\\ 
& + \gamma_{\{2,5,7\}}a_{3,55} + \gamma_{\{2,5,7\}}a_{3,56} + \gamma_{\{2,5\}}a_{3,58} + \gamma_{\{2,5\}}a_{3,59} + \gamma_{\{4,5,6,8\}}a_{3,63}\\ 
& + \gamma_{\{7,8\}}a_{3,65} + \gamma_{\{6,7\}}a_{3,70} + \gamma_{\{2,5\}}a_{3,72} + \gamma_{\{2,5\}}a_{3,73} + \gamma_{\{5,6,7\}}a_{3,75}\\ 
& + \gamma_{\{5,6\}}a_{3,76} + \gamma_{\{5,6,7\}}a_{3,77} + \gamma_{\{5,6\}}a_{3,79} + \gamma_{\{1,2,3,5\}}a_{3,82} + \gamma_{\{1,2,3,5\}}a_{3,83}\\ 
& + \gamma_{\{6,7,8\}}a_{3,84} + \gamma_{\{5,6\}}a_{3,88} + \gamma_{\{5,6\}}a_{3,89} + \gamma_{\{2,4,5,6,7\}}a_{3,90} \equiv_{(3)|^{3}|(2)|^2} 0.
\end{align*}
The above equality implies $\gamma_u = 0$ for $1 \leqslant u \leqslant 9$. So, $[f]_{(3)|^{3}|(2)|^2} = 0$. The proposition is proved for $s = 3$.

For $s \geqslant 4$, by computing $\rho_4(f)+f$ in terms of the admissible monomials, we get
\begin{align*}
\rho_4(f)&+f \equiv_{(3)|^{s}|(2)|^2} \gamma_{\{1,3\}}a_{s,1} + \gamma_{\{1,4\}}a_{s,2} + \gamma_{\{1,4\}}a_{s,3} + \gamma_{\{1,2\}}a_{s,5}\\ 
& + \gamma_{\{1,2\}}a_{s,6} + \gamma_{\{2,5\}}a_{s,13} + \gamma_{\{2,5\}}a_{s,14} + \gamma_{\{2,7\}}a_{s,15} + \gamma_{2}a_{s,20} + \gamma_{2}a_{s,21}\\ 
& + \gamma_{\{3,4,5,6\}}a_{s,26} + \gamma_{\{3,4,5,6\}}a_{s,27} + \gamma_{\{1,2,3,4,5,6\}}a_{s,31} + \gamma_{\{1,2,3,4,5,6\}}a_{s,32}\\ 
& + \gamma_{\{4,8\}}a_{s,35} + \gamma_{\{4,8\}}a_{s,36} + \gamma_{\{1,2,3,5\}}a_{s,37} + \gamma_{\{1,2,3,5\}}a_{s,38} + \gamma_{\{4,7\}}a_{s,39}\\ 
& + \gamma_{\{5,7\}}a_{s,45} + \gamma_{\{5,6\}}a_{s,47} + \gamma_{\{5,6\}}a_{s,48} + \gamma_{\{5,7\}}a_{s,49} + \gamma_{5}a_{s,52}\\ 
& + \gamma_{\{2,5,7\}}a_{s,55} + \gamma_{\{2,5,7\}}a_{s,56} + \gamma_{\{2,5\}}a_{s,58} + \gamma_{\{2,5\}}a_{s,59}\\ 
& + \gamma_{\{4,5,6,7,9\}}a_{s,63} + \gamma_{\{2,4,5,6,7,9\}}a_{s,65} + \gamma_{\{6,7\}}a_{s,70} + \gamma_{\{2,5\}}a_{s,72}\\ 
& + \gamma_{\{2,5\}}a_{s,73} + \gamma_{\{5,6,7\}}a_{s,75} + \gamma_{\{5,6\}}a_{s,76} + \gamma_{\{5,6,7\}}a_{s,77} + \gamma_{\{5,6\}}a_{s,79}\\ 
& + \gamma_{\{1,2,3,5\}}a_{s,82} + \gamma_{\{1,2,3,5\}}a_{s,83} + \gamma_{\{6,9\}}a_{s,84} + \gamma_{\{2,4,5,6,7\}}a_{s,86}\\ 
& + \gamma_{\{5,6\}}a_{s,88} + \gamma_{\{5,6\}}a_{s,89} + \gamma_{\{2,4,5,6,7\}}a_{s,90}  \equiv_{(3)|^{s}|(2)|^2} 0.
\end{align*}
From this equality it implies $\gamma_u = 0$ for $1 \leqslant u \leqslant 9$. Therefore, $f \equiv_{(3)|^s|(2)|^2} 0$. The proposition is completely proved.
\end{proof}

\subsubsection{\textbf{The case $t \geqslant 4$}}\

\medskip
Following \cite{su5,su50}, a basis of $QP_4((3)|^{s}|(2)|^{t-1})$ is the set of all classes represented by the admissible monomials $a_{s,j} = a_{t,s,j}$ which are determined as follows:

\medskip
For $s \geqslant 2$,

\medskip
\centerline{\begin{tabular}{lll} 
$a_{s,1} = x_2^{2^{s}-1}x_3^{2^{s+t-1}-1}x_4^{2^{s+t-1}-1}$\cr  $a_{s,2} = x_2^{2^{s+t-1}-1}x_3^{2^{s}-1}x_4^{2^{s+t-1}-1}$\cr  $a_{s,3} = x_2^{2^{s+t-1}-1}x_3^{2^{s+t-1}-1}x_4^{2^{s}-1}$\cr  $a_{s,4} = x_1^{2^{s}-1}x_3^{2^{s+t-1}-1}x_4^{2^{s+t-1}-1}$\cr  $a_{s,5} = x_1^{2^{s}-1}x_2^{2^{s+t-1}-1}x_4^{2^{s+t-1}-1}$\cr  $a_{s,6} = x_1^{2^{s}-1}x_2^{2^{s+t-1}-1}x_3^{2^{s+t-1}-1}$\cr  $a_{s,7} = x_1^{2^{s+t-1}-1}x_3^{2^{s}-1}x_4^{2^{s+t-1}-1}$\cr  $a_{s,8} = x_1^{2^{s+t-1}-1}x_3^{2^{s+t-1}-1}x_4^{2^{s}-1}$\cr  $a_{s,9} = x_1^{2^{s+t-1}-1}x_2^{2^{s}-1}x_4^{2^{s+t-1}-1}$\cr  $a_{s,10} = x_1^{2^{s+t-1}-1}x_2^{2^{s}-1}x_3^{2^{s+t-1}-1}$\cr  $a_{s,11} = x_1^{2^{s+t-1}-1}x_2^{2^{s+t-1}-1}x_4^{2^{s}-1}$\cr  $a_{s,12} = x_1^{2^{s+t-1}-1}x_2^{2^{s+t-1}-1}x_3^{2^{s}-1}$\cr  $a_{s,13} = x_2^{2^{s+1}-1}x_3^{2^{s+t-1}-2^{s}-1}x_4^{2^{s+t-1}-1}$\cr  $a_{s,14} = x_2^{2^{s+1}-1}x_3^{2^{s+t-1}-1}x_4^{2^{s+t-1}-2^{s}-1}$\cr  $a_{s,15} = x_2^{2^{s+t-1}-1}x_3^{2^{s+1}-1}x_4^{2^{s+t-1}-2^{s}-1}$\cr  $a_{s,16} = x_1^{2^{s+1}-1}x_3^{2^{s+t-1}-2^{s}-1}x_4^{2^{s+t-1}-1}$\cr  $a_{s,17} = x_1^{2^{s+1}-1}x_3^{2^{s+t-1}-1}x_4^{2^{s+t-1}-2^{s}-1}$\cr  $a_{s,18} = x_1^{2^{s+1}-1}x_2^{2^{s+t-1}-2^{s}-1}x_4^{2^{s+t-1}-1}$\cr  $a_{s,19} = x_1^{2^{s+1}-1}x_2^{2^{s+t-1}-2^{s}-1}x_3^{2^{s+t-1}-1}$\cr  $a_{s,20} = x_1^{2^{s+1}-1}x_2^{2^{s+t-1}-1}x_4^{2^{s+t-1}-2^{s}-1}$\cr $a_{s,21} = x_1^{2^{s+1}-1}x_2^{2^{s+t-1}-1}x_3^{2^{s+t-1}-2^{s}-1}$\cr  $a_{s,22} = x_1^{2^{s+t-1}-1}x_3^{2^{s+1}-1}x_4^{2^{s+t-1}-2^{s}-1}$\cr  $a_{s,23} = x_1^{2^{s+t-1}-1}x_2^{2^{s+1}-1}x_4^{2^{s+t-1}-2^{s}-1}$\cr  $a_{s,24} = x_1^{2^{s+t-1}-1}x_2^{2^{s+1}-1}x_3^{2^{s+t-1}-2^{s}-1}$\cr  $a_{s,25} = x_2^{2^{s+2}-1}x_3^{2^{s+t-1}-2^{s+1}-1}x_4^{2^{s+t-1}-2^{s}-1}$\cr  $a_{s,26} = x_1^{2^{s+2}-1}x_3^{2^{s+t-1}-2^{s+1}-1}x_4^{2^{s+t-1}-2^{s}-1}$\cr $a_{s,27} = x_1^{2^{s+2}-1}x_2^{2^{s+t-1}-2^{s+1}-1}x_4^{2^{s+t-1}-2^{s}-1}$\cr  $a_{s,28} = x_1^{2^{s+2}-1}x_2^{2^{s+t-1}-2^{s+1}-1}x_3^{2^{s+t-1}-2^{s}-1}$\cr  $a_{s,29} = x_1x_2^{2^{s}-2}x_3^{2^{s+t-1}-1}x_4^{2^{s+t-1}-1}$\cr  $a_{s,30} = x_1x_2^{2^{s+t-1}-1}x_3^{2^{s}-2}x_4^{2^{s+t-1}-1}$\cr  
\end{tabular}}
\centerline{\begin{tabular}{lll} 
$a_{s,31} = x_1x_2^{2^{s+t-1}-1}x_3^{2^{s+t-1}-1}x_4^{2^{s}-2}$\cr  $a_{s,32} = x_1^{2^{s+t-1}-1}x_2x_3^{2^{s}-2}x_4^{2^{s+t-1}-1}$\cr  $a_{s,33} = x_1^{2^{s+t-1}-1}x_2x_3^{2^{s+t-1}-1}x_4^{2^{s}-2}$\cr  $a_{s,34} = x_1^{2^{s+t-1}-1}x_2^{2^{s+t-1}-1}x_3x_4^{2^{s}-2}$\cr  $a_{s,35} = x_1x_2^{2^{s}-1}x_3^{2^{s+t-1}-2}x_4^{2^{s+t-1}-1}$\cr  $a_{s,36} = x_1x_2^{2^{s}-1}x_3^{2^{s+t-1}-1}x_4^{2^{s+t-1}-2}$\cr $a_{s,37} = x_1x_2^{2^{s+t-1}-2}x_3^{2^{s}-1}x_4^{2^{s+t-1}-1}$\cr  $a_{s,38} = x_1x_2^{2^{s+t-1}-2}x_3^{2^{s+t-1}-1}x_4^{2^{s}-1}$\cr  $a_{s,39} = x_1x_2^{2^{s+t-1}-1}x_3^{2^{s}-1}x_4^{2^{s+t-1}-2}$\cr  $a_{s,40} = x_1x_2^{2^{s+t-1}-1}x_3^{2^{s+t-1}-2}x_4^{2^{s}-1}$\cr  $a_{s,41} = x_1^{2^{s}-1}x_2x_3^{2^{s+t-1}-2}x_4^{2^{s+t-1}-1}$\cr  $a_{s,42} = x_1^{2^{s}-1}x_2x_3^{2^{s+t-1}-1}x_4^{2^{s+t-1}-2}$\cr $a_{s,43} = x_1^{2^{s}-1}x_2^{2^{s+t-1}-1}x_3x_4^{2^{s+t-1}-2}$\cr  $a_{s,44} = x_1^{2^{s+t-1}-1}x_2x_3^{2^{s}-1}x_4^{2^{s+t-1}-2}$\cr  $a_{s,45} = x_1^{2^{s+t-1}-1}x_2x_3^{2^{s+t-1}-2}x_4^{2^{s}-1}$\cr  $a_{s,46} = x_1^{2^{s+t-1}-1}x_2^{2^{s}-1}x_3x_4^{2^{s+t-1}-2}$\cr  $a_{s,47} = x_1x_2^{2^{s+1}-2}x_3^{2^{s+t-1}-2^{s}-1}x_4^{2^{s+t-1}-1}$\cr  $a_{s,48} = x_1x_2^{2^{s+1}-2}x_3^{2^{s+t-1}-1}x_4^{2^{s+t-1}-2^{s}-1}$\cr  $a_{s,49} = x_1x_2^{2^{s+t-1}-1}x_3^{2^{s+1}-2}x_4^{2^{s+t-1}-2^{s}-1}$\cr  $a_{s,50} = x_1^{2^{s+t-1}-1}x_2x_3^{2^{s+1}-2}x_4^{2^{s+t-1}-2^{s}-1}$\cr  $a_{s,51} = x_1x_2^{2^{s+1}-1}x_3^{2^{s+t-1}-2^{s}-2}x_4^{2^{s+t-1}-1}$\cr  $a_{s,52} = x_1x_2^{2^{s+1}-1}x_3^{2^{s+t-1}-1}x_4^{2^{s+t-1}-2^{s}-2}$\cr  $a_{s,53} = x_1x_2^{2^{s+t-1}-1}x_3^{2^{s+1}-1}x_4^{2^{s+t-1}-2^{s}-2}$\cr  $a_{s,54} = x_1^{2^{s+1}-1}x_2x_3^{2^{s+t-1}-2^{s}-2}x_4^{2^{s+t-1}-1}$\cr  $a_{s,55} = x_1^{2^{s+1}-1}x_2x_3^{2^{s+t-1}-1}x_4^{2^{s+t-1}-2^{s}-2}$\cr  $a_{s,56} = x_1^{2^{s+1}-1}x_2^{2^{s+t-1}-1}x_3x_4^{2^{s+t-1}-2^{s}-2}$\cr  $a_{s,57} = x_1^{2^{s+t-1}-1}x_2x_3^{2^{s+1}-1}x_4^{2^{s+t-1}-2^{s}-2}$\cr  $a_{s,58} = x_1^{2^{s+t-1}-1}x_2^{2^{s+1}-1}x_3x_4^{2^{s+t-1}-2^{s}-2}$\cr  $a_{s,59} = x_1^{3}x_2^{2^{s+t-1}-3}x_3^{2^{s}-2}x_4^{2^{s+t-1}-1}$\cr  $a_{s,60} = x_1^{3}x_2^{2^{s+t-1}-3}x_3^{2^{s+t-1}-1}x_4^{2^{s}-2}$\cr  $a_{s,61} = x_1^{3}x_2^{2^{s+t-1}-1}x_3^{2^{s+t-1}-3}x_4^{2^{s}-2}$\cr  $a_{s,62} = x_1^{2^{s+t-1}-1}x_2^{3}x_3^{2^{s+t-1}-3}x_4^{2^{s}-2}$\cr  $a_{s,63} = x_1^{3}x_2^{2^{s+1}-3}x_3^{2^{s+t-1}-2^{s}-2}x_4^{2^{s+t-1}-1}$\cr  $a_{s,64} = x_1^{3}x_2^{2^{s+1}-3}x_3^{2^{s+t-1}-1}x_4^{2^{s+t-1}-2^{s}-2}$\cr  $a_{s,65} = x_1^{3}x_2^{2^{s+t-1}-1}x_3^{2^{s+1}-3}x_4^{2^{s+t-1}-2^{s}-2}$\cr  $a_{s,66} = x_1^{2^{s+t-1}-1}x_2^{3}x_3^{2^{s+1}-3}x_4^{2^{s+t-1}-2^{s}-2}$\cr  $a_{s,67} = x_1x_2^{2^{s+1}-1}x_3^{2^{s+t-1}-2^{s}-1}x_4^{2^{s+t-1}-2}$\cr  $a_{s,68} = x_1x_2^{2^{s+1}-1}x_3^{2^{s+t-1}-2}x_4^{2^{s+t-1}-2^{s}-1}$\cr  $a_{s,69} = x_1x_2^{2^{s+t-1}-2}x_3^{2^{s+1}-1}x_4^{2^{s+t-1}-2^{s}-1}$\cr   $a_{s,70} = x_1^{2^{s+1}-1}x_2x_3^{2^{s+t-1}-2^{s}-1}x_4^{2^{s+t-1}-2}$\cr   $a_{s,71} = x_1^{2^{s+1}-1}x_2x_3^{2^{s+t-1}-2}x_4^{2^{s+t-1}-2^{s}-1}$\cr  
\end{tabular}}
\centerline{\begin{tabular}{lll} 
$a_{s,72} = x_1^{2^{s+1}-1}x_2^{2^{s+t-1}-2^{s}-1}x_3x_4^{2^{s+t-1}-2}$\cr  $a_{s,73} = x_1x_2^{2^{s+2}-2}x_3^{2^{s+t-1}-2^{s+1}-1}x_4^{2^{s+t-1}-2^{s}-1}$\cr  $a_{s,74} = x_1x_2^{2^{s+2}-1}x_3^{2^{s+t-1}-2^{s+1}-2}x_4^{2^{s+t-1}-2^{s}-1}$\cr  $a_{s,75} = x_1^{2^{s+2}-1}x_2x_3^{2^{s+t-1}-2^{s+1}-2}x_4^{2^{s+t-1}-2^{s}-1}$\cr $a_{s,76} = x_1x_2^{2^{s+2}-1}x_3^{2^{s+t-1}-2^{s+1}-1}x_4^{2^{s+t-1}-2^{s}-2}$\cr  $a_{s,77} = x_1^{2^{s+2}-1}x_2x_3^{2^{s+t-1}-2^{s+1}-1}x_4^{2^{s+t-1}-2^{s}-2}$\cr  $a_{s,78} = x_1^{2^{s+2}-1}x_2^{2^{s+t-1}-2^{s+1}-1}x_3x_4^{2^{s+t-1}-2^{s}-2}$\cr  $a_{s,79} = x_1^{3}x_2^{2^{s}-1}x_3^{2^{s+t-1}-3}x_4^{2^{s+t-1}-2}$\cr  $a_{s,80} = x_1^{3}x_2^{2^{s+t-1}-3}x_3^{2^{s}-1}x_4^{2^{s+t-1}-2}$\cr  $a_{s,81} = x_1^{3}x_2^{2^{s+t-1}-3}x_3^{2^{s+t-1}-2}x_4^{2^{s}-1}$\cr  $a_{s,82} = x_1^{3}x_2^{2^{s+1}-3}x_3^{2^{s+t-1}-2^{s}-1}x_4^{2^{s+t-1}-2}$\cr  $a_{s,83} = x_1^{3}x_2^{2^{s+1}-3}x_3^{2^{s+t-1}-2}x_4^{2^{s+t-1}-2^{s}-1}$\cr  $a_{s,84} = x_1^{3}x_2^{2^{s+t-1}-3}x_3^{2^{s+1}-2}x_4^{2^{s+t-1}-2^{s}-1}$\cr  $a_{s,85} = x_1^{3}x_2^{2^{s+1}-1}x_3^{2^{s+t-1}-2^{s}-3}x_4^{2^{s+t-1}-2}$\cr  $a_{s,86} = x_1^{2^{s+1}-1}x_2^{3}x_3^{2^{s+t-1}-2^{s}-3}x_4^{2^{s+t-1}-2}$\cr  $a_{s,87} = x_1^{3}x_2^{2^{s+1}-1}x_3^{2^{s+t-1}-3}x_4^{2^{s+t-1}-2^{s}-2}$\cr  $a_{s,88} = x_1^{3}x_2^{2^{s+t-1}-3}x_3^{2^{s+1}-1}x_4^{2^{s+t-1}-2^{s}-2}$\cr  $a_{s,89} = x_1^{2^{s+1}-1}x_2^{3}x_3^{2^{s+t-1}-3}x_4^{2^{s+t-1}-2^{s}-2}$\cr  $a_{s,90} = x_1^{3}x_2^{2^{s+2}-3}x_3^{2^{s+t-1}-2^{s+1}-2}x_4^{2^{s+t-1}-2^{s}-1}$\cr  $a_{s,91} = x_1^{3}x_2^{2^{s+2}-3}x_3^{2^{s+t-1}-2^{s+1}-1}x_4^{2^{s+t-1}-2^{s}-2}$\cr  $a_{s,92} = x_1^{3}x_2^{2^{s+2}-1}x_3^{2^{s+t-1}-2^{s+1}-3}x_4^{2^{s+t-1}-2^{s}-2}$\cr  $a_{s,93} = x_1^{2^{s+2}-1}x_2^{3}x_3^{2^{s+t-1}-2^{s+1}-3}x_4^{2^{s+t-1}-2^{s}-2}$\cr  $a_{s,94} = x_1^{7}x_2^{2^{s+t-1}-5}x_3^{2^{s+t-1}-3}x_4^{2^{s}-2}$\cr  $a_{s,95} = x_1^{7}x_2^{2^{s+t-1}-5}x_3^{2^{s+1}-3}x_4^{2^{s+t-1}-2^{s}-2}$\cr  $a_{s,96} = x_1^{7}x_2^{2^{s+2}-5}x_3^{2^{s+t-1}-2^{s+1}-3}x_4^{2^{s+t-1}-2^{s}-2}$\cr
\end{tabular}}

\medskip
For $s = 2$,

\medskip
\centerline{\begin{tabular}{lll}
$a_{2,97} = x_1^{3}x_2^{3}x_3^{2^{t+1}-4}x_4^{2^{t+1}-1}$ & &$a_{2,98} = x_1^{3}x_2^{3}x_3^{2^{t+1}-1}x_4^{2^{t+1}-4}$\cr  $a_{2,99} = x_1^{3}x_2^{2^{t+1}-1}x_3^{3}x_4^{2^{t+1}-4}$ & &$a_{2,100} = x_1^{2^{t+1}-1}x_2^{3}x_3^{3}x_4^{2^{t+1}-4}$\cr  $a_{2,101} = x_1^{3}x_2^{7}x_3^{2^{t+1}-5}x_4^{2^{t+1}-4}$ & &$a_{2,102} = x_1^{7}x_2^{3}x_3^{2^{t+1}-5}x_4^{2^{t+1}-4}$\cr  $a_{2,103} = x_1^{7}x_2^{2^{t+1}-5}x_3^{3}x_4^{2^{t+1}-4}$ & &$a_{2,104} = x_1^{7}x_2^{7}x_3^{2^{t+1}-8}x_4^{2^{t+1}-5}$\cr  $a_{2,105} = x_1^{7}x_2^{7}x_3^{2^{t+1}-7}x_4^{2^{t+1}-6}$ & &\cr     	 
\end{tabular}}

\medskip
For $s \geqslant 3$,

\medskip
\centerline{\begin{tabular}{lll}   	 
$a_{s,97} = x_1^{3}x_2^{2^{s}-3}x_3^{2^{s+t-1}-2}x_4^{2^{s+t-1}-1}$\cr  $a_{s,98} = x_1^{3}x_2^{2^{s}-3}x_3^{2^{s+t-1}-1}x_4^{2^{s+t-1}-2}$\cr  $a_{s,99} = x_1^{3}x_2^{2^{s+t-1}-1}x_3^{2^{s}-3}x_4^{2^{s+t-1}-2}$\cr  $a_{s,100} = x_1^{2^{s+t-1}-1}x_2^{3}x_3^{2^{s}-3}x_4^{2^{s+t-1}-2}$\cr  $a_{s,101} = x_1^{2^{s}-1}x_2^{3}x_3^{2^{s+t-1}-3}x_4^{2^{s+t-1}-2}$\cr  $a_{s,102} = x_1^{7}x_2^{2^{s+t-1}-5}x_3^{2^{s}-3}x_4^{2^{s+t-1}-2}$\cr   $a_{s,103} = x_1^{7}x_2^{2^{s+1}-5}x_3^{2^{s+t-1}-2^{s}-3}x_4^{2^{s+t-1}-2}$\cr  $a_{s,104} = x_1^{7}x_2^{2^{s+1}-5}x_3^{2^{s+t-1}-3}x_4^{2^{s+t-1}-2^{s}-2}$\cr
\end{tabular}}

\medskip
For $s = 3$, $a_{3,105} =  x_1^{7}x_2^{7}x_3^{2^{t+2}-7}x_4^{2^{t+2}-2}$.

\medskip
For $s \geqslant 4$, $a_{s,105} = x_1^{7}x_2^{2^{s}-5}x_3^{2^{s+t-1}-3}x_4^{2^{s+t-1}-2}$.

\medskip
By using this basis, we prove the following.

\begin{props}\label{mdt41} For any $s \geqslant 2$ and $t \geqslant 4$, we have 
$$QP_4((3)|^{s}|(2)|^{t-1})^{GL_4} = \langle [\bar\xi_{t,s}]_{(3)|^{s}|(2)|^{t-1}} \rangle,$$
where 
$$
\bar\xi_{t,s} = \begin{cases}
\sum_{1 \leqslant j \leqslant 104,\, j \notin\{35,\, 36,\, 39,\, 41,\, 42,\, 43,\, 44,\, 46,\, 68,\, 71,\, 72\}\atop \hskip1.7cm\cup\{79,\, 80,\, 83,\, 85,\, 86,\, 87,\, 89,\, 101,\, 102 \}}a_{2,j}, &\mbox{if } s = 2,\\
\sum_{1 \leqslant j \leqslant 105,\, j \ne 79,\, 101}a_{3,j}, &\mbox{if } s = 3,\\ 
\sum_{1 \leqslant j \leqslant 105}a_{s,j}, &\mbox{if } s \geqslant 4.
\end{cases}
$$
\end{props}
\begin{rems}
The dimensional result of the proposition is stated in \cite[Prop. 4.1.13]{pp25} but there is no detailed proof. The dimensional result can easily predict when assume that Singer's conjecture is true for $k = 4$.	 
\end{rems}

We note that for a fixed value of $t$, we denote $a_{s,j} = a_{t,s,j}$. By a direct computation, we have
\begin{align*}
&[\Sigma_4(a_{s,1})]_{(3)|^{s}|(2)|^{t-1}} = \langle \{[a_{s,j}]_{(3)|^{s}|(2)|^{t-1}}: 1 \leqslant j \leqslant 12 \}\rangle,\\
&[\Sigma_4(a_{s,13})]_{(3)|^{s}|(2)|^{t-1}} = \langle \{[a_{s,j}]_{(3)|^{s}|(2)|^{t-1}}: 13 \leqslant j \leqslant 24 \}\rangle,\\
&[\Sigma_4(a_{s,25})]_{(3)|^{s}|(2)|^{t-1}} = \langle \{[a_{s,j}]_{(3)|^{s}|(2)|^{t-1}}: 25 \leqslant j \leqslant 28\} \rangle\\
&[\Sigma_4(a_{s,29})]_{(3)|^{s}|(2)|^{t-1}} = \langle \{[a_{s,j}]_{(3)|^{s}|(2)|^{t-1}}: 29 \leqslant j \leqslant 34\} \rangle.
\end{align*}

We set $p_{s,u} = p_{t,s,u}$, $1 \leqslant u \leqslant 4$ with
\begin{equation}\label{ctt1}
p_{s,1} = \sum_{j=1}^{12}a_{s,j}, \ 
p_{s,2} = \sum_{j=13}^{24}a_{s,j}, \ 
p_{s,3} = \sum_{j=25}^{28}a_{s,j},\
p_{s,4} = \sum_{j=29}^{34}a_{s,j}.
\end{equation}

We need the following lemma for the proof of the proposition. 

\begin{lems}\label{bdtts} For $s \geqslant 2$ and $t\geqslant 4$, we have
$$QP_4((3)|^{s}|(2)|^{t-1})^{\Sigma_4} = \langle \{[p_{s,u}]_{(3)|^{s}|(2)|^{t-1}} : 1 \leqslant u \leqslant 13\} \rangle,$$
where $p_{s,u},\, u =1\, 2,\, 3,\, 4$ are defined by \eqref{ctt1} and
\begin{align*}	
&p_{s,5} = 
\begin{cases}
\sum_{j \in \{35,\, 36,\, 37,\, 38,\, 39,\, 40,\, 44,\, 45,\, 63,\, 64,\, 65,\, 66\}}a_{2,j}, &\mbox{if } s=2,\\
\sum_{j \in \{35,\, 36,\, 37,\, 38,\, 39,\, 40,\, 41,\, 42,\, 43,\, 44,\, 45,\, 46,\, 59,\, 60,\, 61,\, 62\}}a_{s,j}, &\mbox{if } s\geqslant 3,
\end{cases}\\
&p_{s,6} = 
\begin{cases}
\sum_{j \in \{37,\, 38,\, 40,\, 41,\, 42,\, 43,\, 45,\, 46,\, 59,\, 60,\, 61,\, 62,\, 97,\, 98,\, 99,\, 100\}}a_{2,j}, &\mbox{if } s=2,\\
\sum_{j \in \{47,\, 48,\, 49,\, 50,\, 51,\, 52,\, 53,\, 54,\, 55,\, 56,\, 57,\, 58,\, 59,\, 60,\, 61,\, 62\}}a_{s,j}, &\mbox{if } s\geqslant 3,
\end{cases}\\
&p_{s,7} = 
\begin{cases}
\sum_{j \in \{41,\, 42,\, 43,\, 46,\, 47,\, 48,\, 49,\, 50,\, 51,\, 52,\, 53\}\atop\hskip0.8cm\cup\{54,\, 55,\, 56,\, 57,\, 58,\, 63,\, 64,\, 65,\, 66\}}a_{2,j}, &\mbox{if } s=2,\\
	\sum_{j \in \{59,\, 60,\, 61,\, 62,\, 63,\, 64,\, 65,\, 66,\, 97,\, 98,\, 99,\, 100\}}a_{s,j}, &\mbox{if } s\geqslant 3,
\end{cases}\\
&p_{s,8} = 
\begin{cases}
\sum_{j \in \{67,\, 68,\, 70,\, 71,\, 81,\, 83,\, 85,\, 86,\, 87,\, 89,\, 95,\, 96,\, 101\}\atop\hskip3cm \cup\{102,\, 103,\, 104,\, 105\}}a_{2,j}, &\mbox{if } s=2,\\
\sum_{j \in \{80,\, 81,\, 105\}}a_{s,j}, &\mbox{if } s = 3,\\
	\sum_{j \in \{94,\, 102,\, 105\}}a_{s,j}, &\mbox{if } s\geqslant 4,
\end{cases}\\
&p_{s,9} = 
\begin{cases}
\sum_{j \in \{68,\, 69,\, 70,\, 71,\, 80,\, 84,\, 85,\, 87,\, 88,\, 96,\, 101,\, 103\}}a_{2,j}, &\mbox{if } s=2,\\
\sum_{j \in \{95,\, 96,\, 103,\, 104\}}a_{s,j}, &\mbox{if } s\geqslant 3,
\end{cases}\\
&p_{s,10} = 
\begin{cases}
\sum_{j \in \{69,\, 70,\, 72,\, 83,\, 84,\, 86,\, 88,\, 89,\, 94,\, 95,\, 101,\, 104\}}a_{2,j}, &\mbox{if } s=2,\\
\sum_{j \in \{79,\, 80,\, 81,\, 94,\, 101,\, 102\}}a_{s,j}, &\mbox{if } s = 3,\\
\sum_{j \in \{79,\, 80,\, 81,\, 101,\, 105\}}a_{s,j}, &\mbox{if } s\geqslant 4,
\end{cases}\\
&p_{s,11} = 
\begin{cases}
\sum_{j \in \{70,\, 73,\, 74,\, 75,\, 76,\, 77,\, 78,\, 81,\, 84,\, 85,\, 87,\, 88,\, 94,\, 96,\, 102\}\atop\hskip4.7cm \cup\{103,\, 105\}}a_{2,j}, &\mbox{if } s=2,\\
\sum_{j \in \{67,\, 68,\, 69,\, 70,\, 71,\, 72,\, 84,\, 87,\, 88,\, 89,\, 96,\, 101,\, 102,\, 104\}}a_{3,j}, &\mbox{if } s = 3,\\
\sum_{j \in \{67,\, 68,\, 69,\, 70,\, 71,\, 72,\, 84,\, 87,\, 88,\, 89,\, 96,\, 102,\, 104,\, 105\}}a_{s,j}, &\mbox{if } s\geqslant 4,
\end{cases}\\
&p_{s,12} = 
\begin{cases}
\sum_{j \in \{79,\, 80,\, 81\}}a_{2,j}, &\mbox{if } s=2,\\
\sum_{j \in \{79,\, 80,\, 81,\, 82,\, 83,\, 84,\, 85,\, 86,\, 87,\, 88,\, 89,\, 90,\, 91,\, 92,\, 93,\, 101\}}a_{s,j}, &\mbox{if } s = 3,\\
\sum_{j \in \{82,\, 83,\, 84,\, 85,\, 86,\, 87,\, 88,\, 89,\, 90,\, 91,\, 92,\, 93,\, 105\}}a_{s,j}, &\mbox{if } s\geqslant 4,
\end{cases}\\
&p_{s,13} = 
\begin{cases}
\sum_{j \in \{80,\, 81,\, 82,\, 83,\, 84,\, 85,\, 86,\, 87,\, 88,\, 89,\, 90,\, 91,\, 92,\, 93\}}a_{2,j}, &\mbox{if } s=2,\\
\sum_{j \in \{73,\, 74,\, 75,\, 76,\, 77,\, 78,\, 84,\, 87,\, 88,\, 89,\, 96,\, 101,\, 102,\, 104\}}a_{3,j}, &\mbox{if } s = 3,\\
\sum_{j \in \{73,\, 74,\, 75,\, 76,\, 77,\, 78,\, 84,\, 87,\, 88,\, 89,\, 96,\, 102,\, 104,\, 105\}}a_{s,j}, &\mbox{if } s\geqslant 4,
\end{cases}
\end{align*}
\end{lems}
\begin{proof} By a simple computation, we see that $[p_{s,u}]_{(3)|^s|(2)|^{t-1}}$ is an $\Sigma_4$ invariant for $1 \leqslant u \leqslant 13$, and there is a direct summand decomposition of $\Sigma_4$-modules:
\begin{align*}QP_4(&(3)|^s|(2)|^{t-1}) = [\Sigma_4(a_{s,1})]_{(3)|^s|(2)|^{t-1}}\bigoplus [\Sigma_4(a_{s,13})]_{(3)|^s|(2)|^{t-1}}\\ &\bigoplus [\Sigma_4(a_{s,25})]_{(3)|^s|(2)|^{t-1}} \bigoplus [\Sigma_4(a_{s,29})]_{(3)|^s|(2)|^{t-1}}\bigoplus \mathcal U_{t,s} \bigoplus \mathcal V_{t,s},
\end{align*}
where 
\begin{align*}
\mathcal U_{t,s} &= \langle [a_{s,j}]_{(3)|^s|(2)|^{t-1}}: j \in \mathbb J_{t,s} = \{35,\ldots 66,\, 97, 98,\, 99,\, 100\}\rangle,\\
\mathcal V_{t,s} &= \langle [a_{s,j}]_{(3)|^s|(2)|^{t-1}}: j \in \mathbb K_{t,s} =  \{67,\ldots 96,\, 101,\ldots ,\, 105\}\rangle.
\end{align*}
By a simple computation, we easily obtain
\begin{align*}
&[\Sigma_4(a_{s,1})]_{(3)|^s|(2)|^{t-1}}^{\Sigma_4} = \langle[p_{s,1}]_{(3)|^s|(2)|^{t-1}}\rangle,\\ &[\Sigma_4(a_{s,13})]_{(3)|^s|(2)|^{t-1}}^{\Sigma_4} = \langle[p_{s,2}]_{(3)|^s|(2)|^{t-1}}\rangle,\\
&[\Sigma_4(a_{s,25})]_{(3)|^s|(2)|^{t-1}}^{\Sigma_4} = \langle[p_{s,3}]_{(3)|^s|(2)|^{t-1}}\rangle,\\
&[\Sigma_4(a_{s,29})]_{(3)|^s|(2)|^{t-1}}^{\Sigma_4} = \langle[p_{s,4}]_{(3)|^s|(2)|^{t-1}}\rangle.
\end{align*}
We prove $\mathcal U_{t,s}^{\Sigma_4} = \langle [p_{s,u}]_{(3)|^{s}|(2)|^{t-1}} : 5 \leqslant u \leqslant 7\}\rangle$, $\mathcal V_{t,s}^{\Sigma_4} = \langle [p_{s,u}]_{(3)|^{s}|(2)|^{t-1}} : 8 \leqslant u \leqslant 13\}\rangle$. The computation of $U_{t,s}^{\Sigma_4}$ is similar to the one of $U_{3,2}^{\Sigma_4}$ in Lemma \ref{bdt32} with some changes. We only present the detailed computations for the case $s = 2$.

We observe that the leading monomials of $p_{2,5}$, $p_{2,6}$, $p_{2,7}$ respectively are $a_{2,66}$, $a_{2,62}$, $a_{2,58}$. Hence, if $[f]_{(3)|^{2}|(2)|^{t-1}} \in \mathcal U_{t,2}^{\Sigma_4}$ with $f \in P_4$, then there are $\gamma_5,\, \gamma_6,\, \gamma_7 \in \mathbb F_2$ such that
$$\bar f := f + \gamma_5p_{2,5}+ \gamma_6p_{2,6} + \gamma_7p_{2,7}  \equiv _{(3)|^{2}|(2)|^{t-1}} \sum_{j\in \mathbb J_{t,2},\, j \ne 58,62,66}a_{2,j}.$$
Then, $[\bar f]_{(3)|^{2}|(2)|^{t-1}}$ is also an $\Sigma_4$-invariant. A direct computation shows
\begin{align*}
\rho_1(\bar f)& + \bar f \equiv_{(3)|^{2}|(2)|^{t-1}} \gamma_{\{35,37,41,47\}}a_{2,35} + \gamma_{\{36,38,42,48\}}a_{2,36} + \gamma_{\{39,44\}}a_{2,39}\\ 
&+ \gamma_{\{40,45\}}a_{2,40} + \gamma_{\{35,37,41,47\}}a_{2,41} + \gamma_{\{36,38,42,48\}}a_{2,42} + \gamma_{\{43,46\}}a_{2,43}\\ 
& + \gamma_{\{39,44\}}a_{2,44} + \gamma_{\{40,45\}}a_{2,45} + \gamma_{\{43,46\}}a_{2,46} + \gamma_{\{49,50\}}a_{2,49}\\ 
& + \gamma_{\{49,50\}}a_{2,50} + \gamma_{\{51,54\}}a_{2,51} + \gamma_{\{52,55\}}a_{2,52} + \gamma_{\{53,57\}}a_{2,53}\\ 
& + \gamma_{\{51,54\}}a_{2,54} + \gamma_{\{52,55\}}a_{2,55} + \gamma_{\{56,58\}}a_{2,56} + \gamma_{\{53,57\}}a_{2,57}\\ 
& + \gamma_{\{56,58\}}a_{2,58} + \gamma_{\{61,62\}}a_{2,61} + \gamma_{\{61,62\}}a_{2,62} + \gamma_{\{65,66\}}a_{2,65}\\ 
& + \gamma_{\{65,66\}}a_{2,66} + \gamma_{\{37,47,59,63\}}a_{2,97} + \gamma_{\{38,48,60,64\}}a_{2,98} + \gamma_{\{99,100\}}a_{2,99}\\ 
& + \gamma_{\{99,100\}}a_{2,100} \equiv_{(3)|^{2}|(2)|^{t-1}} 0,\\
\rho_2(\bar f) &+ \bar f \equiv_{(3)|^{2}|(2)|^{t-1}} \gamma_{\{35,37,97\}}a_{2,35} + \gamma_{\{36,39\}}a_{2,36} + \gamma_{\{35,37,97\}}a_{2,37}\\ 
& + \gamma_{\{38,40\}}a_{2,38} + \gamma_{\{36,39\}}a_{2,39} + \gamma_{\{38,40\}}a_{2,40} + \gamma_{\{41,54,59\}}a_{2,41}\\ 
& + \gamma_{\{42,43\}}a_{2,42} + \gamma_{\{42,43\}}a_{2,43} + \gamma_{\{44,45,46,50\}}a_{2,44} + \gamma_{\{44,45,46,50\}}a_{2,46}\\ 
& + \gamma_{\{47,51\}}a_{2,47} + \gamma_{\{48,49\}}a_{2,48} + \gamma_{\{48,49\}}a_{2,49} + \gamma_{\{47,51\}}a_{2,51}\\ 
& + \gamma_{\{52,53\}}a_{2,52} + \gamma_{\{52,53\}}a_{2,53} + \gamma_{\{55,56\}}a_{2,55} + \gamma_{\{55,56\}}a_{2,56}\\ 
& + \gamma_{\{57,58\}}a_{2,57} + \gamma_{\{57,58\}}a_{2,58} + \gamma_{\{41,54,59\}}a_{2,59} + \gamma_{\{60,61\}}a_{2,60}\\ 
& + \gamma_{\{60,61\}}a_{2,61} + \gamma_{\{64,65\}}a_{2,64} + \gamma_{\{64,65\}}a_{2,65} + \gamma_{\{98,99\}}a_{2,98}\\ 
& + \gamma_{\{98,99\}}a_{2,99} + \gamma_{\{45,50,62,66\}}a_{2,100} \equiv_{(3)|^{2}|(2)|^{t-1}} 0,\\  
\rho_3(\bar f) &+ \bar f \equiv_{(3)|^{2}|(2)|^{t-1}} \gamma_{\{35,36\}}a_{2,35} + \gamma_{\{35,36\}}a_{2,36} + \gamma_{\{37,38\}}a_{2,37}\\ 
& + \gamma_{\{37,38\}}a_{2,38} + \gamma_{\{39,40,99\}}a_{2,39} + \gamma_{\{39,40,99\}}a_{2,40} + \gamma_{\{41,42\}}a_{2,41}\\ 
& + \gamma_{\{41,42\}}a_{2,42} + \gamma_{\{43,56,61\}}a_{2,43} + \gamma_{\{44,45,100\}}a_{2,44} + \gamma_{\{44,45,100\}}a_{2,45}\\ 
& + \gamma_{\{46,58,62\}}a_{2,46} + \gamma_{\{47,48\}}a_{2,47} + \gamma_{\{47,48\}}a_{2,48} + \gamma_{\{49,53\}}a_{2,49}\\ 
& + \gamma_{\{50,57\}}a_{2,50} + \gamma_{\{51,52\}}a_{2,51} + \gamma_{\{51,52\}}a_{2,52} + \gamma_{\{49,53\}}a_{2,53}\\ 
& + \gamma_{\{54,55\}}a_{2,54} + \gamma_{\{54,55\}}a_{2,55} + \gamma_{\{50,57\}}a_{2,57} + \gamma_{\{59,60\}}a_{2,59}\\ 
& + \gamma_{\{59,60\}}a_{2,60} + \gamma_{\{43,56,61\}}a_{2,61} + \gamma_{\{46,58,62\}}a_{2,62} + \gamma_{\{63,64\}}a_{2,63}\\ 
& + \gamma_{\{63,64\}}a_{2,64} + \gamma_{\{97,98\}}a_{2,97} + \gamma_{\{97,98\}}a_{2,98} \equiv_{(3)|^{2}|(2)|^{t-1}} 0.
\end{align*}
From these equalities we get $\gamma_{j} = 0$ for all $j \in \mathbb J_{t,2}\setminus\{58,62,66\}$. Hence, we get $\bar f \equiv_{(3)|^{2}|(2)|^{t-1}} 0$ and
$f \equiv_{(3)|^{2}|(2)|^{t-1}}  \gamma_5p_{2,5}+ \gamma_6p_{2,6} + \gamma_7p_{2,7}.$

\medskip
We now compute $\mathcal V_{t,s}^{\Sigma_4}$.

For $s = 2$, the leading monomials of $p_{2,u}$, $8 \leqslant u \leqslant 13$, respectively are $a_{2,95}$, $a_{2,103}$, $a_{2,94}$, $a_{2,78}$, $a_{2,81}$, $a_{2,93}$. So, if $[f]_{(3)|^{2}|(2)|^{t-1}} \in \mathcal V_{t,2}^{\Sigma_4}$ with $f \in P_4$, then there are $\gamma_u\in \mathbb F_2$, $8 \leqslant u \leqslant 13$, such that
$$g := f + \sum_{8 \leqslant u \leqslant 13}\gamma_up_{2,u} \equiv _{(3)|^{2}|(2)|^{t-1}} \sum_{j\in \mathbb K_{t,2}\setminus\{78,\, 81,\, 93,\, 94,\, 95,\, 103\}}a_{2,j}.$$
Then, $[g]_{3)|^{2}|(2)|^{t-1}}$ is also an $\Sigma_4$-invariant. A direct computation shows
\begin{align*}
\rho_1(g)& + g \equiv_{3)|^{2}|(2)|^{t-1}} \gamma_{\{67,69,70,73,84,88,90,91\}}a_{2,67} + \gamma_{\{68,71,84,88,90,91\}}a_{2,68}\\ 
& + \gamma_{\{67,69,70,73,84,88,90,91\}}a_{2,70} + \gamma_{\{68,71,84,88,90,91\}}a_{2,71} + \gamma_{\{74,75\}}a_{2,74}\\ 
& + \gamma_{\{74,75\}}a_{2,75} + \gamma_{\{76,77\}}a_{2,76} + \gamma_{\{69,73,80,81,82,83,84,88,90,91\}}a_{2,79}\\ 
& + \gamma_{\{76,77\}}a_{2,77} + \gamma_{\{69,73,88,91\}}a_{2,82} + \gamma_{\{69,73,88,91\}}a_{2,83}\\ 
& + \gamma_{\{84,85,86,90\}}a_{2,85} + \gamma_{\{84,85,86,90\}}a_{2,86} + \gamma_{\{87,88,89,91\}}a_{2,87}\\ 
& + \gamma_{\{87,88,89,91\}}a_{2,89} + \gamma_{\{92,93\}}a_{2,92} + \gamma_{\{69,73,84,88,90,91,101,102\}}a_{2,101}\\ 
& + \gamma_{\{92,93\}}a_{2,93} + \gamma_{\{69,73,84,88,90,91,101,102\}}a_{2,102} + \gamma_{\{84,88,90,91\}}a_{2,104}\\ 
& + \gamma_{\{88,91,95\}}a_{2,105} \equiv_{3)|^{2}|(2)|^{t-1}} 0,\\
\rho_2(g) &+ g \equiv_{3)|^{2}|(2)|^{t-1}} \gamma_{\{71,75,94,96,104,105\}}a_{2,67} + \gamma_{\{68,69,71,75,94,96,105\}}a_{2,68}\\ 
& + \gamma_{\{68,69,104\}}a_{2,69} + \gamma_{\{70,72,94,96,104,105\}}a_{2,70} + \gamma_{\{71,75,94,96,104,105\}}a_{2,71}\\ 
& + \gamma_{\{70,71,72,75\}}a_{2,72} + \gamma_{\{73,74\}}a_{2,73} + \gamma_{\{73,74\}}a_{2,74} + \gamma_{\{77,78\}}a_{2,77}\\ 
& + \gamma_{\{77,78\}}a_{2,78} + \gamma_{\{79,80,86,96,104,105\}}a_{2,79} + \gamma_{\{71,75,79,80,86,94\}}a_{2,80}\\ 
& + \gamma_{\{82,85,94,96,104,105\}}a_{2,82} + \gamma_{\{83,84,94,96,104,105\}}a_{2,83} + \gamma_{\{71,75,83,84\}}a_{2,84}\\ 
& + \gamma_{\{71,75,82,85\}}a_{2,85} + \gamma_{\{87,88,94,105\}}a_{2,87} + \gamma_{\{87,88,94,105\}}a_{2,88}\\ 
& + \gamma_{\{89,93,95\}}a_{2,89} + \gamma_{\{91,92\}}a_{2,91} + \gamma_{\{91,92\}}a_{2,92} + \gamma_{\{89,93,95\}}a_{2,95}\\ 
& + \gamma_{\{71,75,94,96,104,105\}}a_{2,101} + \gamma_{\{95,96,102,103,104,105\}}a_{2,102}\\ 
& + \gamma_{\{71,75,89,93,94,102,103\}}a_{2,103} + \gamma_{\{71,75,94,96,104,105\}}a_{2,104} \equiv_{3)|^{2}|(2)|^{t-1}} 0,\\  
\rho_3(g) &+ g \equiv_{3)|^{2}|(2)|^{t-1}} \gamma_{\{67,68,101,104\}}a_{2,67} + \gamma_{\{67,68,101,104\}}a_{2,68}\\ 
& + \gamma_{\{70,71,102,104\}}a_{2,70} + \gamma_{\{70,71,102,104\}}a_{2,71} + \gamma_{\{72,78,94\}}a_{2,72}\\ 
& + \gamma_{\{74,76\}}a_{2,74} + \gamma_{\{75,77\}}a_{2,75} + \gamma_{\{74,76\}}a_{2,76} + \gamma_{\{75,77\}}a_{2,77}\\ 
& + \gamma_{\{80,81,103\}}a_{2,80} + \gamma_{\{80,81,103\}}a_{2,81} + \gamma_{\{82,83,104\}}a_{2,82} + \gamma_{\{82,83,104\}}a_{2,83}\\ 
& + \gamma_{\{84,88\}}a_{2,84} + \gamma_{\{85,87\}}a_{2,85} + \gamma_{\{86,89\}}a_{2,86} + \gamma_{\{85,87\}}a_{2,87}\\ 
& + \gamma_{\{84,88\}}a_{2,88} + \gamma_{\{86,89\}}a_{2,89} + \gamma_{\{90,91\}}a_{2,90} + \gamma_{\{90,91\}}a_{2,91}\\ 
& + \gamma_{\{72,78,94\}}a_{2,94} + \gamma_{\{95,96,103,104\}}a_{2,105} \equiv_{3)|^{2}|(2)|^{t-1}} 0.
\end{align*}
From these equalities we get $\gamma_{j} = 0$ for all $j \in \mathbb K_{t,2}\setminus\{78,\, 81,\, 93,\, 94,\, 95,\, 103\}$. Hence, we get $g \equiv_{3)|^{2}|(2)|^{t-1}} 0$ and
$f \equiv_{3)|^{2}|(2)|^{t-1}}  \sum_{8 \leqslant u \leqslant 13}\gamma_up_{2,u}.$

For $s = 3$, the leading monomials of $p_{3,u}$, $8 \leqslant u \leqslant 13$, respectively are $a_{3,105}$, $a_{3,95}$, $a_{3,94}$, $a_{3,72}$, $a_{3,93}$, $a_{3,78}$. So, if $[f]_{(3)|^{3}|(2)|^{t-1}} \in \mathcal V_{t,3}^{\Sigma_4}$ with $f \in P_4$, then there are $\gamma_u\in \mathbb F_2$, $8 \leqslant u \leqslant 13$, such that
$$g := f + \sum_{8 \leqslant u \leqslant 13}\gamma_up_{3,u} \equiv _{(3)|^{3}|(2)|^{t-1}} \sum_{j\in \mathbb K_{t,3}\setminus\{72,\, 78,\, 93,\, 94,\, 95,\, 105\}}a_{3,j}.$$
Then, $[g]_{(3)|^{3}|(2)|^{t-1}}$ is also an $\Sigma_4$-invariant. A direct computation shows
\begin{align*}
\rho_1(g)& + g \equiv_{(3)|^{3}|(2)|^{t-1}} \gamma_{\{67,70\}}a_{3,67} + \gamma_{\{68,71\}}a_{3,68} + \gamma_{\{67,70\}}a_{3,70}\\ 
& + \gamma_{\{68,71\}}a_{3,71} + \gamma_{\{74,75\}}a_{3,74} + \gamma_{\{74,75\}}a_{3,75} + \gamma_{\{76,77\}}a_{3,76}\\ 
& + \gamma_{\{76,77\}}a_{3,77} + \gamma_{\{69,73,79,80,81,82,83,101\}}a_{3,79} + \gamma_{\{69,73,88,91\}}a_{3,82}\\ 
& + \gamma_{\{69,73,84,90\}}a_{3,83} + \gamma_{\{85,86\}}a_{3,85} + \gamma_{\{85,86\}}a_{3,86} + \gamma_{\{87,89\}}a_{3,87}\\ 
& + \gamma_{\{87,89\}}a_{3,89} + \gamma_{\{92,93\}}a_{3,92} + \gamma_{\{69,73,79,80,81,82,83,84,88,90,91,101\}}a_{3,101}\\ 
& + \gamma_{\{92,93\}}a_{3,93} + \gamma_{\{84,90,95,96\}}a_{3,103} + \gamma_{\{88,91,95,96\}}a_{3,104}\\ 
& + \gamma_{\{69,73,80,81,82,83,94,95,96,102,103,104\}}a_{3,105} \equiv_{(3)|^{3}|(2)|^{t-1}} 0,\\
\rho_2(g) &+ g \equiv_{(3)|^{3}|(2)|^{t-1}} \gamma_{\{68,69\}}a_{3,68} + \gamma_{\{68,69\}}a_{3,69} + \gamma_{\{70,72\}}a_{3,70}\\ 
& + \gamma_{\{70,72\}}a_{3,72} + \gamma_{\{73,74\}}a_{3,73} + \gamma_{\{73,74\}}a_{3,74} + \gamma_{\{77,78\}}a_{3,77}\\ 
& + \gamma_{\{77,78\}}a_{3,78} + \gamma_{\{79,80,105\}}a_{3,79} + \gamma_{\{79,80,105\}}a_{3,80} + \gamma_{\{82,85\}}a_{3,82}\\ 
& + \gamma_{\{71,75,83,84\}}a_{3,83} + \gamma_{\{71,75,83,84\}}a_{3,84} + \gamma_{\{82,85\}}a_{3,85} + \gamma_{\{71,75,89,93\}}a_{3,86}\\ 
& + \gamma_{\{87,88\}}a_{3,87} + \gamma_{\{87,88\}}a_{3,88} + \gamma_{\{91,92\}}a_{3,91} + \gamma_{\{91,92\}}a_{3,92}\\ 
& + \gamma_{\{89,93,95,104\}}a_{3,95} + \gamma_{\{86,101,102\}}a_{3,101} + \gamma_{\{71,75,86,89,93,101,102\}}a_{3,102}\\ 
& + \gamma_{\{71,75,89,93\}}a_{3,103} + \gamma_{\{89,93,95,104\}}a_{3,104} \equiv_{(3)|^{3}|(2)|^{t-1}} 0,\\  
\rho_3(g) &+ g \equiv_{(3)|^{3}|(2)|^{t-1}} \gamma_{\{67,68\}}a_{3,67} + \gamma_{\{67,68\}}a_{3,68} + \gamma_{\{70,71\}}a_{3,70}\\ 
& + \gamma_{\{70,71\}}a_{3,71} + \gamma_{\{74,76\}}a_{3,74} + \gamma_{\{75,77\}}a_{3,75} + \gamma_{\{74,76\}}a_{3,76}\\ 
& + \gamma_{\{75,77\}}a_{3,77} + \gamma_{\{80,81\}}a_{3,80} + \gamma_{\{80,81\}}a_{3,81} + \gamma_{\{82,83\}}a_{3,82}\\ 
& + \gamma_{\{82,83\}}a_{3,83} + \gamma_{\{84,88\}}a_{3,84} + \gamma_{\{72,78,85,87\}}a_{3,85} + \gamma_{\{72,78,86,89\}}a_{3,86}\\ 
& + \gamma_{\{72,78,85,87\}}a_{3,87} + \gamma_{\{84,88\}}a_{3,88} + \gamma_{\{72,78,86,89\}}a_{3,89} + \gamma_{\{90,91\}}a_{3,90}\\ 
& + \gamma_{\{90,91\}}a_{3,91} + \gamma_{\{72,78,94,102\}}a_{3,94} + \gamma_{\{72,78,94,102\}}a_{3,102}\\ 
& + \gamma_{\{72,78,103,104\}}a_{3,103} + \gamma_{\{72,78,103,104\}}a_{3,104} \equiv_{(3)|^{3}|(2)|^{t-1}} 0.
\end{align*}
From these equalities we get $\gamma_{j} = 0$ for all $j \in \mathbb K_{t,3}\setminus\{72,\, 78,\, 93,\, 94,\, 95,\, 105\}$. Hence, we get $g \equiv_{(3)|^{3}|(2)|^{t-1}} 0$ and
$f \equiv_{(3)|^{3}|(2)|^{t-1}}  \sum_{8 \leqslant u \leqslant 13}\gamma_up_{3,u}.$

For $s \geqslant 4$, we see that the leading monomials of $p_{3,u}$, $8 \leqslant u \leqslant 13$, respectively are $a_{3,94}$, $a_{3,95}$, $a_{3,101}$, $a_{3,72}$, $a_{3,93}$, $a_{3,78}$. Therefore, if $[h]_{(3)|^{s}|(2)|^{t-1}} \in \mathcal V_{t,s}^{\Sigma_4}$ with $h \in P_4$, then there are $\gamma_u\in \mathbb F_2$, $8 \leqslant u \leqslant 13$, such that
$$h^* := h + \sum_{8 \leqslant u \leqslant 13}\gamma_up_{s,u}  \equiv _{(3)|^{s}|(2)|^{t-1}} \sum_{j\in \mathbb K_{t,s}\setminus\{72,\, 78,\, 93,\, 94,\, 95,\, 101\}}a_{s,j}.$$
Since $[h]_{(3)|^{s}|(2)|^{t-1}}$ is an $\Sigma_4$-invariant, so is $[h^*]_{(3)|^{s}|(2)|^{t-1}}$. By a direct computation we have
\begin{align*}
\rho_1(h^*)& + h^* \equiv_{(3)|^{s}|(2)|^{t-1}} \gamma_{\{67,70\}}a_{s,67} + \gamma_{\{68,71\}}a_{s,68} + \gamma_{\{67,70\}}a_{s,70} + \gamma_{\{68,71\}}a_{s,71}\\ 
& + \gamma_{\{74,75\}}a_{s,74} + \gamma_{\{74,75\}}a_{s,75} + \gamma_{\{76,77\}}a_{s,76} + \gamma_{\{76,77\}}a_{s,77} + \gamma_{\{79,101\}}a_{s,79}\\ 
& + \gamma_{\{69,73,88,91\}}a_{s,82} + \gamma_{\{69,73,84,90\}}a_{s,83} + \gamma_{\{85,86\}}a_{s,85} + \gamma_{\{85,86\}}a_{s,86}\\ 
& + \gamma_{\{87,89\}}a_{s,87} + \gamma_{\{87,89\}}a_{s,89} + \gamma_{\{92,93\}}a_{s,92} + \gamma_{\{92,93\}}a_{s,93}\\ 
& + \gamma_{\{79,101\}}a_{s,101} + \gamma_{\{84,90,95,96\}}a_{s,103} + \gamma_{\{88,91,95,96\}}a_{s,104}\\ 
& + \gamma_{\{69,73,80,81,82,83,84,88,90,91,94,95,96,102,103,104\}}a_{s,105} \equiv_{(3)|^{s}|(2)|^{t-1}} 0,\\
\rho_2(h^*) &+ h^* \equiv_{(3)|^{s}|(2)|^{t-1}} \gamma_{\{68,69\}}a_{s,68} + \gamma_{\{68,69\}}a_{s,69} + \gamma_{\{70,72\}}a_{s,70} + \gamma_{\{70,72\}}a_{s,72}\\ 
& + \gamma_{\{73,74\}}a_{s,73} + \gamma_{\{73,74\}}a_{s,74} + \gamma_{\{77,78\}}a_{s,77} + \gamma_{\{77,78\}}a_{s,78}\\ 
& + \gamma_{\{79,80\}}a_{s,79} + \gamma_{\{79,80\}}a_{s,80} + \gamma_{\{82,85\}}a_{s,82} + \gamma_{\{71,75,83,84\}}a_{s,83}\\ 
& + \gamma_{\{71,75,83,84\}}a_{s,84} + \gamma_{\{82,85\}}a_{s,85} + \gamma_{\{71,75,89,93\}}a_{s,86} + \gamma_{\{87,88\}}a_{s,87}\\ 
& + \gamma_{\{87,88\}}a_{s,88} + \gamma_{\{91,92\}}a_{s,91} + \gamma_{\{91,92\}}a_{s,92} + \gamma_{\{89,93,95,104\}}a_{s,95}\\ 
& + \gamma_{\{71,75,86,89,93,101,102,105\}}a_{s,102} + \gamma_{\{71,75,89,93\}}a_{s,103} + \gamma_{\{89,93,95,104\}}a_{s,104}\\ 
& + \gamma_{\{86,101,102,105\}}a_{s,105} \equiv_{(3)|^{s}|(2)|^{t-1}} 0,\\  
\rho_3(h^*) &+ h^* \equiv_{(3)|^{s}|(2)|^{t-1}} \gamma_{\{67,68\}}a_{s,67} + \gamma_{\{67,68\}}a_{s,68} + \gamma_{\{70,71\}}a_{s,70} + \gamma_{\{70,71\}}a_{s,71}\\ 
& + \gamma_{\{74,76\}}a_{s,74} + \gamma_{\{75,77\}}a_{s,75} + \gamma_{\{74,76\}}a_{s,76} + \gamma_{\{75,77\}}a_{s,77}\\ 
& + \gamma_{\{80,81\}}a_{s,80} + \gamma_{\{80,81\}}a_{s,81} + \gamma_{\{82,83\}}a_{s,82} + \gamma_{\{82,83\}}a_{s,83}\\ 
& + \gamma_{\{84,88\}}a_{s,84} + \gamma_{\{72,78,85,87\}}a_{s,85} + \gamma_{\{72,78,86,89\}}a_{s,86}\\ 
& + \gamma_{\{72,78,85,87\}}a_{s,87} + \gamma_{\{84,88\}}a_{s,88} + \gamma_{\{72,78,86,89\}}a_{s,89} + \gamma_{\{90,91\}}a_{s,90}\\ 
& + \gamma_{\{90,91\}}a_{s,91} + \gamma_{\{72,78,94,102\}}a_{s,94} + \gamma_{\{72,78,94,102\}}a_{s,102}\\ 
& + \gamma_{\{72,78,103,104\}}a_{s,103} + \gamma_{\{72,78,103,104\}}a_{s,104} \equiv_{(3)|^{s}|(2)|^{t-1}} 0.
\end{align*}
These equalities imply $\gamma_{j} = 0$ for all $j \in \mathbb K_{t,s}\setminus\{72,\, 78,\, 93,\, 94,\, 95,\, 101\}$ and $h^* \equiv_{(3)|^{s}|(2)|^{t-1}} 0$. Hence, we obtain
$h \equiv_{(3)|^{s}|(2)|^{t-1}}  \sum_{8 \leqslant u \leqslant 13}\gamma_up_{s,u}.$
The lemma is completely proved.
\end{proof}

\begin{proof}[Proof of Proposition \ref{mdt41}]
Suppose  $[f]_{(3)|^s|(2)|^{t-1}} \in QP_4((3)|^s|(2)|^{t-1})^{GL_4}$ with $f \in QP_4((3)|^s|(2)|^{t-1})$. Then, $[f]_{(3)|^s|(2)|^2} \in QP_4((3)|^s|(2)|^2)^{\Sigma_4}$. 
	
By Lemma \ref{bdtts}, we have 
$$ f \equiv_{(3)|^s|(2)|^{t-1}} \sum_{1\leqslant u \leqslant 13} \gamma_up_{s,u},$$
where $\gamma_u \in \mathbb F_2$. 
	
For $s = 2$, by computing $\rho_4(f)+f$ in terms of the admissible monomials, we get
\begin{align*}
\rho_4(f)&+f \equiv_{(3)|^2|(2)|^{t-1}} \gamma_{\{1,4\}}a_{2,1} + \gamma_{\{1,5,6\}}a_{2,2} + \gamma_{\{1,5,6\}}a_{2,3} + \gamma_{\{1,2\}}a_{2,5}\\ 
& + \gamma_{\{1,2\}}a_{2,6} + \gamma_{\{2,7\}}a_{2,13} + \gamma_{\{2,7\}}a_{2,14} + \gamma_{\{2,9,10\}}a_{2,15} + \gamma_{\{2,3\}}a_{2,20}\\ 
& + \gamma_{\{2,3\}}a_{2,21} + \gamma_{\{3,11\}}a_{2,25} + \gamma_{\{4,6\}}a_{2,30} + \gamma_{\{4,6\}}a_{2,31} + \gamma_{\{6,7\}}a_{2,35}\\ 
& + \gamma_{\{6,7\}}a_{2,36} + \gamma_{\{5,9,12,13\}}a_{2,39} + \gamma_{\{5,6,8,11,12,13\}}a_{2,40} + \gamma_{\{6,7,10\}}a_{2,43}\\ 
& + \gamma_{\{7,9,10,11,13\}}a_{2,49} + \gamma_{5}a_{2,51} + \gamma_{5}a_{2,52} + \gamma_{\{7,9,10,11,13\}}a_{2,53}\\ 
& + \gamma_{\{7,11\}}a_{2,56} + \gamma_{\{6,10,11\}}a_{2,61} + \gamma_{\{5,7,8,10\}}a_{2,65} + \gamma_{\{2,3,7,8,9,10,13\}}a_{2,67}\\ 
& + \gamma_{\{2,3,7,9,10,11,13\}}a_{2,68} + \gamma_{\{2,3,7,11\}}a_{2,70} + \gamma_{\{2,3,7,11\}}a_{2,71}\\ 
& + \gamma_{\{11,13\}}a_{2,74} + \gamma_{\{11,13\}}a_{2,76} + \gamma_{\{2,3,6,10,11\}}a_{2,79} + \gamma_{\{2,3,7,11\}}a_{2,82}\\ 
& + \gamma_{\{2,3,7,11\}}a_{2,83} + \gamma_{\{8,10,13\}}a_{2,85} + \gamma_{\{8,10,13\}}a_{2,87} + \gamma_{\{8,9,11,13\}}a_{2,92}\\ 
& + \gamma_{\{1,2,4,7\}}a_{2,97} + \gamma_{\{1,2,4,7\}}a_{2,98} + \gamma_{\{6,8,9,11\}}a_{2,99} + \gamma_{\{8,11\}}a_{2,101}\\ 
& + \gamma_{\{5,11,13\}}a_{2,105} \equiv_{(3)|^2|(2)|^{t-1}} 0.
	\end{align*}
This equality implies  $\gamma_u = 0$ for $u = 5,\, 10,\, 12$ and $\gamma_u = \gamma_1$ for $1 \leqslant u \leqslant 13,\, u \ne 5,\, 10,\, 12$. Hence, 
\begin{align*}f &\equiv_{(3)|^s|(2)|^{t-1}} \gamma_1\Big(\sum_{1 \leqslant u \leqslant 13,\, u \ne 5,\, 10,\, 12}p_{2,u}\Big)
= \gamma_1\bar\xi_{t,2}.\end{align*} 
The proposition is proved for $s = 2$.
	
For $s = 3$, by computing $\rho_4(f)+f$ in terms of the admissible monomials, we get
\begin{align*}
\rho_4(f)&+f \equiv_{(3)|^3|(2)|^{t-1}} \gamma_{\{1,4\}}a_{3,1} + \gamma_{\{1,5\}}a_{3,2} + \gamma_{\{1,5\}}a_{3,3} + \gamma_{\{1,2\}}a_{3,5}\\ 
& + \gamma_{\{1,2\}}a_{3,6} + \gamma_{\{2,6\}}a_{3,13} + \gamma_{\{2,6\}}a_{3,14} + \gamma_{\{2,11\}}a_{3,15} + \gamma_{\{2,3\}}a_{3,20}\\ 
& + \gamma_{\{2,3\}}a_{3,21} + \gamma_{\{3,13\}}a_{3,25} + \gamma_{\{4,5,6,7\}}a_{3,30} + \gamma_{\{4,5,6,7\}}a_{3,31}\\ 
& + \gamma_{\{1,2,4,5,6,7\}}a_{3,35} + \gamma_{\{1,2,4,5,6,7\}}a_{3,36} + \gamma_{\{5,8,10,12\}}a_{3,39} + \gamma_{\{5,8,10,12\}}a_{3,40}\\ 
& + \gamma_{\{1,2,4,6\}}a_{3,41} + \gamma_{\{1,2,4,6\}}a_{3,42} + \gamma_{\{5,11\}}a_{3,43} + \gamma_{\{6,11,12,13\}}a_{3,49}\\ 
& + \gamma_{\{6,7\}}a_{3,51} + \gamma_{\{6,7\}}a_{3,52} + \gamma_{\{6,11,12,13\}}a_{3,53} + \gamma_{\{6,13\}}a_{3,56}\\ 
& + \gamma_{\{5,6,7,10\}}a_{3,61} + \gamma_{\{7,9\}}a_{3,65} + \gamma_{\{2,3,6,11,12,13\}}a_{3,67} + \gamma_{\{2,3,6,11,12,13\}}a_{3,68}\\ 
& + \gamma_{\{2,3,6,13\}}a_{3,70} + \gamma_{\{2,3,6,13\}}a_{3,71} + \gamma_{\{12,13\}}a_{3,74} + \gamma_{\{12,13\}}a_{3,76}\\ 
& + \gamma_{\{10,11,12,13\}}a_{3,79} + \gamma_{\{2,3,6,13\}}a_{3,82} + \gamma_{\{2,3,6,13\}}a_{3,83} + \gamma_{\{6,7,9,13\}}a_{3,85}\\ 
& + \gamma_{\{6,7,12,13\}}a_{3,86} + \gamma_{\{6,7,9,13\}}a_{3,87} + \gamma_{\{6,7,12,13\}}a_{3,89} + \gamma_{\{9,11,12,13\}}a_{3,92}\\ 
& + \gamma_{\{1,2,4,6\}}a_{3,97} + \gamma_{\{1,2,4,6\}}a_{3,98} + \gamma_{\{7,10,11,13\}}a_{3,99} + \gamma_{\{6,7,12,13\}}a_{3,103}\\ 
& + \gamma_{\{6,7,12,13\}}a_{3,104} + \gamma_{\{2,3,5,6,7,11,12,13\}}a_{3,105} \equiv_{(3)|^3|(2)|^{t-1}} 0.
\end{align*}
The above equality implies $\gamma_u = \gamma_1$ for $1 \leqslant u \leqslant 13$. Hence, we get 
$$f \equiv_{(3)|^3|(2)|^{t-1}} \gamma_1\Big(\sum_{1 \leqslant u \leqslant 13}p_{3,u}\Big) = \gamma_1\Big(\sum_{1 \leqslant j \leqslant 105,\, j \ne 79,\, 101}a_{3,j}\Big)= \gamma_1\bar \xi_{t,3}.$$ 
The proposition is proved for $s = 3$.
	
For $s \geqslant 4$, by computing $\rho_4(f)+f$ in terms of the admissible monomials, we get
\begin{align*}
\rho_4(f)&+f \equiv_{(3)|^{s}|(2)|^{t-1}} \gamma_{\{1,4\}}a_{s,1} + \gamma_{\{1,5\}}a_{s,2} + \gamma_{\{1,5\}}a_{s,3} + \gamma_{\{1,2\}}a_{s,5} + \gamma_{\{1,2\}}a_{s,6}\\ 
& + \gamma_{\{2,6\}}a_{s,13} + \gamma_{\{2,6\}}a_{s,14} + \gamma_{\{2,11\}}a_{s,15} + \gamma_{\{2,3\}}a_{s,20} + \gamma_{\{2,3\}}a_{s,21}\\ 
& + \gamma_{\{3,13\}}a_{s,25} + \gamma_{\{4,5,6,7\}}a_{s,30} + \gamma_{\{4,5,6,7\}}a_{s,31} + \gamma_{\{1,2,4,5,6,7\}}a_{s,35}\\ 
& + \gamma_{\{1,2,4,5,6,7\}}a_{s,36} + \gamma_{\{5,10\}}a_{s,39} + \gamma_{\{5,10\}}a_{s,40} + \gamma_{\{1,2,4,6\}}a_{s,41}\\ 
& + \gamma_{\{1,2,4,6\}}a_{s,42} + \gamma_{\{5,11\}}a_{s,43} + \gamma_{\{6,11,12,13\}}a_{s,49} + \gamma_{\{6,7\}}a_{s,51} + \gamma_{\{6,7\}}a_{s,52}\\ 
& + \gamma_{\{6,11,12,13\}}a_{s,53} + \gamma_{\{6,13\}}a_{s,56} + \gamma_{\{5,6,7,8\}}a_{s,61} + \gamma_{\{7,9\}}a_{s,65}\\ 
& + \gamma_{\{2,3,6,11,12,13\}}a_{s,67} + \gamma_{\{2,3,6,11,12,13\}}a_{s,68} + \gamma_{\{2,3,6,13\}}a_{s,70} + \gamma_{\{2,3,6,13\}}a_{s,71}\\ 
& + \gamma_{\{12,13\}}a_{s,74} + \gamma_{\{12,13\}}a_{s,76} + \gamma_{\{2,3,5,6,7,8\}}a_{s,79} + \gamma_{\{2,3,6,13\}}a_{s,82}\\ 
& + \gamma_{\{2,3,6,13\}}a_{s,83} + \gamma_{\{6,7,9,13\}}a_{s,85} + \gamma_{\{6,7,12,13\}}a_{s,86} + \gamma_{\{6,7,9,13\}}a_{s,87}\\ 
& + \gamma_{\{6,7,12,13\}}a_{s,89} + \gamma_{\{9,11,12,13\}}a_{s,92} + \gamma_{\{1,2,4,6\}}a_{s,97} + \gamma_{\{1,2,4,6\}}a_{s,98}\\ 
& + \gamma_{\{7,8,11,13\}}a_{s,99} + \gamma_{\{2,3,5,6,7,11,12,13\}}a_{s,101} + \gamma_{\{6,7,12,13\}}a_{s,103}\\ 
& + \gamma_{\{6,7,12,13\}}a_{s,104} + \gamma_{\{2,3,5,6,7,11,12,13\}}a_{s,105} \equiv_{(3)|^{s}|(2)|^{t-1}} 0.
\end{align*}
From this equality it implies $\gamma_u = \gamma_1$ for $1 \leqslant u \leqslant 13$. Therefore, 
$$f \equiv_{(3)|^s|(2)|^{t-1}} \gamma_1\Big(\sum_{1 \leqslant u \leqslant 13}p_{s,u}\Big) = \gamma_1\Big(\sum_{1 \leqslant j \leqslant 105}a_{s,j}\Big)= \gamma_1\bar \xi_{t,s}.$$
The proposition is completely proved.
\end{proof}	

\subsection{Proof of Theorem \ref{thm1} for the case $t = 1$}\

\medskip
Recall that for $t = 1$, $d_{s,1} = 2^{s+1} + 2^{s} - 3$ and 
$$(QP_4)_{d_{s,1}} \cong QP_4((3)|^s) \bigoplus QP_4((3)|^{s-1}|(1)|^2).$$

A basis of $QP_4((3)|^{s-1}|(1)|^2)$ is the set of all the classes represented by the admissible monomials $\widetilde a_{s,j} := \widetilde a_{1,s,j}$ which are determined as follows:

\medskip
For $s \geqslant 2$,

\medskip
\centerline{\begin{tabular}{lll} 
$\widetilde a_{s,1} = x_2^{2^{s-1}-1}x_3^{2^{s-1}-1}x_4^{2^{s+1}-1}$&$\widetilde a_{s,2} = x_2^{2^{s-1}-1}x_3^{2^{s+1}-1}x_4^{2^{s-1}-1}$\cr  $\widetilde a_{s,3} = x_2^{2^{s+1}-1}x_3^{2^{s-1}-1}x_4^{2^{s-1}-1}$&$\widetilde a_{s,4} = x_1^{2^{s-1}-1}x_3^{2^{s-1}-1}x_4^{2^{s+1}-1}$\cr  $\widetilde a_{s,5} = x_1^{2^{s-1}-1}x_3^{2^{s+1}-1}x_4^{2^{s-1}-1}$&$\widetilde a_{s,6} = x_1^{2^{s-1}-1}x_2^{2^{s-1}-1}x_4^{2^{s+1}-1}$\cr  $\widetilde a_{s,7} = x_1^{2^{s-1}-1}x_2^{2^{s-1}-1}x_3^{2^{s+1}-1}$&$\widetilde a_{s,8} = x_1^{2^{s-1}-1}x_2^{2^{s+1}-1}x_4^{2^{s-1}-1}$\cr  $\widetilde a_{s,9} = x_1^{2^{s-1}-1}x_2^{2^{s+1}-1}x_3^{2^{s-1}-1}$&$\widetilde a_{s,10} = x_1^{2^{s+1}-1}x_3^{2^{s-1}-1}x_4^{2^{s-1}-1}$\cr  $\widetilde a_{s,11} = x_1^{2^{s+1}-1}x_2^{2^{s-1}-1}x_4^{2^{s-1}-1}$&$\widetilde a_{s,12} = x_1^{2^{s+1}-1}x_2^{2^{s-1}-1}x_3^{2^{s-1}-1}$\cr  $\widetilde a_{s,13} = x_2^{2^{s-1}-1}x_3^{2^{s}-1}x_4^{3.2^{s-1}-1}$&$\widetilde a_{s,14} = x_2^{2^{s}-1}x_3^{2^{s-1}-1}x_4^{3.2^{s-1}-1}$\cr  $\widetilde a_{s,15} = x_2^{2^{s}-1}x_3^{3.2^{s-1}-1}x_4^{2^{s-1}-1}$&$\widetilde a_{s,16} = x_1^{2^{s-1}-1}x_3^{2^{s}-1}x_4^{3.2^{s-1}-1}$\cr  $\widetilde a_{s,17} = x_1^{2^{s-1}-1}x_2^{2^{s}-1}x_4^{3.2^{s-1}-1}$&$\widetilde a_{s,18} = x_1^{2^{s-1}-1}x_2^{2^{s}-1}x_3^{3.2^{s-1}-1}$\cr  $\widetilde a_{s,19} = x_1^{2^{s}-1}x_3^{2^{s-1}-1}x_4^{3.2^{s-1}-1}$&$\widetilde a_{s,20} = x_1^{2^{s}-1}x_3^{3.2^{s-1}-1}x_4^{2^{s-1}-1}$\cr  $\widetilde a_{s,21} = x_1^{2^{s}-1}x_2^{2^{s-1}-1}x_4^{3.2^{s-1}-1}$&$\widetilde a_{s,22} = x_1^{2^{s}-1}x_2^{2^{s-1}-1}x_3^{3.2^{s-1}-1}$\cr  $\widetilde a_{s,23} = x_1^{2^{s}-1}x_2^{3.2^{s-1}-1}x_4^{2^{s-1}-1}$&$\widetilde a_{s,24} = x_1^{2^{s}-1}x_2^{3.2^{s-1}-1}x_3^{2^{s-1}-1}$\cr  $\widetilde a_{s,25} = x_1x_2^{2^{s-1}-1}x_3^{2^{s-1}-1}x_4^{2^{s+1}-2}$&$\widetilde a_{s,26} = x_1x_2^{2^{s-1}-1}x_3^{2^{s+1}-2}x_4^{2^{s-1}-1}$\cr  $\widetilde a_{s,27} = x_1x_2^{2^{s+1}-2}x_3^{2^{s-1}-1}x_4^{2^{s-1}-1}$&$\widetilde a_{s,28} = x_1x_2^{2^{s-1}-1}x_3^{2^{s}-2}x_4^{3.2^{s-1}-1}$\cr  $\widetilde a_{s,29} = x_1x_2^{2^{s}-2}x_3^{2^{s-1}-1}x_4^{3.2^{s-1}-1}$&$\widetilde a_{s,30} = x_1x_2^{2^{s}-2}x_3^{3.2^{s-1}-1}x_4^{2^{s-1}-1}$\cr   
\end{tabular}}

\medskip
For $s = 2$,

\medskip
\centerline{\begin{tabular}{lll}   	 
$\widetilde a_{2,31} =  x_1x_2x_3^{3}x_4^{4}$& $\widetilde a_{2,32} =  x_1x_2^{3}x_3x_4^{4}$& $\widetilde a_{2,33} =  x_1x_2^{3}x_3^{4}x_4$\cr  $\widetilde a_{2,34} =  x_1^{3}x_2x_3x_4^{4}$& $\widetilde a_{2,35} =  x_1^{3}x_2x_3^{4}x_4$& $\widetilde a_{2,36} =  x_1^{3}x_2^{4}x_3x_4$\cr   
\end{tabular}}

\medskip
For $s \geqslant 3$,

\medskip
\centerline{\begin{tabular}{lll}
$\widetilde a_{s,31} = x_1x_2^{2^{s-1}-2}x_3^{2^{s-1}-1}x_4^{2^{s+1}-1}$&$\widetilde a_{s,32} = x_1x_2^{2^{s-1}-2}x_3^{2^{s+1}-1}x_4^{2^{s-1}-1}$\cr  $\widetilde a_{s,33} = x_1x_2^{2^{s-1}-1}x_3^{2^{s-1}-2}x_4^{2^{s+1}-1}$&$\widetilde a_{s,34} = x_1x_2^{2^{s-1}-1}x_3^{2^{s+1}-1}x_4^{2^{s-1}-2}$\cr  $\widetilde a_{s,35} = x_1x_2^{2^{s+1}-1}x_3^{2^{s-1}-2}x_4^{2^{s-1}-1}$&$\widetilde a_{s,36} = x_1x_2^{2^{s+1}-1}x_3^{2^{s-1}-1}x_4^{2^{s-1}-2}$\cr  $\widetilde a_{s,37} = x_1^{2^{s-1}-1}x_2x_3^{2^{s-1}-2}x_4^{2^{s+1}-1}$&$\widetilde a_{s,38} = x_1^{2^{s-1}-1}x_2x_3^{2^{s+1}-1}x_4^{2^{s-1}-2}$\cr  $\widetilde a_{s,39} = x_1^{2^{s-1}-1}x_2^{2^{s+1}-1}x_3x_4^{2^{s-1}-2}$&$\widetilde a_{s,40} = x_1^{2^{s+1}-1}x_2x_3^{2^{s-1}-2}x_4^{2^{s-1}-1}$\cr  $\widetilde a_{s,41} = x_1^{2^{s+1}-1}x_2x_3^{2^{s-1}-1}x_4^{2^{s-1}-2}$&$\widetilde a_{s,42} = x_1^{2^{s+1}-1}x_2^{2^{s-1}-1}x_3x_4^{2^{s-1}-2}$\cr  $\widetilde a_{s,43} = x_1x_2^{2^{s-1}-2}x_3^{2^{s}-1}x_4^{3.2^{s-1}-1}$&$\widetilde a_{s,44} = x_1x_2^{2^{s}-1}x_3^{2^{s-1}-2}x_4^{3.2^{s-1}-1}$\cr  $\widetilde a_{s,45} = x_1x_2^{2^{s}-1}x_3^{3.2^{s-1}-1}x_4^{2^{s-1}-2}$&$\widetilde a_{s,46} = x_1^{2^{s}-1}x_2x_3^{2^{s-1}-2}x_4^{3.2^{s-1}-1}$\cr  $\widetilde a_{s,47} = x_1^{2^{s}-1}x_2x_3^{3.2^{s-1}-1}x_4^{2^{s-1}-2}$&$\widetilde a_{s,48} = x_1^{2^{s}-1}x_2^{3.2^{s-1}-1}x_3x_4^{2^{s-1}-2}$\cr  $\widetilde a_{s,49} = x_1^{2^{s-1}-1}x_2x_3^{2^{s-1}-1}x_4^{2^{s+1}-2}$&$\widetilde a_{s,50} = x_1^{2^{s-1}-1}x_2x_3^{2^{s+1}-2}x_4^{2^{s-1}-1}$\cr  $\widetilde a_{s,51} = x_1^{2^{s-1}-1}x_2^{2^{s-1}-1}x_3x_4^{2^{s+1}-2}$&$\widetilde a_{s,52} = x_1^{2^{s-1}-1}x_2x_3^{2^{s}-2}x_4^{3.2^{s-1}-1}$\cr  $\widetilde a_{s,53} = x_1x_2^{2^{s-1}-1}x_3^{2^{s}-1}x_4^{3.2^{s-1}-2}$&$\widetilde a_{s,54} = x_1x_2^{2^{s}-1}x_3^{2^{s-1}-1}x_4^{3.2^{s-1}-2}$\cr  $\widetilde a_{s,55} = x_1x_2^{2^{s}-1}x_3^{3.2^{s-1}-2}x_4^{2^{s-1}-1}$&$\widetilde a_{s,56} = x_1^{2^{s-1}-1}x_2x_3^{2^{s}-1}x_4^{3.2^{s-1}-2}$\cr  $\widetilde a_{s,57} = x_1^{2^{s-1}-1}x_2^{2^{s}-1}x_3x_4^{3.2^{s-1}-2}$&$\widetilde a_{s,58} = x_1^{2^{s}-1}x_2x_3^{2^{s-1}-1}x_4^{3.2^{s-1}-2}$\cr  $\widetilde a_{s,59} = x_1^{2^{s}-1}x_2x_3^{3.2^{s-1}-2}x_4^{2^{s-1}-1}$&$\widetilde a_{s,60} = x_1^{2^{s}-1}x_2^{2^{s-1}-1}x_3x_4^{3.2^{s-1}-2}$\cr  $\widetilde a_{s,61} = x_1^{3}x_2^{2^{s-1}-1}x_3^{2^{s+1}-3}x_4^{2^{s-1}-2}$&$\widetilde a_{s,62} = x_1^{3}x_2^{2^{s+1}-3}x_3^{2^{s-1}-2}x_4^{2^{s-1}-1}$\cr  $\widetilde a_{s,63} = x_1^{3}x_2^{2^{s+1}-3}x_3^{2^{s-1}-1}x_4^{2^{s-1}-2}$&$\widetilde a_{s,64} = x_1^{3}x_2^{2^{s}-3}x_3^{2^{s-1}-2}x_4^{3.2^{s-1}-1}$\cr  $\widetilde a_{s,65} = x_1^{3}x_2^{2^{s}-3}x_3^{3.2^{s-1}-1}x_4^{2^{s-1}-2}$&$\widetilde a_{s,66} = x_1^{3}x_2^{2^{s}-1}x_3^{3.2^{s-1}-3}x_4^{2^{s-1}-2}$\cr  $\widetilde a_{s,67} = x_1^{2^{s}-1}x_2^{3}x_3^{3.2^{s-1}-3}x_4^{2^{s-1}-2}$&$\widetilde a_{s,68} = x_1^{3}x_2^{2^{s}-3}x_3^{2^{s-1}-1}x_4^{3.2^{s-1}-2}$\cr  $\widetilde a_{s,69} = x_1^{3}x_2^{2^{s}-3}x_3^{3.2^{s-1}-2}x_4^{2^{s-1}-1}$&$\widetilde a_{s,70} = x_1^{2^{s-1}-1}x_2^{3}x_3^{2^{s}-3}x_4^{3.2^{s-1}-2}$\cr  
\end{tabular}}

\medskip
For $s = 3$,

\medskip
\centerline{\begin{tabular}{lll}  
$\widetilde a_{3,71} =  x_1^{3}x_2^{3}x_3^{3}x_4^{12}$& $\widetilde a_{3,72} =  x_1^{3}x_2^{3}x_3^{12}x_4^{3}$& $\widetilde a_{3,73} =  x_1^{3}x_2^{3}x_3^{4}x_4^{11}$\cr  $\widetilde a_{3,74} =  x_1^{3}x_2^{3}x_3^{7}x_4^{8}$& $\widetilde a_{3,75} =  x_1^{3}x_2^{7}x_3^{3}x_4^{8}$& $\widetilde a_{3,76} =  x_1^{3}x_2^{7}x_3^{8}x_4^{3}$\cr  $\widetilde a_{3,77} =  x_1^{7}x_2^{3}x_3^{3}x_4^{8}$& $\widetilde a_{3,78} =  x_1^{7}x_2^{3}x_3^{8}x_4^{3}$& $\widetilde a_{3,79} =  x_1^{7}x_2^{9}x_3^{2}x_4^{3}$\cr  $\widetilde a_{3,80} =  x_1^{7}x_2^{9}x_3^{3}x_4^{2}$& \cr  
\end{tabular}}

\medskip
For $s \geqslant 4$,

\medskip
\centerline{\begin{tabular}{lll}   	 
$\widetilde a_{s,71} = x_1^{3}x_2^{2^{s-1}-3}x_3^{2^{s-1}-2}x_4^{2^{s+1}-1}$&$\widetilde a_{s,72} = x_1^{3}x_2^{2^{s-1}-3}x_3^{2^{s+1}-1}x_4^{2^{s-1}-2}$\cr  $\widetilde a_{s,73} = x_1^{3}x_2^{2^{s+1}-1}x_3^{2^{s-1}-3}x_4^{2^{s-1}-2}$&$\widetilde a_{s,74} = x_1^{2^{s+1}-1}x_2^{3}x_3^{2^{s-1}-3}x_4^{2^{s-1}-2}$\cr  $\widetilde a_{s,75} = x_1^{3}x_2^{2^{s-1}-3}x_3^{2^{s-1}-1}x_4^{2^{s+1}-2}$&$\widetilde a_{s,76} = x_1^{3}x_2^{2^{s-1}-3}x_3^{2^{s+1}-2}x_4^{2^{s-1}-1}$\cr  $\widetilde a_{s,77} = x_1^{3}x_2^{2^{s-1}-1}x_3^{2^{s-1}-3}x_4^{2^{s+1}-2}$&$\widetilde a_{s,78} = x_1^{2^{s-1}-1}x_2^{3}x_3^{2^{s-1}-3}x_4^{2^{s+1}-2}$\cr  $\widetilde a_{s,79} = x_1^{3}x_2^{2^{s-1}-3}x_3^{2^{s}-2}x_4^{3.2^{s-1}-1}$&$\widetilde a_{s,80} = x_1^{3}x_2^{2^{s-1}-3}x_3^{2^{s}-1}x_4^{3.2^{s-1}-2}$\cr  $\widetilde a_{s,81} = x_1^{3}x_2^{2^{s}-1}x_3^{2^{s-1}-3}x_4^{3.2^{s-1}-2}$&$\widetilde a_{s,82} = x_1^{2^{s}-1}x_2^{3}x_3^{2^{s-1}-3}x_4^{3.2^{s-1}-2}$\cr  $\widetilde a_{s,83} = x_1^{2^{s-1}-1}x_2^{3}x_3^{2^{s+1}-3}x_4^{2^{s-1}-2}$&$\widetilde a_{s,84} = x_1^{3}x_2^{2^{s-1}-1}x_3^{2^{s}-3}x_4^{3.2^{s-1}-2}$\cr  $\widetilde a_{s,85} = x_1^{7}x_2^{2^{s+1}-5}x_3^{2^{s-1}-3}x_4^{2^{s-1}-2}$&$\widetilde a_{s,86} = x_1^{7}x_2^{2^{s}-5}x_3^{2^{s-1}-3}x_4^{3.2^{s-1}-2}$\cr  $\widetilde a_{s,87} = x_1^{7}x_2^{2^{s}-5}x_3^{3.2^{s-1}-3}x_4^{2^{s-1}-2}$&\cr  
\end{tabular}}

\medskip
For $s = 4$, $\widetilde a_{4,88} =  x_1^{7}x_2^{7}x_3^{7}x_4^{24}$, $\widetilde a_{4,89} =  x_1^{7}x_2^{7}x_3^{9}x_4^{22}$, $\widetilde a_{4,90} =  x_1^{7}x_2^{7}x_3^{25}x_4^{6}$

\medskip
For $s \geqslant 5$,

\medskip
\centerline{\begin{tabular}{lll}   	 
$\widetilde a_{s,88} = x_1^{7}x_2^{2^{s-1}-5}x_3^{2^{s-1}-3}x_4^{2^{s+1}-2}$&$\widetilde a_{s,89} = x_1^{7}x_2^{2^{s-1}-5}x_3^{2^{s+1}-3}x_4^{2^{s-1}-2}$\cr  $\widetilde a_{s,90} = x_1^{7}x_2^{2^{s-1}-5}x_3^{2^{s}-3}x_4^{3.2^{s-1}-2}$&\cr  
\end{tabular}}

\medskip
We prove the following theorem which is equivalent to Theorem \ref{thm1} for $t = 1$.
\begin{thms}\label{mdct1} For $s \geqslant 2$ and $d_{s,1} = 2^{s+1}+2^s-3$, we have
$$(QP_4)_{d_{s,1}}^{GL_4} = \langle [\xi_{1,s}]\rangle,$$ 
where
$$\xi_{1,s} = \begin{cases}
\sum_{j \in \{25,\, 26,\, 27,\, 34,\, 35,\, 36\}}\widetilde a_{2,j}, &\mbox{if } s = 2,\\
0, &\mbox{if } s = 3,\\
\bar \xi_{1,4} + \sum_{1 \leqslant j \leqslant 90,\, j \notin \{25,\, 49,\, 51,\, 61,\, 70,\, 75,\, 77,\, 78,\, 83,\, 84\}}\widetilde a_{4,j}, &\mbox{if } s = 4,\\
\bar \xi_{1,s} + \sum_{1 \leqslant j \leqslant 90}\widetilde a_{s,j}, &\mbox{if } s \geqslant 5.
\end{cases} $$
Here, for $s \geqslant 4$, $\bar \xi_{1,s}$ is defined as in Proposition \ref{mdt11}.
\end{thms}
\begin{rems} In \cite[Page 458]{pp25}, the author stated that $(QP_4)_{d_{s,1}}^{GL_4} = \langle [\bar \xi_{1,s}] \rangle$ but it is false because $\bar \xi_{1,s}$ is not an $GL_4$-invariant in $(QP_4)_{d_{s,1}}$. The case $s = 2$ of this theorem is also presented in \cite[Prop. 4.1.4]{pp25} but this result is due to Singer \cite{si1}. In \cite[Page 1534]{p231}, the author stated that $\dim(QP_4)_{d_{3,1}}^{GL_4} = 1$ but this is false.
\end{rems}
We observe that
\begin{equation}\label{ctct1}
\begin{cases}
[\Sigma_4(\widetilde a_{s,1})] = \{[\widetilde a_{s,j}]: 1 \leqslant j \leqslant 12 \},\ [\Sigma_4(\widetilde a_{s,1})]^{\Sigma_4} = \langle[\widetilde p_{s,1}]\rangle\\
[\Sigma_4(\widetilde a_{s,13})] = \{[\widetilde a_{s,j}]: 13 \leqslant j \leqslant 24\},\ [\Sigma_4(\widetilde a_{s,13})]^{\Sigma_4} = \langle[\widetilde p_{s,2}]\rangle,
\end{cases}
\end{equation}
where 
$$\widetilde p_{s,1} = \widetilde p_{1,s,1} = \sum_{1 \leqslant j \leqslant 12}\widetilde a_{s,j}, \ \widetilde p_{s,2} = \widetilde p_{1,s,2} = \sum_{13 \leqslant j \leqslant 24}\widetilde a_{s,j}.$$

\medskip
For $s = 2$, the space $(QP_4)_9^{GL_4}$ had been determined by Singer \cite{si1}. By Proposition \ref{mdt11}, $QP_4(3,3)^{GL_4} = 0$, so $(QP_4)_9^{GL_4} = QP_4(3,1,1)^{GL_4}$.  

\begin{lems}\label{bdct11} We have
$$QP_4(3,1,1)^{\Sigma_4} = \langle[\widetilde p_{2,u}]: 1\leqslant u \leqslant 3\rangle,$$ 
where $\widetilde p_{2,1}$, $\widetilde p_{2,2}$ are defined by \eqref{ctct1} and
\begin{align*}
\widetilde p_{2,3} &= \sum_{j \in \{25,\, 26,\, 27,\, 34,\, 35,\, 36\}}\widetilde a_{2,j} \equiv x_1x_2x_3x_4^{6} + x_1x_2x_3^{6}x_4 + x_1x_2^{6}x_3x_4 + x_1^{6}x_2x_3x_4.
\end{align*}
\end{lems}

\begin{proof} Note that $x_1^{6}x_2x_3x_4 \equiv x_1^{3}x_2x_3x_4^{4} + x_1^{3}x_2x_3^{4}x_4 + x_1^{3}x_2^{4}x_3x_4$. It is clear that there is a direct summand decomposition of $\Sigma_4$-submodules
$$QP_4(3,1,1) = [\Sigma_4(\widetilde a_{2,1})]\bigoplus [\Sigma_4(\widetilde a_{2,13})]\bigoplus QP_4^+(3,1,1).$$ 
	
By using \eqref{ctct1}, we need only to prove $QP_4^+(3,1,1)^{\Sigma_4} = \langle[\widetilde p_{2,3}]\rangle$.
	
If $[f] \in QP_4^+(3,1,1)^{\Sigma_4}$ with $f \in QP_4^+(3,1,1) = \langle [\widetilde a_{2,j}]: 25 \leqslant j \leqslant 36\rangle$, then $f \equiv \sum_{25 \leqslant j \leqslant 36}\widetilde a_{2,j}$. By a simple computation, we have
\begin{align*}
\rho_1(f) + f &\equiv \gamma_{29}\widetilde a_{s,25} + \gamma_{30}\widetilde a_{s,26} + \gamma_{\{27,36\}}\widetilde a_{s,27} + \gamma_{29}\widetilde a_{s,28} + \gamma_{30}\widetilde a_{s,31}\\ 
&\quad + \gamma_{\{32,34,36\}}\widetilde a_{s,32} + \gamma_{\{33,35,36\}}\widetilde a_{s,33} + \gamma_{\{27,32,34\}}\widetilde a_{s,34}\\ 
&\quad + \gamma_{\{27,33,35\}}\widetilde a_{s,35} + \gamma_{\{27,36\}}\widetilde a_{s,36} \equiv 0,\\
\rho_2(f) + f &\equiv \gamma_{\{26,27\}}\widetilde a_{s,26} + \gamma_{\{26,27\}}\widetilde a_{s,27} + \gamma_{\{28,29\}}\widetilde a_{s,28} + \gamma_{\{28,29\}}\widetilde a_{s,29}\\ 
&\quad  + \gamma_{\{30,33\}}\widetilde a_{s,30} + \gamma_{\{31,32\}}\widetilde a_{s,31} + \gamma_{\{31,32\}}\widetilde a_{s,32} + \gamma_{\{30,33\}}\widetilde a_{s,33}\\ 
&\quad  + \gamma_{\{35,36\}}\widetilde a_{s,35} + \gamma_{\{35,36\}}\widetilde a_{s,36} \equiv 0,\\  
\rho_3(f) + f &\equiv \gamma_{\{25,26\}}\widetilde a_{s,25} + \gamma_{\{25,26\}}\widetilde a_{s,26} + \gamma_{\{28,31\}}\widetilde a_{s,28} + \gamma_{\{29,30\}}\widetilde a_{s,29}\\ 
&\quad + \gamma_{\{29,30\}}\widetilde a_{s,30} + \gamma_{\{28,31\}}\widetilde a_{s,31} + \gamma_{\{32,33\}}\widetilde a_{s,32} + \gamma_{\{32,33\}}\widetilde a_{s,33}\\ 
&\quad + \gamma_{\{34,35\}}\widetilde a_{s,34} + \gamma_{\{34,35\}}\widetilde a_{s,35} \equiv 0.
\end{align*}
These equalities imply $\gamma_j = 0$ for $28 \leqslant j \leqslant 33$ and $\gamma_j = \gamma_{25}$ for $j \in \{$26, 27, 34, 35, 36$\}$. Hence $ \equiv \gamma_{25}\widetilde p_{2,3}$. The lemma is proved.
\end{proof}
It is easy to see that for $s \geqslant 3$, we have
\begin{equation}\label{ctct12}
\begin{cases}
[\Sigma_4(\widetilde a_{s,31})] = \{[\widetilde a_{s,j}]: 31 \leqslant j \leqslant 42 \},\ [\Sigma_4(\widetilde a_{s,31})]^{\Sigma_4} = \langle[\widetilde p_{s,3}]\rangle\\
[\Sigma_4(\widetilde a_{s,43})] = \{[\widetilde a_{s,j}]: 43 \leqslant j \leqslant 48\},\ [\Sigma_4(\widetilde a_{s,43})]^{\Sigma_4} = \langle[\widetilde p_{s,4}]\rangle,
\end{cases}
\end{equation}
where 
$$\widetilde p_{s,3} = \widetilde p_{1,s,3} = \sum_{31 \leqslant j \leqslant 42}\widetilde a_{s,j}, \ \widetilde p_{s,4} = \widetilde p_{1,s,4} = \sum_{43 \leqslant j \leqslant 48}\widetilde a_{s,j}.$$

\begin{lems}\label{bdct12} For $s =3$, we have
$$QP_4((3)|^2|(1)|^2)^{\Sigma_4} = \langle \{[\widetilde p_{3,u}] : 1 \leqslant u \leqslant 6\}\rangle, $$ 
where  $\widetilde p_{3,u}$, $ 1 \leqslant u \leqslant 4$ are determined by \eqref{ctct1}, \eqref{ctct12} and  
\begin{align*}
\widetilde p_{3,5} &= \sum_{j \in \{27,\, 54,\, 55,\, 68,\, 69,\, 70,\, 71,\, 72,\, 75,\, 76\}}\widetilde a_{3,j},\\
\widetilde p_{3,6} &= \sum_{j \in \{25,\, 26,\, 27,\, 49,\, 50,\, 51,\, 58,\, 59,\, 60,\, 61,\, 62,\, 63,\, 67,\, 68,\, 69,\, 70,\, 79,\, 80\}}\widetilde a_{3,j}.
\end{align*}
\end{lems}
\begin{proof} We observe that there is a direct summand decomposition of $\Sigma_4$-submodules
\begin{align*}QP_4((3)|^2|(1)|^2) = [\Sigma_4(\widetilde a_{3,1})]&\bigoplus [\Sigma_4(\widetilde a_{3,13})]\bigoplus [\Sigma_4(\widetilde a_{3,31})]\\ &\hskip 1.5cm\bigoplus [\Sigma_4(\widetilde a_{3,43})]\bigoplus \widetilde {\mathcal U}_{1,3},\end{align*} 
where  $\widetilde {\mathcal U}_{1,3} = \langle [\widetilde a_{3,j}] : j \in \{25, \ldots ,\, 30,\, 49,\ldots 80\}\rangle $. From \eqref{ctct1} and \eqref{ctct12} we need only to prove  $\widetilde {\mathcal U}_{1,3}^{\Sigma_4} = \langle[\widetilde p_{3,5}],\, [\widetilde p_{3,6}]\rangle$. 

We see that $\widetilde p_{3,5}$, $\widetilde p_{3,6}$ are $\Sigma_4$-invariants and the leading monomials respectively are $\widetilde a_{3,76}$ and $\widetilde a_{3,80}$. Hence, if $[f] \in \widetilde{\mathcal U}_{1,3}^{\Sigma_4}$ with $f \in P_4$, then there are $\gamma_5,\, \gamma_6 \in \mathbb F_2$ such that
$$g := f + \gamma_5\widetilde p_{3,5} + \gamma_6\widetilde p_{3,6} \equiv  \sum_{j\in \{25,\ldots,30,\, 49, \ldots, 80\}\setminus\{76,\, 80\}}\gamma_j\widetilde a_{3,j}.$$
Then, $[g]$ is also an $\Sigma_4$-invariant. A direct computation shows
\begin{align*}
\rho_1(g)& + g \equiv \gamma_{\{25,29,49\}}\widetilde a_{3,25} + \gamma_{\{26,30,50\}}\widetilde a_{3,26} + \gamma_{\{28,29,52\}}\widetilde a_{3,28}\\ 
& + \gamma_{\{25,29,49\}}\widetilde a_{3,49} + \gamma_{\{26,30,50\}}\widetilde a_{3,50} + \gamma_{64}\widetilde a_{3,51} + \gamma_{\{28,29,52\}}\widetilde a_{3,52}\\ 
& + \gamma_{\{30,53,56\}}\widetilde a_{3,53} + \gamma_{\{27,54,58\}}\widetilde a_{3,54} + \gamma_{\{27,55,59\}}\widetilde a_{3,55} + \gamma_{\{30,53,56\}}\widetilde a_{3,56}\\ 
& + \gamma_{\{57,60,79\}}\widetilde a_{3,57} + \gamma_{\{27,54,58\}}\widetilde a_{3,58} + \gamma_{\{27,55,59\}}\widetilde a_{3,59} + \gamma_{\{57,60,62\}}\widetilde a_{3,60}\\ 
& + \gamma_{65}\widetilde a_{3,61} + \gamma_{\{62,79\}}\widetilde a_{3,62} + \gamma_{63}\widetilde a_{3,63} + \gamma_{\{66,67\}}\widetilde a_{3,66} + \gamma_{\{63,66,67\}}\widetilde a_{3,67}\\ 
& + \gamma_{\{68,69\}}\widetilde a_{3,70} + \gamma_{\{27,29,68\}}\widetilde a_{3,71} + \gamma_{\{27,30,69\}}\widetilde a_{3,72} + \gamma_{\{29,64\}}\widetilde a_{3,73}\\ 
& + \gamma_{\{30,65\}}\widetilde a_{3,74} + \gamma_{\{27,75,77\}}\widetilde a_{3,75} + \gamma_{\{27,78,79\}}\widetilde a_{3,76} + \gamma_{\{27,63,75,77\}}\widetilde a_{3,77}\\ 
& + \gamma_{\{27,62,78\}}\widetilde a_{3,78} + \gamma_{\{62,79\}}\widetilde a_{3,79} + \gamma_{63}\widetilde a_{3,80} \equiv 0,\\
\rho_2(g) &+ g \equiv \gamma_{73}\widetilde a_{3,25} + \gamma_{\{26,27\}}\widetilde a_{3,26} + \gamma_{\{26,27,72,78\}}\widetilde a_{3,27} + \gamma_{\{28,29,73\}}\widetilde a_{3,28}\\ 
& + \gamma_{\{28,29,73\}}\widetilde a_{3,29} + \gamma_{\{30,55\}}\widetilde a_{3,30} + \gamma_{\{49,51\}}\widetilde a_{3,49} + \gamma_{\{50,62\}}\widetilde a_{3,50}\\ 
& + \gamma_{\{49,51\}}\widetilde a_{3,51} + \gamma_{\{52,64\}}\widetilde a_{3,52} + \gamma_{\{53,54\}}\widetilde a_{3,53} + \gamma_{\{53,54,72,78\}}\widetilde a_{3,54}\\ 
& + \gamma_{\{30,55,72,78\}}\widetilde a_{3,55} + \gamma_{\{56,57\}}\widetilde a_{3,56} + \gamma_{\{56,57\}}\widetilde a_{3,57} + \gamma_{\{58,60\}}\widetilde a_{3,58}\\ 
& + \gamma_{\{59,79\}}\widetilde a_{3,59} + \gamma_{\{58,60\}}\widetilde a_{3,60} + \gamma_{\{61,63\}}\widetilde a_{3,61} + \gamma_{\{50,62\}}\widetilde a_{3,62}\\ 
& + \gamma_{\{61,63\}}\widetilde a_{3,63} + \gamma_{\{52,64\}}\widetilde a_{3,64} + \gamma_{\{65,66\}}\widetilde a_{3,65} + \gamma_{\{65,66\}}\widetilde a_{3,66}\\ 
& + \gamma_{67}\widetilde a_{3,67} + \gamma_{\{68,70\}}\widetilde a_{3,68} + \gamma_{\{68,70\}}\widetilde a_{3,70} + \gamma_{\{73,78\}}\widetilde a_{3,71} + \gamma_{\{72,78\}}\widetilde a_{3,72}\\ 
& + \gamma_{\{74,75\}}\widetilde a_{3,74} + \gamma_{\{72,74,75,78\}}\widetilde a_{3,75} + \gamma_{\{72,78\}}\widetilde a_{3,76} + \gamma_{78}\widetilde a_{3,77}\\ 
& + \gamma_{\{59,79\}}\widetilde a_{3,79} + \gamma_{67}\widetilde a_{3,80} \equiv 0,\\  
\rho_3(g) &+ g \equiv \gamma_{\{25,26\}}\widetilde a_{3,25} + \gamma_{\{25,26\}}\widetilde a_{3,26} + \gamma_{\{28,53\}}\widetilde a_{3,28} + \gamma_{\{29,30\}}\widetilde a_{3,29}\\ 
& + \gamma_{\{29,30\}}\widetilde a_{3,30} + \gamma_{\{49,50\}}\widetilde a_{3,49} + \gamma_{\{49,50\}}\widetilde a_{3,50} + \gamma_{\{51,61\}}\widetilde a_{3,51}\\ 
& + \gamma_{\{52,56\}}\widetilde a_{3,52} + \gamma_{\{28,53\}}\widetilde a_{3,53} + \gamma_{\{54,55\}}\widetilde a_{3,54} + \gamma_{\{54,55\}}\widetilde a_{3,55}\\ 
& + \gamma_{\{52,56\}}\widetilde a_{3,56} + \gamma_{\{57,66\}}\widetilde a_{3,57} + \gamma_{\{58,59\}}\widetilde a_{3,58} + \gamma_{\{58,59\}}\widetilde a_{3,59}\\ 
& + \gamma_{\{60,67\}}\widetilde a_{3,60} + \gamma_{\{51,61\}}\widetilde a_{3,61} + \gamma_{\{62,63\}}\widetilde a_{3,62} + \gamma_{\{62,63\}}\widetilde a_{3,63}\\ 
& + \gamma_{\{64,65\}}\widetilde a_{3,64} + \gamma_{\{64,65\}}\widetilde a_{3,65} + \gamma_{\{57,66\}}\widetilde a_{3,66} + \gamma_{\{60,67\}}\widetilde a_{3,67}\\ 
& + \gamma_{\{68,69\}}\widetilde a_{3,68} + \gamma_{\{68,69\}}\widetilde a_{3,69} + \gamma_{\{71,72\}}\widetilde a_{3,71} + \gamma_{\{71,72\}}\widetilde a_{3,72}\\ 
& + \gamma_{\{73,74\}}\widetilde a_{3,73} + \gamma_{\{73,74\}}\widetilde a_{3,74} + \gamma_{75}\widetilde a_{3,75} + \gamma_{75}\widetilde a_{3,76} + \gamma_{\{77,78\}}\widetilde a_{3,77}\\ 
& + \gamma_{\{77,78\}}\widetilde a_{3,78} + \gamma_{79}\widetilde a_{3,79} + \gamma_{79}\widetilde a_{3,80} \equiv 0.
\end{align*}
From these equalities we get $\gamma_{j} = 0$ for all $j \in\{25,\ldots,30,\, 49,\ldots,80\},, j \ne 76,\, 80$. Hence, we get $g \equiv 0$ and
$f \equiv \gamma_5\widetilde p_{3,5} + \gamma_6\widetilde p_{3,6}.$ The lemma follows.
\end{proof}

\begin{lems}\label{bdct13} For $s \geqslant 4$, we have
$$QP_4((3)|^{s-1}|(1)|^2)^{\Sigma_4} = \langle \{[\widetilde p_{s,u}] : 1 \leqslant u \leqslant 9\}\rangle, $$ 
where  $\widetilde p_{s,u}$, $ 1 \leqslant u \leqslant 4$ are determined by \eqref{ctct1}, \eqref{ctct12} and  
\begin{align*}
\widetilde p_{s,5} &= \mbox{$\sum_{71 \leqslant j \leqslant 74}\widetilde a_{3,j}$},\\
\widetilde p_{3,6} &= \begin{cases}\sum_{j \in \{25,\, 26,\, 27,\, 49,\, 50,\, 68,\, 69,\, 75,\, 76,\, 89\}}\widetilde a_{s,j}, &\mbox{if } s = 4,\\
\sum_{j \in \{25,\, 26,\, 27,\, 49,\, 50,\, 51,\, 68,\, 69,\, 70,\, 75,\, 76,\, 77,\, 78,\, 84,\, 88,\, 90\}}\widetilde a_{s,j}, &\mbox{if } s \geqslant 5, \end{cases}\\
\widetilde p_{3,7} &= \begin{cases}\sum_{j \in \{26,\, 27,\, 50,\, 51,\, 62,\, 63,\, 75,\, 78,\, 83,\, 85,\, 88,\, 90\}}\widetilde a_{s,j}, &\mbox{if } s = 4,\\
\sum_{j \in \{25,\, 26,\, 27,\, 49,\, 50,\, 51,\, 68,\, 69,\, 70,\, 75,\, 76,\, 77,\, 78,\, 84,\, 88,\, 90\}}\widetilde a_{s,j}, &\mbox{if } s \geqslant 5, \end{cases}\\
\widetilde p_{3,8} &= \begin{cases}\sum_{j \in \{28,\, 29,\, 30,\, 52,\, 53,\, 54,\, 55,\, 56,\, 57,\, 58,\, 59,\, 60,\, 61,\, 62,\, 63,\, 64,\, 65\}\atop\hskip5.2cm \cup\{66,\, 67,\, 83\}}\widetilde a_{s,j}, &\mbox{if } s = 4,\\
\sum_{j \in \{61,\, 62,\, 63,\, 68,\, 69,\, 70,\, 75,\, 76,\, 77,\, 78,\, 83,\, 84,\, 86,\, 87,\, 88,\, 89\}}\widetilde a_{s,j}, &\mbox{if } s \geqslant 5, \end{cases}\\
\widetilde p_{3,9} &= \begin{cases}\sum_{j \in \{61,\, 62,\, 63,\, 68,\, 69,\, 70,\, 75,\, 76,\, 77,\, 84,\, 86,\, 87\}}\widetilde a_{s,j}, &\mbox{if } s = 4,\\
\sum_{j \in \{85,\, 86,\, 87,\, 88,\, 89,\, 90\}}\widetilde a_{s,j}, &\mbox{if } s \geqslant 5. \end{cases}
\end{align*}
\end{lems}
	 
\begin{proof} By a direct computation we can easily verify that $\widetilde p_{3,u}$ is an $\Sigma_4$-invariant for $1 \leqslant u \leqslant 9$. We also observe that there is a direct summand decomposition of $\Sigma_4$-submodules
\begin{align*}QP_4((3)|^{s-1}|(1)|^2) = [\Sigma_4(\widetilde a_{s,1})]&\bigoplus [\Sigma_4(\widetilde a_{s,13})]\bigoplus [\Sigma_4(\widetilde a_{s,31})]\\ &\bigoplus [\Sigma_4(\widetilde a_{s,43})]\bigoplus [\Sigma_4(\widetilde a_{s,71})]\bigoplus \widetilde {\mathcal U}_{1,s},\end{align*} 
where  $\widetilde {\mathcal U}_{1,s} = \langle [\widetilde a_{s,j}] : j \in \{25, \ldots ,\, 30,\, 49,\ldots 70,\, 75,\ldots,90\}\rangle $. It easy to verify that $[\Sigma_4(\widetilde a_{s,71})]^{\Sigma_4} =\langle[\widetilde p_{s,5}]\rangle$. Hence, from \eqref{ctct1} and \eqref{ctct12} we need only to prove  $\widetilde {\mathcal U}_{1,s}^{\Sigma_4} = \langle[\widetilde p_{s,u}] : 6 \leqslant u \leqslant 9\rangle$. 
	
For $s =4$, the leading monomials of $\widetilde p_{4,u}$, for  $6 \leqslant u \leqslant 9\rangle$, respectively are $\widetilde a_{4,89}$, $\widetilde a_{4,85}$, $\widetilde a_{4,67}$, $\widetilde a_{4,87}$. Hence, if $[f] \in \widetilde{\mathcal U}_{1,4}^{\Sigma_4}$ with $f \in P_4$, then there are $\gamma_u \in \mathbb F_2$, $6 \leqslant u \leqslant 9$, such that
$$g := f + \sum_{6 \leqslant u \leqslant 9}\gamma_u\widetilde p_{4,u} \equiv  \sum_{j\in \{25, \ldots ,\, 30,\, 49,\ldots 70,\, 75,\ldots,90\}\setminus\{67,\, 85,\, 87,\, 89\}}\gamma_j\widetilde a_{4,j}.$$
Then, $[g]$ is also an $\Sigma_4$-invariant. By a direct computation we get
\begin{align*}
\rho_1(g) &+ g \equiv \gamma_{\{25,49,62,63,64,65,68,69\}}\widetilde a_{4,25} + \gamma_{\{26,50\}}\widetilde a_{4,26} + \gamma_{\{28,52\}}\widetilde a_{4,28}\\ 
& + \gamma_{\{25,49,62,63,64,65,68,69\}}\widetilde a_{4,49} + \gamma_{\{26,50\}}\widetilde a_{4,50} + \gamma_{\{62,64,86\}}\widetilde a_{4,51}\\ 
& + \gamma_{\{28,52\}}\widetilde a_{4,52} + \gamma_{\{53,56\}}\widetilde a_{4,53} + \gamma_{\{54,58\}}\widetilde a_{4,54} + \gamma_{\{55,59\}}\widetilde a_{4,55}\\ 
& + \gamma_{\{53,56\}}\widetilde a_{4,56} + \gamma_{\{57,60\}}\widetilde a_{4,57} + \gamma_{\{54,58\}}\widetilde a_{4,58} + \gamma_{\{55,59\}}\widetilde a_{4,59}\\ 
& + \gamma_{\{57,60\}}\widetilde a_{4,60} + \gamma_{\{61,63,65,83\}}\widetilde a_{4,61} + \gamma_{66}\widetilde a_{4,66} + \gamma_{66}\widetilde a_{4,67}\\ 
& + \gamma_{\{68,69,70,84\}}\widetilde a_{4,70} + \gamma_{\{27,29,62,64,65,69\}}\widetilde a_{4,75} + \gamma_{\{27,30,62,69\}}\widetilde a_{4,76}\\ 
& + \gamma_{\{62,64,77,78\}}\widetilde a_{4,77} + \gamma_{\{62,64,77,78\}}\widetilde a_{4,78} + \gamma_{\{29,64\}}\widetilde a_{4,79} + \gamma_{\{30,65\}}\widetilde a_{4,80}\\ 
& + \gamma_{\{81,82\}}\widetilde a_{4,81} + \gamma_{\{81,82\}}\widetilde a_{4,82} + \gamma_{\{61,63,65,83\}}\widetilde a_{4,83} + \gamma_{\{68,69,70,84\}}\widetilde a_{4,84}\\ 
& + \gamma_{\{62,63,64,65,68,69\}}\widetilde a_{4,88} + \gamma_{\{68,69,86\}}\widetilde a_{4,89} + \gamma_{\{63,65\}}\widetilde a_{4,90} \equiv 0,\\
\rho_2(g) &+ g \equiv \gamma_{\{50,52,59,90\}}\widetilde a_{4,25} + \gamma_{\{26,27\}}\widetilde a_{4,26} + \gamma_{\{26,27\}}\widetilde a_{4,27} + \gamma_{\{28,29\}}\widetilde a_{4,28}\\ 
& + \gamma_{\{28,29\}}\widetilde a_{4,29} + \gamma_{\{30,55\}}\widetilde a_{4,30} + \gamma_{\{49,50,51,52,59\}}\widetilde a_{4,49} + \gamma_{\{49,50,51,52,59\}}\widetilde a_{4,51}\\ 
& + \gamma_{\{53,54\}}\widetilde a_{4,53} + \gamma_{\{53,54\}}\widetilde a_{4,54} + \gamma_{\{30,55\}}\widetilde a_{4,55} + \gamma_{\{56,57\}}\widetilde a_{4,56} + \gamma_{\{56,57\}}\widetilde a_{4,57}\\ 
& + \gamma_{\{58,60\}}\widetilde a_{4,58} + \gamma_{\{58,60\}}\widetilde a_{4,60} + \gamma_{\{61,63,90\}}\widetilde a_{4,61} + \gamma_{\{50,59,62,76\}}\widetilde a_{4,62}\\ 
& + \gamma_{\{61,63,90\}}\widetilde a_{4,63} + \gamma_{\{52,64,79\}}\widetilde a_{4,64} + \gamma_{\{65,66\}}\widetilde a_{4,65} + \gamma_{\{65,66\}}\widetilde a_{4,66}\\ 
& + \gamma_{\{68,84\}}\widetilde a_{4,68} + \gamma_{\{70,86\}}\widetilde a_{4,70} + \gamma_{\{50,52,59,75,77\}}\widetilde a_{4,75} + \gamma_{\{50,59,62,76\}}\widetilde a_{4,76}\\ 
& + \gamma_{\{50,52,59,75,77\}}\widetilde a_{4,77} + \gamma_{59}\widetilde a_{4,78} + \gamma_{\{52,64,79\}}\widetilde a_{4,79} + \gamma_{\{80,81\}}\widetilde a_{4,80}\\ 
& + \gamma_{\{80,81\}}\widetilde a_{4,81} + \gamma_{59}\widetilde a_{4,82} + \gamma_{83}\widetilde a_{4,83} + \gamma_{\{68,84\}}\widetilde a_{4,84} + \gamma_{83}\widetilde a_{4,85}\\ 
& + \gamma_{\{70,86\}}\widetilde a_{4,86} + \gamma_{\{50,52,59,90\}}\widetilde a_{4,88} \equiv 0,\\  
\rho_3(g) &+ g \equiv \gamma_{\{25,26,88\}}\widetilde a_{4,25} + \gamma_{\{25,26,88\}}\widetilde a_{4,26} + \gamma_{\{28,53\}}\widetilde a_{4,28} + \gamma_{\{29,30\}}\widetilde a_{4,29}\\ 
& + \gamma_{\{29,30\}}\widetilde a_{4,30} + \gamma_{\{49,50,88\}}\widetilde a_{4,49} + \gamma_{\{49,50,88\}}\widetilde a_{4,50} + \gamma_{\{51,57,60,90\}}\widetilde a_{4,51}\\ 
& + \gamma_{\{52,56\}}\widetilde a_{4,52} + \gamma_{\{28,53\}}\widetilde a_{4,53} + \gamma_{\{54,55\}}\widetilde a_{4,54} + \gamma_{\{54,55\}}\widetilde a_{4,55}\\ 
& + \gamma_{\{52,56\}}\widetilde a_{4,56} + \gamma_{\{58,59\}}\widetilde a_{4,58} + \gamma_{\{58,59\}}\widetilde a_{4,59} + \gamma_{\{57,61,77\}}\widetilde a_{4,61}\\ 
& + \gamma_{\{62,63\}}\widetilde a_{4,62} + \gamma_{\{62,63\}}\widetilde a_{4,63} + \gamma_{\{64,65\}}\widetilde a_{4,64} + \gamma_{\{64,65\}}\widetilde a_{4,65}\\ 
& + \gamma_{\{57,66,81\}}\widetilde a_{4,66} + \gamma_{\{60,82\}}\widetilde a_{4,67} + \gamma_{\{68,69\}}\widetilde a_{4,68} + \gamma_{\{68,69\}}\widetilde a_{4,69}\\ 
& + \gamma_{\{75,76,88\}}\widetilde a_{4,75} + \gamma_{\{75,76,88\}}\widetilde a_{4,76} + \gamma_{\{57,61,77\}}\widetilde a_{4,77} + \gamma_{\{60,78,83\}}\widetilde a_{4,78}\\ 
& + \gamma_{\{79,80\}}\widetilde a_{4,79} + \gamma_{\{79,80\}}\widetilde a_{4,80} + \gamma_{\{57,66,81\}}\widetilde a_{4,81} + \gamma_{\{60,82\}}\widetilde a_{4,82}\\ 
& + \gamma_{\{60,78,83\}}\widetilde a_{4,83} + \gamma_{86}\widetilde a_{4,86} + \gamma_{86}\widetilde a_{4,87} + \gamma_{\{51,57,60,90\}}\widetilde a_{4,90} \equiv 0.
\end{align*}
From these equalities we obtain $\gamma_{j} = 0$ for all $j \in\{25,\ldots,30,\, 49,\ldots,70,75,\ldots,90\}\setminus\{67,\, 85,\, 87,\, 89\}$. Hence, we get $g \equiv 0$ and
$f   \equiv \sum_{6 \leqslant u \leqslant 9}\gamma_u\widetilde p_{4,u}.$ The lemma is proved for $s = 4$.

For $s \geqslant 5$, the leading monomials of $\widetilde p_{s,u}$, for  $6 \leqslant u \leqslant 9\rangle$, respectively are $\widetilde a_{s,51}$, $\widetilde a_{s,60}$, $\widetilde a_{s,83}$, $\widetilde a_{s,85}$. Hence, if $[f] \in \widetilde{\mathcal U}_{1,s}^{\Sigma_4}$ with $f \in P_4$, then there are $\gamma_u \in \mathbb F_2$, $6 \leqslant u \leqslant 9$, such that
$$g := f + \sum_{6 \leqslant u \leqslant 9}\gamma_u\widetilde p_{s,u} \equiv  \sum_{j\in \{25, \ldots ,\, 30,\, 49,\ldots 70,\, 75,\ldots,90\}\setminus\{51,\, 60,\, 83,\, 85\}}\gamma_j\widetilde a_{s,j}.$$
Then, $[g]$ is also an $\Sigma_4$-invariant. By a direct computation we get
\begin{align*}
\rho_1(g)& + g \equiv \gamma_{\{25,49\}}\widetilde a_{s,25} + \gamma_{\{26,50\}}\widetilde a_{s,26} + \gamma_{\{28,52\}}\widetilde a_{s,28} + \gamma_{\{25,49\}}\widetilde a_{s,49}\\ 
& + \gamma_{\{26,50\}}\widetilde a_{s,50} + \gamma_{\{28,52\}}\widetilde a_{s,52} + \gamma_{\{53,56\}}\widetilde a_{s,53} + \gamma_{\{54,58\}}\widetilde a_{s,54}\\ 
& + \gamma_{\{55,59\}}\widetilde a_{s,55} + \gamma_{\{53,56\}}\widetilde a_{s,56} + \gamma_{57}\widetilde a_{s,57} + \gamma_{\{54,58\}}\widetilde a_{s,58} + \gamma_{\{55,59\}}\widetilde a_{s,59}\\ 
& + \gamma_{57}\widetilde a_{s,60} + \gamma_{61}\widetilde a_{s,61} + \gamma_{\{66,67\}}\widetilde a_{s,66} + \gamma_{\{66,67\}}\widetilde a_{s,67} + \gamma_{\{70,84\}}\widetilde a_{s,70}\\ 
& + \gamma_{\{27,29,63,68\}}\widetilde a_{s,75} + \gamma_{\{27,30,62,69\}}\widetilde a_{s,76} + \gamma_{\{77,78\}}\widetilde a_{s,77} + \gamma_{\{77,78\}}\widetilde a_{s,78}\\ 
& + \gamma_{\{29,64\}}\widetilde a_{s,79} + \gamma_{\{30,65\}}\widetilde a_{s,80} + \gamma_{\{81,82\}}\widetilde a_{s,81} + \gamma_{\{81,82\}}\widetilde a_{s,82} + \gamma_{61}\widetilde a_{s,83}\\ 
& + \gamma_{\{70,84\}}\widetilde a_{s,84} + \gamma_{\{62,64,86\}}\widetilde a_{s,88} + \gamma_{\{63,65,87\}}\widetilde a_{s,89} + \gamma_{\{68,69,86,87\}}\widetilde a_{s,90} \equiv 0,\\
\rho_2(g) &+ g \equiv \gamma_{\{26,27\}}\widetilde a_{s,26} + \gamma_{\{26,27\}}\widetilde a_{s,27} + \gamma_{\{28,29\}}\widetilde a_{s,28} + \gamma_{\{28,29\}}\widetilde a_{s,29}\\ 
& + \gamma_{\{30,55\}}\widetilde a_{s,30} + \gamma_{49}\widetilde a_{s,49} + \gamma_{49}\widetilde a_{s,51} + \gamma_{\{53,54\}}\widetilde a_{s,53} + \gamma_{\{53,54\}}\widetilde a_{s,54}\\ 
& + \gamma_{\{30,55\}}\widetilde a_{s,55} + \gamma_{\{56,57\}}\widetilde a_{s,56} + \gamma_{\{56,57\}}\widetilde a_{s,57} + \gamma_{58}\widetilde a_{s,58} + \gamma_{58}\widetilde a_{s,60}\\ 
& + \gamma_{\{61,63\}}\widetilde a_{s,61} + \gamma_{\{50,59,62,76\}}\widetilde a_{s,62} + \gamma_{\{61,63\}}\widetilde a_{s,63} + \gamma_{\{52,64,79\}}\widetilde a_{s,64}\\ 
& + \gamma_{\{65,66\}}\widetilde a_{s,65} + \gamma_{\{65,66\}}\widetilde a_{s,66} + \gamma_{\{68,84\}}\widetilde a_{s,68} + \gamma_{\{75,77\}}\widetilde a_{s,75} + \gamma_{\{75,77\}}\widetilde a_{s,77}\\ 
& + \gamma_{\{50,59,62,76\}}\widetilde a_{s,76} + \gamma_{\{50,52,70\}}\widetilde a_{s,78} + \gamma_{\{52,64,79\}}\widetilde a_{s,79} + \gamma_{\{80,81\}}\widetilde a_{s,80}\\ 
& + \gamma_{\{80,81\}}\widetilde a_{s,81} + \gamma_{\{59,67\}}\widetilde a_{s,82} + \gamma_{\{68,84\}}\widetilde a_{s,84} + \gamma_{\{67,89\}}\widetilde a_{s,85} + \gamma_{\{70,86,90\}}\widetilde a_{s,86}\\ 
& + \gamma_{\{50,52,59,67,70\}}\widetilde a_{s,88} + \gamma_{\{67,89\}}\widetilde a_{s,89} + \gamma_{\{70,86,90\}}\widetilde a_{s,90} \equiv 0,\\  
\rho_3(g) &+ g \equiv \gamma_{\{25,26\}}\widetilde a_{s,25} + \gamma_{\{25,26\}}\widetilde a_{s,26} + \gamma_{\{28,53\}}\widetilde a_{s,28} + \gamma_{\{29,30\}}\widetilde a_{s,29}\\ 
& + \gamma_{\{29,30\}}\widetilde a_{s,30} + \gamma_{\{49,50\}}\widetilde a_{s,49} + \gamma_{\{49,50\}}\widetilde a_{s,50} + \gamma_{\{52,56\}}\widetilde a_{s,52} + \gamma_{\{28,53\}}\widetilde a_{s,53}\\ 
& + \gamma_{\{54,55\}}\widetilde a_{s,54} + \gamma_{\{54,55\}}\widetilde a_{s,55} + \gamma_{\{52,56\}}\widetilde a_{s,56} + \gamma_{\{58,59\}}\widetilde a_{s,58} + \gamma_{\{58,59\}}\widetilde a_{s,59}\\ 
& + \gamma_{\{61,77\}}\widetilde a_{s,61} + \gamma_{\{62,63\}}\widetilde a_{s,62} + \gamma_{\{62,63\}}\widetilde a_{s,63} + \gamma_{\{64,65\}}\widetilde a_{s,64} + \gamma_{\{64,65\}}\widetilde a_{s,65}\\ 
& + \gamma_{\{57,66,81\}}\widetilde a_{s,66} + \gamma_{\{67,82\}}\widetilde a_{s,67} + \gamma_{\{68,69\}}\widetilde a_{s,68} + \gamma_{\{68,69\}}\widetilde a_{s,69}\\ 
& + \gamma_{\{75,76\}}\widetilde a_{s,75} + \gamma_{\{75,76\}}\widetilde a_{s,76} + \gamma_{\{61,77\}}\widetilde a_{s,77} + \gamma_{\{57,78\}}\widetilde a_{s,78} + \gamma_{\{79,80\}}\widetilde a_{s,79}\\ 
& + \gamma_{\{79,80\}}\widetilde a_{s,80} + \gamma_{\{57,66,81\}}\widetilde a_{s,81} + \gamma_{\{67,82\}}\widetilde a_{s,82} + \gamma_{\{57,78\}}\widetilde a_{s,83}\\ 
& + \gamma_{\{86,87\}}\widetilde a_{s,86} + \gamma_{\{86,87\}}\widetilde a_{s,87} + \gamma_{\{57,88,89\}}\widetilde a_{s,88} + \gamma_{\{57,88,89\}}\widetilde a_{s,89} \equiv 0.
\end{align*}
From these equalities we get $\gamma_{j} = 0$ for all $j \in\{25,\ldots,30,\, 49,\ldots,70,75,\ldots,90\}\setminus\{51,\, 60,\, 83,\, 85\}$. Hence, we have $g \equiv 0$ and
$f \equiv \sum_{6 \leqslant u \leqslant 9}\gamma_u\widetilde p_{s,u}.$ This completes the proof of the lemma.
\end{proof}	
\begin{proof}[Proof of Theorem \ref{mdct1}] Let $[f] \in (QP_4)_{d_{s,1}}^{GL_4}$ with $f \in (P_4)_{d_{s,1}}$.
	
For $s=2$, by Proposition \ref{mdt21}, $[f]_{(3)|^2} =0$, hence $[f]\in QP_4(3,1,1)$. By using Lemma \ref{bdct11}, we have $f \equiv \gamma_1\widetilde p_{2,1} + \gamma_2\widetilde p_{2,2} + \gamma_3\widetilde p_{2,3},$ with $\gamma_1,\, \gamma_2,\, \gamma_3 \in \mathbb F_2$. By computing $\rho_4(f) + f$ in terms of the admissible monomials we get
\begin{align*}
\rho_4(f)&+f \equiv \gamma_{1}\widetilde a_{2,1} + \gamma_{1}\widetilde a_{2,2} + \gamma_{1}\widetilde a_{2,3} + \gamma_{\{1,2\}}\widetilde a_{2,8} + \gamma_{\{1,2\}}\widetilde a_{2,9} + \gamma_{2}\widetilde a_{2,13}\\ 
& + \gamma_{2}\widetilde a_{2,14} + \gamma_{2}\widetilde a_{2,15} + \gamma_{2}\widetilde a_{2,17} + \gamma_{2}\widetilde a_{2,18} + \gamma_{2}\widetilde a_{2,25} + \gamma_{2}\widetilde a_{2,26} + \gamma_{1}\widetilde a_{2,27}\\ 
& + \gamma_{2}\widetilde a_{2,28} + \gamma_{2}\widetilde a_{2,31} + \gamma_{1}\widetilde a_{2,34} + \gamma_{1}\widetilde a_{2,35} + \gamma_{1}\widetilde a_{2,36} \equiv 0.
\end{align*} 
This equality implies $\gamma_1 = \gamma_2 = 0$ and $f \equiv \gamma_3\widetilde p_{2,3} = \gamma_3\xi_{1,2}$. The theorem is proved for $s = 2$. 

For $s = 3$, applying Proposition \ref{mdt21} and Lemma \ref{bdct11} to get $f \equiv \sum_{1\leqslant u \leqslant 6}\gamma_u\widetilde p_{3,u}$
with $\gamma_u \in \mathbb F_2$. By computing $\rho_4(f) + f$ in terms of the admissible monomials we obtain
\begin{align*}
\rho_4(f)&+f \equiv \gamma_{\{1,3\}}\widetilde a_{3,1} + \gamma_{\{1,3\}}\widetilde a_{3,2} + \gamma_{\{1,5,6\}}\widetilde a_{3,3} + \gamma_{\{1,2\}}\widetilde a_{3,8} + \gamma_{\{1,2\}}\widetilde a_{3,9}\\ 
& + \gamma_{\{2,4\}}\widetilde a_{3,13} + \gamma_{2}\widetilde a_{3,14} + \gamma_{2}\widetilde a_{3,15} + \gamma_{2}\widetilde a_{3,17} + \gamma_{2}\widetilde a_{3,18} + \gamma_{6}\widetilde a_{3,25} + \gamma_{6}\widetilde a_{3,26}\\ 
& + \gamma_{3}\widetilde a_{3,33} + \gamma_{3}\widetilde a_{3,34} + \gamma_{3}\widetilde a_{3,35} + \gamma_{3}\widetilde a_{3,36} + \gamma_{\{3,4\}}\widetilde a_{3,39} + \gamma_{4}\widetilde a_{3,44} + \gamma_{4}\widetilde a_{3,45}\\ 
& + \gamma_{4}\widetilde a_{3,51} + \gamma_{5}\widetilde a_{3,54} + \gamma_{5}\widetilde a_{3,55} + \gamma_{3}\widetilde a_{3,60} + \gamma_{4}\widetilde a_{3,61} + \gamma_{3}\widetilde a_{3,62} + \gamma_{3}\widetilde a_{3,63}\\ 
& + \gamma_{3}\widetilde a_{3,67} + \gamma_{\{1,2,6\}}\widetilde a_{3,71} + \gamma_{\{1,2,6\}}\widetilde a_{3,72} + \gamma_{\{2,4\}}\widetilde a_{3,73} + \gamma_{\{2,4\}}\widetilde a_{3,74}\\ 
& + \gamma_{6}\widetilde a_{3,75} + \gamma_{6}\widetilde a_{3,76} + \gamma_{3}\widetilde a_{3,77} + \gamma_{3}\widetilde a_{3,78} + \gamma_{3}\widetilde a_{3,79} + \gamma_{3}\widetilde a_{3,80} \equiv 0.
\end{align*} 	
This equality implies $\gamma_u = 0$ for $1\leqslant u \leqslant 6$ and $f \equiv 0$. The theorem is proved for $s = 3$.

For $s = 4$, by Proposition \ref{mdt21} we have $f \equiv_{(3)|^4} \gamma_0\bar \xi_{1,4}$ with $\gamma_0 \in \mathbb F_2$. Hence, $[f + \gamma_0\bar \xi_{1,4}] \in QP_4((3)|^3|(1)|^2)$. By a direct computation we obtain
\begin{equation}\label{cthr1}
\begin{cases}
\rho_1(\bar \xi_{1,4}) + \bar \xi_{1,4} \equiv 0,\\
\rho_2(\bar \xi_{1,4}) + \bar \xi_{1,4} \equiv \sum_{j \in \{25,\, 49,\, 51,\, 64,\, 70,\, 75,\, 77,\, 79,\, 86,\, 88\}}\widetilde a_{4,j},\\
\rho_3(\bar \xi_{1,4}) + \bar \xi_{1,4} \equiv \sum_{j \in \{51,\, 61,\, 66,\, 67,\, 77,\, 78,\, 81,\, 82,\, 83,\, 90\}}\widetilde a_{4,j},\\
\rho_4(\bar \xi_{1,4}) + \bar \xi_{1,4} \equiv \sum_{j \in \{17,\, 18,\, 57\}}\widetilde a_{4,j}.
\end{cases}
\end{equation}
By using \eqref{cthr1} we see that $\rho_i(f) + f \equiv 0$ for $i=1,\, 2,\, 3$ if and only if
$$f \equiv \gamma_0(\bar \xi_{1,4} + \widetilde p_{4,0}) + \sum _{1 \leqslant u \leqslant 9} \gamma_u\widetilde p_{4,u},$$
where $\gamma_u \in \mathbb F_2$ and 
$$\widetilde p_{4,0} = \mbox{$\sum_{j \in \{25,\, 26,\, 27,\, 49,\, 50,\, 51,\, 68,\, 69,\, 70,\, 75,\, 76,\, 77,\, 78,\, 79,\, 80,\, 81,\, 82,\, 84\}}\widetilde a_{4,j}$}.$$
Now by computing $\rho_4(f) + f$ in terms of the admissible monomials we get
\begin{align*}
\rho_4(f)&+f \equiv \gamma_{\{1,3\}}\widetilde a_{4,1} + \gamma_{\{1,3\}}\widetilde a_{4,2} + \gamma_{\{0,1,6,7\}}\widetilde a_{4,3} + \gamma_{\{1,2\}}\widetilde a_{4,8} + \gamma_{\{1,2\}}\widetilde a_{4,9}\\ 
& + \gamma_{\{2,4\}}\widetilde a_{4,13} + \gamma_{\{2,8\}}\widetilde a_{4,14} + \gamma_{\{2,8\}}\widetilde a_{4,15} + \gamma_{\{0,2\}}\widetilde a_{4,17} + \gamma_{\{0,2\}}\widetilde a_{4,18}\\ 
& + \gamma_{\{1,2,3,7,8,9\}}\widetilde a_{4,25} + \gamma_{\{1,2,3,7,8,9\}}\widetilde a_{4,26} + \gamma_{\{0,2,4,8\}}\widetilde a_{4,28} + \gamma_{\{3,5\}}\widetilde a_{4,33}\\ 
& + \gamma_{\{3,5\}}\widetilde a_{4,34} + \gamma_{\{3,7,8,9\}}\widetilde a_{4,35} + \gamma_{\{3,7,8,9\}}\widetilde a_{4,36} + \gamma_{\{3,4\}}\widetilde a_{4,39} + \gamma_{\{4,8\}}\widetilde a_{4,44}\\ 
& + \gamma_{\{4,8\}}\widetilde a_{4,45} + \gamma_{\{1,2,3,8\}}\widetilde a_{4,49} + \gamma_{\{1,2,3,8\}}\widetilde a_{4,50} + \gamma_{\{0,3,4,5\}}\widetilde a_{4,51}\\ 
& + \gamma_{\{2,4\}}\widetilde a_{4,52} + \gamma_{\{0,2,4,8\}}\widetilde a_{4,53} + \gamma_{\{0,6,8,9\}}\widetilde a_{4,54} + \gamma_{\{0,6,8,9\}}\widetilde a_{4,55}\\ 
& + \gamma_{\{2,4\}}\widetilde a_{4,56} + \gamma_{\{0,8\}}\widetilde a_{4,57} + \gamma_{\{7,8\}}\widetilde a_{4,61} + \gamma_{\{8,9\}}\widetilde a_{4,66} + \gamma_{\{5,7\}}\widetilde a_{4,73}\\ 
& + \gamma_{\{1,2,3,8\}}\widetilde a_{4,75} + \gamma_{\{1,2,3,8\}}\widetilde a_{4,76} + \gamma_{\{0,7\}}\widetilde a_{4,77} + \gamma_{\{2,4\}}\widetilde a_{4,79} + \gamma_{\{2,4\}}\widetilde a_{4,80}\\ 
& + \gamma_{\{0,9\}}\widetilde a_{4,81} + \gamma_{\{0,9\}}\widetilde a_{4,84} + \gamma_{\{0,8\}}\widetilde a_{4,89} + \gamma_{\{3,4,5,8\}}\widetilde a_{4,90} \equiv 0.
\end{align*} 
From this equality we obtain $\gamma_u = \gamma_0$ for $1 \leqslant u \leqslant 9$ and
\begin{align*}
f &\equiv \gamma_0\Big(\bar \xi_{1,4} + \sum _{0 \leqslant u \leqslant 9}\widetilde p_{4,u}\Big) = \gamma_0\xi_{1,4}.
\end{align*}
The theorem is proved for $s= 4$.

For $s \geqslant 5$, by Proposition \ref{mdt21} we have $f \equiv_{(3)|^s} \gamma_0\bar \xi_{1,s}$ with $\gamma_0 \in \mathbb F_2$. Hence, $[f + \gamma_0\bar \xi_{1,4}] \in QP_4((3)|^{S-1}|(1)|^2)$. By a direct computation we obtain
\begin{equation}\label{cthr2}
\begin{cases}
\rho_1(\bar \xi_{1,s}) + \bar \xi_{1,s} \equiv 0,\\
\rho_2(\bar \xi_{1,s}) + \bar \xi_{1,s} \equiv \sum_{j \in \{64,\, 79,\, 86,\, 90\}}\widetilde a_{s,j},\\
\rho_3(\bar \xi_{1,s}) + \bar \xi_{1,s} \equiv \sum_{j \in \{66,\, 67,\, 81,\, 82,\, 88,\, 89\}}\widetilde a_{s,j},\\
\rho_4(\bar \xi_{1,s}) + \bar \xi_{1,s} \equiv \sum_{j \in \{17,\, 18,\, 57\}}\widetilde a_{s,j}.
\end{cases}
\end{equation}
By using \eqref{cthr1} we see that $\rho_i(f) + f \equiv 0$ for $i=1,\, 2,\, 3$ if and only if
$$f \equiv \gamma_0(\bar \xi_{1,s} + \widetilde p_{s,0}) + \sum _{1 \leqslant u \leqslant 9} \gamma_u\widetilde p_{s,u},$$
where $\gamma_u \in \mathbb F_2$ and 
$$\widetilde p_{s,0} = \mbox{$\sum_{j \in \{79,\, 80,\, 81,\, 82,\, 88,\, 90\}}\widetilde a_{s,j}$}.$$
Now by computing $\rho_4(f) + f$ in terms of the admissible monomials we get
\begin{align*}
\rho_4(f)&+f \equiv \gamma_{\{1,3\}}\widetilde a_{s,1} + \gamma_{\{1,3\}}\widetilde a_{s,2} + \gamma_{\{1,6\}}\widetilde a_{s,3} + \gamma_{\{1,2\}}\widetilde a_{s,8} + \gamma_{\{1,2\}}\widetilde a_{s,9}\\ 
& + \gamma_{\{2,4\}}\widetilde a_{s,13} + \gamma_{\{2,7\}}\widetilde a_{s,14} + \gamma_{\{2,7\}}\widetilde a_{s,15} + \gamma_{\{0,2\}}\widetilde a_{s,17} + \gamma_{\{0,2\}}\widetilde a_{s,18}\\ 
& + \gamma_{\{1,2,3,8\}}\widetilde a_{s,25} + \gamma_{\{1,2,3,8\}}\widetilde a_{s,26} + \gamma_{\{0,2,4,7\}}\widetilde a_{s,28} + \gamma_{\{3,5\}}\widetilde a_{s,33}\\ 
& + \gamma_{\{3,5\}}\widetilde a_{s,34} + \gamma_{\{3,8\}}\widetilde a_{s,35} + \gamma_{\{3,8\}}\widetilde a_{s,36} + \gamma_{\{3,4\}}\widetilde a_{s,39} + \gamma_{\{4,7\}}\widetilde a_{s,44}\\ 
& + \gamma_{\{4,7\}}\widetilde a_{s,45} + \gamma_{\{1,2,3,7\}}\widetilde a_{s,49} + \gamma_{\{1,2,3,7\}}\widetilde a_{s,50} + \gamma_{\{2,4\}}\widetilde a_{s,52}\\ 
& + \gamma_{\{0,2,4,7\}}\widetilde a_{s,53} + \gamma_{\{6,8\}}\widetilde a_{s,54} + \gamma_{\{6,8\}}\widetilde a_{s,55} + \gamma_{\{2,4\}}\widetilde a_{s,56} + \gamma_{\{0,7\}}\widetilde a_{s,57}\\ 
& + \gamma_{\{3,4,5,9\}}\widetilde a_{s,61} + \gamma_{\{8,9\}}\widetilde a_{s,66} + \gamma_{\{0,7\}}\widetilde a_{s,70} + \gamma_{\{5,9\}}\widetilde a_{s,73} + \gamma_{\{1,2,3,7\}}\widetilde a_{s,75}\\ 
& + \gamma_{\{1,2,3,7\}}\widetilde a_{s,76} + \gamma_{\{3,4,5,9\}}\widetilde a_{s,77} + \gamma_{\{0,3,4,5\}}\widetilde a_{s,78} + \gamma_{\{2,4\}}\widetilde a_{s,79}\\ 
& + \gamma_{\{2,4\}}\widetilde a_{s,80} + \gamma_{\{0,7,8,9\}}\widetilde a_{s,81} + \gamma_{\{3,4,5,7\}}\widetilde a_{s,83} + \gamma_{\{8,9\}}\widetilde a_{s,84}\\ 
& + \gamma_{\{0,3,4,5\}}\widetilde a_{s,88} + \gamma_{\{3,4,5,7\}}\widetilde a_{s,89} + \gamma_{\{0,7\}}\widetilde a_{s,90} \equiv 0.
\end{align*} 
This equality imply $\gamma_u = \gamma_0$ for $1 \leqslant u \leqslant 9$ and
\begin{align*}
f &\equiv \gamma_0\Big(\bar \xi_{1,s} + \sum _{0 \leqslant u \leqslant 9}\widetilde p_{s,u}\Big) = \gamma_0\Big(\bar \xi_{1,s} + \sum _{1 \leqslant j \leqslant 90}\widetilde a_{s,j}\Big) = \gamma_0\xi_{1,s}.
\end{align*}
The theorem is completely proved.
\end{proof}

\subsection{Computation of $QP_4((3)|^{s-1}|(1)|^{t+1})^{GL_4}$ for $t \geqslant 2$}\

\medskip
Following our work \cite{su5}, a basis of $QP_4((3)|^{s-1}|(1)|^{t+1})$ is the set of all the classes represented by the admissible monomials $\widetilde a_{s,j} := \widetilde a_{t,s,j}$ which are determined as follows:

\medskip
For $s \geqslant 2$,

\medskip
\centerline{\begin{tabular}{lll} 
$\widetilde a_{s,1} = x_2^{2^{s-1}-1}x_3^{2^{s-1}-1}x_4^{2^{s+t}-1}$\cr  $\widetilde a_{s,2} = x_2^{2^{s-1}-1}x_3^{2^{s+t}-1}x_4^{2^{s-1}-1}$\cr  $\widetilde a_{s,3} = x_2^{2^{s+t}-1}x_3^{2^{s-1}-1}x_4^{2^{s-1}-1}$\cr  $\widetilde a_{s,4} = x_1^{2^{s-1}-1}x_3^{2^{s-1}-1}x_4^{2^{s+t}-1}$\cr  $\widetilde a_{s,5} = x_1^{2^{s-1}-1}x_3^{2^{s+t}-1}x_4^{2^{s-1}-1}$\cr  $\widetilde a_{s,6} = x_1^{2^{s-1}-1}x_2^{2^{s-1}-1}x_4^{2^{s+t}-1}$\cr  $\widetilde a_{s,7} = x_1^{2^{s-1}-1}x_2^{2^{s-1}-1}x_3^{2^{s+t}-1}$\cr  $\widetilde a_{s,8} = x_1^{2^{s-1}-1}x_2^{2^{s+t}-1}x_4^{2^{s-1}-1}$\cr  $\widetilde a_{s,9} = x_1^{2^{s-1}-1}x_2^{2^{s+t}-1}x_3^{2^{s-1}-1}$\cr  $\widetilde a_{s,10} = x_1^{2^{s+t}-1}x_3^{2^{s-1}-1}x_4^{2^{s-1}-1}$\cr  $\widetilde a_{s,11} = x_1^{2^{s+t}-1}x_2^{2^{s-1}-1}x_4^{2^{s-1}-1}$\cr  $\widetilde a_{s,12} = x_1^{2^{s+t}-1}x_2^{2^{s-1}-1}x_3^{2^{s-1}-1}$\cr  $\widetilde a_{s,13} = x_2^{2^{s-1}-1}x_3^{2^{s}-1}x_4^{2^{s+t}-2^{s-1}-1}$\cr  $\widetilde a_{s,14} = x_2^{2^{s}-1}x_3^{2^{s-1}-1}x_4^{2^{s+t}-2^{s-1}-1}$\cr  $\widetilde a_{s,15} = x_2^{2^{s}-1}x_3^{2^{s+t}-2^{s-1}-1}x_4^{2^{s-1}-1}$\cr  $\widetilde a_{s,16} = x_1^{2^{s-1}-1}x_3^{2^{s}-1}x_4^{2^{s+t}-2^{s-1}-1}$\cr  $\widetilde a_{s,17} = x_1^{2^{s-1}-1}x_2^{2^{s}-1}x_4^{2^{s+t}-2^{s-1}-1}$\cr  $\widetilde a_{s,18} = x_1^{2^{s-1}-1}x_2^{2^{s}-1}x_3^{2^{s+t}-2^{s-1}-1}$\cr  $\widetilde a_{s,19} = x_1^{2^{s}-1}x_3^{2^{s-1}-1}x_4^{2^{s+t}-2^{s-1}-1}$\cr  $\widetilde a_{s,20} = x_1^{2^{s}-1}x_3^{2^{s+t}-2^{s-1}-1}x_4^{2^{s-1}-1}$\cr  $\widetilde a_{s,21} = x_1^{2^{s}-1}x_2^{2^{s-1}-1}x_4^{2^{s+t}-2^{s-1}-1}$\cr  $\widetilde a_{s,22} = x_1^{2^{s}-1}x_2^{2^{s-1}-1}x_3^{2^{s+t}-2^{s-1}-1}$\cr  $\widetilde a_{s,23} = x_1^{2^{s}-1}x_2^{2^{s+t}-2^{s-1}-1}x_4^{2^{s-1}-1}$\cr  $\widetilde a_{s,24} = x_1^{2^{s}-1}x_2^{2^{s+t}-2^{s-1}-1}x_3^{2^{s-1}-1}$\cr  $\widetilde a_{s,25} = x_2^{2^{s}-1}x_3^{3.2^{s-1}-1}x_4^{2^{s+t}-3.2^{s-1}-1}$\cr  $\widetilde a_{s,26} = x_1^{2^{s}-1}x_3^{3.2^{s-1}-1}x_4^{2^{s+t}-3.2^{s-1}-1}$\cr  $\widetilde a_{s,27} = x_1^{2^{s}-1}x_2^{3.2^{s-1}-1}x_4^{2^{s+t}-3.2^{s-1}-1}$\cr  $\widetilde a_{s,28} = x_1^{2^{s}-1}x_2^{3.2^{s-1}-1}x_3^{2^{s+t}-3.2^{s-1}-1}$\cr  $\widetilde a_{s,29} = x_1x_2^{2^{s-1}-1}x_3^{2^{s-1}-1}x_4^{2^{s+t}-2}$\cr  $\widetilde a_{s,30} = x_1x_2^{2^{s-1}-1}x_3^{2^{s+t}-2}x_4^{2^{s-1}-1}$\cr  $\widetilde a_{s,31} = x_1x_2^{2^{s+t}-2}x_3^{2^{s-1}-1}x_4^{2^{s-1}-1}$\cr  $\widetilde a_{s,32} = x_1x_2^{2^{s-1}-1}x_3^{2^{s}-2}x_4^{2^{s+t}-2^{s-1}-1}$\cr  $\widetilde a_{s,33} = x_1x_2^{2^{s}-2}x_3^{2^{s-1}-1}x_4^{2^{s+t}-2^{s-1}-1}$\cr  $\widetilde a_{s,34} = x_1x_2^{2^{s}-2}x_3^{2^{s+t}-2^{s-1}-1}x_4^{2^{s-1}-1}$\cr  $\widetilde a_{s,35} = x_1x_2^{2^{s}-2}x_3^{3.2^{s-1}-1}x_4^{2^{s+t}-3.2^{s-1}-1}$\cr    	 
\end{tabular}}

\medskip
For $s = 2$,

\medskip
\centerline{\begin{tabular}{lll}  
$\widetilde a_{2,36} = x_1x_2x_3^{3}x_4^{2^{t+2}-4}$ & &$\widetilde a_{2,37} = x_1x_2^{3}x_3x_4^{2^{t+2}-4}$\cr  $\widetilde a_{2,38} = x_1x_2^{3}x_3^{2^{t+2}-4}x_4$ & &$\widetilde a_{2,39} = x_1^{3}x_2x_3x_4^{2^{t+2}-4}$\cr  $\widetilde a_{2,40} = x_1^{3}x_2x_3^{2^{t+2}-4}x_4$ & &$\widetilde a_{2,41} = x_1x_2^{3}x_3^{4}x_4^{2^{t+2}-7}$\cr  $\widetilde a_{2,42} = x_1^{3}x_2x_3^{4}x_4^{2^{t+2}-7}$ & &$\widetilde a_{2,43} = x_1x_2^{3}x_3^{5}x_4^{2^{t+2}-8}$\cr  $\widetilde a_{2,44} = x_1^{3}x_2x_3^{5}x_4^{2^{t+2}-8}$ & &$\widetilde a_{2,45} = x_1^{3}x_2^{5}x_3x_4^{2^{t+2}-8}$\cr 	 
\end{tabular}}

\medskip
For $s = t = 2$,\ $\widetilde a_{2,46} = x_1^{3}x_2^{5}x_3^8x_4$.

\medskip
For $s \geqslant 3$,

\medskip
\centerline{\begin{tabular}{lll}  
$\widetilde a_{s,36} = x_1x_2^{2^{s-1}-2}x_3^{2^{s-1}-1}x_4^{2^{s+t}-1}$\cr  $\widetilde a_{s,37} = x_1x_2^{2^{s-1}-2}x_3^{2^{s+t}-1}x_4^{2^{s-1}-1}$\cr  $\widetilde a_{s,38} = x_1x_2^{2^{s-1}-1}x_3^{2^{s-1}-2}x_4^{2^{s+t}-1}$\cr  $\widetilde a_{s,39} = x_1x_2^{2^{s-1}-1}x_3^{2^{s+t}-1}x_4^{2^{s-1}-2}$\cr  $\widetilde a_{s,40} = x_1x_2^{2^{s+t}-1}x_3^{2^{s-1}-2}x_4^{2^{s-1}-1}$\cr  $\widetilde a_{s,41} = x_1x_2^{2^{s+t}-1}x_3^{2^{s-1}-1}x_4^{2^{s-1}-2}$\cr  $\widetilde a_{s,42} = x_1^{2^{s-1}-1}x_2x_3^{2^{s-1}-2}x_4^{2^{s+t}-1}$\cr  $\widetilde a_{s,43} = x_1^{2^{s-1}-1}x_2x_3^{2^{s+t}-1}x_4^{2^{s-1}-2}$\cr  $\widetilde a_{s,44} = x_1^{2^{s-1}-1}x_2^{2^{s+t}-1}x_3x_4^{2^{s-1}-2}$\cr  $\widetilde a_{s,45} = x_1^{2^{s+t}-1}x_2x_3^{2^{s-1}-2}x_4^{2^{s-1}-1}$\cr  $\widetilde a_{s,46} = x_1^{2^{s+t}-1}x_2x_3^{2^{s-1}-1}x_4^{2^{s-1}-2}$\cr  $\widetilde a_{s,47} = x_1^{2^{s+t}-1}x_2^{2^{s-1}-1}x_3x_4^{2^{s-1}-2}$\cr  $\widetilde a_{s,48} = x_1x_2^{2^{s-1}-2}x_3^{2^{s}-1}x_4^{2^{s+t}-2^{s-1}-1}$\cr  $\widetilde a_{s,49} = x_1x_2^{2^{s}-1}x_3^{2^{s-1}-2}x_4^{2^{s+t}-2^{s-1}-1}$\cr  $\widetilde a_{s,50} = x_1x_2^{2^{s}-1}x_3^{2^{s+t}-2^{s-1}-1}x_4^{2^{s-1}-2}$\cr  $\widetilde a_{s,51} = x_1^{2^{s}-1}x_2x_3^{2^{s-1}-2}x_4^{2^{s+t}-2^{s-1}-1}$\cr  $\widetilde a_{s,52} = x_1^{2^{s}-1}x_2x_3^{2^{s+t}-2^{s-1}-1}x_4^{2^{s-1}-2}$\cr  $\widetilde a_{s,53} = x_1^{2^{s}-1}x_2^{2^{s+t}-2^{s-1}-1}x_3x_4^{2^{s-1}-2}$\cr  $\widetilde a_{s,54} = x_1^{2^{s-1}-1}x_2x_3^{2^{s-1}-1}x_4^{2^{s+t}-2}$\cr  $\widetilde a_{s,55} = x_1^{2^{s-1}-1}x_2x_3^{2^{s+t}-2}x_4^{2^{s-1}-1}$\cr  $\widetilde a_{s,56} = x_1^{2^{s-1}-1}x_2^{2^{s-1}-1}x_3x_4^{2^{s+t}-2}$\cr  $\widetilde a_{s,57} = x_1^{2^{s-1}-1}x_2x_3^{2^{s}-2}x_4^{2^{s+t}-2^{s-1}-1}$\cr  $\widetilde a_{s,58} = x_1x_2^{2^{s-1}-1}x_3^{2^{s}-1}x_4^{2^{s+t}-2^{s-1}-2}$\cr  $\widetilde a_{s,59} = x_1x_2^{2^{s}-1}x_3^{2^{s-1}-1}x_4^{2^{s+t}-2^{s-1}-2}$\cr  $\widetilde a_{s,60} = x_1x_2^{2^{s}-1}x_3^{2^{s+t}-2^{s-1}-2}x_4^{2^{s-1}-1}$\cr  $\widetilde a_{s,61} = x_1^{2^{s-1}-1}x_2x_3^{2^{s}-1}x_4^{2^{s+t}-2^{s-1}-2}$\cr  $\widetilde a_{s,62} = x_1^{2^{s-1}-1}x_2^{2^{s}-1}x_3x_4^{2^{s+t}-2^{s-1}-2}$\cr  $\widetilde a_{s,63} = x_1^{2^{s}-1}x_2x_3^{2^{s-1}-1}x_4^{2^{s+t}-2^{s-1}-2}$\cr  $\widetilde a_{s,64} = x_1^{2^{s}-1}x_2x_3^{2^{s+t}-2^{s-1}-2}x_4^{2^{s-1}-1}$\cr  $\widetilde a_{s,65} = x_1^{2^{s}-1}x_2^{2^{s-1}-1}x_3x_4^{2^{s+t}-2^{s-1}-2}$\cr  $\widetilde a_{s,66} = x_1x_2^{2^{s}-1}x_3^{3.2^{s-1}-2}x_4^{2^{s+t}-3.2^{s-1}-1}$\cr  $\widetilde a_{s,67} = x_1^{2^{s}-1}x_2x_3^{3.2^{s-1}-2}x_4^{2^{s+t}-3.2^{s-1}-1}$\cr  $\widetilde a_{s,68} = x_1x_2^{2^{s}-1}x_3^{3.2^{s-1}-1}x_4^{2^{s+t}-3.2^{s-1}-2}$\cr  $\widetilde a_{s,69} = x_1^{2^{s}-1}x_2x_3^{3.2^{s-1}-1}x_4^{2^{s+t}-3.2^{s-1}-2}$\cr  $\widetilde a_{s,70} = x_1^{2^{s}-1}x_2^{3.2^{s-1}-1}x_3x_4^{2^{s+t}-3.2^{s-1}-2}$\cr  $\widetilde a_{s,71} = x_1^{3}x_2^{2^{s-1}-1}x_3^{2^{s+t}-3}x_4^{2^{s-1}-2}$\cr  $\widetilde a_{s,72} = x_1^{3}x_2^{2^{s+t}-3}x_3^{2^{s-1}-2}x_4^{2^{s-1}-1}$\cr  $\widetilde a_{s,73} = x_1^{3}x_2^{2^{s+t}-3}x_3^{2^{s-1}-1}x_4^{2^{s-1}-2}$\cr  
\end{tabular}}
\centerline{\begin{tabular}{lll}   	 
$\widetilde a_{s,74} = x_1^{3}x_2^{2^{s}-3}x_3^{2^{s-1}-2}x_4^{2^{s+t}-2^{s-1}-1}$\cr  $\widetilde a_{s,75} = x_1^{3}x_2^{2^{s}-3}x_3^{2^{s+t}-2^{s-1}-1}x_4^{2^{s-1}-2}$\cr  $\widetilde a_{s,76} = x_1^{3}x_2^{2^{s}-1}x_3^{2^{s+t}-2^{s-1}-3}x_4^{2^{s-1}-2}$\cr  $\widetilde a_{s,77} = x_1^{2^{s}-1}x_2^{3}x_3^{2^{s+t}-2^{s-1}-3}x_4^{2^{s-1}-2}$\cr  $\widetilde a_{s,78} = x_1^{3}x_2^{2^{s}-3}x_3^{2^{s-1}-1}x_4^{2^{s+t}-2^{s-1}-2}$\cr  $\widetilde a_{s,79} = x_1^{3}x_2^{2^{s}-3}x_3^{2^{s+t}-2^{s-1}-2}x_4^{2^{s-1}-1}$\cr  $\widetilde a_{s,80} = x_1^{3}x_2^{2^{s}-3}x_3^{3.2^{s-1}-2}x_4^{2^{s+t}-3.2^{s-1}-1}$\cr  $\widetilde a_{s,81} = x_1^{3}x_2^{2^{s}-3}x_3^{3.2^{s-1}-1}x_4^{2^{s+t}-3.2^{s-1}-2}$\cr  $\widetilde a_{s,82} = x_1^{3}x_2^{2^{s-1}-1}x_3^{2^{s}-3}x_4^{2^{s+t}-2^{s-1}-2}$\cr  $\widetilde a_{s,83} = x_1^{3}x_2^{2^{s}-1}x_3^{3.2^{s-1}-3}x_4^{2^{s+t}-3.2^{s-1}-2}$\cr  $\widetilde a_{s,84} = x_1^{2^{s}-1}x_2^{3}x_3^{3.2^{s-1}-3}x_4^{2^{s+t}-3.2^{s-1}-2}$\cr   
\end{tabular}}

For $s = 3$,

\medskip
\centerline{\begin{tabular}{lll}   	 
$\widetilde a_{3,85} = x_1^{3}x_2^{3}x_3^{3}x_4^{2^{t+3}-4}$ & &$\widetilde a_{3,86} = x_1^{3}x_2^{3}x_3^{2^{t+3}-4}x_4^{3}$\cr  $\widetilde a_{3,87} = x_1^{3}x_2^{3}x_3^{4}x_4^{2^{t+3}-5}$ & &$\widetilde a_{3,88} = x_1^{3}x_2^{3}x_3^{7}x_4^{2^{t+3}-8}$\cr  $\widetilde a_{3,89} = x_1^{3}x_2^{7}x_3^{3}x_4^{2^{t+3}-8}$ & &$\widetilde a_{3,90} = x_1^{7}x_2^{3}x_3^{3}x_4^{2^{t+3}-8}$\cr  
\end{tabular}}

\medskip
For $s \geqslant 4$,

\medskip
\centerline{\begin{tabular}{lll}
$\widetilde a_{s,85} = x_1^{3}x_2^{2^{s-1}-3}x_3^{2^{s-1}-2}x_4^{2^{s+t}-1}$\cr  $\widetilde a_{s,86} = x_1^{3}x_2^{2^{s-1}-3}x_3^{2^{s+t}-1}x_4^{2^{s-1}-2}$\cr  $\widetilde a_{s,87} = x_1^{3}x_2^{2^{s+t}-1}x_3^{2^{s-1}-3}x_4^{2^{s-1}-2}$\cr  $\widetilde a_{s,88} = x_1^{2^{s+t}-1}x_2^{3}x_3^{2^{s-1}-3}x_4^{2^{s-1}-2}$\cr  $\widetilde a_{s,89} = x_1^{3}x_2^{2^{s-1}-3}x_3^{2^{s-1}-1}x_4^{2^{s+t}-2}$\cr  $\widetilde a_{s,90} = x_1^{3}x_2^{2^{s-1}-3}x_3^{2^{s+t}-2}x_4^{2^{s-1}-1}$\cr  $\widetilde a_{s,91} = x_1^{3}x_2^{2^{s-1}-1}x_3^{2^{s-1}-3}x_4^{2^{s+t}-2}$\cr  $\widetilde a_{s,92} = x_1^{2^{s-1}-1}x_2^{3}x_3^{2^{s-1}-3}x_4^{2^{s+t}-2}$\cr  $\widetilde a_{s,93} = x_1^{3}x_2^{2^{s-1}-3}x_3^{2^{s}-2}x_4^{2^{s+t}-2^{s-1}-1}$\cr  $\widetilde a_{s,94} = x_1^{3}x_2^{2^{s-1}-3}x_3^{2^{s}-1}x_4^{2^{s+t}-2^{s-1}-2}$\cr  $\widetilde a_{s,95} = x_1^{3}x_2^{2^{s}-1}x_3^{2^{s-1}-3}x_4^{2^{s+t}-2^{s-1}-2}$\cr  $\widetilde a_{s,96} = x_1^{2^{s}-1}x_2^{3}x_3^{2^{s-1}-3}x_4^{2^{s+t}-2^{s-1}-2}$\cr  $\widetilde a_{s,97} = x_1^{2^{s-1}-1}x_2^{3}x_3^{2^{s+t}-3}x_4^{2^{s-1}-2}$\cr  $\widetilde a_{s,98} = x_1^{2^{s-1}-1}x_2^{3}x_3^{2^{s}-3}x_4^{2^{s+t}-2^{s-1}-2}$\cr  $\widetilde a_{s,99} = x_1^{7}x_2^{2^{s-1}-5}x_3^{2^{s-1}-3}x_4^{2^{s+t}-2}$\cr  $\widetilde a_{s,100} = x_1^{7}x_2^{2^{s-1}-5}x_3^{2^{s+t}-3}x_4^{2^{s-1}-2}$\cr  $\widetilde a_{s,101} = x_1^{7}x_2^{2^{s}-5}x_3^{3.2^{s-1}-3}x_4^{2^{s+t}-3.2^{s-1}-2}$\cr  $\widetilde a_{s,102} = x_1^{7}x_2^{2^{s+t}-5}x_3^{2^{s-1}-3}x_4^{2^{s-1}-2}$\cr  
\end{tabular}}

\medskip
For $s = 4$, 
$$\widetilde a_{4,103} = x_1^{7}x_2^{7}x_3^{7}x_4^{2^{t+4}-8},\ \widetilde a_{4,104} = x_1^{7}x_2^{7}x_3^{9}x_4^{2^{t+4}-10},\ \widetilde a_{4,105} = x_1^{7}x_2^{7}x_3^{2^{t+4}-7}x_4^{6}.$$

For $s \geqslant 5$,

\medskip
\centerline{\begin{tabular}{lll}   	 
$\widetilde a_{s,103} = x_1^{7}x_2^{2^{s-1}-5}x_3^{2^{s}-3}x_4^{2^{s+t}-2^{s-1}-2}$\cr  $\widetilde a_{s,104} = x_1^{7}x_2^{2^{s}-5}x_3^{2^{s-1}-3}x_4^{2^{s+t}-2^{s-1}-2}$\cr  $\widetilde a_{s,105} = x_1^{7}x_2^{2^{s}-5}x_3^{2^{s+t}-2^{s-1}-3}x_4^{2^{s-1}-2}$\cr   
\end{tabular}}

\medskip
\centerline{\begin{tabular}{lll}   	 
\end{tabular}}

We note that for a fixed value of $t \geqslant 2$, we denote $\widetilde a_{s,j} = \widetilde a_{t,s,j}$ and $\widetilde p_{s,u} = \widetilde p_{t,s,u}$. We prove the following.

\begin{thms}\label{dlt21} For any $t \geqslant 2$, we have $QP_4((3)|^{s-1}|(1)|^{t+1})^{GL_4} = \langle [\xi_{t,s}]\rangle$, where
$$\xi_{t,s} = \begin{cases}
\sum_{j\in \{29,\, 30,\, 37,\, 38,\, 39,\, 40,\, 45,\, 46\}}\widetilde a_{2,j}, &\mbox{if } s = t = 2,\\
0, &\mbox{if } s = 2,\, t>2; s= 3,\\
\sum_{1 \leqslant j \leqslant 105,\, j \notin \{29,\, 54,\, 56,\, 71,\, 82,\, 89,\, 91,\, 92,\, 97,\, 98\}}\widetilde a_{4,j}, &\mbox{if } s = 4,\\
\sum_{1 \leqslant j \leqslant 105}\widetilde a_{s,j}, &\mbox{if } s \geqslant 5.
 \end{cases}$$
\end{thms}
\begin{rems} The dimensional result of this theorem is also stated in \cite[Props. 4.1.6, 4.1.13]{pp25} but the detailed proof is not provided. So, it is only a prediction with the assumption that Singer's conjecture is true for $k = 4$.
\end{rems}
By a direct computation, we easily observe that
\begin{align}
\begin{cases}\label{ctt21}
[\Sigma_4(\widetilde a_{s,1})] = \langle \{[\widetilde a_{s,j}]: 1 \leqslant j \leqslant 12 \}\rangle,\ [\Sigma_4(\widetilde a_{s,1})]^{\Sigma_4} = \langle [\widetilde p_{s,1}]\rangle,\\
[\Sigma_4(\widetilde a_{s,13})] = \langle \{[\widetilde a_{s,j}]: 13 \leqslant j \leqslant 24\} \rangle,\ [\Sigma_4(\widetilde a_{s,13})]^{\Sigma_4} = \langle [\widetilde p_{s,2}]\rangle,\\
[\Sigma_4(\widetilde a_{s,25})] = \langle \{[\widetilde a_{s,j} : 25 \leqslant j \leqslant 28 \}\rangle, \ [\Sigma_4(\widetilde a_{s,25})]^{\Sigma_4} = \langle [\widetilde p_{s,3}]\rangle,\\
[\Sigma_4(\widetilde a_{2,36})] = \langle \{[\widetilde a_{s,j}]: 36 \leqslant j \leqslant 47\} \rangle, \ [\Sigma_4(\widetilde a_{s,36})]^{\Sigma_4} = \langle [\widetilde p_{s,4}]\rangle,\ s \geqslant 3,\\
[\Sigma_4(a_{s,48})] = \langle \{[\widetilde a_{s,j}]: 48 \leqslant j \leqslant 53\} \rangle, \ [\Sigma_4(\widetilde a_{s,48})]^{\Sigma_4} = \langle [\widetilde p_{s,5}]\rangle,\  s \geqslant 3,\\
[\Sigma_4(a_{s,85})] = \langle \{[\widetilde a_{s,j}]: 85 \leqslant j \leqslant 88\} \rangle, \ [\Sigma_4(\widetilde a_{s,85})]^{\Sigma_4} = \langle [\widetilde p_{s,6}]\rangle,\  s \geqslant 4,
\end{cases}
\end{align}
where
\begin{align*}
&\widetilde p_{s,1} = \sum_{1 \leqslant j \leqslant 12}\widetilde a_{s,j}, \ \widetilde p_{s,2} = \sum_{13 \leqslant j \leqslant 24}\widetilde a_{s,j}, \ \widetilde p_{s,3} = \sum_{25\leqslant j \leqslant 28}\widetilde a_{s,j},\\
&\widetilde p_{s,4} = \sum_{36 \leqslant j \leqslant 47}\widetilde a_{s,j},\,  \widetilde p_{s,5} = \sum_{48 \leqslant j \leqslant 53}\widetilde a_{s,j},\
\widetilde p_{s,6} = \sum_{85 \leqslant j \leqslant 88}\widetilde a_{s,j}.
\end{align*}

We begin with the case $s=t=2$.
\begin{lems}\label{bdt22} We have
$QP_4((3)|(1)|^3)^{\Sigma_4} = \langle \{[\widetilde p_{2,u}] : 1 \leqslant u \leqslant 4\} \rangle,$ where $\widetilde p_{2,u}$ is defined by \eqref{ctt21} for $1 \leqslant u \leqslant 3$ and 
$$\widetilde p_{2,4} = \begin{cases} \sum_{j\in \{29,\, 30,\, 37,\, 38,\, 39,\, 40,\, 45,\, 46\}}\widetilde a_{2,j}, &\mbox{if } t = 2,\\
\sum_{j\in \{29,\, 31,\, 33,\, 34,\, 35,\, 38,\, 40,\, 41,\, 42,\, 43,\, 44,\, 45\}}\widetilde a_{2,j}, &\mbox{if } t \geqslant 3,
\end{cases}$$
\end{lems}
\begin{proof} For $s = t= 2$, $QP_4((3)|(1)|^3) = \langle \{[\widetilde a_{2,j}] : 1 \leqslant j \leqslant 46\}\rangle$. It is easy to see that there is a direct summand decomposition of $\Sigma_4$-modules:
$$QP_4((3)|(1)|^3) = [\Sigma_4(\widetilde a_{2,1})]\bigoplus [\Sigma_4(\widetilde a_{2,13})] \bigoplus [\Sigma_4(\widetilde a_{2,25})]\bigoplus QP_4^+((3)|(1)|^3).$$
	
We easily observe $QP_4^+((3)|(1)|^3) = \langle[\widetilde a_{2,j}] : 29 \leqslant j \leqslant 46\rangle$.
 Suppose $f \in P_4((3)|(1)|^3)$ and $[f] \in QP_4^+((3)|(1)|^3)^{\Sigma_4}$. Then, we have $f \equiv\sum_{j=29}^{46}\gamma_j\widetilde a_{2,j}$, with $\gamma_j \in \mathbb F_2$. By computing $\rho_i(f)+f$ in terms of the admissible monomials, we get 
\begin{align*}
\rho_1(f) &+ f \equiv \gamma_{33}\widetilde a_{2,29} + \gamma_{34}\widetilde a_{2,30} + \gamma_{\{33,35\}}\widetilde a_{2,32} + \gamma_{\{34,35\}}\widetilde a_{2,36}\\ 
& + \gamma_{\{31,37,39\}}\widetilde a_{2,37} + \gamma_{\{31,38,40\}}\widetilde a_{2,38} + \gamma_{\{31,37,39\}}\widetilde a_{2,39}\\ 
& + \gamma_{\{31,38,40\}}\widetilde a_{2,40} + \gamma_{\{41,42\}}\widetilde a_{2,41} + \gamma_{\{41,42\}}\widetilde a_{2,42} + \gamma_{\{43,44\}}\widetilde a_{2,43}\\ 
& + \gamma_{\{43,44\}}\widetilde a_{2,44} + \gamma_{31}\widetilde a_{2,45} + \gamma_{31}\widetilde a_{2,46} \equiv 0,\\ 
\rho_2(f) &+ f \equiv \gamma_{42}\widetilde a_{2,29} + \gamma_{\{30,31,46\}}\widetilde a_{2,30} + \gamma_{\{30,31,40\}}\widetilde a_{2,31} + \gamma_{\{32,33,42\}}\widetilde a_{2,32}\\ 
& + \gamma_{\{32,33,42\}}\widetilde a_{2,33} + \gamma_{\{34,38,46\}}\widetilde a_{2,34} + \gamma_{\{35,41\}}\widetilde a_{2,35} + \gamma_{\{36,37,46\}}\widetilde a_{2,36}\\ 
& + \gamma_{\{36,37,40\}}\widetilde a_{2,37} + \gamma_{\{34,38,40\}}\widetilde a_{2,38} + \gamma_{42}\widetilde a_{2,39} + \gamma_{\{40,46\}}\widetilde a_{2,40}\\ 
& + \gamma_{\{35,41\}}\widetilde a_{2,41} + \gamma_{\{44,45,46\}}\widetilde a_{2,44} + \gamma_{\{40,44,45\}}\widetilde a_{2,45} + \gamma_{\{40,46\}}\widetilde a_{2,46} \equiv 0,\\ 
\rho_3(f) &+ f \equiv \gamma_{\{29,30\}}\widetilde a_{2,29} + \gamma_{\{29,30\}}\widetilde a_{2,30} + \gamma_{\{32,36\}}\widetilde a_{2,32} + \gamma_{\{33,34\}}\widetilde a_{2,33}\\ 
& + \gamma_{\{33,34\}}\widetilde a_{2,34} + \gamma_{\{32,36\}}\widetilde a_{2,36} + \gamma_{\{37,38\}}\widetilde a_{2,37} + \gamma_{\{37,38\}}\widetilde a_{2,38}\\ 
& + \gamma_{\{39,40\}}\widetilde a_{2,39} + \gamma_{\{39,40\}}\widetilde a_{2,40} + \gamma_{\{41,43\}}\widetilde a_{2,41} + \gamma_{\{42,44\}}\widetilde a_{2,42}\\ 
& + \gamma_{\{41,43\}}\widetilde a_{2,43} + \gamma_{\{42,44\}}\widetilde a_{2,44} + \gamma_{\{45,46\}}\widetilde a_{2,45} + \gamma_{\{45,46\}}\widetilde a_{2,46} \equiv 0.
\end{align*}
These equalities imply $\gamma_j = 0$ for either $31\leqslant j \leqslant 36$ or $41\leqslant j \leqslant 44$ and $\gamma_j = \gamma_{29}$ for $j \in \{29,\, 30,\, 37,\, 38,\, 39,\, 40,\, 45,\, 46\}$. Hence, $f \equiv \gamma_{29}\widetilde p_{2,4}$ and 
$QP_4^+((3)|(1)|^3)^{\Sigma_4} = \langle[\widetilde p_{2,4}]\rangle.$
The lemma follows from this and \eqref{ctt21}.

For $t > 2$, we have $QP_4((3)|(1)|^{t+1}) = \langle \{[\widetilde a_{2,j}] : 1 \leqslant j \leqslant 45\}\rangle$. It is easy to see that there is a direct summand decomposition of $\Sigma_4$-modules:
$$QP_4((3)|(1)|^{t+1}) = [\Sigma_4(\widetilde a_{2,1})]\bigoplus [\Sigma_4(\widetilde a_{2,13})] \bigoplus [\Sigma_4(\widetilde a_{2,25})]\bigoplus QP_4^+((3)|(1)|^{t+1}).$$
	
We easily observe $QP_4^+((3)|(1)|^{t+1}) = \langle[\widetilde a_{2,j}] : 29 \leqslant j \leqslant 45\rangle$. We can easily verify that $\widetilde p_{2,4}$ is an $\Sigma_4$-invariant and the leading monomial is $\widetilde a_{2,45}$. Hence, if $f \in P_4((3)|(1)|^{t+1})$ and $[f] \in QP_4^+((3)|(1)|^{t+1})^{\Sigma_4}$, then there are $\gamma \in \mathbb F_2$ such that
 $f^* = f + \gamma\widetilde p_{2,4} \equiv\sum_{29\leqslant j \leqslant 44}\gamma_j\widetilde a_{2,j}$, with $\gamma_j \in \mathbb F_2$. By computing $\rho_i(f^*)+f^*$ in terms of the admissible monomials, we get 
\begin{align*}
\rho_1(f^*) &+ f^* \equiv \gamma_{\{31,33\}}\widetilde a_{2,29} + \gamma_{\{31,34\}}\widetilde a_{2,30} + \gamma_{\{33,35\}}\widetilde a_{2,32} + \gamma_{\{34,35\}}\widetilde a_{2,35}\\ 
& + \gamma_{\{37,39\}}\widetilde a_{2,36} + \gamma_{\{38,40\}}\widetilde a_{2,37} + \gamma_{\{37,39\}}\widetilde a_{2,38} + \gamma_{\{38,40\}}\widetilde a_{2,39}\\ 
& + \gamma_{\{41,42\}}\widetilde a_{2,41} + \gamma_{\{41,42\}}\widetilde a_{2,42} + \gamma_{\{43,44\}}\widetilde a_{2,43} + \gamma_{\{43,44\}}\widetilde a_{2,44} \equiv 0,\\ 
\rho_2(f^*) &+ f^* \equiv \gamma_{\{40,42\}}\widetilde a_{2,29} + \gamma_{\{30,31,40\}}\widetilde a_{2,30} + \gamma_{\{30,31,40\}}\widetilde a_{2,31}\\ 
& + \gamma_{\{32,33,42\}}\widetilde a_{2,32} + \gamma_{\{32,33,42\}}\widetilde a_{2,33} + \gamma_{\{34,38\}}\widetilde a_{2,34} + \gamma_{\{36,37\}}\widetilde a_{2,35}\\ 
& + \gamma_{\{36,37\}}\widetilde a_{2,36} + \gamma_{\{34,38\}}\widetilde a_{2,37} + \gamma_{\{40,42\}}\widetilde a_{2,38} + \gamma_{\{35,41\}}\widetilde a_{2,40}\\ 
& + \gamma_{\{35,41\}}\widetilde a_{2,41} + \gamma_{44}\widetilde a_{2,44} + \gamma_{44}\widetilde a_{2,45} \equiv 0,\\ 
\rho_3(f^*) &+ f^* \equiv \gamma_{\{29,30\}}\widetilde a_{2,29} + \gamma_{\{29,30\}}\widetilde a_{2,30} + \gamma_{\{32,36\}}\widetilde a_{2,32} + \gamma_{\{33,34\}}\widetilde a_{2,33}\\ 
& + \gamma_{\{33,34\}}\widetilde a_{2,34} + \gamma_{\{32,36\}}\widetilde a_{2,35} + \gamma_{\{37,38\}}\widetilde a_{2,36} + \gamma_{\{37,38\}}\widetilde a_{2,37}\\ 
& + \gamma_{\{39,40\}}\widetilde a_{2,38} + \gamma_{\{39,40\}}\widetilde a_{2,39} + \gamma_{\{41,43\}}\widetilde a_{2,41} + \gamma_{\{42,44\}}\widetilde a_{2,42}\\ 
& + \gamma_{\{41,43\}}\widetilde a_{2,43} + \gamma_{\{42,44\}}\widetilde a_{2,44} \equiv 0.
\end{align*}
Computing from these equalities gives $\gamma_j = 0$ for $29\leqslant j \leqslant 44$. Hence, $f^* \equiv 0$ and $f \equiv  \gamma\widetilde p_{2,4}$. This implies
$QP_4^+((3)|(1)|^{t+1})^{\Sigma_4} = \langle[\widetilde p_{2,4}]\rangle.$
The lemma follows from this and \eqref{ctt21}.
\end{proof}

\begin{lems}\label{bdt23} We have
$QP_4((3)|^{2}|(1)|^{t+1})^{\Sigma_4} = \langle \{[\widetilde p_{3,u}] : 1 \leqslant u \leqslant 9\} \rangle,$ where $\widetilde p_{3,u}$ is defined by \eqref{ctt21} for $1 \leqslant u \leqslant 5$ and 
\begin{align*}
\widetilde p_{3,6} &= \mbox{$\sum_{j\in \mathbb J_{3,6} = \{29,\, 30,\, 31,\, 78,\, 79,\, 82\}}$}\widetilde a_{3,j},\\ 
\widetilde p_{3,7} &=\mbox{$\sum_{j\in \mathbb J_{3,7} = \{29,\, 30,\, 31,\, 32,\, 33,\, 34,\, 54,\, 55,\, 58,\, 59,\, 60,\, 63,\, 64,\, 80,\, 81,\, 82,\, 83,\, 84\}}$}\widetilde a_{3,j},\\ 
\widetilde p_{3,8} &=\mbox{$\sum_{j\in \mathbb J_{3,8} = \{30,\, 33,\, 34,\, 54,\, 57,\, 60,\, 61,\, 62,\, 64,\, 65,\, 72,\, 73,\, 74,\, 75,\, 76,\, 85,\, 86,\, 87,\, 88,\, 89,\, 90\}}$}\widetilde a_{3,j},\\ 
\widetilde p_{3,9} &=\mbox{$\sum_{j\in\mathbb J_{3,9} = \{29,\, 30,\, 31,\, 32,\, 33,\, 34,\, 35,\, 54,\, 55,\, 56,\, 57,\, 58,\, 59,\, 60,\, 61,\, 62,\, 63,\, 64,\, 65\}\atop\hskip5.8cm\cup \{66,\, 67,\, 68,\, 69,\, 70\}}$}\widetilde a_{3,j}. 
\end{align*}
Consequently, $\dim	QP_4((3)|^{2}|(1)|^{t+1})^{\Sigma_4} = 9$ for any $t \geqslant 2$.
\end{lems}

\begin{proof} 
For $s = 3$, we have $QP_4((3)|^2|(1)|^{t+1}) = \langle \{[\widetilde a_{3,j}]: 1 \leqslant j \leqslant 90\}\rangle$. We have a direct summand decomposition of $\Sigma_4$-modules:
\begin{align*}QP_4((3)|^2|(1)|^{t+1}) = [\Sigma_4(\widetilde a_{3,1})]&\bigoplus [\Sigma_4(\widetilde a_{3,13})] \bigoplus [\Sigma_4(\widetilde a_{3,25})]\\ &\bigoplus [\Sigma_4(\widetilde a_{3,36})]\bigoplus [\Sigma_4(\widetilde a_{3,48})]\bigoplus \mathcal M_{t,3},
\end{align*}
where $\mathcal M_{t,3} = \langle \{[\widetilde a_{3,j}]: 29\leqslant j \leqslant 35 \mbox{ or } 54\leqslant j \leqslant 90\}\rangle$. By using \eqref{ctt21} we need only to prove $\mathcal M_{t,3}^{\Sigma_4} = \langle \{[\widetilde p_{3,u}] : 6 \leqslant u \leqslant 9\}\rangle$. By a direct computation we can verify that $[\widetilde p_{3,u}] \in \mathcal M_{t,3}^{\Sigma_4}$ for $6 \leqslant u \leqslant 9$. The leading monomials of
$\widetilde p_{3,6}$, $\widetilde p_{3,7}$, $\widetilde p_{3,8}$, $\widetilde p_{3,9}$ are $\widetilde a_{3,79}$, $\widetilde a_{3,84}$, $\widetilde a_{3,77}$, $\widetilde a_{3,70}$ respectively. 
	
Suppose $f \in P_4((3)|^2|(1)|^3)$ and $[f] \in \mathcal M_{t,3}^{\Sigma_4}$. Then, there are $\gamma_u \in \mathbb F_2$ with $6 \leqslant u \leqslant 9$ such that 
$$h := f + \sum_{6 \leqslant u \leqslant 9}\gamma_{u}\widetilde p_{3,u} \equiv\sum_{29 \leqslant j \leqslant 35}\gamma_j\widetilde a_{3,j}+ \sum_{54 \leqslant j \leqslant 90, j \ne 70, 77, 79, 84}\gamma_j\widetilde a_{3,j},$$
with suitable $\gamma_j \in \mathbb F_2$. Since $[f],\, [\widetilde p_{3,u}] \in \mathcal M_{t,3}^{\Sigma_4}$ for  $6 \leqslant u \leqslant 9$, we have $[h] \in \mathcal M_{t,3}^{\Sigma_4}$.
A direct computation shows
\begin{align*}
\rho_1(h)+h &\equiv \gamma_{\{29,31,33,54\}}\widetilde a_{3,29} + \gamma_{\{30,31,34,55\}}\widetilde a_{3,30} + \gamma_{\{32,33,35,57\}}\widetilde a_{3,32}\\ 
&\quad \ + \gamma_{\{29,31,33,54\}}\widetilde a_{3,54} + \gamma_{\{30,31,34,55\}}\widetilde a_{3,55} + \gamma_{\{72,74\}}\widetilde a_{3,56}\\ 
&\quad \ + \gamma_{\{32,33,35,57\}}\widetilde a_{3,57} + \gamma_{\{34,35,58,61\}}\widetilde a_{3,58} + \gamma_{\{59,63\}}\widetilde a_{3,59}\\ 
&\quad \ + \gamma_{\{60,64\}}\widetilde a_{3,60} + \gamma_{\{34,35,58,61\}}\widetilde a_{3,61} + \gamma_{\{62,65\}}\widetilde a_{3,62} + \gamma_{\{59,63\}}\widetilde a_{3,63}\\ 
&\quad \ + \gamma_{\{60,64\}}\widetilde a_{3,64} + \gamma_{\{62,65\}}\widetilde a_{3,65} + \gamma_{\{66,67\}}\widetilde a_{3,66} + \gamma_{\{66,67\}}\widetilde a_{3,67}\\ 
&\quad \ + \gamma_{\{68,69\}}\widetilde a_{3,68} + \gamma_{\{68,69\}}\widetilde a_{3,69} + \gamma_{\{73,75\}}\widetilde a_{3,71} + \gamma_{76}\widetilde a_{3,76} + \gamma_{76}\widetilde a_{3,77}\\ 
&\quad \ + \gamma_{\{78,80,81\}}\widetilde a_{3,82} + \gamma_{83}\widetilde a_{3,83} + \gamma_{83}\widetilde a_{3,84} + \gamma_{\{31,33,73,78\}}\widetilde a_{3,85}\\ 
&\quad \ + \gamma_{\{31,34,72\}}\widetilde a_{3,86} + \gamma_{\{33,35,74,80\}}\widetilde a_{3,87} + \gamma_{\{34,35,75,81\}}\widetilde a_{3,88}\\ 
&\quad \ + \gamma_{\{89,90\}}\widetilde a_{3,89} + \gamma_{\{89,90\}}\widetilde a_{3,90} \equiv 0.
\end{align*}
By computing from this equality we get $\gamma_{76} = \gamma_{83} = 0$  and $\gamma_i = \gamma_j $ for $(i,j) \in \{(59,63),\, (60,64),\, (62,65),\, (66,67),\, (68,69),\, (72,74),\, (73,75),\, (89,90)\}$. By using these results we get
\begin{align*}
\rho_2(h)+h &\equiv \gamma_{\{86,87\}}\widetilde a_{3,29} + \gamma_{\{30,31,86\}}\widetilde a_{3,30} + \gamma_{\{30,31,86\}}\widetilde a_{3,31} + \gamma_{\{32,33,87\}}\widetilde a_{3,32}\\ 
&\quad \ + \gamma_{\{32,33,87\}}\widetilde a_{3,33} + \gamma_{\{34,60\}}\widetilde a_{3,34} + \gamma_{\{35,66\}}\widetilde a_{3,35} + \gamma_{\{54,56,60,66\}}\widetilde a_{3,54}\\ 
&\quad \ + \gamma_{\{55,60,72\}}\widetilde a_{3,55} + \gamma_{\{54,56,60,66\}}\widetilde a_{3,56} + \gamma_{\{57,66,72\}}\widetilde a_{3,57}\\ 
&\quad \ + \gamma_{\{58,59\}}\widetilde a_{3,58} + \gamma_{\{58,59\}}\widetilde a_{3,59} + \gamma_{\{34,60\}}\widetilde a_{3,60} + \gamma_{\{61,62\}}\widetilde a_{3,61}\\ 
&\quad \ + \gamma_{\{61,62\}}\widetilde a_{3,62} + \gamma_{\{59,60,62,66\}}\widetilde a_{3,63} + \gamma_{\{59,60,62,66\}}\widetilde a_{3,65}\\ 
&\quad \ + \gamma_{\{35,66\}}\widetilde a_{3,66} + \gamma_{68}\widetilde a_{3,69} + \gamma_{68}\widetilde a_{3,70} + \gamma_{\{71,73\}}\widetilde a_{3,71} + \gamma_{\{55,60,72\}}\widetilde a_{3,72}\\ 
&\quad \ + \gamma_{\{71,73\}}\widetilde a_{3,73} + \gamma_{\{57,66,72\}}\widetilde a_{3,74} + \gamma_{73}\widetilde a_{3,75} + \gamma_{73}\widetilde a_{3,76} + \gamma_{\{78,82\}}\widetilde a_{3,78}\\ 
&\quad \ + \gamma_{81}\widetilde a_{3,81} + \gamma_{\{78,82\}}\widetilde a_{3,82} + \gamma_{81}\widetilde a_{3,83} + \gamma_{\{60,66,86,87\}}\widetilde a_{3,85}\\ 
&\quad \ + \gamma_{\{88,89\}}\widetilde a_{3,88} + \gamma_{\{88,89\}}\widetilde a_{3,89} + \gamma_{\{60,66\}}\widetilde a_{3,90} \equiv 0.
\end{align*}
By computing from the above equalities we get $\gamma_{j} = 0$ for $j \in \{$68,\, 69,\, 71,\, 73,\, 75,\, 81$\}$  and $\gamma_i = \gamma_j $ for $(i,j) \in \{$(34,35),\, (34,60),\, (34,66),\, (54,56),\, (58,61),\, (58,59),\, (58,62),\, (78,80),\, (78,82),\, (86,87),\, (88,89)$\}$. By using the above results we get
\begin{align*}
\rho_3(h)+h &\equiv \gamma_{\{29,30,88\}}\widetilde a_{3,29} + \gamma_{\{29,30,88\}}\widetilde a_{3,30} + \gamma_{\{32,58\}}\widetilde a_{3,32} + \gamma_{\{33,34\}}\widetilde a_{3,33}\\ 
&\quad \ + \gamma_{\{33,34\}}\widetilde a_{3,34} + \gamma_{\{54,55,88\}}\widetilde a_{3,54} + \gamma_{\{54,55,88\}}\widetilde a_{3,55} + \gamma_{54}\widetilde a_{3,56}\\ 
&\quad \ + \gamma_{\{57,58\}}\widetilde a_{3,57} + \gamma_{\{32,58\}}\widetilde a_{3,58} + \gamma_{\{34,58,88\}}\widetilde a_{3,59} + \gamma_{\{34,58,88\}}\widetilde a_{3,60}\\ 
&\quad \ + \gamma_{\{57,58\}}\widetilde a_{3,61} + \gamma_{58}\widetilde a_{3,62} + \gamma_{\{34,58,88\}}\widetilde a_{3,63} + \gamma_{\{34,58,88\}}\widetilde a_{3,64}\\ 
&\quad \ + \gamma_{58}\widetilde a_{3,65} + \gamma_{34}\widetilde a_{3,66} + \gamma_{34}\widetilde a_{3,67} + \gamma_{34}\widetilde a_{3,68} + \gamma_{34}\widetilde a_{3,69} + \gamma_{54}\widetilde a_{3,71}\\ 
&\quad \ + \gamma_{72}\widetilde a_{3,72} + \gamma_{72}\widetilde a_{3,73} + \gamma_{72}\widetilde a_{3,74} + \gamma_{72}\widetilde a_{3,75} + \gamma_{58}\widetilde a_{3,76} + \gamma_{58}\widetilde a_{3,77}\\ 
&\quad \ + \gamma_{78}\widetilde a_{3,78} + \gamma_{78}\widetilde a_{3,79} + \gamma_{78}\widetilde a_{3,80} + \gamma_{78}\widetilde a_{3,81} + \gamma_{\{85,86\}}\widetilde a_{3,85}\\ 
&\quad \ + \gamma_{\{85,86\}}\widetilde a_{3,86} + \gamma_{\{86,88\}}\widetilde a_{3,87} + \gamma_{\{86,88\}}\widetilde a_{3,88} \equiv 0.
\end{align*}

From the above equalities we get $\gamma_j =0$ for all $j$, $h \equiv 0$ and 
$$f \equiv \sum_{6 \leqslant u \leqslant 9}\gamma_{u}\widetilde p_{3,u}.$$ 
The lemma is proved.
\end{proof}

\begin{lems}\label{bdt24} For $s \geqslant 4$ and $t \geqslant 2$, we have
$QP_4((3)|^{s-1}|(1)|^{t+1})^{\Sigma_4} = \langle \{[\widetilde p_{s,u}] : 1 \leqslant u \leqslant 13\} \rangle$,
where $\widetilde p_{s,u}$ is defined by \eqref{ctt21} for $1 \leqslant u \leqslant 6$, $\widetilde p_{s,7} = \widetilde a_{s,101}$ and
\begin{align*}
\widetilde p_{s,8} &= \begin{cases}
\sum_{j \in\mathbb J_{4,8}= \{29,\, 30,\, 31,\, 54,\, 55,\, 78,\, 79,\, 89,\, 90,\, 104\}}\widetilde a_{4,j},& \mbox{if } s =4,\\ 
\sum_{j \in\mathbb J_{s,8}=\{99,\, 100,\, 102,\, 103,\, 104,\, 105\}}\widetilde a_{s,j}, & \mbox{if } s > 4, \end{cases}\\ 
\widetilde p_{s,9} &= \begin{cases}
\sum_{j \in\mathbb J_{4,9}=\{71,\, 72,\, 73,\, 78,\, 79,\, 82,\, 89,\, 90,\, 91,\, 98,\, 99,\, 100\}}\widetilde a_{4,j}, & \mbox{if } s =4,\\
\sum_{j \in\mathbb J_{s,9}=\{71,\, 72,\, 73,\, 78,\, 79,\, 82,\, 89,\, 90,\, 91,\, 92,\, 97,\, 98,\, 99,\, 100,\, 104,\, 105\}}\widetilde a_{s,j}, & \mbox{if } s > 4,
\end{cases}\\ 
\widetilde p_{s,10} &= \begin{cases}
\sum_{j \in\mathbb J_{4,10}= \{30,\, 31,\, 55,\, 56,\, 72,\, 73,\, 89,\, 92,\, 97,\, 102,\, 103,\, 105\}}\widetilde a_{4,j},  & \mbox{if } s =4,\\
\sum_{j \in\mathbb J_{s,10}= \{29,\, 30,\, 31,\, 54,\, 55,\, 56,\, 78,\, 79,\, 82,\, 89,\, 90,\, 91,\, 92,\, 98,\, 99,\, 103\}}\widetilde a_{s,j}, & \mbox{if } s > 4,\end{cases}\\ 
\widetilde p_{s,11} &= \begin{cases}
\sum_{j \in\mathbb J_{4,11}= \{71,\, 72,\ldots,\, 83,\, 84,\, 89,\, 90,\ldots ,\, 97,\, 98\}}\widetilde a_{4,j},& \mbox{if } s =4,\\ 
\sum_{j \in\mathbb J_{s,11}= \{74,\, 75,\, 76,\, 77,\, 80,\, 81,\, 83,\, 84,\, 93,\, 94,\, 95,\, 96,\, 99,\, 100,\, 104,\, 105\}}\widetilde a_{s,j}, & \mbox{if } s > 4,\end{cases}\\
\widetilde p_{s,12} &= \mbox{$\sum_{j \in \mathbb J_{s,12} = \{32,\, 33,\, 34,\, 57,\ldots,\, 65,\, 78,\ldots,\, 84,\, 89,\ldots,\, 96,\, 98\}}$}\widetilde a_{s,j},\\
\widetilde p_{s,13} &= \begin{cases}\sum_{j \in\mathbb J_{4,13}= \{29,\, 30,\, 31,\, 35,\, 54,\, 55,\, 56,\, 66,\, 67,\, 68,\, 69,\, 70,\, 78,\, 79,\, 80,\, 81,\, 82\}\atop\hskip3cm\cup \{83,\, 84,\, 89,\, 90,\, 91,\, 92,\, 93,\, 94,\, 95,\, 96,\, 98\}}\widetilde a_{4,j},& \mbox{if } s =4,\\
\sum_{j \in\mathbb J_{s,13}= \{35,\, 66,\ldots,\, 70,\, 80,\, 81,\, 83,\, 84,\, 93,\ldots,\, 96,\, 99,\, 103\}}\widetilde a_{s,j}, & \mbox{if } s > 4.\end{cases} 
\end{align*}		
\end{lems}
\begin{proof} 
For $s > 3$, we have $QP_4((3)|^{s-1}|(1)|^{t+1}) = \langle \{[a_{s,j}] : 1 \leqslant j \leqslant 105\}\rangle$. We have a direct summand decomposition of $\Sigma_4$-modules:
\begin{align*}QP_4((3)|^{s-1}|(1)|^{t+1}) &= [\Sigma_4(\widetilde a_{s,1})]\bigoplus [\Sigma_4(\widetilde a_{s,13})] \bigoplus [\Sigma_4(\widetilde a_{s,25})]\\ &\quad \bigoplus [\Sigma_4(\widetilde a_{s,36})]\bigoplus [\Sigma_4(\widetilde a_{s,48})]\bigoplus [\Sigma_4(\widetilde a_{s,85})]\bigoplus \mathcal M_{t,s},
\end{align*}
where $\mathcal M_{t,s} = \langle \{[\widetilde a_{s,j}]: 29\leqslant j \leqslant 35 \mbox{ or } 54\leqslant j \leqslant 84 \mbox{ or } 89\leqslant j \leqslant 105\}\rangle$. By using \eqref{ctt21} we need only to prove $\mathcal M_{t,s}^{\Sigma_4} = \langle \{ [\widetilde p_{s,u}] : 7 \leqslant t \leqslant 13\}\rangle$. We prove this for $s > 4$. The case $s = 4$ is realized by a similar computation. 

By a direct computation we can easily verify that $[\widetilde p_{s,u}] \in \mathcal M_{t,s}^{\Sigma_4}$ for $7 \leqslant t \leqslant 13$. 
The leading monomials of $\widetilde p_{s,u}$, for $8 \leqslant u \leqslant 13$ are respectively $\widetilde a_{s,102}$, $\widetilde a_{s,97}$, $\widetilde a_{s,56}$, $\widetilde a_{s,77}$, $\widetilde a_{s,65}$, $\widetilde a_{s,70}$. 
	
Suppose $g \in P_4((3)|^{s-1}|(1)|^{t+1})$ and $[g] \in \mathcal M_{t,s}^{\Sigma_4}$. 
Then, there are $\gamma_u \in \mathbb F_2$ with $7 \leqslant u \leqslant 13$ such that 
\begin{align*}h = g &+ \sum_{7 \leqslant u \leqslant 13}\gamma_{u}\widetilde p_{s,u} \equiv\sum_{29 \leqslant j \leqslant 35}\gamma_j\widetilde a_{s,j}\\ &+ \sum_{54 \leqslant j \leqslant 84, j \ne 56, 65, 70, 77}\gamma_j\widetilde a_{s,j} + \sum_{89 \leqslant j \leqslant 105, j \ne 97,101,102}\gamma_j\widetilde a_{s,j},
\end{align*}
with suitable $\gamma_j \in \mathbb F_2$. Since $[g],\, [\widetilde p_{s,u}] \in \mathcal M_{t,s}^{\Sigma_4}$ for  $7 \leqslant u \leqslant 13$, we have $[h] \in \mathcal M_{t,s}^{\Sigma_4}$.
By computing $\rho_i(h)+h$ in terms of the admissible monomials we obtain	
\begin{align*}
\rho_1(h)+h &\equiv \gamma_{\{29,54\}}\widetilde a_{s,29} + \gamma_{\{30,55\}}\widetilde a_{s,30} + \gamma_{\{32,57\}}\widetilde a_{s,32} + \gamma_{\{29,54\}}\widetilde a_{s,54}\\
&\quad + \gamma_{\{30,55\}}\widetilde a_{s,55} + \gamma_{\{32,57\}}\widetilde a_{s,57} + \gamma_{\{58,61\}}\widetilde a_{s,58} + \gamma_{\{59,63\}}\widetilde a_{s,59}\\
&\quad + \gamma_{\{60,64\}}\widetilde a_{s,60} + \gamma_{\{58,61\}}\widetilde a_{s,61} + \gamma_{62}\widetilde a_{s,62} + \gamma_{\{59,63\}}\widetilde a_{s,63}\\
&\quad + \gamma_{\{60,64\}}\widetilde a_{s,64} + \gamma_{62}\widetilde a_{s,65} + \gamma_{\{66,67\}}\widetilde a_{s,66} + \gamma_{\{66,67\}}\widetilde a_{s,67}\\
&\quad + \gamma_{\{68,69\}}\widetilde a_{s,68} + \gamma_{\{68,69\}}\widetilde a_{s,69} + \gamma_{71}\widetilde a_{s,71} + \gamma_{76}\widetilde a_{s,76} + \gamma_{76}\widetilde a_{s,77}\\
&\quad + \gamma_{\{82,98\}}\widetilde a_{s,82} + \gamma_{\{83,84\}}\widetilde a_{s,83} + \gamma_{\{83,84\}}\widetilde a_{s,84} + \gamma_{\{31,33,73,78\}}\widetilde a_{s,89}\\
&\quad + \gamma_{\{31,34,72,79\}}\widetilde a_{s,90} + \gamma_{\{91,92\}}\widetilde a_{s,91} + \gamma_{\{91,92\}}\widetilde a_{s,92}\\
&\quad + \gamma_{\{33,35,74,80\}}\widetilde a_{s,93} + \gamma_{\{34,35,75,81\}}\widetilde a_{s,94} + \gamma_{\{95,96\}}\widetilde a_{s,95}\\
&\quad + \gamma_{\{95,96\}}\widetilde a_{s,96} + \gamma_{71}\widetilde a_{s,97} + \gamma_{\{82,98\}}\widetilde a_{s,98} + \gamma_{\{72,74,104\}}\widetilde a_{s,99}\\
&\quad + \gamma_{\{73,75,105\}}\widetilde a_{s,100} + \gamma_{\{78,79,80,81,104,105\}}\widetilde a_{s,103} \equiv 0.
\end{align*}
By computing from this equality we get $\gamma_{62} = \gamma_{71} = \gamma_{76} = 0$  and $\gamma_i = \gamma_j $ for $(i,j) \in \{$(29,54),\, (30,55),\, (32,57),\, (58,61),\, (59,63),\, (60,64),\, (66,67),\, (68,69),\, (82,98),\, (83,84),\, (91,92),\, (95,96)$\}$. By using these results we get
\begin{align*}
\rho_2(h)+h &\equiv \gamma_{\{30,31\}}\widetilde a_{s,30} + \gamma_{\{30,31\}}\widetilde a_{s,31} + \gamma_{\{32,33\}}\widetilde a_{s,32} + \gamma_{\{32,33\}}\widetilde a_{s,33}\\
&\quad + \gamma_{\{34,60\}}\widetilde a_{s,34} + \gamma_{\{35,66\}}\widetilde a_{s,35} + \gamma_{29}\widetilde a_{s,54} + \gamma_{29}\widetilde a_{s,56} + \gamma_{\{58,59\}}\widetilde a_{s,58}\\
&\quad + \gamma_{\{58,59\}}\widetilde a_{s,59} + \gamma_{\{34,60\}}\widetilde a_{s,60} + \gamma_{58}\widetilde a_{s,61} + \gamma_{58}\widetilde a_{s,62} + \gamma_{59}\widetilde a_{s,63}\\
&\quad + \gamma_{59}\widetilde a_{s,65} + \gamma_{\{35,66\}}\widetilde a_{s,66} + \gamma_{68}\widetilde a_{s,69} + \gamma_{68}\widetilde a_{s,70} + \gamma_{73}\widetilde a_{s,71}\\
&\quad + \gamma_{\{30,60,72,90\}}\widetilde a_{s,72} + \gamma_{73}\widetilde a_{s,73} + \gamma_{\{32,66,74,93\}}\widetilde a_{s,74} + \gamma_{75}\widetilde a_{s,75}\\
&\quad + \gamma_{75}\widetilde a_{s,76} + \gamma_{\{78,82\}}\widetilde a_{s,78} + \gamma_{\{81,83\}}\widetilde a_{s,81} + \gamma_{\{78,82\}}\widetilde a_{s,82}\\
&\quad + \gamma_{\{81,83\}}\widetilde a_{s,83} + \gamma_{\{89,91\}}\widetilde a_{s,89} + \gamma_{\{30,60,72,90\}}\widetilde a_{s,90} + \gamma_{\{89,91\}}\widetilde a_{s,91}\\
&\quad + \gamma_{\{30,32,82\}}\widetilde a_{s,92} + \gamma_{\{32,66,74,93\}}\widetilde a_{s,93} + \gamma_{\{94,95\}}\widetilde a_{s,94} + \gamma_{\{94,95\}}\widetilde a_{s,95}\\
&\quad + \gamma_{\{60,66,83\}}\widetilde a_{s,96} + \gamma_{\{30,32,60,66,82,83\}}\widetilde a_{s,99} + \gamma_{100}\widetilde a_{s,100} + \gamma_{100}\widetilde a_{s,102}\\
&\quad + \gamma_{\{82,83,103,104\}}\widetilde a_{s,103} + \gamma_{\{82,83,103,104\}}\widetilde a_{s,104} \equiv 0.
\end{align*}

Computing from the above equalities gives $\gamma_{j} = 0$ for $j \in \{$29,\, 54,\, 58,\, 59,\, 61,\, 63,\, 68,\, 69,\, 73,\, 75,\, 100,\, 105$\}$  and $\gamma_i = \gamma_j $ for $(i,j) \in \{$(30,31),\, (32,33),\, (34,60),\, (35,66),\, (78,82),\, (81,83),\, (89,91),\, (94,95)$\}$. By using the above results we get
\begin{align*}
\rho_3(h)+h &\equiv \gamma_{30}\widetilde a_{s,29} + \gamma_{30}\widetilde a_{s,30} + \gamma_{32}\widetilde a_{s,32} + \gamma_{\{32,34\}}\widetilde a_{s,33} + \gamma_{\{32,34\}}\widetilde a_{s,34}\\
&\quad + \gamma_{30}\widetilde a_{s,54} + \gamma_{30}\widetilde a_{s,55} + \gamma_{32}\widetilde a_{s,57} + \gamma_{32}\widetilde a_{s,58} + \gamma_{34}\widetilde a_{s,59} + \gamma_{34}\widetilde a_{s,60}\\
&\quad + \gamma_{32}\widetilde a_{s,61} + \gamma_{34}\widetilde a_{s,63} + \gamma_{34}\widetilde a_{s,64} + \gamma_{35}\widetilde a_{s,66} + \gamma_{35}\widetilde a_{s,67} + \gamma_{35}\widetilde a_{s,68}\\
&\quad + \gamma_{35}\widetilde a_{s,69} + \gamma_{89}\widetilde a_{s,71} + \gamma_{72}\widetilde a_{s,72} + \gamma_{72}\widetilde a_{s,73} + \gamma_{74}\widetilde a_{s,74} + \gamma_{74}\widetilde a_{s,75}\\
&\quad + \gamma_{94}\widetilde a_{s,76} + \gamma_{94}\widetilde a_{s,77} + \gamma_{\{78,79\}}\widetilde a_{s,78} + \gamma_{\{78,79\}}\widetilde a_{s,79} + \gamma_{\{80,81\}}\widetilde a_{s,80}\\
&\quad + \gamma_{\{80,81\}}\widetilde a_{s,81} + \gamma_{\{89,90\}}\widetilde a_{s,89} + \gamma_{\{89,90\}}\widetilde a_{s,90} + \gamma_{89}\widetilde a_{s,91} + \gamma_{89}\widetilde a_{s,92}\\
&\quad + \gamma_{\{93,94\}}\widetilde a_{s,93} + \gamma_{\{93,94\}}\widetilde a_{s,94} + \gamma_{94}\widetilde a_{s,95} + \gamma_{94}\widetilde a_{s,96} + \gamma_{89}\widetilde a_{s,97}\\
&\quad + \gamma_{99}\widetilde a_{s,99} + \gamma_{99}\widetilde a_{s,100} + \gamma_{104}\widetilde a_{s,104} + \gamma_{104}\widetilde a_{s,105} \equiv 0.
\end{align*}
	
By computing from the above equalities we get $\gamma_j =0$ for all $j$. This implies $h \equiv 0$ and 
$ g \equiv  \sum_{7 \leqslant u \leqslant 13}\gamma_{u}\widetilde p_{s,u}.$
The lemma is proved.
\end{proof}

\begin{proof}[Proof of Theorem \ref{dlt21}] 
Suppose $f \in QP_4((3)|^{s-1}|(1)|^{t+1})$ such that the class $[f] \in QP_4((3)|^{s-1}|(1)|^{t+1})^{GL_4}$. Then, $[f] \in QP_4((3)|^{s-1}|(1)|^{t+1})^{\Sigma_4}$. 

For $s = t = 2$, by Lemma \ref{bdt22}, we have 
$f \equiv \sum_{1 \leqslant u \leqslant 4}\gamma_u\widetilde p_{2,u},$
where $\gamma_u \in \mathbb F_2$. By computing $\rho_4(f)+f$ in terms of the admissible monomials, we get
\begin{align*}
\rho_4(f)+f &\equiv \gamma_{1}\widetilde a_{2,1} + \gamma_{1}\widetilde a_{2,2} + \gamma_{1}\widetilde a_{2,3} + \gamma_{\{1,2\}}\widetilde a_{2,8} + \gamma_{\{1,2\}}\widetilde a_{2,9} + \gamma_{2}\widetilde a_{2,13}\\ 
&\quad + \gamma_{2}\widetilde a_{2,14} + \gamma_{2}\widetilde a_{2,15} + \gamma_{\{2,3\}}\widetilde a_{2,17} + \gamma_{\{2,3\}}\widetilde a_{2,18} + \gamma_{3}\widetilde a_{2,25}\\ 
&\quad + \gamma_{2}\widetilde a_{2,29} + \gamma_{2}\widetilde a_{2,30} + \gamma_{\{2,3\}}\widetilde a_{2,32} + \gamma_{\{2,3\}}\widetilde a_{2,36} + \gamma_{1}\widetilde a_{2,37}\\ 
&\quad + \gamma_{1}\widetilde a_{2,38} + \gamma_{1}\widetilde a_{2,39} + \gamma_{1}\widetilde a_{2,40} + \gamma_{1}\widetilde a_{2,45} + \gamma_{1}\widetilde a_{2,46} \equiv 0.
\end{align*}
This equality implies $\gamma_1 = \gamma_2 = \gamma_3 = 0$. Hence, $f \equiv \gamma_4\widetilde p_{2,4} = \gamma_4\xi_{2,2}$. The theorem is proved for $s = t =2$.

For $s = 2$ and $t \geqslant 3$, by Lemma \ref{bdt22}, we have 
$f \equiv \sum_{1 \leqslant u \leqslant 4}\gamma_u\widetilde p_{2,u},$
where $\gamma_u \in \mathbb F_2$. By computing $\rho_4(f)+f$ in terms of the admissible monomials, we get
\begin{align*}
\rho_4(f)&+f \equiv \gamma_{1}\widetilde a_{2,1} + \gamma_{1}\widetilde a_{2,2} + \gamma_{\{1,4\}}\widetilde a_{2,3} + \gamma_{\{1,2\}}\widetilde a_{2,8} + \gamma_{\{1,2\}}\widetilde a_{2,9} + \gamma_{2}\widetilde a_{2,13}\\ 
& + \gamma_{\{2,4\}}\widetilde a_{2,14} + \gamma_{\{2,4\}}\widetilde a_{2,15} + \gamma_{\{2,3\}}\widetilde a_{2,17} + \gamma_{\{2,3\}}\widetilde a_{2,18} + \gamma_{\{3,4\}}\widetilde a_{2,25}\\ 
& + \gamma_{\{1,2\}}\widetilde a_{2,29} + \gamma_{\{1,2\}}\widetilde a_{2,30} + \gamma_{\{2,3\}}\widetilde a_{2,32} + \gamma_{\{2,3\}}\widetilde a_{2,35} + \gamma_{4}\widetilde a_{2,36}\\ 
& + \gamma_{4}\widetilde a_{2,37} + \gamma_{4}\widetilde a_{2,41} + \gamma_{4}\widetilde a_{2,43} \equiv 0.
\end{align*}
This equality implies $\gamma_1 = \gamma_2 = \gamma_3 = \gamma_4 = 0$. Hence, $f \equiv 0$. The theorem is proved for $s =2$ and $t \geqslant 3$.

For $s = 3$, by using Lemma \ref{bdt23}, we have 
$ f \equiv \sum_{1 \leqslant u \leqslant 9}\gamma_u\widetilde p_{3,u},$
where $\gamma_u \in \mathbb F_2$. By computing $\rho_4(f)+f$ in terms of the admissible monomials we get
\begin{align*}
\rho_4(f)+f &\equiv \gamma_{\{1,4\}}\widetilde a_{3,1} + \gamma_{\{1,4\}}\widetilde a_{3,2} + \gamma_{\{1,6,7,9\}}\widetilde a_{3,3} + \gamma_{\{1,2\}}\widetilde a_{3,8} + \gamma_{\{1,2\}}\widetilde a_{3,9}\\ 
&\quad + \gamma_{\{2,5\}}\widetilde a_{3,13} + \gamma_{\{2,7,8,9\}}\widetilde a_{3,14} + \gamma_{\{2,7,8,9\}}\widetilde a_{3,15} + \gamma_{\{2,3\}}\widetilde a_{3,17}\\ 
&\quad + \gamma_{\{2,3\}}\widetilde a_{3,18} + \gamma_{\{3,9\}}\widetilde a_{3,25} + \gamma_{\{7,8,9\}}\widetilde a_{3,29} + \gamma_{\{7,9\}}\widetilde a_{3,30} + \gamma_{\{8,9\}}\widetilde a_{3,32}\\ 
&\quad + \gamma_{4}\widetilde a_{3,38} + \gamma_{4}\widetilde a_{3,39} + \gamma_{\{4,8\}}\widetilde a_{3,40} + \gamma_{\{4,8\}}\widetilde a_{3,41} + \gamma_{\{4,5\}}\widetilde a_{3,44}\\ 
&\quad + \gamma_{\{5,8\}}\widetilde a_{3,49} + \gamma_{\{5,8\}}\widetilde a_{3,50} + \gamma_{\{4,5\}}\widetilde a_{3,56} + \gamma_{\{8,9\}}\widetilde a_{3,58} + \gamma_{\{6,7,9\}}\widetilde a_{3,59}\\ 
&\quad + \gamma_{\{6,7,8,9\}}\widetilde a_{3,60} + \gamma_{8}\widetilde a_{3,62} + \gamma_{\{7,9\}}\widetilde a_{3,66} + \gamma_{\{7,9\}}\widetilde a_{3,68} + \gamma_{\{4,5\}}\widetilde a_{3,71}\\ 
&\quad + \gamma_{8}\widetilde a_{3,76} + \gamma_{8}\widetilde a_{3,82} + \gamma_{7}\widetilde a_{3,83} + \gamma_{\{1,2,4,7,9\}}\widetilde a_{3,85} + \gamma_{\{1,2,4,7,8,9\}}\widetilde a_{3,86}\\ 
&\quad + \gamma_{\{2,3,5,9\}}\widetilde a_{3,87} + \gamma_{\{2,3,5,9\}}\widetilde a_{3,88} + \gamma_{8}\widetilde a_{3,89} \equiv 0.
\end{align*}
From this equality we get $\gamma_u = 0$ for $1 \leqslant u \leqslant 9$. Hence, $[f] = 0$. The theorem is proved for $s = 3$.

For $s = 4$, by using Lemma \ref{bdt24}, we have 
$ f \equiv \sum_{1 \leqslant u \leqslant 13}\gamma_u\widetilde p_{4,u},$
where $\gamma_u \in \mathbb F_2$. By computing $\rho_4(f)+f$ in terms of the admissible monomials, we obtain
\begin{align*}
\rho_4(f)+f &\equiv \gamma_{\{1,4\}}\widetilde a_{4,1} + \gamma_{\{1,4\}}\widetilde a_{4,2} + \gamma_{\{1,8,10,13\}}\widetilde a_{4,3} + \gamma_{\{1,2\}}\widetilde a_{4,8} + \gamma_{\{1,2\}}\widetilde a_{4,9}\\ 
&\quad + \gamma_{\{2,5\}}\widetilde a_{4,13} + \gamma_{\{2,12\}}\widetilde a_{4,14} + \gamma_{\{2,12\}}\widetilde a_{4,15} + \gamma_{\{2,3\}}\widetilde a_{4,17} + \gamma_{\{2,3\}}\widetilde a_{4,18}\\ 
&\quad + \gamma_{\{3,13\}}\widetilde a_{4,25} + \gamma_{\{1,2,4,9,10,11\}}\widetilde a_{4,29} + \gamma_{\{1,2,4,9,10,11\}}\widetilde a_{4,30}\\ 
&\quad + \gamma_{\{2,3,5,11\}}\widetilde a_{4,32} + \gamma_{\{4,6\}}\widetilde a_{4,38} + \gamma_{\{4,6\}}\widetilde a_{4,39} + \gamma_{\{4,9,10,11\}}\widetilde a_{4,40}\\ 
&\quad + \gamma_{\{4,9,10,11\}}\widetilde a_{4,41} + \gamma_{\{4,5\}}\widetilde a_{4,44} + \gamma_{\{5,11\}}\widetilde a_{4,49} + \gamma_{\{5,11\}}\widetilde a_{4,50}\\ 
&\quad + \gamma_{\{1,2,4,12\}}\widetilde a_{4,54} + \gamma_{\{1,2,4,12\}}\widetilde a_{4,55} + \gamma_{\{4,5,6,11,12,13\}}\widetilde a_{4,56}\\ 
&\quad + \gamma_{\{2,3,5,13\}}\widetilde a_{4,57} + \gamma_{\{2,3,5,11\}}\widetilde a_{4,58} + \gamma_{\{8,9,11,13\}}\widetilde a_{4,59}\\ 
&\quad + \gamma_{\{8,9,11,13\}}\widetilde a_{4,60} + \gamma_{\{2,3,5,13\}}\widetilde a_{4,61} + \gamma_{\{12,13\}}\widetilde a_{4,62} + \gamma_{\{11,12\}}\widetilde a_{4,66}\\ 
&\quad + \gamma_{\{11,12\}}\widetilde a_{4,68} + \gamma_{\{10,11\}}\widetilde a_{4,71} + \gamma_{\{9,11\}}\widetilde a_{4,76} + \gamma_{\{9,11,12,13\}}\widetilde a_{4,82}\\ 
&\quad + \gamma_{\{7,11,12,13\}}\widetilde a_{4,83} + \gamma_{\{6,10\}}\widetilde a_{4,87} + \gamma_{\{1,2,4,12\}}\widetilde a_{4,89} + \gamma_{\{1,2,4,12\}}\widetilde a_{4,90}\\ 
&\quad + \gamma_{\{10,11,12,13\}}\widetilde a_{4,91} + \gamma_{\{2,3,5,13\}}\widetilde a_{4,93} + \gamma_{\{2,3,5,13\}}\widetilde a_{4,94}\\ 
&\quad + \gamma_{\{9,11,12,13\}}\widetilde a_{4,95} + \gamma_{\{12,13\}}\widetilde a_{4,104} + \gamma_{\{4,5,6,11\}}\widetilde a_{4,105} \equiv 0.
\end{align*}
From this equality we get $\gamma_u = \gamma_1$ for $1 \leqslant u \leqslant 13$. Hence, 
$$f \equiv \gamma_1\Big(\sum_{1 \leqslant u \leqslant 13}\widetilde p_{4,u}\Big) = \gamma_1 \xi_{t,4}.$$
The theorem is proved for $s = 4$.

By using Lemma \ref{bdt24} for $s\geqslant 5$, we have 
$f \equiv \sum_{1 \leqslant u \leqslant 13}\gamma_u\widetilde p_{s,u},$
where $\gamma_u \in \mathbb F_2$. By computing $\rho_4(f)+f$ in terms of the admissible monomials, we obtain
\begin{align*}
\rho_4(f)+f &\equiv \gamma_{\{1,4\}}\widetilde a_{s,1} + \gamma_{\{1,4\}}\widetilde a_{s,2} + \gamma_{\{1,10\}}\widetilde a_{s,3} + \gamma_{\{1,2\}}\widetilde a_{s,8} + \gamma_{\{1,2\}}\widetilde a_{s,9}\\ 
&\quad + \gamma_{\{2,5\}}\widetilde a_{s,13} + \gamma_{\{2,12\}}\widetilde a_{s,14} + \gamma_{\{2,12\}}\widetilde a_{s,15} + \gamma_{\{2,3\}}\widetilde a_{s,17} + \gamma_{\{2,3\}}\widetilde a_{s,18}\\ 
&\quad + \gamma_{\{3,13\}}\widetilde a_{s,25} + \gamma_{\{1,2,4,9\}}\widetilde a_{s,29} + \gamma_{\{1,2,4,9\}}\widetilde a_{s,30} + \gamma_{\{2,3,5,11\}}\widetilde a_{s,32}\\ 
&\quad + \gamma_{\{4,6\}}\widetilde a_{s,38} + \gamma_{\{4,6\}}\widetilde a_{s,39} + \gamma_{\{4,9\}}\widetilde a_{s,40} + \gamma_{\{4,9\}}\widetilde a_{s,41} + \gamma_{\{4,5\}}\widetilde a_{s,44}\\ 
&\quad + \gamma_{\{5,11\}}\widetilde a_{s,49} + \gamma_{\{5,11\}}\widetilde a_{s,50} + \gamma_{\{1,2,4,12\}}\widetilde a_{s,54} + \gamma_{\{1,2,4,12\}}\widetilde a_{s,55}\\ 
&\quad + \gamma_{\{2,3,5,13\}}\widetilde a_{s,57} + \gamma_{\{2,3,5,11\}}\widetilde a_{s,58} + \gamma_{\{9,10\}}\widetilde a_{s,59} + \gamma_{\{9,10\}}\widetilde a_{s,60}\\ 
&\quad + \gamma_{\{2,3,5,13\}}\widetilde a_{s,61} + \gamma_{\{12,13\}}\widetilde a_{s,62} + \gamma_{\{11,12\}}\widetilde a_{s,66} + \gamma_{\{11,12\}}\widetilde a_{s,68}\\ 
&\quad + \gamma_{\{4,5,6,8\}}\widetilde a_{s,71} + \gamma_{\{8,9\}}\widetilde a_{s,76} + \gamma_{\{8,9\}}\widetilde a_{s,82} + \gamma_{\{7,11,12,13\}}\widetilde a_{s,83}\\ 
&\quad + \gamma_{\{6,8\}}\widetilde a_{s,87} + \gamma_{\{1,2,4,12\}}\widetilde a_{s,89} + \gamma_{\{1,2,4,12\}}\widetilde a_{s,90} + \gamma_{\{4,5,6,8\}}\widetilde a_{s,91}\\ 
&\quad + \gamma_{\{4,5,6,11,12,13\}}\widetilde a_{s,92} + \gamma_{\{2,3,5,13\}}\widetilde a_{s,93} + \gamma_{\{2,3,5,13\}}\widetilde a_{s,94}\\ 
&\quad + \gamma_{\{8,9,12,13\}}\widetilde a_{s,95} + \gamma_{\{4,5,6,11\}}\widetilde a_{s,97} + \gamma_{\{12,13\}}\widetilde a_{s,98}\\ 
&\quad + \gamma_{\{4,5,6,11,12,13\}}\widetilde a_{s,99} + \gamma_{\{4,5,6,11\}}\widetilde a_{s,100} + \gamma_{\{12,13\}}\widetilde a_{s,103} \equiv 0.
\end{align*}
From this equality we get $\gamma_u = \gamma_1$ for $1 \leqslant u \leqslant 13$. Hence, 
$$f \equiv \gamma_1\Big(\sum_{1 \leqslant u \leqslant 13}\widetilde p_{4,u}\Big) = \gamma_1\Big(\sum_{1 \leqslant u \leqslant 105}\widetilde a_{s,u}\Big) =\gamma_1 \xi_{t,s}.$$
Theorem \ref{dlt21} is completely proved. 
\end{proof}

Since $QP_4((3)|^{s}|(2))^{GL_4} =0$, we have $(QP_4)_{d_{s,2}}^{GL_4} = QP_4((3)|^{s-1}|(1)|^3)^{GL_4}$. So, we get the following.

\begin{corls} Theorem \ref{thm1} holds for $t=2$.
\end{corls}	

\begin{rems} In \cite{pp25}, the author proved that $(QP_4)_{d_{2,2}}^{GL_4}= \langle [\xi_{2,2}]\rangle$ but there are some mistakes in the proof. However, the $GL_4$-invariant $[\xi_{2,2}]$ was first found in \cite[Remark 6.3.9]{su50}. The author also stated that $(QP_4)_{d_{s,2}}^{GL_4}= 0$ for $s > 2$ but this result is false and the proof in \cite{pp25} is also false. 
\end{rems}

\subsection{Proof of Theorem \ref{thm1} for $t \geqslant 3$}\

\medskip
For $t = 3$ and $s \geqslant 3$, we have $QP_4((3)|^3|(2)|^2)^{GL_4} = 0$, hence Theorem \ref{thm1} follows from Theorem \ref{dlt21}. So, we need only to prove the theorem for $s = 2$.
\begin{props}\label{mdt32} We have $(QP_4)_{d_{2,3}}^{GL_4} =  \langle [\xi_{3,2}]\rangle,$ where
$$\xi_{3,2} = \bar\xi_{3,2} + \sum_{j\in \{29,\, 36,\, 37,\, 39,\, 41,\, 42,\, 45\}}\widetilde a_{2,j}$$
\end{props}

\begin{proof}
By Proposition \ref{mdt31}, for $s = 2$, we have $QP_4((3)|^2|(2)|^2)^{GL_4} = \langle [\bar\xi_{3,2}]\rangle$. Hence, if $[f] \in (QP_4)_{33}^{GL_4}$ with $f \in P_4$, then there is $\lambda_0 \in \mathbb F_2$ such that $[f] + \lambda_0[\bar\xi_{3,2}] \in QP_4((3)|(1)|^4)$. By a direct computation we obtain
\begin{equation}\label{cthrt3}
\begin{cases}
\rho_1(\bar \xi_{3,2}) + \bar \xi_{3,2} \equiv 0,\\
\rho_2(\bar \xi_{3,2}) + \bar \xi_{3,2} \equiv \sum_{j \in \{29,\, 32,\, 33,\, 35,\, 39,\, 41,\, 44,\, 45\}}\widetilde a_{2,j},\\
\rho_3(\bar \xi_{3,2}) + \bar \xi_{3,2} \equiv \sum_{j \in \{32,\, 36,\, 41,\, 42,\, 43,\, 44\}}\widetilde a_{2,j},\\
\rho_4(\bar \xi_{3,2}) + \bar \xi_{3,2} \equiv \widetilde a_{2,41}.
\end{cases}
\end{equation}
By a computation using \eqref{cthrt3} and Lemma \ref{bdt24} we see that $\rho_i(f) + f \equiv 0$ for $i=1,\, 2,\, 3$ if and only if
$$f \equiv \gamma_0(\bar \xi_{3,2} + \widetilde p_{2,0}) + \sum _{1 \leqslant u \leqslant 4} \gamma_u\widetilde p_{2,u},$$
where $\gamma_u \in \mathbb F_2$ and 
$\widetilde p_{2,0} = \sum_{j \in \{29,\, 36,\, 37,\, 39,\, 41,\, 42,\, 45\}}\widetilde a_{2,j}.$

Now by computing $\rho_4(f) + f$ in terms of the admissible monomials we get
\begin{align*}
\rho_4(f)&+f \equiv \gamma_{1}\widetilde aa_{2,1} + \gamma_{1}\widetilde aa_{2,2} + \gamma_{\{1,4\}}\widetilde aa_{2,3} + \gamma_{\{1,2\}}\widetilde aa_{2,8} + \gamma_{\{1,2\}}\widetilde aa_{2,9}\\ 
& + \gamma_{2}\widetilde aa_{2,13} + \gamma_{\{2,4\}}\widetilde aa_{2,14} + \gamma_{\{2,4\}}\widetilde aa_{2,15} + \gamma_{\{2,3\}}\widetilde aa_{2,17} + \gamma_{\{2,3\}}\widetilde aa_{2,18}\\ 
& + \gamma_{\{3,4\}}\widetilde aa_{2,25} + \gamma_{\{1,2\}}\widetilde aa_{2,29} + \gamma_{\{1,2\}}\widetilde aa_{2,30} + \gamma_{\{2,3\}}\widetilde aa_{2,32}\\ 
& + \gamma_{\{2,3\}}\widetilde aa_{2,35} + \gamma_{4}\widetilde aa_{2,36} + \gamma_{4}\widetilde aa_{2,37} + \gamma_{4}\widetilde aa_{2,41} + \gamma_{4}\widetilde aa_{2,43} \equiv 0.
\end{align*} 
From this equality we obtain $\gamma_u = 0$ for $1 \leqslant u \leqslant 4$ and
\begin{align*}
f &\equiv \gamma_0(\bar \xi_{3,2} + \widetilde p_{2,0}) = \gamma_0\xi_{3,2}.
\end{align*}
The proposition is proved.
\end{proof}

For $t \geqslant 4$, Theorem \ref{thm1} equivalent to the following.

\begin{thms}\label{d1t4} For $t \geqslant 4$, we have
$(QP_4)_{d_{s,t}}^{GL_4} = \langle [\xi_{t,s}],\,[\theta_{t,s}]\rangle$, where $\xi_{t,s}$ is determined as in Theorem \ref{dlt21} and
$$\theta_{t,s} = \begin{cases}
\bar \xi_{t,2} + \sum_{j\in \{25,\, 26,\, 27,\, 28,\, 29,\, 32,\, 35,\, 37,\, 39,\, 43,\, 44,\, 45\}}\widetilde a_{2,j}, &\mbox{if } s=2,\\
\bar \xi_{t,3} + \sum_{j\in \{25,\, 26,\, 27,\, 28,\, 35,\, 56,\, 62,\, 65,\, 66,\, 67,\, 68,\, 69,\, 78,\, 81\}\atop \hskip4.1cm\cup\{82,\, 83,\, 84\}}\widetilde a_{3,j}, &\mbox{if } s=3,\\
\bar \xi_{t,4} + \sum_{j\in \{25,\, 26,\, 27,\, 28,\, 29,\, 35,\, 54,\, 56,\, 66,\, 67,\, 68,\, 69,\, 70,\, 72,\, 73\}\atop \hskip0.5cm\cup\{79,\, 80,\, 81,\, 82,\, 83,\, 84,\, 90,\, 92,\, 97,\, 99,\, 103,\, 104,\, 105\}}\widetilde a_{4,j}, &\mbox{if } s=4,\\
\bar \xi_{t,s} + \sum_{j\in \{25,\, 26,\, 27,\, 28,\, 35,\, 66,\, 67,\, 68,\, 69,\, 70,\, 81,\, 82,\, 83,\, 84\}}\widetilde a_{s,j}, &\mbox{if } s\geqslant 5.\\
\end{cases} $$
Here, $\bar \xi_{t,s}$ is determined as in Proposition \ref{mdt41} and $\widetilde a_{s,j} = \widetilde a_{t,s,j}$.	
\end{thms}
\begin{proof}
By Proposition \ref{mdt41}, for $t \geqslant 4$, we have $QP_4((3)|^s|(2)|^{t-1})^{GL_4} = \langle [\bar\xi_{t,s}]\rangle$. Hence, if $[f] \in (QP_4)_{d_{s,t}}^{GL_4}$ with $f \in P_4$, then there is $\lambda_0 \in \mathbb F_2$ such that $[f] + \lambda_0[\bar\xi_{t,s}] \in QP_4((3)|^{s-1}|(1)|^{t+1})$. 

For $s = 2$, we have
\begin{equation}\label{cthrt4}
\begin{cases}
\rho_1(\bar \xi_{t,2}) + \bar \xi_{t,2} \equiv \widetilde a_{2,32} + \widetilde a_{2,36},\\
\rho_2(\bar \xi_{t,2}) + \bar \xi_{t,2} \equiv \sum_{j \in \{32,\, 33,\, 35,\, 36,\, 37,\, 41\}}\widetilde a_{2,j},\\
\rho_3(\bar \xi_{t,2}) + \bar \xi_{t,2} \equiv \sum_{j \in \{32,\, 36,\, 41,\, 42,\, 43,\, 44\}}\widetilde a_{2,j},\\
\rho_4(\bar \xi_{t,2}) + \bar \xi_{t,2} \equiv \sum_{j \in \{17,\, 18,\, 32,\, 36,\, 43\}}\widetilde a_{2,j}.
\end{cases}
\end{equation}
By a computation using \eqref{cthrt4} and Lemma \ref{bdt24} we see that $\rho_i(f) + f \equiv 0$ for $i=1,\, 2,\, 3$ if and only if
$$f \equiv \gamma_0(\bar \xi_{t,2} + \widetilde p_{2,0}) + \sum _{1 \leqslant u \leqslant 4} \gamma_u\widetilde p_{2,u},$$
where $\gamma_u \in \mathbb F_2$ and 
$\widetilde p_{2,0} = \sum_{j \in \{29,\, 32,\, 35,\, 37,\, 39,\, 43,\, 44,\, 45\}}\widetilde a_{2,j}.$

By computing $\rho_4(f) + f$ in terms of the admissible monomials we get
\begin{align*}
\rho_4(f)&+f \equiv \gamma_{1}\widetilde a_{2,1} + \gamma_{1}\widetilde a_{2,2} + \gamma_{\{1,4\}}\widetilde a_{2,3} + \gamma_{\{1,2\}}\widetilde a_{2,8} + \gamma_{\{1,2\}}\widetilde a_{2,9} + \gamma_{2}\widetilde a_{2,13}\\ 
& + \gamma_{\{2,4\}}\widetilde a_{2,14} + \gamma_{\{2,4\}}\widetilde a_{2,15} + \gamma_{\{0,2,3\}}\widetilde a_{2,17} + \gamma_{\{0,2,3\}}\widetilde a_{2,18} + \gamma_{\{0,3,4\}}\widetilde a_{2,25}\\ 
& + \gamma_{\{1,2\}}\widetilde a_{2,29} + \gamma_{\{1,2\}}\widetilde a_{2,30} + \gamma_{\{0,2,3\}}\widetilde a_{2,32} + \gamma_{\{0,2,3\}}\widetilde a_{2,36} + \gamma_{4}\widetilde a_{2,37}\\ 
& + \gamma_{4}\widetilde a_{2,38} + \gamma_{4}\widetilde a_{2,41} + \gamma_{4}\widetilde a_{2,43}  \equiv 0.
\end{align*} 
From this equality we obtain $\gamma_u = 0$ for $u = 1,\, 2,\, 4$ and $\gamma_5 = \gamma_3$. hence, we get
\begin{align*}
f &\equiv \gamma_0(\bar \xi_{3,2} + \widetilde p_{2,0} + \widetilde p_{2,3}) = \gamma_0\theta_{t,2}.
\end{align*}
The theorem is proved for $s = 2$.

For $s = 3$, we have
\begin{equation}\label{cthrt5}
\begin{cases}
\rho_1(\bar \xi_{t,3}) + \bar \xi_{t,3} \equiv \sum_{j \in \{32,\, 57,\, 58,\, 61\}}\widetilde a_{3,j},\\
\rho_2(\bar \xi_{t,3}) + \bar \xi_{t,3} \equiv \sum_{j \in \{57,\, 61,\, 62,\, 69,\, 70,\, 74\}}\widetilde a_{3,j},\\
\rho_3(\bar \xi_{t,3}) + \bar \xi_{t,3} \equiv \sum_{j \in \{56,\, 62,\, 65,\, 71,\, 76,\, 77\}}\widetilde a_{3,j},\\
\rho_4(\bar \xi_{t,3}) + \bar \xi_{t,3} \equiv \sum_{j \in \{17,\, 18,\, 62,\, 83\}}\widetilde a_{3,j}.
\end{cases}
\end{equation}
By a computation using \eqref{cthrt5} and Lemma \ref{bdt24} we see that $\rho_i(f) + f \equiv 0$ for $i=1,\, 2,\, 3$ if and only if
$$f \equiv \gamma_0(\bar \xi_{t,3} + \widetilde p_{3,0}) + \sum _{1 \leqslant u \leqslant 9} \gamma_u\widetilde p_{3,u},$$
where $\gamma_u \in \mathbb F_2$ and 
$\widetilde p_{3,0} = \sum_{j \in \{29,\, 30,\, 31,\, 32,\, 33,\, 34,\, 35,\, 54,\, 55,\, 56,\, 58\}\atop\hskip0.5cm\cup\{59,\, 60,\, 62,\, 63,\, 64,\, 65,\, 66,\, 67,\, 68,\, 69\}}\widetilde a_{3,j}.$

Computing $\rho_4(f) + f$ in terms of the admissible monomials gives
\begin{align*}
\rho_4(f)&+f \equiv \gamma_{\{1,4\}}\widetilde a_{3,1} + \gamma_{\{1,4\}}\widetilde a_{3,2} + \gamma_{\{0,1,6\}}\widetilde a_{3,3} + \gamma_{\{1,2\}}\widetilde a_{3,8} + \gamma_{\{1,2\}}\widetilde a_{3,9}\\ 
& + \gamma_{\{2,5\}}\widetilde a_{3,13} + \gamma_{\{0,2,8\}}\widetilde a_{3,14} + \gamma_{\{0,2,8\}}\widetilde a_{3,15} + \gamma_{\{0,2,3\}}\widetilde a_{3,17} + \gamma_{\{0,2,3\}}\widetilde a_{3,18}\\ 
& + \gamma_{\{0,3,9\}}\widetilde a_{3,25} + \gamma_{\{0,8\}}\widetilde a_{3,29} + \gamma_{\{0,7,8\}}\widetilde a_{3,30} + \gamma_{\{7,9\}}\widetilde a_{3,32} + \gamma_{4}\widetilde a_{3,38} + \gamma_{4}\widetilde a_{3,39}\\ 
& + \gamma_{\{4,7\}}\widetilde a_{3,40} + \gamma_{\{4,7\}}\widetilde a_{3,41} + \gamma_{\{4,5\}}\widetilde a_{3,44} + \gamma_{\{5,7\}}\widetilde a_{3,49} + \gamma_{\{5,7\}}\widetilde a_{3,50}\\ 
& + \gamma_{\{4,5\}}\widetilde a_{3,56} + \gamma_{\{7,9\}}\widetilde a_{3,58} + \gamma_{\{0,6\}}\widetilde a_{3,59} + \gamma_{\{0,6,7\}}\widetilde a_{3,60} + \gamma_{7}\widetilde a_{3,62}\\ 
& + \gamma_{\{0,7,8\}}\widetilde a_{3,66} + \gamma_{\{0,7,8\}}\widetilde a_{3,68} + \gamma_{\{4,5\}}\widetilde a_{3,71} + \gamma_{7}\widetilde a_{3,76} + \gamma_{7}\widetilde a_{3,78}\\ 
& + \gamma_{\{0,7,8,9\}}\widetilde a_{3,83} + \gamma_{\{0,1,2,4,7,8\}}\widetilde a_{3,85} + \gamma_{\{0,1,2,4,8\}}\widetilde a_{3,86} + \gamma_{\{0,2,3,5,9\}}\widetilde a_{3,87}\\ 
& + \gamma_{\{0,2,3,5,9\}}\widetilde a_{3,88} + \gamma_{7}\widetilde a_{3,89} \equiv 0.
\end{align*} 
From this equality we obtain $\gamma_u = \gamma_0$ for $u = 3,\, 6,\, 8$, and $\gamma_u = 0$ for $u \ne 3,\, 6,\, 8$. So, we get
\begin{align*}
f &\equiv \gamma_0(\bar \xi_{t,3} + \widetilde p_{3,0} + \widetilde p_{3,3} + \widetilde p_{3,6} + \widetilde p_{3,8}) = \gamma_0\theta_{t,3}.
\end{align*}
The theorem is proved for $s = 3$.

For $s = 4$, we have
\begin{equation}\label{cthrt6}
\begin{cases}
\rho_1(\bar \xi_{t,4}) + \bar \xi_{t,4} \equiv 0,\\
\rho_2(\bar \xi_{t,4}) + \bar \xi_{t,4} \equiv \sum_{j \in \{29,\, 54,\, 56,\, 74,\, 89,\, 91,\, 93,\, 98,\, 100,\, 103\}}\widetilde a_{4,j},\\
\rho_3(\bar \xi_{t,4}) + \bar \xi_{t,4} \equiv \sum_{j \in \{56,\, 71,\, 76,\, 77,\, 91,\, 92,\, 95,\, 96,\, 97,\, 105\}}\widetilde a_{4,j},\\
\rho_4(\bar \xi_{t,4}) + \bar \xi_{t,4} \equiv \sum_{j \in \{17,\, 18,\, 29,\, 30,\, 40,\, 41,\, 59,\, 60,\, 62,\, 71,\, 83,\, 87,\, 91\}}\widetilde a_{4,j}.
\end{cases}
\end{equation}
By a computation using \eqref{cthrt6} and Lemma \ref{bdt24} we see that $\rho_i(f) + f \equiv 0$ for $i=1,\, 2,\, 3$ if and only if
$$f \equiv \gamma_0(\bar \xi_{t,4} + \widetilde p_{4,0}) + \sum _{1 \leqslant u \leqslant 13} \gamma_u\widetilde p_{4,u},$$
where $\gamma_u \in \mathbb F_2$ and 
$\widetilde p_{4,0} = \sum_{\scriptscriptstyle j \in \{29,\, 30,\, 31,\, 54,\, 55,\, 56,\, 78,\, 79,\, 80,\, 89,\, 90,\, 91,\, 92,\, 93,\, 94,\, 95,\, 96,\, 98\}}\widetilde a_{4,j}.$

By computing $\rho_4(f) + f$ in terms of the admissible monomials we have
\begin{align*}
\rho_4(f)&+f \equiv \gamma_{\{1,4\}}\widetilde a_{4,1} + \gamma_{\{1,4\}}\widetilde a_{4,2} + \gamma_{\{0,1,8,10,13\}}\widetilde a_{4,3} + \gamma_{\{1,2\}}\widetilde a_{4,8} + \gamma_{\{1,2\}}\widetilde a_{4,9}\\
& + \gamma_{\{2,5\}}\widetilde a_{4,13} + \gamma_{\{2,12\}}\widetilde a_{4,14} + \gamma_{\{2,12\}}\widetilde a_{4,15} + \gamma_{\{0,2,3\}}\widetilde a_{4,17} + \gamma_{\{0,2,3\}}\widetilde a_{4,18}\\
& + \gamma_{\{3,13\}}\widetilde a_{4,25} + \gamma_{\{0,1,2,4,9,10,11\}}\widetilde a_{4,29} + \gamma_{\{0,1,2,4,9,10,11\}}\widetilde a_{4,30}\\
& + \gamma_{\{0,2,3,5,11\}}\widetilde a_{4,32} + \gamma_{\{4,6\}}\widetilde a_{4,38} + \gamma_{\{4,6\}}\widetilde a_{4,39} + \gamma_{\{0,4,9,10,11\}}\widetilde a_{4,40}\\
& + \gamma_{\{0,4,9,10,11\}}\widetilde a_{4,41} + \gamma_{\{4,5\}}\widetilde a_{4,44} + \gamma_{\{5,11\}}\widetilde a_{4,49} + \gamma_{\{5,11\}}\widetilde a_{4,50}\\
& + \gamma_{\{1,2,4,12\}}\widetilde a_{4,54} + \gamma_{\{1,2,4,12\}}\widetilde a_{4,55} + \gamma_{\{0,4,5,6,11,12,13\}}\widetilde a_{4,56} + \gamma_{\{2,3,5,13\}}\widetilde a_{4,57}\\
& + \gamma_{\{0,2,3,5,11\}}\widetilde a_{4,58} + \gamma_{\{8,9,11,13\}}\widetilde a_{4,59} + \gamma_{\{8,9,11,13\}}\widetilde a_{4,60} + \gamma_{\{2,3,5,13\}}\widetilde a_{4,61}\\
& + \gamma_{\{0,12,13\}}\widetilde a_{4,62} + \gamma_{\{11,12\}}\widetilde a_{4,66} + \gamma_{\{11,12\}}\widetilde a_{4,68} + \gamma_{\{0,10,11\}}\widetilde a_{4,71} + \gamma_{\{9,11\}}\widetilde a_{4,76}\\
& + \gamma_{\{0,9,11,12,13\}}\widetilde a_{4,78} + \gamma_{\{0,7,11,12,13\}}\widetilde a_{4,83} + \gamma_{\{0,6,10\}}\widetilde a_{4,87} + \gamma_{\{1,2,4,12\}}\widetilde a_{4,89}\\
& + \gamma_{\{1,2,4,12\}}\widetilde a_{4,90} + \gamma_{\{10,11,12,13\}}\widetilde a_{4,91} + \gamma_{\{2,3,5,13\}}\widetilde a_{4,93} + \gamma_{\{2,3,5,13\}}\widetilde a_{4,94}\\
& + \gamma_{\{0,9,11,12,13\}}\widetilde a_{4,95} + \gamma_{\{0,12,13\}}\widetilde a_{4,104} + \gamma_{\{4,5,6,11\}}\widetilde a_{4,105} \equiv 0.
\end{align*} 
From this equality we obtain $\gamma_u = \gamma_1$ for $1 \leqslant u \leqslant 13,\, u \ne 3,\, 8,\, 10, 13$, and $\gamma_u = \gamma_{0} + \gamma_1$ for $u = 3,\, 8,\, 10, 13$. So, we get
\begin{align*}
f &\equiv \gamma_0\Big(\bar \xi_{t,4} + \sum_{u\in \{0,\, 3,\, 8,\, 10, 13\}}\widetilde p_{4,u}\Big) + \gamma_1\Big(\sum_{1 \leqslant u \leqslant 13}\widetilde p_{4,u}\Big) = \gamma_0\theta_{t,4} + \gamma_1\xi_{t,4}.
\end{align*}
The theorem is proved for $s = 4$.

For $s \geqslant 5$, we have
\begin{equation}\label{cthrt7}
\begin{cases}
\rho_1(\bar \xi_{t,s}) + \bar \xi_{t,s} \equiv 0,\\
\rho_2(\bar \xi_{t,s}) + \bar \xi_{t,s} \equiv \sum_{j \in \{74,\, 93,\, 100,\, 104\}}\widetilde a_{s,j},\\
\rho_3(\bar \xi_{t,s}) + \bar \xi_{t,s} \equiv \sum_{j \in \{76,\, 77,\, 95,\, 96,\, 103,\, 105\}}\widetilde a_{s,j},\\
\rho_4(\bar \xi_{t,s}) + \bar \xi_{t,s} \equiv \sum_{j \in \{17,\, 18,\, 62,\, 83\}}\widetilde a_{s,j}.
\end{cases}
\end{equation}
By a computation using \eqref{cthrt7} and Lemma \ref{bdt24} we see that $\rho_i(f) + f \equiv 0$ for $i=1,\, 2,\, 3$ if and only if
$$f \equiv \gamma_0(\bar \xi_{t,s} + \widetilde p_{s,0}) + \sum _{1 \leqslant u \leqslant 13} \gamma_u\widetilde p_{s,u},$$
where $\gamma_u \in \mathbb F_2$ and 
$\widetilde p_{s,0} = \sum_{j \in \{93,\, 94,\, 95,\, 96,\, 103,\, 104\}}\widetilde a_{s,j}.$

Computing $\rho_4(f) + f$ in terms of the admissible monomials gives
\begin{align*}
\rho_4(f)&+f \equiv \gamma_{\{1,4\}}\widetilde a_{s,1} + \gamma_{\{1,4\}}\widetilde a_{s,2} + \gamma_{\{1,10\}}\widetilde a_{s,3} + \gamma_{\{1,2\}}\widetilde a_{s,8} + \gamma_{\{1,2\}}\widetilde a_{s,9}\\ 
& + \gamma_{\{2,5\}}\widetilde a_{s,13} + \gamma_{\{2,12\}}\widetilde a_{s,14} + \gamma_{\{2,12\}}\widetilde a_{s,15} + \gamma_{\{0,2,3\}}\widetilde a_{s,17} + \gamma_{\{0,2,3\}}\widetilde a_{s,18}\\ 
& + \gamma_{\{3,13\}}\widetilde a_{s,25} + \gamma_{\{1,2,4,9\}}\widetilde a_{s,29} + \gamma_{\{1,2,4,9\}}\widetilde a_{s,30} + \gamma_{\{0,2,3,5,11\}}\widetilde a_{s,32}\\ 
& + \gamma_{\{4,6\}}\widetilde a_{s,38} + \gamma_{\{4,6\}}\widetilde a_{s,39} + \gamma_{\{4,9\}}\widetilde a_{s,40} + \gamma_{\{4,9\}}\widetilde a_{s,41} + \gamma_{\{4,5\}}\widetilde a_{s,44}\\ 
& + \gamma_{\{5,11\}}\widetilde a_{s,49} + \gamma_{\{5,11\}}\widetilde a_{s,50} + \gamma_{\{1,2,4,12\}}\widetilde a_{s,54} + \gamma_{\{1,2,4,12\}}\widetilde a_{s,55}\\ 
& + \gamma_{\{2,3,5,13\}}\widetilde a_{s,57} + \gamma_{\{0,2,3,5,11\}}\widetilde a_{s,58} + \gamma_{\{9,10\}}\widetilde a_{s,59} + \gamma_{\{9,10\}}\widetilde a_{s,60}\\ 
& + \gamma_{\{2,3,5,13\}}\widetilde a_{s,61} + \gamma_{\{0,12,13\}}\widetilde a_{s,62} + \gamma_{\{11,12\}}\widetilde a_{s,66} + \gamma_{\{11,12\}}\widetilde a_{s,68}\\ 
& + \gamma_{\{4,5,6,8\}}\widetilde a_{s,71} + \gamma_{\{8,9\}}\widetilde a_{s,76} + \gamma_{\{8,9\}}\widetilde a_{s,78} + \gamma_{\{0,7,11,12,13\}}\widetilde a_{s,83}\\ 
& + \gamma_{\{6,8\}}\widetilde a_{s,87} + \gamma_{\{1,2,4,12\}}\widetilde a_{s,89} + \gamma_{\{1,2,4,12\}}\widetilde a_{s,90} + \gamma_{\{4,5,6,8\}}\widetilde a_{s,91}\\ 
& + \gamma_{\{0,4,5,6,11,12,13\}}\widetilde a_{s,92} + \gamma_{\{2,3,5,13\}}\widetilde a_{s,93} + \gamma_{\{2,3,5,13\}}\widetilde a_{s,94}\\ 
& + \gamma_{\{0,8,9,12,13\}}\widetilde a_{s,95} + \gamma_{\{4,5,6,11\}}\widetilde a_{s,97} + \gamma_{\{0,12,13\}}\widetilde a_{s,98}\\ 
& + \gamma_{\{0,4,5,6,11,12,13\}}\widetilde a_{s,103} + \gamma_{\{0,12,13\}}\widetilde a_{s,104} + \gamma_{\{4,5,6,11\}}\widetilde a_{s,105} \equiv 0.
\end{align*} 
From this equality we obtain $\gamma_u = \gamma_1$ for $2 \leqslant u \leqslant 12,\, u \ne 3$ and $\gamma_u = \gamma_0 + \gamma_1$ for $u = 3,\, 13$. Hence, we get
\begin{align*}
f &\equiv \gamma_0(\bar \xi_{t,s} + \widetilde p_{s,0} + \widetilde p_{s,3} + \widetilde p_{s,13}) + \gamma_1\Big(\sum_{1 \leqslant u \leqslant 13}\widetilde p_{s,u}\Big)= \gamma_0\theta_{t,s} + \gamma_1\xi_{t,s}.
\end{align*}
Theorem \ref{d1t4} is completely proved.
\end{proof}

For $t \geqslant 4$, we have $\xi_{t,2} = 0$ and $\xi_{t,3} = 0$. Hence, $\dim (QP_4)_{d_{s,t}} = 1$ for $s = 2,\, 3$ and $\dim (QP_4)_{d_{s,t}} = 2$ for $s \geqslant 4$. This completes the proof of Theorem \ref{thm1}.


\section{Proof of Theorem \ref{thm2}}\label{s4}
\setcounter{equation}{0}

Let $m = n_{s,t} = 2^{s+t} + 2^s - 2$ with $s,\, t$ positive integers. 
Consider Kameko's homomorphism
$$(\widetilde{Sq}^0_*)_{(4,n_{s,t})} : (QP_{4})_{n_{s,t}}\to (QP_4)_{d_{s-1,t}} .$$
It is well-known that $(\widetilde{Sq}^0_*)_{(4,n_{s,t})}$ is a homomorphism of $GL_4$-modules, hence if $[f] \in (QP_4)_{n_{s,t}}^{GL_4}$, then  $(\widetilde{Sq}^0_*)_{(4,n_{s,t})}([f])\in (QP_4)_{d_{s-1,t}}^{GL_4}$.  Hence, there is $[\xi] \in (QP_4)_{d_{s-1,t}}^{GL_4}$ with $\xi \in (P_4)_{d_{s-1,t}}$ such that $f \equiv \psi_{s,t}(\xi) + h$, where $\psi_{s,t}: (P_4)_{d_{s-1,t}} \to (P_4)_{n_{s,t}}$ is determined by $\psi_{s,t}(y) = x_1x_2x_3x_4y^2$ for all $y \in (P_4)_{d_{s-1,t}}$ and $[h] \in \mbox{Ker}\big((\widetilde{Sq}^0_*)_{(4,n_{s,t})}\big)$. If $[h_0] \in \mbox{Ker}\big((\widetilde{Sq}^0_*)_{(4,n_{s,t})}\big)$ such that $[\psi_{s,t}(\xi) + h_0] \in (QP_4)_{n_{s,t}}^{GL_4}$, then $[g]$ is an element of $\mbox{Ker}(\widetilde{Sq}^0_*)_{(4,n_{s,t})}^{GL_4}$ with $g = h-h_0$. Thus, we need to explicitly determine the space $\mbox{Ker}\big((\widetilde{Sq}^0_*)_{(4,n_{s,t})}\big)^{GL_4}$.

\medskip
\subsection{The case $s= 1$}\

\medskip
For $s= 1$, we have $n_{1,t} = 2^{t+1}$ and $d_{0,t} = 2^{t} - 2$. The space $(QP_4)_{2^t-2}$ had been explicitly determined in our work \cite{su2}.
\begin{thms}[See \cite{su2}]\label{dlmd1}
For $t \geqslant 2$, we have
$$(QP_4)_{2^t-2}^{GL_4} = \big(\mbox{\rm Ker}(\widetilde{Sq}^0_* )_{(4,2^{t}-2)}\big)^{GL_4} = \begin{cases} 0, &\mbox{if } t \leqslant 3,\\
\langle [\xi_{t,0}]\rangle, &\text{if } t \geqslant 4,\end{cases}
	$$
where $\xi_{4,0} = x_1x_2x_3^6x_4^6 + x_1^3x_2^3x_3^4x_4^4$ and for $t \geqslant 5$,
\begin{align*}
\xi_{t,0} &=  x_3^{2^t-1}x_4^{2^t-1} +  x_2^{2^t-1}x_4^{2^t-1} +  x_2^{2^t-1}x_3^{2^t-1} +  x_1^{2^t-1}x_4^{2^t-1} +  x_1^{2^t-1}x_3^{2^t-1}\\ 
&\quad +  x_1^{2^t-1}x_2^{2^t-1} +  x_2x_3^{2^t-2}x_4^{2^t-1} +  x_2x_3^{2^t-1}x_4^{2^t-2} +  x_2^{2^t-1}x_3x_4^{2^t-2}\\ 
&\quad +  x_1x_3^{2^t-2}x_4^{2^t-1} +  x_1x_3^{2^t-1}x_4^{2^t-2} +  x_1x_2^{2^t-2}x_4^{2^t-1} +  x_1x_2^{2^t-2}x_3^{2^t-1}\\ 
&\quad +  x_1x_2^{2^t-1}x_4^{2^t-2} +  x_1x_2^{2^t-1}x_3^{2^t-2} +  x_1^{2^t-1}x_3x_4^{2^t-2} +  x_1^{2^t-1}x_2x_4^{2^t-2}\\ 
&\quad +  x_1^{2^t-1}x_2x_3^{2^t-2} +  x_2^{3}x_3^{2^t-3}x_4^{2^t-2} +  x_1^{3}x_3^{2^t-3}x_4^{2^t-2} +  x_1^{3}x_2^{2^t-3}x_4^{2^t-2}\\ 
&\quad +  x_1^{3}x_2^{2^t-3}x_3^{2^t-2} +  x_1x_2x_3^{2^t-2}x_4^{2^t-2} +  x_1x_2^{2^t-2}x_3x_4^{2^t-2} +  x_1x_2^{2}x_3^{2^t-4}x_4^{2^t-1}\\ 
&\quad +  x_1x_2^{2}x_3^{2^t-1}x_4^{2^t-4} +  x_1x_2^{2^t-1}x_3^{2}x_4^{2^t-4} +  x_1^{2^t-1}x_2x_3^{2}x_4^{2^t-4}\\ 
&\quad +  x_1x_2^{2}x_3^{2^t-3}x_4^{2^t-2} +  x_1x_2^{3}x_3^{2^t-4}x_4^{2^t-2} +  x_1x_2^{3}x_3^{2^t-2}x_4^{2^t-4}\\ 
&\quad +  x_1^{3}x_2x_3^{2^t-4}x_4^{2^t-2} +  x_1^{3}x_2x_3^{2^t-2}x_4^{2^t-4} +  x_1^{3}x_2^{2^t-3}x_3^{2}x_4^{2^t-4}\\ 
&\quad +  x_1^{3}x_2^{5}x_3^{2^t-6}x_4^{2^t-4}. 
\end{align*}
\end{thms}

\begin{rems}
The space $(QP_4)_{2^{t+1}}^{GL_4}$ had been determined by Singer \cite{si1} for $t = 1,\, 2$. He proved that $(QP_4)_{4}^{GL_4}=0$ and $(QP_4)_{8}^{GL_4} = 0$. In \cite{p24}, the author also proved $(QP_4)_{8}^{GL_4} = 0$ but there is a mistake in his proof. 
\end{rems}

By our works \cite{su50,su5}, for $t \geqslant 2$, a basis of $\mbox{\rm Ker}(\widetilde{Sq}^0_* )_{(4,2^{t+1})}$ is the set of all classes represented by the admissible monomials $b_{t,j} = b_{t,1,j}$ which are determined as follows:

\medskip
For $t \geqslant 2$,

\medskip  
\centerline{\begin{tabular}{lll}
$b_{t,1} = x_3x_4^{2^{t+1}-1}$&$b_{t,2} = x_3^{2^{t+1}-1}x_4$&$b_{t,3} = x_2x_4^{2^{t+1}-1}$\cr  $b_{t,4} = x_2x_3^{2^{t+1}-1}$&$b_{t,5} = x_2^{2^{t+1}-1}x_4$&$b_{t,6} = x_2^{2^{t+1}-1}x_3$\cr  $b_{t,7} = x_1x_4^{2^{t+1}-1}$&$b_{t,8} = x_1x_3^{2^{t+1}-1}$&$b_{t,9} = x_1x_2^{2^{t+1}-1}$\cr  $b_{t,10} = x_1^{2^{t+1}-1}x_4$&$b_{t,11} = x_1^{2^{t+1}-1}x_3$&$b_{t,12} = x_1^{2^{t+1}-1}x_2$\cr  $b_{t,13} = x_3^{3}x_4^{2^{t+1}-3}$&$b_{t,14} = x_2^{3}x_4^{2^{t+1}-3}$&$b_{t,15} = x_2^{3}x_3^{2^{t+1}-3}$\cr  $b_{t,16} = x_1^{3}x_4^{2^{t+1}-3}$&$b_{t,17} = x_1^{3}x_3^{2^{t+1}-3}$&$b_{t,18} = x_1^{3}x_2^{2^{t+1}-3}$\cr  $b_{t,19} = x_2x_3x_4^{2^{t+1}-2}$&$b_{t,20} = x_2x_3^{2^{t+1}-2}x_4$&$b_{t,21} = x_1x_3x_4^{2^{t+1}-2}$\cr  $b_{t,22} = x_1x_3^{2^{t+1}-2}x_4$&$b_{t,23} = x_1x_2x_4^{2^{t+1}-2}$&$b_{t,24} = x_1x_2x_3^{2^{t+1}-2}$\cr  $b_{t,25} = x_1x_2^{2^{t+1}-2}x_4$&$b_{t,26} = x_1x_2^{2^{t+1}-2}x_3$&$b_{t,27} = x_2x_3^{2}x_4^{2^{t+1}-3}$\cr  $b_{t,28} = x_1x_3^{2}x_4^{2^{t+1}-3}$&$b_{t,29} = x_1x_2^{2}x_4^{2^{t+1}-3}$&$b_{t,30} = x_1x_2^{2}x_3^{2^{t+1}-3}$\cr  $b_{t,31} = x_2x_3^{3}x_4^{2^{t+1}-4}$&$b_{t,32} = x_2^{3}x_3x_4^{2^{t+1}-4}$&$b_{t,33} = x_1x_3^{3}x_4^{2^{t+1}-4}$\cr  $b_{t,34} = x_1x_2^{3}x_4^{2^{t+1}-4}$&$b_{t,35} = x_1x_2^{3}x_3^{2^{t+1}-4}$&$b_{t,36} = x_1^{3}x_3x_4^{2^{t+1}-4}$\cr  $b_{t,37} = x_1^{3}x_2x_4^{2^{t+1}-4}$&$b_{t,38} = x_1^{3}x_2x_3^{2^{t+1}-4}$&\cr  
\end{tabular}}

\medskip
For $t = 2$,

\medskip  
\centerline{\begin{tabular}{lll}
$b_{2,39} =  x_1^{3}x_2^{4}x_4$& $b_{2,40} =  x_1^{3}x_2^{4}x_3$& $b_{2,41} =  x_2^{3}x_3^{4}x_4$\cr  $b_{2,42} =  x_1^{3}x_3^{4}x_4$& $b_{2,43} =  x_1x_2x_3^{2}x_4^{4}$& $b_{2,44} =  x_1x_2^{2}x_3x_4^{4}$\cr  $b_{2,45} =  x_1x_2^{2}x_3^{4}x_4$& $b_{2,46} =  x_1x_2^{2}x_3^{2}x_4^{3}$& $b_{2,47} =  x_1x_2^{2}x_3^{3}x_4^{2}$\cr  $b_{2,48} =  x_1x_2^{3}x_3^{2}x_4^{2}$& $b_{2,49} =  x_1^{3}x_2x_3^{2}x_4^{2}$& \cr  
\end{tabular}}

\medskip
For $t \geqslant 3$,

\medskip  
\centerline{\begin{tabular}{lll}
$b_{t,39} = x_1x_2x_3^{2}x_4^{2^{t+1}-4}$&$b_{t,40} = x_1x_2^{2}x_3x_4^{2^{t+1}-4}$&$b_{t,41} = x_1x_2^{2}x_3^{2^{t+1}-4}x_4$\cr  $b_{t,42} = x_1x_2^{2}x_3^{4}x_4^{2^{t+1}-7}$&$b_{t,43} = x_1x_2^{2}x_3^{5}x_4^{2^{t+1}-8}$&$b_{t,44} = x_1x_2^{3}x_3^{4}x_4^{2^{t+1}-8}$\cr  $b_{t,45} = x_1^{3}x_2x_3^{4}x_4^{2^{t+1}-8}$& &\cr  
\end{tabular}}

\medskip
For $t \geqslant 3$,

\medskip  
\centerline{\begin{tabular}{lll}
$b_{2,46} =  x_1x_2^{3}x_3^{6}x_4^{6}$& $b_{2,47} =  x_1^{3}x_2x_3^{6}x_4^{6}$\cr  $b_{2,48} =  x_1^{3}x_2^{5}x_3^{2}x_4^{6}$& $b_{2,49} =  x_1^{3}x_2^{5}x_3^{6}x_4^{2}$\cr 
\end{tabular}}

\medskip
We prove the following.
\begin{thms}\label{dlt1} We have
$$(QP_4)_{2^{t+1}}^{GL_5} = \begin{cases}0, &\mbox{if } t \leqslant 3,\\
\langle [\zeta_{t,1}] \rangle, &\mbox{if } t \geqslant 4, 
 \end{cases} $$
where 
$$\zeta_{t,1} = \begin{cases}\psi_{1,4}(\xi_{4,0}), &\mbox{if } t = 4,\\
\psi_{1,t}(\xi_{t,0}) + b_{t,30} + b_{t,40} + b_{t,41}, &\mbox{if } t \geqslant 5. 
\end{cases}$$
Here, $\xi_{t,0}$ is determined as in Theorem \ref{dlmd1}.	
\end{thms}
\begin{rems}
In \cite[Page 471, Line 15$\uparrow$]{pp25}, the author stated that $(QP_4)_{n_{1,t}}^{GL_4}$ is generated by $[\psi_{1,t}(\xi_{t,0})]$ but this is false because the class	$[\psi_{1,t}(\xi_{t,0})]$ is not an $GL_4$-invariant in $(QP_4)_{n_{1,t}}$ for $t \geqslant 5$.
\end{rems}
\medskip
By a simple computation we observe that
\begin{equation}
\begin{cases}
[\Sigma_4(b_{t,1})] = \langle \{[b_{t,j}]: 1 \leqslant j \leqslant 12 \}\rangle,\ [\Sigma_4(b_{t,1})]^{\Sigma_4} = \langle [q_{t,1}] \rangle,\\
[\Sigma_4(b_{t,13})] = \langle \{[b_{t,j}] : 13 \leqslant j \leqslant 18 \}\rangle,\ [\Sigma_4(b_{t,13})]^{\Sigma_4} = \langle [q_{t,2}]\rangle\\
[\Sigma_4(b_{t,35})] = \begin{cases}
\langle \{[b_{2,j}]: 19 \leqslant j \leqslant 42\} \rangle, &\mbox{if } t = 2,\\
\langle \{[b_{t,j}]: 19 \leqslant j \leqslant 38\} \rangle, &\mbox{if } t > 2.
\end{cases}
\end{cases}
\end{equation}
where 
$q_{t,1}  = \mbox{$\sum_{1 \leqslant j \leqslant 12}$}b_{t,j},\,
 q_{t,2}  = \mbox{$\sum_{13 \leqslant j \leqslant 18}$}b_{t,j}.$

\begin{lems}\label{bdts21} For $t\geqslant 2$, $[\Sigma_4(b_{t,35})]^{\Sigma_4} = \langle [q_{t,3}]\rangle$, where
$$q_{t,3} = \begin{cases}
\sum_{j \in \{19,\, 20,\, 21,\, 22,\, 23,\, 24,\, 25,\, 26,\, 32,\, 36,\, 37,\, 38,\, 39,\, 40,\, 41,\, 42\}}b_{2,j}, &\mbox{if } t = 2,\\
\sum_{j \in \{20,\, 22,\, 25,\, 26,\, 27,\, 28,\, 29,\, 30,\, 31,\, 32,\, 33,\, 34,\, 35,\, 36,\, 37,\, 38\}}b_{t,j}, &\mbox{if } t > 2.
\end{cases}$$
\end{lems}

\begin{proof} Let $[f] \in [\Sigma_4(b_{t,35})]$ with $f \in P_4$. 
	
For $t = 2$, we have $f \equiv \sum_{19\leqslant j \leqslant 42}b_{2,j}$. By computing $\rho_i(f) + f$ in terms of the admissible monomials we have
\begin{align*}
\rho_1(f) &+ f \equiv \gamma_{\{19,21\}}b_{2,19} + \gamma_{\{20,22\}}b_{2,20} + \gamma_{\{19,21\}}b_{2,21} + \gamma_{\{20,22\}}b_{2,22} + \gamma_{29}b_{2,23}\\ 
& + \gamma_{30}b_{2,24} + \gamma_{\{25,39\}}b_{2,25} + \gamma_{\{26,40\}}b_{2,26} + \gamma_{\{27,28\}}b_{2,27} + \gamma_{\{27,28\}}b_{2,28}\\ 
& + \gamma_{\{31,33\}}b_{2,31} + \gamma_{\{32,36\}}b_{2,32} + \gamma_{\{31,33\}}b_{2,33} + \gamma_{\{34,37,39\}}b_{2,34}\\ 
& + \gamma_{\{35,38,40\}}b_{2,35} + \gamma_{\{32,36\}}b_{2,36} + \gamma_{\{25,34,37\}}b_{2,37} + \gamma_{\{26,35,38\}}b_{2,38}\\ 
& + \gamma_{\{25,39\}}b_{2,39} + \gamma_{\{26,40\}}b_{2,40} + \gamma_{\{41,42\}}b_{2,41} + \gamma_{\{41,42\}}b_{2,42} \equiv 0,\\
\rho_2(f) &+ f \equiv \gamma_{27}b_{2,19} + \gamma_{\{20,41\}}b_{2,20} + \gamma_{\{21,23\}}b_{2,21} + \gamma_{\{22,25\}}b_{2,22} + \gamma_{\{21,23\}}b_{2,23}\\ 
& + \gamma_{\{24,26\}}b_{2,24} + \gamma_{\{22,25\}}b_{2,25} + \gamma_{\{24,26\}}b_{2,26} + \gamma_{\{28,29\}}b_{2,28} + \gamma_{\{28,29\}}b_{2,29}\\ 
& + \gamma_{\{30,35\}}b_{2,30} + \gamma_{\{31,32,41\}}b_{2,31} + \gamma_{\{20,31,32\}}b_{2,32} + \gamma_{\{33,34\}}b_{2,33}\\ 
& + \gamma_{\{33,34\}}b_{2,34} + \gamma_{\{30,35\}}b_{2,35} + \gamma_{\{36,37\}}b_{2,36} + \gamma_{\{36,37\}}b_{2,37} + \gamma_{\{38,40\}}b_{2,38}\\ 
& + \gamma_{\{39,42\}}b_{2,39} + \gamma_{\{38,40\}}b_{2,40} + \gamma_{\{20,41\}}b_{2,41} + \gamma_{\{39,42\}}b_{2,42} \equiv 0,\\
\rho_3(f) &+ f \equiv \gamma_{\{19,20\}}b_{2,19} + \gamma_{\{19,20\}}b_{2,20} + \gamma_{\{21,22\}}b_{2,21} + \gamma_{\{21,22\}}b_{2,22}\\ 
& + \gamma_{\{23,24\}}b_{2,23} + \gamma_{\{23,24\}}b_{2,24} + \gamma_{\{25,26\}}b_{2,25} + \gamma_{\{25,26\}}b_{2,26}\\ 
& + \gamma_{\{27,31\}}b_{2,27} + \gamma_{\{28,33\}}b_{2,28} + \gamma_{\{29,30\}}b_{2,29} + \gamma_{\{29,30\}}b_{2,30}\\ 
& + \gamma_{\{27,31\}}b_{2,31} + \gamma_{\{32,41\}}b_{2,32} + \gamma_{\{28,33\}}b_{2,33} + \gamma_{\{34,35\}}b_{2,34}\\ 
& + \gamma_{\{34,35\}}b_{2,35} + \gamma_{\{36,42\}}b_{2,36} + \gamma_{\{37,38\}}b_{2,37} + \gamma_{\{37,38\}}b_{2,38}\\ 
& + \gamma_{\{39,40\}}b_{2,39} + \gamma_{\{39,40\}}b_{2,40} + \gamma_{\{32,41\}}b_{2,41} + \gamma_{\{36,42\}}b_{2,42} \equiv 0.
\end{align*}
From these equalities we get $\gamma_j = \gamma_{19}$ for $j \in \mathbb J_1 = \{$19,\, 20,\, 21,\, 22,\, 23,\, 24,\, 25,\, 26,\, 32,\, 36,\, 37,\, 38,\, 39,\, 40,\, 41,\, 42$\}$ and $\gamma_j = 0$ for $j \notin \mathbb J_1$. The lemma is proved for $t = 2$.

For $t \geqslant 3$, $f \equiv \sum_{19\leqslant j \leqslant 38}b_{t,j}$. By computing $\rho_i(f) + f$ in terms of the admissible monomials we have
\begin{align*}
\rho_1(f) &+ f \equiv \gamma_{\{19,21\}}b_{t,19} + \gamma_{\{20,22\}}b_{t,20} + \gamma_{\{19,21\}}b_{t,21} + \gamma_{\{20,22\}}b_{t,22}\\ 
& + \gamma_{\{25,29\}}b_{t,23} + \gamma_{\{26,30\}}b_{t,24} + \gamma_{\{27,28\}}b_{t,27} + \gamma_{\{27,28\}}b_{t,28} + \gamma_{\{31,33\}}b_{t,31}\\ 
& + \gamma_{\{32,36\}}b_{t,32} + \gamma_{\{31,33\}}b_{t,33} + \gamma_{\{34,37\}}b_{t,34} + \gamma_{\{35,38\}}b_{t,35} + \gamma_{\{32,36\}}b_{t,36}\\ 
& + \gamma_{\{34,37\}}b_{t,37} + \gamma_{\{35,38\}}b_{t,38} \equiv 0,\\
\rho_2(f) &+ f \equiv \gamma_{\{20,27\}}b_{t,19} + \gamma_{\{21,23\}}b_{t,21} + \gamma_{\{22,25\}}b_{t,22} + \gamma_{\{21,23\}}b_{t,23}\\ 
& + \gamma_{\{24,26,38\}}b_{t,24} + \gamma_{\{22,25\}}b_{t,25} + \gamma_{\{24,26,38\}}b_{t,26} + \gamma_{\{28,29\}}b_{t,28}\\ 
& + \gamma_{\{28,29\}}b_{t,29} + \gamma_{\{30,35\}}b_{t,30} + \gamma_{\{31,32\}}b_{t,31} + \gamma_{\{31,32\}}b_{t,32} + \gamma_{\{33,34\}}b_{t,33}\\ 
& + \gamma_{\{33,34\}}b_{t,34} + \gamma_{\{30,35\}}b_{t,35} + \gamma_{\{36,37\}}b_{t,36} + \gamma_{\{36,37\}}b_{t,37} \equiv 0,\\
\rho_3(f) &+ f \equiv \gamma_{\{19,20,32\}}b_{t,19} + \gamma_{\{19,20,32\}}b_{t,20} + \gamma_{\{21,22,36\}}b_{t,21} + \gamma_{\{21,22,36\}}b_{t,22}\\ 
& + \gamma_{\{23,24\}}b_{t,23} + \gamma_{\{23,24\}}b_{t,24} + \gamma_{\{25,26\}}b_{t,25} + \gamma_{\{25,26\}}b_{t,26} + \gamma_{\{27,31\}}b_{t,27}\\ 
& + \gamma_{\{28,33\}}b_{t,28} + \gamma_{\{29,30\}}b_{t,29} + \gamma_{\{29,30\}}b_{t,30} + \gamma_{\{27,31\}}b_{t,31} + \gamma_{\{28,33\}}b_{t,33}\\ 
& + \gamma_{\{34,35\}}b_{t,34} + \gamma_{\{34,35\}}b_{t,35} + \gamma_{\{37,38\}}b_{t,37} + \gamma_{\{37,38\}}b_{t,38} \equiv 0.
\end{align*}
From these equalities we get $\gamma_j = \gamma_{20}$ for $j \in \mathbb J_2 = \{$20,\, 22,\, 25,\, 26,\, 27,\, 28,\, 29,\, 30,\, 31,\, 32,\, 33,\, 34,\, 35,\, 36,\, 37,\, 38$\}$ and $\gamma_j = 0$ for $j \notin \mathbb J_2$. The lemma is proved.
\end{proof}

We set
$$\mathcal W_t = \begin{cases} 
\langle [b_{2,j}]: 43\leqslant j \leqslant 45\rangle, &\mbox{if } t = 2,\\
\langle [b_{t,j}]: 39\leqslant j \leqslant 45\rangle, &\mbox{if } t \geqslant 3.
\end{cases}$$
It is easy to verify that $\mathcal W_t$ is an $\Sigma_4$-submodule of $\mbox{Ker}((\widetilde{Sq}^0_*)_{(4,2^{t+1})})$.
\begin{lems}\label{bdts22} We have 

\medskip
{\rm i)} $QP_4((2,3))^{GL_4} = 0$ and $QP_4((2,3,2))^{GL_4} = 0$.

{\rm ii)} $\mathcal W_t^{\Sigma_4} = \begin{cases} 
\langle [q_{2,4}]\rangle &\mbox{if } t = 2,\\
\langle [q_{t,4}],\, [q_{t,5}] \rangle, &\mbox{if } t \geqslant 3.\end{cases}$

\medskip\noindent
where $q_{2,4}= \sum_{43\leqslant j \leqslant 45}b_{2,j}$, $q_{t,4} = \sum_{39\leqslant j \leqslant 41}b_{t,j}$, $q_{t,5} = \sum_{42\leqslant j \leqslant 45}b_{t,j}$, for $t \geqslant 3$.
\end{lems}
\begin{proof} It is easy to see that $QP_4((2,3)) = \langle [b_{2,j}]_{(2,3)}: 46\leqslant j \leqslant 49 \rangle$, $QP_4((2,3))^{\Sigma_4} = \langle [q_{2,5}]_{(2,3)}\rangle$ with $q_{2,5} = \sum_{46\leqslant j \leqslant 49}b_{2,j}$. So, if $[f]_{(2,3)} \in QP_4((2,3))^{GL_4}$ with $f \in P_4$, then $f \equiv_{(2,3)} \gamma q_{2,5}$ with $\gamma \in \mathbb F_2$. Then we have $\rho_4(f) + f \equiv_{(2,3)} \gamma b_{2,48} \equiv_{(2,3)} 0$. This implies $\gamma =0$ and $f \equiv_{(2,3)} 0$. The relation $QP_4((2,3,2))^{GL_4} = 0$ is easily proved by a similar computation.
	
We prove Part ii). The case $t = 2$ is easy, so we assume $t > 2$.  If $[f] \in \mathcal W_t$ with $f \in P_4$, then $f \equiv \sum_{39\leqslant j \leqslant 45} \gamma_jb_{t,j}$ with $\gamma_j \in \mathbb F_2$. By a simple computation we have
\begin{align*}
\rho_1(f) &+ f \equiv \gamma_{\{40,41,42,43\}}b_{t,39} + \gamma_{\{44,45\}}b_{t,44} + \gamma_{\{44,45\}}b_{t,45} \equiv 0,\\
\rho_2(f) &+ f \equiv \gamma_{\{39,40,45\}}b_{t,39} + \gamma_{\{39,40,45\}}b_{t,40} + \gamma_{\{43,44\}}b_{t,43} + \gamma_{\{43,44\}}b_{t,44} \equiv 0,\\
\rho_3(f) &+ f \equiv \gamma_{\{40,41\}}b_{t,40} + \gamma_{\{40,41\}}b_{t,41} + \gamma_{\{42,43\}}b_{t,42} + \gamma_{\{42,43\}}b_{t,43} \equiv 0.
\end{align*}
These equalities imply $\gamma_{40} = \gamma_{41}$, $\gamma_{j} = \gamma_{42}$ for $j = 43,\, 44,\, 45$ and $\gamma_{39} = \gamma_{40} + \gamma_{42}$. Hence, we obtain
$f \equiv \gamma_{40}(b_{t,39} + b_{t,40} + b_{t,41}) + \gamma_{42}(b_{t,39} + b_{t,42} + b_{t,43} + b_{t,44} + b_{t,45}) = \gamma_{40}q_{t,4} + \gamma_{42}q_{t,5}$. The lemma is proved.	
\end{proof}

\begin{proof}[Proof of Theorem \ref{dlt1}] Let $[f] \in (QP_4)_{2^{t+1}}^{GL_4}$ with $f \in (P_4)_{2^{t+1}}$. Consider Kameko's homomorphism
$$(\widetilde{Sq}^0_*)_{(4,2^{t+1})} : (QP_{4})_{2^{t+1}}\to (QP_4)_{2^{t}-2}.$$

For $t = 2$, we have a direct summand decomposition of $\Sigma_4$-submodules
$$\mbox{Ker}((\widetilde{Sq}^0_*)_{(4,8)}) = \Sigma_4(b_{2,1})\bigoplus\Sigma_4(b_{2,13})\bigoplus \Sigma_4(b_{2,35})\bigoplus\mathcal W_2 \bigoplus QP_4((2,3)).$$
Since $(QP_4)_2^{GL_4} = 0$, by applying Lemmas \ref{bdts21} and \ref{bdts22} we have $f \equiv \sum_{1\leqslant u \leqslant 4}q_{2,u}$ with $\gamma_u \in \mathbb F_2$. By computing $\rho_4(f) + f$ in terms of the admissible monomials, we get
\begin{align*}
\rho_4(f) &+ f \equiv \gamma_{1}b_{2,3} + \gamma_{1}b_{2,4} + \gamma_{1}b_{2,5} + \gamma_{1}b_{2,6} + \gamma_{\{1,2\}}b_{2,9} + \gamma_{2}b_{2,14} + \gamma_{2}b_{2,15}\\ 
& + \gamma_{3}b_{2,19} + \gamma_{3}b_{2,20} + \gamma_{2}b_{2,23} + \gamma_{2}b_{2,24} + \gamma_{1}b_{2,25} + \gamma_{1}b_{2,26} + \gamma_{\{3,4\}}b_{2,32}\\ 
& + \gamma_{1}b_{2,37} + \gamma_{1}b_{2,38} + \gamma_{1}b_{2,39} + \gamma_{1}b_{2,40} + \gamma_{\{3,4\}}b_{2,41} \equiv 0. 
\end{align*}
This equality implies $\gamma_u = 0$ for $1\leqslant u \leqslant 4$. The theorem is proved for $t = 2$.
	
For $t = 3$, we have a direct summand decomposition of $\Sigma_4$-submodules
$$\mbox{Ker}((\widetilde{Sq}^0_*)_{(4,16)}) = \Sigma_4(b_{3,1})\bigoplus\Sigma_4(b_{3,13})\bigoplus \Sigma_4(b_{3,35})\bigoplus\mathcal W_3 \bigoplus QP_4((2,3,2)).$$
By Theorem \ref{dlmd1}, $(QP_4)_6^{GL_4} = 0$, so by applying Lemmas \ref{bdts21} and \ref{bdts22} we have $f \equiv \sum_{1\leqslant u \leqslant 5}q_{3,u}$ with $\gamma_u \in \mathbb F_2$. By computing $\rho_4(f) + f$ in terms of the admissible monomials, we get
\begin{align*}
\rho_4(f) &+ f \equiv \gamma_{1}b_{3,3} + \gamma_{1}b_{3,4} + \gamma_{\{1,3\}}b_{3,5} + \gamma_{\{1,3\}}b_{3,6} + \gamma_{\{1,2\}}b_{3,9} + \gamma_{\{2,3\}}b_{3,14}\\ 
& + \gamma_{\{2,3\}}b_{3,15} + \gamma_{4}b_{3,19} + \gamma_{\{3,4\}}b_{3,20} + \gamma_{\{1,2\}}b_{3,23} + \gamma_{\{1,2\}}b_{3,24} + \gamma_{\{3,5\}}b_{3,27}\\ 
& + \gamma_{\{3,5\}}b_{3,31} + \gamma_{3}b_{3,32} + \gamma_{3}b_{3,34} + \gamma_{3}b_{3,35} + \gamma_{3}b_{3,39} + \gamma_{5}b_{3,44} \equiv 0. 
\end{align*}
This equality implies $\gamma_u = 0$ for $1\leqslant u \leqslant 5$. The theorem is proved for $t = 3$.

For $t \geqslant 4$, we have a direct summand decomposition of $\Sigma_4$-submodules
$$\mbox{Ker}((\widetilde{Sq}^0_*)_{(4,2^{t+1})}) = \Sigma_4(b_{t,1})\bigoplus\Sigma_4(b_{t,13})\bigoplus \Sigma_4(b_{t,35})\bigoplus\mathcal W_t.$$
Since $(\widetilde{Sq}^0_*)_{(4,2^{t+1})}([f]) \in (QP_4)_{2^{t}-2}^{GL_4}$, applying Theorem \ref{dlmd1} we see that there is $\gamma_0 \in \mathbb F_2$ such that $g = f + \lambda_0\psi_{1,t}(\xi_{t,0}) \in \mbox{Ker}((\widetilde{Sq}^0_*)_{(4,2^{t+1})})$. 

For $t = 4$, we have $\rho_i(\xi_{4,0}) + \xi_{4,0} \equiv 0$ for $1 \leqslant i \leqslant 4$, hence by applying Lemmas \ref{bdts21} and \ref{bdts22} we have $g \equiv \sum_{1\leqslant u \leqslant 5}q_{4,u}$ with $\gamma_u \in \mathbb F_2$. By an argument analogous to the one for the case $t = 3$, from the relation $\rho_4(g) + g \equiv 0$ we obtain $\gamma_u = 0$ for $1\leqslant u \leqslant 5$. Hence, $f \equiv \lambda_0\psi_{1,4}(\xi_{4,0}) = \lambda_0\zeta_{4,1}$ and the theorem holds for $t = 4$.

For $t \geqslant 5$, we have $\rho_i(\xi_{t,0}) + \xi_{t,0} \equiv 0$ for $1 \leqslant i \leqslant 3$ and 
$$\rho_4(\xi_{t,0}) + \xi_{t,0} \equiv x_2x_3x_4^{2^{t+1}-2} + x_2x_3^{2^{t+1}-2}x_4 = b_{t,19} + b_{t,20}.$$
Hence, by applying Lemmas \ref{bdts21} and \ref{bdts22} we have $g \equiv \sum_{1\leqslant u \leqslant 5}q_{t,u}$ with $\gamma_u \in \mathbb F_2$. By computing $\rho_4(g) + g$ in terms of the admissible monomials we obtain
\begin{align*}
\rho_4(g) &+ g \equiv \gamma_{1}b_{t,3} + \gamma_{1}b_{t,4} + \gamma_{\{1,3\}}b_{t,5} + \gamma_{\{1,3\}}b_{t,6} + \gamma_{\{1,2\}}b_{t,9}\\ 
& + \gamma_{\{2,3\}}b_{t,14} + \gamma_{\{2,3\}}b_{t,15} + \gamma_{\{0,4\}}b_{t,19} + \gamma_{\{0,3,4\}}b_{t,20}\\ 
& + \gamma_{\{1,2\}}b_{t,23} + \gamma_{\{1,2\}}b_{t,24} + \gamma_{\{3,5\}}b_{t,27} + \gamma_{\{3,5\}}b_{t,31} + \gamma_{3}b_{t,32}\\ 
& + \gamma_{3}b_{t,34} + \gamma_{3}b_{t,35} + \gamma_{3}b_{t,39} + \gamma_{5}b_{t,44} \equiv 0. 
\end{align*} 
This equality implies $\gamma_u = 0$ for $u = 1,\, 2,\, 3,\, 5$ and $\gamma_4 = \gamma_0$. Hence, $f \equiv \lambda_0(\psi_{1,t}(\xi_{t,0})+q_{t,4}) = \lambda_0\zeta_{t,1}$ and the theorem is completely proved.
\end{proof}

\subsection{The case $t = 1$ and $s \geqslant 2$}\

\medskip
We have $n_{s,1} = 2^{s+1} + 2^s - 2$. By our works \cite{su50,su5}, $$\mbox{Ker}\Big((\widetilde{Sq}^0_*)_{(4,n_{s,1})}\Big) = QP_4((2)|^s|(1)|) \cong  QP_4^0((2)|^s|(1)|)\bigoplus QP_4^+((2)|^s|(1)|).$$ 

A basis of $QP_4^0((2)|^{s}|(1))$ is the set of all classes represented by the admissible monomials $b_{s,j} := b_{1,s,j}$, which are determined as follows:

\medskip
For $s \geqslant 2$,

\medskip  
\centerline{\begin{tabular}{lll}
$b_{s,1} = x_3^{2^{s}-1}x_4^{2^{s+1}-1}$&$b_{s,2} = x_3^{2^{s+1}-1}x_4^{2^{s}-1}$\cr $b_{s,3} = x_2^{2^{s}-1}x_4^{2^{s+1}-1}$&$b_{s,4} = x_2^{2^{s}-1}x_3^{2^{s+1}-1}$\cr $b_{s,5} = x_2^{2^{s+1}-1}x_4^{2^{s}-1}$&$b_{s,6} = x_2^{2^{s+1}-1}x_3^{2^{s}-1}$\cr $b_{s,7} = x_1^{2^{s}-1}x_4^{2^{s+1}-1}$&$b_{s,8} = x_1^{2^{s}-1}x_3^{2^{s+1}-1}$\cr $b_{s,9} = x_1^{2^{s}-1}x_2^{2^{s+1}-1}$&$b_{s,10} = x_1^{2^{s+1}-1}x_4^{2^{s}-1}$\cr $b_{s,11} = x_1^{2^{s+1}-1}x_3^{2^{s}-1}$&$b_{s,12} = x_1^{2^{s+1}-1}x_2^{2^{s}-1}$\cr $b_{s,13} = x_2x_3^{2^{s}-2}x_4^{2^{s+1}-1}$&$b_{s,14} = x_2x_3^{2^{s+1}-1}x_4^{2^{s}-2}$\cr $b_{s,15} = x_2^{2^{s+1}-1}x_3x_4^{2^{s}-2}$&$b_{s,16} = x_1x_3^{2^{s}-2}x_4^{2^{s+1}-1}$\cr $b_{s,17} = x_1x_3^{2^{s+1}-1}x_4^{2^{s}-2}$&$b_{s,18} = x_1x_2^{2^{s}-2}x_4^{2^{s+1}-1}$\cr $b_{s,19} = x_1x_2^{2^{s}-2}x_3^{2^{s+1}-1}$&$b_{s,20} = x_1x_2^{2^{s+1}-1}x_4^{2^{s}-2}$\cr $b_{s,21} = x_1x_2^{2^{s+1}-1}x_3^{2^{s}-2}$&$b_{s,22} = x_1^{2^{s+1}-1}x_3x_4^{2^{s}-2}$\cr $b_{s,23} = x_1^{2^{s+1}-1}x_2x_4^{2^{s}-2}$&$b_{s,24} = x_1^{2^{s+1}-1}x_2x_3^{2^{s}-2}$\cr $b_{s,25} = x_2x_3^{2^{s}-1}x_4^{2^{s+1}-2}$&$b_{s,26} = x_2x_3^{2^{s+1}-2}x_4^{2^{s}-1}$\cr $b_{s,27} = x_2^{2^{s}-1}x_3x_4^{2^{s+1}-2}$&$b_{s,28} = x_1x_3^{2^{s}-1}x_4^{2^{s+1}-2}$\cr $b_{s,29} = x_1x_3^{2^{s+1}-2}x_4^{2^{s}-1}$&$b_{s,30} = x_1x_2^{2^{s}-1}x_4^{2^{s+1}-2}$\cr $b_{s,31} = x_1x_2^{2^{s}-1}x_3^{2^{s+1}-2}$&$b_{s,32} = x_1x_2^{2^{s+1}-2}x_4^{2^{s}-1}$\cr 
\end{tabular}}
\centerline{\begin{tabular}{lll}
$b_{s,33} = x_1x_2^{2^{s+1}-2}x_3^{2^{s}-1}$&$b_{s,34} = x_1^{2^{s}-1}x_3x_4^{2^{s+1}-2}$\cr $b_{s,35} = x_1^{2^{s}-1}x_2x_4^{2^{s+1}-2}$&$b_{s,36} = x_1^{2^{s}-1}x_2x_3^{2^{s+1}-2}$\cr $b_{s,37} = x_2^{3}x_3^{2^{s+1}-3}x_4^{2^{s}-2}$&$b_{s,38} = x_1^{3}x_3^{2^{s+1}-3}x_4^{2^{s}-2}$\cr $b_{s,39} = x_1^{3}x_2^{2^{s+1}-3}x_4^{2^{s}-2}$&$b_{s,40} = x_1^{3}x_2^{2^{s+1}-3}x_3^{2^{s}-2}$\cr 
\end{tabular}}

\medskip
For $s = 2$,

\medskip  
\centerline{\begin{tabular}{lll}
$b_{2,41} =  x_2^{3}x_3^{3}x_4^{4}$& $b_{2,42} =  x_1^{3}x_3^{3}x_4^{4}$\cr $b_{2,43} =  x_1^{3}x_2^{3}x_4^{4}$&  $b_{2,44} =  x_1^{3}x_2^{3}x_3^{4}$\cr
\end{tabular}}

\medskip
For $s \geqslant 3$,

\medskip  
\centerline{\begin{tabular}{lll}
$b_{s,41} = x_2^{3}x_3^{2^{s}-3}x_4^{2^{s+1}-2}$&$b_{s,42} = x_1^{3}x_3^{2^{s}-3}x_4^{2^{s+1}-2}$\cr $b_{s,43} = x_1^{3}x_2^{2^{s}-3}x_4^{2^{s+1}-2}$&$b_{s,44} = x_1^{3}x_2^{2^{s}-3}x_3^{2^{s+1}-2}$\cr
\end{tabular}}

\medskip
A basis of $QP_4^+((2)|^{s}|(1))$ is the set of all classes represented by the admissible monomials $b_{s,j} := b_{1,s,j}$, which are determined as follows:

\medskip
For $s \geqslant 2$,

\medskip  
\centerline{\begin{tabular}{lll}
$b_{s,45} = x_1x_2x_3^{2^{s}-2}x_4^{2^{s+1}-2}$&$b_{s,46} = x_1x_2x_3^{2^{s+1}-2}x_4^{2^{s}-2}$\cr $b_{s,47} = x_1x_2^{2^{s}-2}x_3x_4^{2^{s+1}-2}$&$b_{s,48} = x_1x_2^{2^{s+1}-2}x_3x_4^{2^{s}-2}$\cr $b_{s,49} = x_1x_2^{2}x_3^{2^{s+1}-3}x_4^{2^{s}-2}$&\cr 
\end{tabular}}

\medskip
For $s = 2$,

\medskip  
\centerline{\begin{tabular}{lll}
$b_{2,50} =  x_1x_2^{2}x_3^{3}x_4^{4}$& $b_{2,51} =  x_1x_2^{2}x_3^{4}x_4^{3}$& $b_{2,52} =  x_1x_2^{3}x_3^{2}x_4^{4}$\cr  $b_{2,53} =  x_1x_2^{3}x_3^{4}x_4^{2}$& $b_{2,54} =  x_1^{3}x_2x_3^{2}x_4^{4}$& $b_{2,55} =  x_1^{3}x_2x_3^{4}x_4^{2}$\cr  $b_{2,56} =  x_1^{3}x_2^{4}x_3x_4^{2}$& &\cr
\end{tabular}}

\medskip
For $s \geqslant 3$,

\medskip  
\centerline{\begin{tabular}{lll}
$b_{s,50} = x_1x_2^{2}x_3^{2^{s}-4}x_4^{2^{s+1}-1}$&$b_{s,51} = x_1x_2^{2}x_3^{2^{s+1}-1}x_4^{2^{s}-4}$\cr $b_{s,52} = x_1x_2^{2^{s+1}-1}x_3^{2}x_4^{2^{s}-4}$&$b_{s,53} = x_1^{2^{s+1}-1}x_2x_3^{2}x_4^{2^{s}-4}$\cr $b_{s,54} = x_1x_2^{2}x_3^{2^{s}-3}x_4^{2^{s+1}-2}$&$b_{s,55} = x_1x_2^{2}x_3^{2^{s}-1}x_4^{2^{s+1}-4}$\cr $b_{s,56} = x_1x_2^{2}x_3^{2^{s+1}-4}x_4^{2^{s}-1}$&$b_{s,57} = x_1x_2^{2^{s}-1}x_3^{2}x_4^{2^{s+1}-4}$\cr $b_{s,58} = x_1^{2^{s}-1}x_2x_3^{2}x_4^{2^{s+1}-4}$&$b_{s,59} = x_1x_2^{3}x_3^{2^{s}-4}x_4^{2^{s+1}-2}$\cr $b_{s,60} = x_1x_2^{3}x_3^{2^{s+1}-2}x_4^{2^{s}-4}$&$b_{s,61} = x_1^{3}x_2x_3^{2^{s}-4}x_4^{2^{s+1}-2}$\cr $b_{s,62} = x_1^{3}x_2x_3^{2^{s+1}-2}x_4^{2^{s}-4}$&$b_{s,63} = x_1x_2^{3}x_3^{2^{s}-2}x_4^{2^{s+1}-4}$\cr $b_{s,64} = x_1x_2^{3}x_3^{2^{s+1}-4}x_4^{2^{s}-2}$&$b_{s,65} = x_1^{3}x_2x_3^{2^{s}-2}x_4^{2^{s+1}-4}$\cr $b_{s,66} = x_1^{3}x_2x_3^{2^{s+1}-4}x_4^{2^{s}-2}$&$b_{s,67} = x_1^{3}x_2^{2^{s+1}-3}x_3^{2}x_4^{2^{s}-4}$\cr $b_{s,68} = x_1^{3}x_2^{2^{s}-3}x_3^{2}x_4^{2^{s+1}-4}$&$b_{s,69} = x_1^{3}x_2^{5}x_3^{2^{s+1}-6}x_4^{2^{s}-4}$\cr 
\end{tabular}}

\medskip
For $s = 3$,\, $b_{3,70} =  x_1^{3}x_2^{5}x_3^{6}x_4^{8}$.

\medskip
For $s \geqslant 4$, $b_{s,70} = x_1^{3}x_2^{5}x_3^{2^{s}-6}x_4^{2^{s+1}-4}$.

\medskip
The following implies Theorem \ref{thm2} for $t = 1$.

\begin{thms}\label{dl2t1} For $n_{s,1} = 2^{s+1} + 2^s -2$, we have $(QP_4)_{n_{s,1}}^{GL_4} = \langle [\zeta_{1,s}]\rangle$, where
$$\zeta_{1,s} = \begin{cases} 
0, &\mbox{if } s = 2,\, 4,\\
\psi_{3,1}(\xi_{1,2}) + \sum_{j \in \{45,\, 49,\, 54,\, 58,\, 61,\, 64,\, 70\}}b_{3,j}, &\mbox{if } s = 3,\\
\psi_{s,1}(\xi_{1,s-1}) + \sum_{j \in \{45,\, 46,\, 50,\, 51,\, 52,\, 53,\, 54,\, 55,\, 56\}\atop\hskip0.4cm\cup\{58,\, 59,\, 60,\, 62,\, 63,\, 64,\, 67,\, 70\}}b_{s,j}, &\mbox{if } s \geqslant 5.
\end{cases} $$
Here, $\xi_{1,s-1}$ is defined in Theorem \ref{mdct1}.	
\end{thms}

By a direct computation, we observe that
\begin{align}
\begin{cases}\label{ctdl2t1}
[\Sigma_4(b_{s,1})] = \langle \{[b_{s,j}]: 1 \leqslant j \leqslant 12 \}\rangle,\ [\Sigma_4(b_{s,1})]^{\Sigma_4} = \langle [q_{s,1}]\rangle,\\
[\Sigma_4(b_{s,13})] = \langle \{[b_{s,j}]: 13 \leqslant j \leqslant 24\} \rangle,\ [\Sigma_4(b_{s,13})]^{\Sigma_4} = \langle [q_{s,2}]\rangle,\\
\end{cases}
\end{align}
where
$q_{s,1} = \sum_{1 \leqslant j \leqslant 12}b_{s,j}, \ q_{s,2} = \sum_{13 \leqslant j \leqslant 24}b_{s,j}.$

\begin{lems}\label{bdd2t11} We have

\medskip
{\rm i)} $[\Sigma_4(b_{2,25})]^{\Sigma_4} = \langle [q_{2,3}] \rangle$, $[\Sigma_4(b_{2,50})]^{\Sigma_4} = \langle [q_{2,4}] \rangle$, where $$q_{2,3} = b_{2,26} + b_{2,27} + b_{2,29} + \sum_{32\leqslant j\leqslant 44}b_{2,j},\ q_{2,4} =\sum_{45\leqslant j\leqslant 56,\, j \ne 47,\, 49,\, 53}b_{2,j}.$$ 

{\rm ii)} For $s \geqslant 3$, 
$[\Sigma_4(b_{s,25}, b _{s,37})]^{\Sigma_4} = \langle [q_{s,3}] \rangle$,
where $ q_{s,3} =\sum_{25\leqslant j\leqslant 40}b_{s,j}.$ 

\medskip
Consequently, $QP_4(2,2,1)^{\Sigma_4} = \langle [q_{2,u}] : 1 \leqslant u \leqslant 4\rangle$ and $QP_4^0((2)|^s|(1))^{\Sigma_4} = \langle [q_{s,u}] : 1 \leqslant u \leqslant 3\rangle$ for $s \geqslant 3$.
\end{lems}
\begin{proof} We have $[\Sigma_4(b_{2,25})] = \langle [b_{2,j}] : 25\leqslant j\leqslant 44\rangle$. If $[f] \in [\Sigma_4(b_{2,25})]^{\Sigma_4}$, then $f \equiv \gamma_j\sum_{25\leqslant j\leqslant 44}b_{2,j}$. A direct computation gives
\begin{align*}
\rho_1(f) &+ f \equiv \gamma_{\{25,28\}}b_{2,25} + \gamma_{\{26,29\}}b_{2,26} + \gamma_{\{27,34\}}b_{2,27} + \gamma_{\{25,28\}}b_{2,28}\\ 
& + \gamma_{\{26,29\}}b_{2,29} + \gamma_{\{30,32,35\}}b_{2,30} + \gamma_{\{31,33,36\}}b_{2,31} + \gamma_{\{27,34\}}b_{2,34}\\ 
& + \gamma_{\{30,32,35\}}b_{2,35} + \gamma_{\{31,33,36\}}b_{2,36} + \gamma_{\{37,38\}}b_{2,37} + \gamma_{\{37,38\}}b_{2,38}\\ 
& + \gamma_{\{41,42\}}b_{2,41} + \gamma_{\{41,42\}}b_{2,42} + \gamma_{\{32,39\}}b_{2,43} + \gamma_{\{33,40\}}b_{2,44} \equiv 0,\\
\rho_2(f) &+ f \equiv \gamma_{\{25,26,27\}}b_{2,25} + \gamma_{\{25,26,27\}}b_{2,27} + \gamma_{\{28,30\}}b_{2,28} + \gamma_{\{29,32\}}b_{2,29}\\ 
& + \gamma_{\{28,30\}}b_{2,30} + \gamma_{\{31,33,44\}}b_{2,31} + \gamma_{\{29,32\}}b_{2,32} + \gamma_{\{31,33,44\}}b_{2,33}\\ 
& + \gamma_{\{34,35\}}b_{2,34} + \gamma_{\{34,35\}}b_{2,35} + \gamma_{\{36,40\}}b_{2,36} + \gamma_{\{38,39\}}b_{2,38}\\ 
& + \gamma_{\{38,39\}}b_{2,39} + \gamma_{\{36,40\}}b_{2,40} + \gamma_{\{26,37\}}b_{2,41} + \gamma_{\{42,43\}}b_{2,42}\\ 
& + \gamma_{\{42,43\}}b_{2,43} \equiv 0,\\
\rho_3(f) &+ f \equiv \gamma_{\{25,26,41\}}b_{2,25} + \gamma_{\{25,26,41\}}b_{2,26} + \gamma_{\{27,37\}}b_{2,27}\\ 
& + \gamma_{\{28,29,42\}}b_{2,28} + \gamma_{\{28,29,42\}}b_{2,29} + \gamma_{\{30,31\}}b_{2,30} + \gamma_{\{30,31\}}b_{2,31}\\ 
& + \gamma_{\{32,33\}}b_{2,32} + \gamma_{\{32,33\}}b_{2,33} + \gamma_{\{34,38\}}b_{2,34} + \gamma_{\{35,36\}}b_{2,35}\\ 
& + \gamma_{\{35,36\}}b_{2,36} + \gamma_{\{27,37\}}b_{2,37} + \gamma_{\{34,38\}}b_{2,38} + \gamma_{\{39,40\}}b_{2,39}\\ 
& + \gamma_{\{39,40\}}b_{2,40} + \gamma_{\{43,44\}}b_{2,43} + \gamma_{\{43,44\}}b_{2,44} \equiv 0.
\end{align*}
By computing from these equalities we get $\gamma_j = 0$ for $j = 25,\, 28,\, 30,\, 31$ and $\gamma_j = \gamma_{26}$ for $j \ne 25,\, 28,\, 30,\, 31$. Hence, $f \equiv \gamma_{26}q_{2,3}$.

We have $[\Sigma_4(b_{2,50})] = QP_4^+((2)|^2|(1)) = \langle [b_{2,j}] : 45 \leqslant j \leqslant 56 \rangle$. If $[f] \in [\Sigma_4(b_{2,50})]^{\Sigma_4}$, then $f \equiv \gamma_j\sum_{45\leqslant j\leqslant 56}b_{2,j}$. A direct computation gives
\begin{align*}
\rho_1(f) &+ f \equiv \gamma_{47}b_{2,45} + \gamma_{49}b_{2,46} + \gamma_{\{48,56\}}b_{2,48} + \gamma_{\{52,54\}}b_{2,52}\\ 
& + \gamma_{\{53,55,56\}}b_{2,53} + \gamma_{\{52,54\}}b_{2,54} + \gamma_{\{48,53,55\}}b_{2,55} + \gamma_{\{48,56\}}b_{2,56} \equiv 0,\\
\rho_2(f) &+ f \equiv \gamma_{\{45,47,51\}}b_{2,45} + \gamma_{\{46,48\}}b_{2,46} + \gamma_{\{45,47,51\}}b_{2,47} + \gamma_{\{46,48\}}b_{2,48}\\ 
& + \gamma_{\{49,53\}}b_{2,49} + \gamma_{\{50,52,53\}}b_{2,50} + \gamma_{\{49,50,52\}}b_{2,52} + \gamma_{\{49,53\}}b_{2,53}\\ 
& + \gamma_{\{55,56\}}b_{2,55} + \gamma_{\{55,56\}}b_{2,56} \equiv 0,\\
\rho_3(f) &+ f \equiv \gamma_{\{45,46,49\}}b_{2,45} + \gamma_{\{45,46,47\}}b_{2,46} + \gamma_{\{47,49\}}b_{2,47} + \gamma_{\{47,49\}}b_{2,49}\\ 
& + \gamma_{\{50,51\}}b_{2,50} + \gamma_{\{50,51\}}b_{2,51} + \gamma_{\{48,52,53\}}b_{2,52} + \gamma_{\{48,52,53\}}b_{2,53}\\ 
& + \gamma_{\{54,55\}}b_{2,54} + \gamma_{\{54,55\}}b_{2,55} \equiv 0.
\end{align*}
From these equalities we get $\gamma_j = 0$ for $j = 47,\, 49,\, 53$ and $\gamma_j = \gamma_{45}$ for $j \ne 47,\, 49,\, 53$. Hence, $f \equiv \gamma_{45}q_{2,4}$.

For $ t \geqslant 3$, $[\Sigma_4(b_{s,25}, b _{s,37})]^{\Sigma_4} = \langle [b_{s,j}]: 25\leqslant j\leqslant 44 \rangle$. If $[f] \in [\Sigma_4(b_{s,25},b _{s,37})]^{\Sigma_4}$, then $f \equiv \gamma_j\sum_{25\leqslant j\leqslant 44}b_{s,j}$. A direct computation gives
\begin{align*}
\rho_1(f) &+ f \equiv \gamma_{\{25,28\}}b_{s,25} + \gamma_{\{26,29\}}b_{s,26} + \gamma_{\{27,34\}}b_{s,27} + \gamma_{\{25,28\}}b_{s,28}\\
& + \gamma_{\{26,29\}}b_{s,29} + \gamma_{\{30,35\}}b_{s,30} + \gamma_{\{31,36\}}b_{s,31} + \gamma_{\{27,34\}}b_{s,34}\\
& + \gamma_{\{30,35\}}b_{s,35} + \gamma_{\{31,36\}}b_{s,36} + \gamma_{\{37,38\}}b_{s,37} + \gamma_{\{37,38\}}b_{s,38}\\
& + \gamma_{\{41,42\}}b_{s,41} + \gamma_{\{41,42\}}b_{s,42} + \gamma_{\{32,39\}}b_{s,43} + \gamma_{\{33,40\}}b_{s,44} \equiv 0,\\
\rho_2(f) &+ f \equiv \gamma_{\{25,27\}}b_{s,25} + \gamma_{\{25,27\}}b_{s,27} + \gamma_{\{28,30\}}b_{s,28} + \gamma_{\{29,32\}}b_{s,29}\\
& + \gamma_{\{28,30\}}b_{s,30} + \gamma_{\{31,33\}}b_{s,31} + \gamma_{\{29,32\}}b_{s,32} + \gamma_{\{31,33\}}b_{s,33}\\
& + \gamma_{\{34,35\}}b_{s,34} + \gamma_{\{34,35\}}b_{s,35} + \gamma_{\{38,39\}}b_{s,38} + \gamma_{\{38,39\}}b_{s,39}\\
& + \gamma_{\{36,40,44\}}b_{s,40} + \gamma_{\{26,37\}}b_{s,41} + \gamma_{\{42,43\}}b_{s,42} + \gamma_{\{42,43\}}b_{s,43}\\
& + \gamma_{\{36,40,44\}}b_{s,44} \equiv 0,\\
\rho_3(f) &+ f \equiv \gamma_{\{25,26\}}b_{s,25} + \gamma_{\{25,26\}}b_{s,26} + \gamma_{\{28,29\}}b_{s,28} + \gamma_{\{28,29\}}b_{s,29}\\
& + \gamma_{\{30,31\}}b_{s,30} + \gamma_{\{30,31\}}b_{s,31} + \gamma_{\{32,33\}}b_{s,32} + \gamma_{\{32,33\}}b_{s,33}\\
& + \gamma_{\{35,36\}}b_{s,35} + \gamma_{\{35,36\}}b_{s,36} + \gamma_{\{27,37,41\}}b_{s,37} + \gamma_{\{34,38,42\}}b_{s,38}\\
& + \gamma_{\{39,40\}}b_{s,39} + \gamma_{\{39,40\}}b_{s,40} + \gamma_{\{27,37,41\}}b_{s,41} + \gamma_{\{34,38,42\}}b_{s,42}\\
& + \gamma_{\{43,44\}}b_{s,43} + \gamma_{\{43,44\}}b_{s,44} \equiv 0.
\end{align*}
From these equalities we get $\gamma_j = 0$ for $j = 41,\, 42,\, 43,\, 44$ and $\gamma_j = \gamma_{25}$ for $26\leqslant j \leqslant 40$. Hence, $f \equiv \gamma_{25}q_{s,3}$.
\end{proof}

\begin{lems}\label{bdd2t12} For $s \geqslant 3$, $QP_4^+((2)|^s|(1))^{\Sigma_4} = \langle [q_{s,u}]: 4 \leqslant j \leqslant 7\rangle$, where
\begin{align*}
q_{s,4} &= \mbox{$\sum_{50\leqslant j \leqslant 53}$}b_{s,j},\\
q_{s,5} &= \begin{cases}
\sum_{j \in \{46,\, 49,\, 54,\, 64,\, 66,\, 68,\, 69,\, 70\}}b_{3,j}, &\mbox{if } s = 3,\\
\sum_{j \in \{45,\, 46,\, 49,\, 54,\, 63,\, 64,\, 65,\, 66,\, 68,\, 69\}}b_{s,j}, &\mbox{if } s \geqslant 4,
\end{cases}\\
q_{s,6} &= \mbox{$\sum_{j \in \{47,\, 48,\, 54,\, 59,\, 60,\, 61,\, 62,\, 67\}}$}b_{s,j},\\
q_{s,7} &= \begin{cases}
\sum_{j \in \{45,\, 55,\, 56,\, 57,\, 58,\, 59,\, 61,\, 70\}}b_{3,j}, &\mbox{if } s = 3,\\
\sum_{j \in \{55,\, 56,\, 57,\, 58,\, 59,\, 61,\, 63,\, 65,\, 68,\, 70\}}b_{s,j}, &\mbox{if } s \geqslant 4,
\end{cases}
\end{align*}
\end{lems}
\begin{proof} For $s\geqslant 3$, we have a direct summand decomposition of $\Sigma_4$-submodule
$$QP_4^+((2)|^s|(1)) = [\Sigma_4(b_{s,50})]\bigoplus \mathcal W_{1,s},$$
where $\mathcal W_{1,s} = \langle [b_{s,j}]: 45 \leqslant j \leqslant 70,\, j \ne 50,\, 51,\, 52,\, 53\rangle$. It is easy to see that $[\Sigma_4(b_{s,50})] : \langle [b_{s,j}]: 50 \leqslant j \leqslant 53\rangle$ and $[\Sigma_4(b_{s,50})]^{\Sigma_4} = \langle [q_{s,4}]\rangle$ and $q_{s,u}$, $4 \leqslant j \leqslant 7$, are $\Sigma_4$-invariants.  

For $s = 3$, the leading monomials of $q_{3,5}$, $q_{3,6}$, $q_{3,7}$ respectively are $b_{3,69}$, $b_{3,67}$, $b_{3,58}$. Hence, if $[f] \in \mathcal W_{1,3}^{\Sigma_4}$, then there are $\gamma_u$, $u = 5,\, 6,\, 7$, such that  $h = f + \gamma_5q_{3,5} + \gamma_5q_{3,6} + \gamma_5q_{3,7} \equiv \sum_{j \in \{45,\ldots,49,54,\ldots, 70\}\setminus \{58,67,69\}}\gamma_jb_{3,j}$. By computing $\rho_i(h) + h$ in terms of the admissible monomials one gets 
\begin{align*}
\rho_1(h) &+ h \equiv \gamma_{\{47,54,68\}}b_{3,45} + \gamma_{\{48,49\}}b_{3,46} + \gamma_{57}b_{3,57} + \gamma_{57}b_{3,58} + \gamma_{\{59,61\}}b_{3,59}\\
& + \gamma_{\{60,62\}}b_{3,60} + \gamma_{\{59,61\}}b_{3,61} + \gamma_{\{60,62\}}b_{3,62} + \gamma_{\{63,65\}}b_{3,63}\\
& + \gamma_{\{64,66\}}b_{3,64} + \gamma_{\{63,65\}}b_{3,65} + \gamma_{\{64,66\}}b_{3,66} + \gamma_{48}b_{3,70} \equiv 0,\\
\rho_2(h) &+ h \equiv \gamma_{\{45,47,49,61,63,64\}}b_{3,45} + \gamma_{\{46,48,60,66\}}b_{3,46} + \gamma_{\{45,47,61,63\}}b_{3,47}\\
& + \gamma_{\{46,48,60,66\}}b_{3,48} + \gamma_{\{49,64\}}b_{3,49} + \gamma_{\{54,56,59,64\}}b_{3,54} + \gamma_{\{55,57\}}b_{3,55}\\
& + \gamma_{\{55,57\}}b_{3,57} + \gamma_{\{49,54,56,59\}}b_{3,59} + \gamma_{62}b_{3,62} + \gamma_{\{49,64\}}b_{3,63} + \gamma_{\{49,64\}}b_{3,64}\\
& + \gamma_{\{65,66,68\}}b_{3,65} + \gamma_{62}b_{3,67} + \gamma_{\{65,66,68\}}b_{3,68} + \gamma_{66}b_{3,70} \equiv 0,\\
\rho_3(h) &+ h \equiv \gamma_{\{45,46,49,57,70\}}b_{3,45} + \gamma_{\{45,46,47,54,57,68,70\}}b_{3,46} + \gamma_{\{47,49,54\}}b_{3,49}\\
& + \gamma_{\{47,49,54\}}b_{3,54} + \gamma_{\{55,56\}}b_{3,55} + \gamma_{\{55,56\}}b_{3,56} + \gamma_{\{57,59,60\}}b_{3,59}\\
& + \gamma_{\{57,59,60\}}b_{3,60} + \gamma_{\{61,62\}}b_{3,61} + \gamma_{\{61,62\}}b_{3,62} + \gamma_{\{48,57,63,64,70\}}b_{3,63}\\
& + \gamma_{\{48,57,63,64,70\}}b_{3,64} + \gamma_{\{65,66,70\}}b_{3,65} + \gamma_{\{65,66,70\}}b_{3,66} + \gamma_{68}b_{3,68}\\
& + \gamma_{68}b_{3,69} \equiv 0.
\end{align*}
From these equalities we get $\gamma_j = 0$ for $j \in \{45,\ldots,49,54,\ldots, 70\}\setminus \{58,67,69\}$. Hence, $h \equiv 0$ and $f \equiv \gamma_5q_{3,5} + \gamma_5q_{3,6} + \gamma_5q_{3,7}$.

For $s \geqslant 4$, the leading monomials of $q_{s,5}$, $q_{s,6}$, $q_{s,7}$ respectively are $b_{s,68}$, $b_{s,67}$, $b_{s,58}$. Hence, if $[f] \in \mathcal W_{1,s}^{\Sigma_4}$, then there are $\gamma_u$, $u = 5,\, 6,\, 7$, such that  $h = f + \gamma_5q_{s,5} + \gamma_5q_{s,6} + \gamma_5q_{s,7} \equiv \sum_{j \in \{45,\ldots,49,54,\ldots, 70\}\setminus \{58,67,68\}}\gamma_jb_{s,j}$. By computing $\rho_i(h) + h$ in terms of the admissible monomials one gets 
\begin{align*}
\rho_1(h) &+ h \equiv \gamma_{\{47,48,54,70\}}b_{s,45} + \gamma_{\{48,49,69\}}b_{s,46} + \gamma_{57}b_{s,57} + \gamma_{57}b_{s,58}\\
& + \gamma_{\{59,61\}}b_{s,59} + \gamma_{\{60,62\}}b_{s,60} + \gamma_{\{59,61\}}b_{s,61} + \gamma_{\{60,62\}}b_{s,62}\\
& + \gamma_{\{48,63,65\}}b_{s,63} + \gamma_{\{64,66\}}b_{s,64} + \gamma_{\{48,63,65\}}b_{s,65} + \gamma_{\{64,66\}}b_{s,66}\\
& + \gamma_{48}b_{s,68} + \gamma_{48}b_{s,70} \equiv 0,\\
\rho_2(h) &+ h \equiv \gamma_{\{45,47,49,61,63,64,66,69\}}b_{s,45} + \gamma_{\{46,48,60,66\}}b_{s,46} + \gamma_{\{45,47,61,63\}}b_{s,47}\\
& + \gamma_{\{46,48,60,66\}}b_{s,48} + \gamma_{\{49,64\}}b_{s,49} + \gamma_{\{54,56,59,64\}}b_{s,54} + \gamma_{\{55,57\}}b_{s,55}\\
& + \gamma_{\{55,57\}}b_{s,57} + \gamma_{\{49,54,56,59\}}b_{s,59} + \gamma_{62}b_{s,62} + \gamma_{\{49,64,66,69\}}b_{s,63}\\
& + \gamma_{\{49,64\}}b_{s,64} + \gamma_{65}b_{s,65} + \gamma_{62}b_{s,67} + \gamma_{65}b_{s,68} + \gamma_{\{66,69\}}b_{s,70} \equiv 0,\\
\rho_3(h) &+ h \equiv \gamma_{\{45,46,49,57,69\}}b_{s,45} + \gamma_{\{45,46,47,54,57,70\}}b_{s,46} + \gamma_{\{47,49,54\}}b_{s,49}\\
& + \gamma_{\{47,49,54\}}b_{s,54} + \gamma_{\{55,56\}}b_{s,55} + \gamma_{\{55,56\}}b_{s,56} + \gamma_{\{57,59,60\}}b_{s,59}\\
& + \gamma_{\{57,59,60\}}b_{s,60} + \gamma_{\{61,62\}}b_{s,61} + \gamma_{\{61,62\}}b_{s,62} + \gamma_{\{48,57,63,64\}}b_{s,63}\\
& + \gamma_{\{48,57,63,64\}}b_{s,64} + \gamma_{\{65,66\}}b_{s,65} + \gamma_{\{65,66\}}b_{s,66} + \gamma_{\{69,70\}}b_{s,69}\\
& + \gamma_{\{69,70\}}b_{s,70} \equiv 0.
\end{align*}
From these equalities we get $\gamma_j = 0$ for $j \in \{45,\ldots,49,54,\ldots, 70\}\setminus \{58,67,68\}$. Hence, $h \equiv 0$ and $f \equiv \gamma_5q_{3,5} + \gamma_5q_{3,6} + \gamma_5q_{3,7}$. The lemma is proved.
\end{proof}
\begin{proof}[Proof of Theorem \ref{dl2t1}] Let $[f] \in (QP_4)_{n_{s,1}}^{GL_4}$ with $f \in (P_4)_{n_{s,1}}$. Since $(\widetilde{Sq}^0_*)_{(4,n_{s,1})}$ is a homomorphism of $GL_4$-module, we have $(\widetilde{Sq}^0_*)_{(4,n_{s,1})}([f]) \in (QP_4)_{d_{s-1,1}}^{GL_4}$.

For $s = 2$, $d_{1,1} = 3$ and $(QP_4)_{3}^{GL_4} = 0$, hence $[f] \in \mbox{Ker}\big((\widetilde{Sq}^0_*)_{(4,10)}\big)$. By Lemma \ref{bdd2t11}, $f \equiv \sum_{1 \leqslant u\leqslant 4}\gamma_uq_{2,4}$. By computing $\rho_4(f) + f$ in terms of the admissible monomials, we have
\begin{align*}	
\rho_4(f) &+ f \equiv \gamma_{\{1,2\}}b_{2,3} + \gamma_{\{1,2\}}b_{2,4} + \gamma_{\{1,3\}}b_{2,5} + \gamma_{\{1,3\}}b_{2,6} + \gamma_{1}b_{2,9} + \gamma_{2}b_{2,13}\\
& + \gamma_{2}b_{2,14} + \gamma_{2}b_{2,15} + \gamma_{\{2,3\}}b_{2,20} + \gamma_{\{2,3\}}b_{2,21} + \gamma_{4}b_{2,25} + \gamma_{\{3,4\}}b_{2,26}\\
& + \gamma_{3}b_{2,27} + \gamma_{3}b_{2,30} + \gamma_{3}b_{2,31} + \gamma_{3}b_{2,37} + \gamma_{3}b_{2,41} + \gamma_{\{1,2\}}b_{2,43} + \gamma_{\{1,2\}}b_{2,44}\\
& + \gamma_{3}b_{2,45} + \gamma_{3}b_{2,46} + \gamma_{2}b_{2,48} + \gamma_{2}b_{2,52} + \gamma_{2}b_{2,54} + \gamma_{2}b_{2,55} + \gamma_{2}b_{2,56} \equiv 0.
\end{align*}
This equality implies $\gamma_u =0$ for $1 \leqslant u\leqslant 4$. Hence, $[f] = 0$ and the theorem holds for $s= 2$.

For $s = 3$, by using Theorem \ref{mdct1}, we see that there is $\lambda_0 \in \mathbb F_2$ such that $f + \lambda_0\psi_{1,3}(\xi_{1,2}) \in \mbox{Ker}\big((\widetilde{Sq}^0_*)_{(4,22)}\big)$. By a direct computation we have
\begin{equation}\label{cthrt11}
\begin{cases}
\rho_1(\psi_{3,1}(\xi_{1,2}))+\psi_{3,1}(\xi_{1,2}) \equiv \sum_{j \in \{45,\, 46,\, 57,\, 58,\, 59,\, 61,\, 64,\, 66\}}b_{3,j},\\
\rho_2(\psi_{3,1}(\xi_{1,2}))+\psi_{3,1}(\xi_{1,2}) \equiv 0,\\
\rho_3(\psi_{3,1}(\xi_{1,2}))+\psi_{3,1}(\xi_{1,2}) \equiv 0,\\
\rho_4(\psi_{3,1}(\xi_{1,2}))+\psi_{3,1}(\xi_{1,2}) \equiv b_{3,37} + b_{3,41},\\
\end{cases}
\end{equation}
By a computation using \eqref{cthrt11} and Lemmas \ref{bdd2t11}, \ref{bdd2t12} we see that $\rho_i(f) + f \equiv 0$ for $i = 1,\, 2,\, 3$, if and only if $f \equiv \lambda_{0}(\psi_{3,1}(\xi_{1,2}) + q_{3,0}) + \sum_{1 \leqslant u \leqslant 7}\gamma_uq_{3,u}$ with $\gamma_u \in \mathbb F_2$ and $q_{3,0} = \sum_{j \in \{49,\, 54,\, 55,\, 56,\, 57,\, 59,\, 64\}}b_{3,j}$. By computing $\rho_4(f) + f$ in terms of the admissible monomials we obtain
\begin{align*}	
\rho_4(f) &+ f \equiv \gamma_{\{1,2\}}b_{3,3} + \gamma_{\{1,2\}}b_{3,4} + \gamma_{\{1,3\}}b_{3,5} + \gamma_{\{1,3\}}b_{3,6} + \gamma_{1}b_{3,9} + \gamma_{\{2,4\}}b_{3,13}\\
& + \gamma_{\{2,4\}}b_{3,14} + \gamma_{\{2,6\}}b_{3,15} + \gamma_{\{2,3\}}b_{3,20} + \gamma_{\{2,3\}}b_{3,21} + \gamma_{\{0,3,7\}}b_{3,25}\\
& + \gamma_{\{0,3,7\}}b_{3,26} + \gamma_{\{3,6\}}b_{3,27} + \gamma_{\{1,2,3\}}b_{3,30} + \gamma_{\{1,2,3\}}b_{3,31} + \gamma_{\{1,2\}}b_{3,35}\\
& + \gamma_{\{1,2\}}b_{3,36} + \gamma_{\{3,5\}}b_{3,37} + \gamma_{\{5,6\}}b_{3,41} + \gamma_{\{1,2\}}b_{3,43} + \gamma_{\{1,2\}}b_{3,44} + \gamma_{3}b_{3,45}\\
& + \gamma_{\{2,3,4,5\}}b_{3,46} + \gamma_{\{4,6\}}b_{3,52} + \gamma_{\{2,4,5\}}b_{3,57} + \gamma_{\{2,4\}}b_{3,58} + \gamma_{\{2,4,6\}}b_{3,59}\\
& + \gamma_{\{5,6\}}b_{3,60} + \gamma_{\{2,4\}}b_{3,61} + \gamma_{5}b_{3,63} + \gamma_{\{2,4\}}b_{3,70} \equiv 0.
\end{align*}
This equality implies $\gamma_u =0$ for $1 \leqslant u\leqslant 6$ and $\gamma_7 = \gamma_0$. Hence, $f \equiv \lambda_{0}(\psi_{3,1}(\xi_{1,2}) + q_{3,0} + q_{3,7}) = \lambda_0\zeta_{1,3}$ and the theorem holds for $s= 3$.

For $s = 4$, by Theorem \ref{mdct1}, $(QP_4)_{21}^{GL_4} = 0$, hence $f \in \mbox{Ker}\big((\widetilde{Sq}^0_*)_{(4,46)}\big)$. By using Lemmas \ref{bdd2t11}, \ref{bdd2t12} we have $f \equiv \sum_{1 \leqslant u \leqslant 7}\gamma_uq_{4,u}$ with $\gamma_u \in \mathbb F_2$. By computing $\rho_4(f) + f$ in terms of the admissible monomials we obtain
\begin{align*}	
\rho_4(f) &+ f \equiv \gamma_{\{1,2\}}b_{4,3} + \gamma_{\{1,2\}}b_{4,4} + \gamma_{\{1,3\}}b_{4,5} + \gamma_{\{1,3\}}b_{4,6} + \gamma_{1}b_{4,9} + \gamma_{\{2,4\}}b_{4,13}\\
& + \gamma_{\{2,4\}}b_{4,14} + \gamma_{\{2,6\}}b_{4,15} + \gamma_{\{2,3\}}b_{4,20} + \gamma_{\{2,3\}}b_{4,21} + \gamma_{\{3,7\}}b_{4,25}\\
& + \gamma_{\{3,7\}}b_{4,26} + \gamma_{\{3,6\}}b_{4,27} + \gamma_{\{1,2,3\}}b_{4,30} + \gamma_{\{1,2,3\}}b_{4,31} + \gamma_{\{1,2\}}b_{4,35}\\
& + \gamma_{\{1,2\}}b_{4,36} + \gamma_{\{3,5\}}b_{4,37} + \gamma_{\{5,6\}}b_{4,41} + \gamma_{\{1,2\}}b_{4,43} + \gamma_{\{1,2\}}b_{4,44}\\
& + \gamma_{\{2,3,4\}}b_{4,45} + \gamma_{\{2,3,4,5\}}b_{4,46} + \gamma_{\{4,6\}}b_{4,52} + \gamma_{\{2,4,5\}}b_{4,57} + \gamma_{\{2,4\}}b_{4,58}\\
& + \gamma_{\{2,4,6\}}b_{4,59} + \gamma_{\{5,6\}}b_{4,60} + \gamma_{\{2,4\}}b_{4,61} + \gamma_{\{2,4,5\}}b_{4,63} + \gamma_{\{2,4\}}b_{4,65}\\
& + \gamma_{2}b_{4,66} + \gamma_{\{2,4\}}b_{4,68} + \gamma_{\{2,4\}}b_{4,70} \equiv 0.
\end{align*}
This equality implies $\gamma_u =0$ for $1 \leqslant u\leqslant 7$. Hence, $f \equiv 0$ and the theorem holds for $s = 4$.

For $s \geqslant 5$, by using Theorem \ref{mdct1}, we see that there is $\lambda_0 \in \mathbb F_2$ such that $f + \lambda_0\psi_{1,s}(\xi_{1,s-1}) \in \mbox{Ker}\big((\widetilde{Sq}^0_*)_{(4,n_{s,1})}\big)$. By a direct computation we have
\begin{equation}\label{cthrt12}
\begin{cases}
\rho_1(\psi_{s,1}(\xi_{1,s-1}))+\psi_{s,1}(\xi_{1,s-1}) \equiv \sum_{j \in \{\scriptscriptstyle 45,\, 46,\, 57,\, 58,\, 59,\, 61,\, 64,\, 66,\, 68,\, 70\}}b_{s,j},\\
\rho_2(\psi_{s,1}(\xi_{1,s-1}))+\psi_{s,1}(\xi_{1,s-1}) \equiv \sum_{j \in \{\scriptscriptstyle 45,\, 49,\, 55,\, 57,\, 59,\, 63,\, 64\}}b_{s,j},\\
\rho_3(\psi_{s,1}(\xi_{1,s-1}))+\psi_{s,1}(\xi_{1,s-1}) \equiv \sum_{j \in \{\scriptscriptstyle 45,\, 46,\, 49,\, 54,\, 63,\, 64,\, 65,\, 66,\, 69,\, 70\}}b_{s,j},\\
\rho_4(\psi_{s,1}(\xi_{1,s-1}))+\psi_{s,1}(\xi_{1,s-1}) \equiv \sum_{j \in \{13,\, 14,\, 25,\, 26,\, 41,\, 45,\, 46,\, 58,\, 60,\, 61\}\atop\hskip2.6cm\cup\{65,\, 68,\, 70\}}b_{s,j},\\
\end{cases}
\end{equation}
By a computation using \eqref{cthrt12} and Lemmas \ref{bdd2t11}, \ref{bdd2t12} we see that $\rho_i(f) + f \equiv 0$ for $i = 1,\, 2,\, 3$, if and only if $f \equiv \lambda_{0}(\psi_{s,1}(\xi_{1,s-1}) + q_{s,0}) + \sum_{1 \leqslant u \leqslant 7}\gamma_uq_{s,u}$ with $\gamma_u \in \mathbb F_2$ and $q_{s,0} = \sum_{j \in \{47,\, 48,\, 49,\, 54,\, 57,\, 59,\, 63,\, 66,\, 69\}}b_{s,j}$. By computing $\rho_4(f) + f$ in terms of the admissible monomials we obtain
\begin{align*}	
\rho_4(f) &+ f \equiv \gamma_{\{1,2\}}b_{s,3} + \gamma_{\{1,2\}}b_{s,4} + \gamma_{\{1,3\}}b_{s,5} + \gamma_{\{1,3\}}b_{s,6} + \gamma_{1}b_{s,9} + \gamma_{\{0,2,4\}}b_{s,13}\\
& + \gamma_{\{0,2,4\}}b_{s,14} + \gamma_{\{0,2,6\}}b_{s,15} + \gamma_{\{2,3\}}b_{s,20} + \gamma_{\{2,3\}}b_{s,21} + \gamma_{\{0,3,7\}}b_{s,25}\\
& + \gamma_{\{0,3,7\}}b_{s,26} + \gamma_{\{0,3,6\}}b_{s,27} + \gamma_{\{1,2,3\}}b_{s,30} + \gamma_{\{1,2,3\}}b_{s,31} + \gamma_{\{1,2\}}b_{s,35}\\
& + \gamma_{\{1,2\}}b_{s,36} + \gamma_{\{0,3,5\}}b_{s,37} + \gamma_{\{5,6\}}b_{s,41} + \gamma_{\{1,2\}}b_{s,43} + \gamma_{\{1,2\}}b_{s,44}\\
& + \gamma_{\{0,2,3,4\}}b_{s,45} + \gamma_{\{2,3,4,5\}}b_{s,46} + \gamma_{\{4,6\}}b_{s,52} + \gamma_{\{2,4,5\}}b_{s,57} + \gamma_{\{0,2,4\}}b_{s,58}\\
& + \gamma_{\{2,4,6\}}b_{s,59} + \gamma_{\{5,6\}}b_{s,60} + \gamma_{\{0,2,4\}}b_{s,61} + \gamma_{\{2,4,5\}}b_{s,63} + \gamma_{\{0,2,4\}}b_{s,65}\\
& + \gamma_{2}b_{s,66} + \gamma_{\{0,2,4\}}b_{s,68} + \gamma_{\{0,2,4\}}b_{s,70} \equiv 0.
\end{align*}
This equality implies $\gamma_u =0$ for $1 \leqslant u\leqslant 3$ and $\gamma_u = \gamma_0$ for $4 \leqslant u\leqslant 7$. Hence, $f \equiv \lambda_{0}(\psi_{s,1}(\xi_{1,s-1}) + q_{s,0} + \sum_{4\leqslant u \leqslant 7}q_{s,u}) = \lambda_0\zeta_{1,s}$ and the theorem is completely proved.
\end{proof}

\begin{rems}
The dimensional result of Theorem \ref{dl2t1} is also stated in \cite{pp25} but its proof is seriously false because Lemma 4.2.1 on Page 467 and Lemma 4.2.2 on Page 468 are false.
\end{rems}

\subsection{The case $s \geqslant 2$, $t \geqslant 2$}\

\medskip

By our works \cite{su50,su5}, $\mbox{Ker}\Big((\widetilde{Sq}^0_*)_{(4,n_{s,t})}\Big) = QP_4((2)|^s|(1)|^t)$ and we have $$QP_4((2)|^s|(1)|^t) = QP_4^0((2)|^s|(1)|^t) \bigoplus QP_4^+((2)|^s|(1)|^t).$$

A basis of $QP_4^0((2)|^{s}|(1)|^t)$ is the set of all classes represented by the admissible monomials $b_{s,j} = b_{t,s,j}$ which are determined as follows:

\medskip
For $s \geqslant 2$,

\medskip
\centerline{\begin{tabular}{lll}
$b_{s,1} = x_3^{2^{s}-1}x_4^{2^{s+t}-1}$&$b_{s,2} = x_3^{2^{s+t}-1}x_4^{2^{s}-1}$\cr  $b_{s,3} = x_2^{2^{s}-1}x_4^{2^{s+t}-1}$&$b_{s,4} = x_2^{2^{s}-1}x_3^{2^{s+t}-1}$\cr  $b_{s,5} = x_2^{2^{s+t}-1}x_4^{2^{s}-1}$&$b_{s,6} = x_2^{2^{s+t}-1}x_3^{2^{s}-1}$\cr  $b_{s,7} = x_1^{2^{s}-1}x_4^{2^{s+t}-1}$&$b_{s,8} = x_1^{2^{s}-1}x_3^{2^{s+t}-1}$\cr  $b_{s,9} = x_1^{2^{s}-1}x_2^{2^{s+t}-1}$&$b_{s,10} = x_1^{2^{s+t}-1}x_4^{2^{s}-1}$\cr  $b_{s,11} = x_1^{2^{s+t}-1}x_3^{2^{s}-1}$&$b_{s,12} = x_1^{2^{s+t}-1}x_2^{2^{s}-1}$\cr  $b_{s,13} = x_3^{2^{s+1}-1}x_4^{2^{s+t}-2^{s}-1}$&$b_{s,14} = x_2^{2^{s+1}-1}x_4^{2^{s+t}-2^{s}-1}$\cr  $b_{s,15} = x_2^{2^{s+1}-1}x_3^{2^{s+t}-2^{s}-1}$&$b_{s,16} = x_1^{2^{s+1}-1}x_4^{2^{s+t}-2^{s}-1}$\cr  $b_{s,17} = x_1^{2^{s+1}-1}x_3^{2^{s+t}-2^{s}-1}$&$b_{s,18} = x_1^{2^{s+1}-1}x_2^{2^{s+t}-2^{s}-1}$\cr  
\end{tabular}}
\centerline{\begin{tabular}{lll}
$b_{s,19} = x_2x_3^{2^{s}-2}x_4^{2^{s+t}-1}$&$b_{s,20} = x_2x_3^{2^{s+t}-1}x_4^{2^{s}-2}$\cr  $b_{s,21} = x_2^{2^{s+t}-1}x_3x_4^{2^{s}-2}$&$b_{s,22} = x_1x_3^{2^{s}-2}x_4^{2^{s+t}-1}$\cr  $b_{s,23} = x_1x_3^{2^{s+t}-1}x_4^{2^{s}-2}$&$b_{s,24} = x_1x_2^{2^{s}-2}x_4^{2^{s+t}-1}$\cr  $b_{s,25} = x_1x_2^{2^{s}-2}x_3^{2^{s+t}-1}$&$b_{s,26} = x_1x_2^{2^{s+t}-1}x_4^{2^{s}-2}$\cr  $b_{s,27} = x_1x_2^{2^{s+t}-1}x_3^{2^{s}-2}$&$b_{s,28} = x_1^{2^{s+t}-1}x_3x_4^{2^{s}-2}$\cr  $b_{s,29} = x_1^{2^{s+t}-1}x_2x_4^{2^{s}-2}$&$b_{s,30} = x_1^{2^{s+t}-1}x_2x_3^{2^{s}-2}$\cr  $b_{s,31} = x_2x_3^{2^{s}-1}x_4^{2^{s+t}-2}$&$b_{s,32} = x_2x_3^{2^{s+t}-2}x_4^{2^{s}-1}$\cr  $b_{s,33} = x_2^{2^{s}-1}x_3x_4^{2^{s+t}-2}$&$b_{s,34} = x_1x_3^{2^{s}-1}x_4^{2^{s+t}-2}$\cr  $b_{s,35} = x_1x_3^{2^{s+t}-2}x_4^{2^{s}-1}$&$b_{s,36} = x_1x_2^{2^{s}-1}x_4^{2^{s+t}-2}$\cr  $b_{s,37} = x_1x_2^{2^{s}-1}x_3^{2^{s+t}-2}$&$b_{s,38} = x_1x_2^{2^{s+t}-2}x_4^{2^{s}-1}$\cr  $b_{s,39} = x_1x_2^{2^{s+t}-2}x_3^{2^{s}-1}$&$b_{s,40} = x_1^{2^{s}-1}x_3x_4^{2^{s+t}-2}$\cr  $b_{s,41} = x_1^{2^{s}-1}x_2x_4^{2^{s+t}-2}$&$b_{s,42} = x_1^{2^{s}-1}x_2x_3^{2^{s+t}-2}$\cr  $b_{s,43} = x_2x_3^{2^{s+1}-2}x_4^{2^{s+t}-2^{s}-1}$&$b_{s,44} = x_1x_3^{2^{s+1}-2}x_4^{2^{s+t}-2^{s}-1}$\cr  $b_{s,45} = x_1x_2^{2^{s+1}-2}x_4^{2^{s+t}-2^{s}-1}$&$b_{s,46} = x_1x_2^{2^{s+1}-2}x_3^{2^{s+t}-2^{s}-1}$\cr  $b_{s,47} = x_2x_3^{2^{s+1}-1}x_4^{2^{s+t}-2^{s}-2}$&$b_{s,48} = x_2^{2^{s+1}-1}x_3x_4^{2^{s+t}-2^{s}-2}$\cr  $b_{s,49} = x_1x_3^{2^{s+1}-1}x_4^{2^{s+t}-2^{s}-2}$&$b_{s,50} = x_1x_2^{2^{s+1}-1}x_4^{2^{s+t}-2^{s}-2}$\cr  $b_{s,51} = x_1x_2^{2^{s+1}-1}x_3^{2^{s+t}-2^{s}-2}$&$b_{s,52} = x_1^{2^{s+1}-1}x_3x_4^{2^{s+t}-2^{s}-2}$\cr  $b_{s,53} = x_1^{2^{s+1}-1}x_2x_4^{2^{s+t}-2^{s}-2}$&$b_{s,54} = x_1^{2^{s+1}-1}x_2x_3^{2^{s+t}-2^{s}-2}$\cr  $b_{s,55} = x_2^{3}x_3^{2^{s+t}-3}x_4^{2^{s}-2}$&$b_{s,56} = x_1^{3}x_3^{2^{s+t}-3}x_4^{2^{s}-2}$\cr  $b_{s,57} = x_1^{3}x_2^{2^{s+t}-3}x_4^{2^{s}-2}$&$b_{s,58} = x_1^{3}x_2^{2^{s+t}-3}x_3^{2^{s}-2}$\cr  $b_{s,59} = x_2^{3}x_3^{2^{s+1}-3}x_4^{2^{s+t}-2^{s}-2}$&$b_{s,60} = x_1^{3}x_3^{2^{s+1}-3}x_4^{2^{s+t}-2^{s}-2}$\cr  $b_{s,61} = x_1^{3}x_2^{2^{s+1}-3}x_4^{2^{s+t}-2^{s}-2}$&$b_{s,62} = x_1^{3}x_2^{2^{s+1}-3}x_3^{2^{s+t}-2^{s}-2}$\cr  	 
\end{tabular}}

\medskip
For $s = 2$,

\medskip
\centerline{\begin{tabular}{lll}
$b_{2,63} = x_2^{3}x_3^{3}x_4^{2^{t+2}-4}$  &$b_{2,64} = x_1^{3}x_3^{3}x_4^{2^{t+2}-4}$\cr  $b_{2,65} = x_1^{3}x_2^{3}x_4^{2^{t+2}-4}$  &$b_{2,66} = x_1^{3}x_2^{3}x_3^{2^{t+2}-4}$\cr         
\end{tabular}}

\medskip
For $s \geqslant 3$,

\medskip
\centerline{\begin{tabular}{lll}
$b_{s,63} = x_2^{3}x_3^{2^{s}-3}x_4^{2^{s+t}-2}$&$b_{s,64} = x_1^{3}x_3^{2^{s}-3}x_4^{2^{s+t}-2}$\cr  $b_{s,65} = x_1^{3}x_2^{2^{s}-3}x_4^{2^{s+t}-2}$&$b_{s,66} = x_1^{3}x_2^{2^{s}-3}x_3^{2^{s+t}-2}$\cr       
\end{tabular}}

\medskip
A basis of $QP_4^+((2)|^{s}|(1)|^t)$ is the set of all classes represented by the admissible monomials $b_{s,j} = b_{t,s,j}$ which are determined as follows:

\medskip
For $s \geqslant 2$,

\medskip
\centerline{\begin{tabular}{lll}
$b_{s,67} = x_1x_2x_3^{2^{s}-2}x_4^{2^{s+t}-2}$&$b_{s,68} = x_1x_2x_3^{2^{s+t}-2}x_4^{2^{s}-2}$\cr  $b_{s,69} = x_1x_2^{2^{s}-2}x_3x_4^{2^{s+t}-2}$&$b_{s,70} = x_1x_2^{2^{s+t}-2}x_3x_4^{2^{s}-2}$\cr  $b_{s,71} = x_1x_2x_3^{2^{s+1}-2}x_4^{2^{s+t}-2^{s}-2}$&$b_{s,72} = x_1x_2^{2^{s+1}-2}x_3x_4^{2^{s+t}-2^{s}-2}$\cr  $b_{s,73} = x_1x_2^{2}x_3^{2^{s+t}-3}x_4^{2^{s}-2}$&$b_{s,74} = x_1x_2^{2}x_3^{2^{s}-1}x_4^{2^{s+t}-4}$\cr  $b_{s,75} = x_1x_2^{2}x_3^{2^{s+t}-4}x_4^{2^{s}-1}$&$b_{s,76} = x_1x_2^{2^{s}-1}x_3^{2}x_4^{2^{s+t}-4}$\cr  $b_{s,77} = x_1^{2^{s}-1}x_2x_3^{2}x_4^{2^{s+t}-4}$&$b_{s,78} = x_1x_2^{2}x_3^{2^{s+1}-4}x_4^{2^{s+t}-2^{s}-1}$\cr  $b_{s,79} = x_1x_2^{2}x_3^{2^{s+1}-3}x_4^{2^{s+t}-2^{s}-2}$&$b_{s,80} = x_1x_2^{2}x_3^{2^{s+1}-1}x_4^{2^{s+t}-2^{s}-4}$\cr  
\end{tabular}}
\centerline{\begin{tabular}{lll}
$b_{s,81} = x_1x_2^{2^{s+1}-1}x_3^{2}x_4^{2^{s+t}-2^{s}-4}$&$b_{s,82} = x_1^{2^{s+1}-1}x_2x_3^{2}x_4^{2^{s+t}-2^{s}-4}$\cr  $b_{s,83} = x_1x_2^{3}x_3^{2^{s+t}-4}x_4^{2^{s}-2}$&$b_{s,84} = x_1^{3}x_2x_3^{2^{s+t}-4}x_4^{2^{s}-2}$\cr  $b_{s,85} = x_1x_2^{3}x_3^{2^{s+1}-4}x_4^{2^{s+t}-2^{s}-2}$&$b_{s,86} = x_1^{3}x_2x_3^{2^{s+1}-4}x_4^{2^{s+t}-2^{s}-2}$\cr  $b_{s,87} = x_1x_2^{3}x_3^{2^{s+1}-2}x_4^{2^{s+t}-2^{s}-4}$&$b_{s,88} = x_1^{3}x_2x_3^{2^{s+1}-2}x_4^{2^{s+t}-2^{s}-4}$\cr  $b_{s,89} = x_1^{3}x_2^{2^{s}-3}x_3^{2}x_4^{2^{s+t}-4}$&\cr  
\end{tabular}}

\medskip
For $s = 2$, $b_{2,90} =  x_1^{3}x_2^{3}x_3^{4}x_4^{2^{t+2}-8}$.

\medskip
For $s = t = 2$, $b_{2,91} =  x_1^{3}x_2^{5}x_3^{8}x_4^{2}.$

\medskip
For $s \geqslant 3$,

\medskip
\centerline{\begin{tabular}{lll}	
$b_{s,90} = x_1x_2^{2}x_3^{2^{s}-4}x_4^{2^{s+t}-1}$&$b_{s,91} = x_1x_2^{2}x_3^{2^{s+t}-1}x_4^{2^{s}-4}$\cr  $b_{s,92} = x_1x_2^{2^{s+t}-1}x_3^{2}x_4^{2^{s}-4}$&$b_{s,93} = x_1^{2^{s+t}-1}x_2x_3^{2}x_4^{2^{s}-4}$\cr  $b_{s,94} = x_1x_2^{2}x_3^{2^{s}-3}x_4^{2^{s+t}-2}$&$b_{s,95} = x_1x_2^{3}x_3^{2^{s}-4}x_4^{2^{s+t}-2}$\cr  $b_{s,96} = x_1x_2^{3}x_3^{2^{s+t}-2}x_4^{2^{s}-4}$&$b_{s,97} = x_1^{3}x_2x_3^{2^{s}-4}x_4^{2^{s+t}-2}$\cr  $b_{s,98} = x_1^{3}x_2x_3^{2^{s+t}-2}x_4^{2^{s}-4}$&$b_{s,99} = x_1x_2^{3}x_3^{2^{s}-2}x_4^{2^{s+t}-4}$\cr  $b_{s,100} = x_1^{3}x_2x_3^{2^{s}-2}x_4^{2^{s+t}-4}$&$b_{s,101} = x_1^{3}x_2^{2^{s+t}-3}x_3^{2}x_4^{2^{s}-4}$\cr  $b_{s,102} = x_1^{3}x_2^{2^{s+1}-3}x_3^{2}x_4^{2^{s+t}-2^{s}-4}$&$b_{s,103} = x_1^{3}x_2^{5}x_3^{2^{s+1}-6}x_4^{2^{s+t}-2^{s}-4}$\cr  $b_{s,104} = x_1^{3}x_2^{5}x_3^{2^{s+t}-6}x_4^{2^{s}-4}$&\cr   
\end{tabular}}

\medskip
For $s = 3$,\ $b_{s,105} = x_1^{3}x_2^{5}x_3^{6}x_4^{2^{t+3}-8}$.

\medskip
For $s \geqslant 4$,\ $b_{s,105} = x_1^{3}x_2^{5}x_3^{2^{s}-6}x_4^{2^{s+t}-4}$.

\subsubsection{\textbf{Computation of $\mbox{\rm Ker}\big((\widetilde{Sq}^0_*)_{(4,n_{s,t})}\big)^{GL_4}$}}\

\medskip
We prove the following.

\begin{thms}\label{dln2} For any $s \geqslant 2$ and $t \geqslant 2$,
$$\mbox{\rm Ker}\big((\widetilde{Sq}^0_*)_{(4,n_{s,t})}\big)^{GL_4} = QP_4((2)|^s|(1)|^t)^{GL_4} = \langle [\zeta_{t,s}] \rangle,$$
where 
\begin{align*}
\zeta_{t,s} = \begin{cases}
\sum_{j\in \{67,\, 68,\, 71,\, 76,\, 77,\, 83,\, 84,\, 89,\, 91\}}b_{2,j}, &\mbox{if } s = t = 2,\\
0, &\mbox{if } s = 2,\, t > 2,\\
\sum_{1 \leqslant j \leqslant 105,\, j \notin \{67,\, 89,\, 99,\, 100\}}b_{3,j}, &\mbox{if } s = 3,\\
\sum_{1 \leqslant j \leqslant 105}b_{s,j}, &\mbox{if } s \geqslant 4.
\end{cases}
\end{align*}
Here, we denote $b_{s,j} = b_{t,s,j}$.

\end{thms}
\begin{rems} The $GL_4$-invariant $[\xi_{2,2}]$ was first found in our work \cite[Remark 7.4.12]{su50}. The dimensional result of the theorem is stated in \cite[Lemmas 3.1.6]{p24} for $t=2$ and \cite[Lemma 4.2.3]{pp25} for $t > 4$. The detailed proof is only provided for $n_{2,2} = 18$ but it is false. No detailed proof is given for other cases.
	
In \cite[Lemma 3.1.9]{p24}, the author stated that $\dim\mbox{\rm Ker}\big((\widetilde{Sq}^0_*)_{(4,n_{s,4})}\big)^{GL_4} =0$ for $s \geqslant 4$, but this result is false because $\dim\mbox{\rm Ker}\big((\widetilde{Sq}^0_*)_{(4,n_{s,4})}\big)^{GL_4} = 1$ for $s \geqslant 4$.
\end{rems}

We set  $q_{s,u} = q_{t,s,u}$ by
\begin{equation}\label{ctn211}
q_{s,1} = \mbox{$\sum_{1 \leqslant j \leqslant 12}$}b_{s,j}, \ q_{s,2} = \mbox{$\sum_{13 \leqslant j \leqslant 18}$}b_{s,j}, \ q_{s,3} = \mbox{$\sum_{19\leqslant j \leqslant 30}$}b_{s,j}.	 
\end{equation}

For $s = 2$, we set  $q_{2,u} = q_{t,2,u}$ by
\begin{equation}\label{ctn212}
\begin{cases}
q_{2,4} = \sum_{j \in \mathbb K_{2,4} = \{31,\, 32,\, 34,\, 35,\, 36,\, 37,\, 38,\, 39,\, 59,\, 60,\, 61,\, 62\}}b_{s,j},\\
q_{2,5} = \sum_{j \in \mathbb K_{2,5} = \{32,\, 33,\, 35,\, 38,\, 39,\, 40,\, 41,\, 42,\, 55,\, 56,\, 57,\, 58,\, 63,\, 64,\, 65,\, 66\}}b_{s,j},\\ 
q_{2,6} = \sum_{j \in \mathbb K_{2,6} = \{31,\, 32,\, 33,\, 34,\, 35,\, 36,\, 37,\, 38,\, 39,\, 40,\, 41,\, 42,\, 43,\, 44,\, 45,\, 46,\, 47,\, 48\}\atop \hskip5cm\cup\{49,\, 50,\, 51,\, 52,\, 53,\, 54\}}b_{s,j},\\
q_{2,7} = \sum_{j \in \mathbb K_{2,7} = \{71,\, 74,\, 75,\, 90\}}b_{s,j},\\
q_{2,8} = \sum_{j \in \mathbb K_{2,9} = \{74,\, 75,\, 76,\, 77,\, 78,\, 80,\, 81,\, 82\}}b_{s,j}.	 
\end{cases}
\end{equation}

By a direct computation, we easily verify that
\begin{align}
\begin{cases}\label{ctn21}
[\Sigma_4(b_{s,1})] = \langle \{[b_{s,j}]: 1 \leqslant j \leqslant 12 \}\rangle,\ [\Sigma_4(b_{s,1})]^{\Sigma_4} = \langle [q_{s,1}]\rangle,\\
[\Sigma_4(b_{s,13})] = \langle \{[b_{s,j}]: 13 \leqslant j \leqslant 18\} \rangle,\ [\Sigma_4(b_{s,13})]^{\Sigma_4} = \langle [q_{s,2}]\rangle,\\
[\Sigma_4(\widetilde a_{s,19})] = \langle \{[b_{s,j} : 19 \leqslant j \leqslant 30 \}\rangle, \ [\Sigma_4(\widetilde a_{s,19})]^{\Sigma_4} = \langle [q_{s,3}]\rangle,\\
\end{cases}
\end{align}

We prepare some lemmas for the proof of Theorem \ref{dln2}.

\begin{lems}\label{bdn21} For $t \geqslant 2$, we have
$QP_4((2)|^2|(1)|^t)^{\Sigma_4} = \langle \{[q_{2,u}] : 1 \leqslant u \leqslant 9\}\rangle,$
where $q_{2,u} = q_{t,2,u}$ for $1 \leqslant u \leqslant 8$, are defined by \eqref{ctn211}, \eqref{ctn212} and
\begin{align*} 
q_{2,9} =\begin{cases}
\sum_{j \in \mathbb K_{2,9} := \mathbb K_{2,2,9} = \{67,\, 68,\, 71,\, 76,\, 77,\, 83,\, 84,\, 89,\, 91\}}b_{2,j} &\mbox{if } s = 2,\\
\sum_{j \in \mathbb K_{2,9} := \mathbb K_{t,2,9} =\{69,\, 70,\, 73,\, 76,\, 79,\, 83,\, 84,\, 85,\, 86,\, 87,\, 88,\, 89\}}b_{2,j} &\mbox{if } s > 2.
\end{cases}
\end{align*}
Consequently, $\dim QP_4((2)|^2|(1)|^t)^{\Sigma_4} = 9$.
\end{lems}
\begin{proof} 
For $s = 2$, we have 
$$QP_4((2)|^2|(1)|^t) = \begin{cases}
\langle \{[b_{2,j}] = b_{2,2,j}: 1 \leqslant j \leqslant 91\}\rangle, &\mbox{if } t = 2\\
\langle \{[b_{2,j} = b_{t,2,j}]: 1 \leqslant j \leqslant 90\}\rangle, &\mbox{if } t > 2.
\end{cases}$$ 
We have a direct summand decomposition of $\Sigma_4$-modules:
\begin{align*}QP_4((2)|^2|(1)|^t) = [\Sigma_4(b_{2,1})]&\bigoplus [\Sigma_4(b_{2,13})] \bigoplus [\Sigma_4(b_{2,19})]\bigoplus \mathcal U_{t,2}\bigoplus \mathcal V_{t,2},
\end{align*}
where $\mathcal U_{t,2} = \langle \{[b_{2,j}]: 31\leqslant j \leqslant 66\}\rangle$ and 
$$\mathcal V_{t,2} = QP_4^+((2)|^2|(1)|^t) = \begin{cases}
\langle \{[b_{2,j} = b_{2,2,j}]: 67\leqslant j \leqslant 91\}\rangle, &\mbox{if } t = 2\\
\langle \{[b_{2,j}= b_{t,2,j}]: 67\leqslant j \leqslant 90\}\rangle, &\mbox{if } t > 2.\end{cases}$$ 
By using \eqref{ctn21} we need only to compute $\mathcal U_{t,2}^{\Sigma_4}$ and $\mathcal V_{t,2}^{\Sigma_4}$. We prove that 
\begin{align*}
\mathcal U_{t,2}^{\Sigma_4} = \langle \{ [q_{2,u}] : 4 \leqslant u \leqslant 6\}\rangle,\ \
\mathcal V_{t,2}^{\Sigma_4}	= \langle \{[q_{2,u}] : 7 \leqslant u \leqslant 9\}\rangle.
\end{align*}
This lemma for $t = 2$ is also proved in \cite{p24} but the calculation for $\mathcal U_{2,2}^{\Sigma_4}$ is false. 

By a direct computation we can verify that $[q_{2,u}] \in \mathcal U_{t,2}^{\Sigma_4}$ for $4 \leqslant u \leqslant 6$. The leading monomials of
$q_{2,4}$, $q_{2,5}$, $q_{2,6}$ respectively are $b_{2,62}$, $b_{2,58}$, $b_{2,54}$. Hence, if $[f] \in \mathcal U_{t,2}^{\Sigma_4}$ with $f \in P_4((2)|^2|(1)|^t)$, then there are $\bar \gamma_u \in \mathbb F_2,\, 4 \leqslant u \leqslant 6$ such that 
$$h = f +\sum_{4 \leqslant u \leqslant 6}\bar\gamma_uq_{2,u} \equiv \sum_{31 \leqslant j \leqslant 66,\, j \ne 54, 58, 62}\gamma_jb_{2,j},$$
with suitable $\gamma_j \in \mathbb F_2$. Since $[f],\, [q_{2,u}] \in \mathcal U_2^{\Sigma_4}$, we have $[h] \in \mathcal U_2^{\Sigma_4}$. 
	
Consider the homomorphisms $\rho_j: P_4 \to P_4$ as defined in Section \ref{s2} with $k =4$ and $1 \leqslant j \leqslant 4$. A direct computation shows
\begin{align*}
\rho_1(h)+h &\equiv \gamma_{\{31,34\}}b_{2,31} + \gamma_{\{32,35\}}b_{2,32} + \gamma_{\{33,40\}}b_{2,33} + \gamma_{\{31,34\}}b_{2,34}\\ 
&\quad + \gamma_{\{32,35\}}b_{2,35} + \gamma_{\{36,38,41,45\}}b_{2,36} + \gamma_{\{37,39,42,46\}}b_{2,37}\\ 
&\quad + \gamma_{\{33,40\}}b_{2,40} + \gamma_{\{36,38,41,45\}}b_{2,41} + \gamma_{\{37,39,42,46\}}b_{2,42}\\ 
&\quad + \gamma_{\{43,44\}}b_{2,43} + \gamma_{\{43,44\}}b_{2,44} + \gamma_{\{47,49\}}b_{2,47} + \gamma_{\{48,52\}}b_{2,48}\\ 
&\quad + \gamma_{\{47,49\}}b_{2,49} + \gamma_{\{50,53\}}b_{2,50} + \gamma_{51}b_{2,51} + \gamma_{\{48,52\}}b_{2,52}\\ 
&\quad + \gamma_{\{50,53\}}b_{2,53} + \gamma_{51}b_{2,54} + \gamma_{\{55,56\}}b_{2,55} + \gamma_{\{55,56\}}b_{2,56}\\ 
&\quad + \gamma_{\{59,60\}}b_{2,59} + \gamma_{\{59,60\}}b_{2,60} + \gamma_{\{63,64\}}b_{2,63} + \gamma_{\{63,64\}}b_{2,64}\\ 
&\quad + \gamma_{\{38,45,57,61\}}b_{2,65} + \gamma_{\{39,46\}}b_{2,66} \equiv 0,\\
\rho_2(h)+h &\equiv \gamma_{\{31,32,33,43\}}b_{2,31} + \gamma_{\{31,32,33,43\}}b_{2,33} + \gamma_{\{34,36\}}b_{2,34}\\ 
&\quad + \gamma_{\{35,38\}}b_{2,35} + \gamma_{\{34,36\}}b_{2,36} + \gamma_{\{37,39,66\}}b_{2,37} + \gamma_{\{35,38\}}b_{2,38}\\ 
&\quad + \gamma_{\{37,39,66\}}b_{2,39} + \gamma_{\{40,41\}}b_{2,40} + \gamma_{\{40,41\}}b_{2,41} + \gamma_{42}b_{2,42}\\ 
&\quad + \gamma_{\{44,45\}}b_{2,44} + \gamma_{\{44,45\}}b_{2,45} + \gamma_{\{46,51\}}b_{2,46} + \gamma_{\{47,48\}}b_{2,47}\\ 
&\quad + \gamma_{\{47,48\}}b_{2,48} + \gamma_{\{49,50\}}b_{2,49} + \gamma_{\{49,50\}}b_{2,50} + \gamma_{\{46,51\}}b_{2,51}\\ 
&\quad + \gamma_{\{52,53\}}b_{2,52} + \gamma_{\{52,53\}}b_{2,53} + \gamma_{\{56,57\}}b_{2,56} + \gamma_{\{56,57\}}b_{2,57}\\ 
&\quad + \gamma_{42}b_{2,58} + \gamma_{\{60,61\}}b_{2,60} + \gamma_{\{60,61\}}b_{2,61} + \gamma_{\{32,43,55,59\}}b_{2,63}\\ 
&\quad + \gamma_{\{64,65\}}b_{2,64} + \gamma_{\{64,65\}}b_{2,65}  \equiv 0,\\
\rho_3(h)+h &\equiv \gamma_{\{31,32,63\}}b_{2,31} + \gamma_{\{31,32,63\}}b_{2,32} + \gamma_{\{33,48,55\}}b_{2,33} + \gamma_{\{34,35,64\}}b_{2,34}\\ 
&\quad + \gamma_{\{34,35,64\}}b_{2,35} + \gamma_{\{36,37\}}b_{2,36} + \gamma_{\{36,37\}}b_{2,37} + \gamma_{\{38,39\}}b_{2,38}\\ 
&\quad + \gamma_{\{38,39\}}b_{2,39} + \gamma_{\{40,52,56\}}b_{2,40} + \gamma_{\{41,42\}}b_{2,41} + \gamma_{\{41,42\}}b_{2,42}\\ 
&\quad + \gamma_{\{43,47\}}b_{2,43} + \gamma_{\{44,49\}}b_{2,44} + \gamma_{\{45,46\}}b_{2,45} + \gamma_{\{45,46\}}b_{2,46}\\ 
&\quad + \gamma_{\{43,47\}}b_{2,47} + \gamma_{\{44,49\}}b_{2,49} + \gamma_{\{50,51\}}b_{2,50} + \gamma_{\{50,51\}}b_{2,51}\\ 
&\quad + \gamma_{53}b_{2,53} + \gamma_{53}b_{2,54} + \gamma_{\{33,48,55\}}b_{2,55} + \gamma_{\{40,52,56\}}b_{2,56} + \gamma_{57}b_{2,57}\\ 
&\quad + \gamma_{57}b_{2,58} + \gamma_{61}b_{2,61} + \gamma_{61}b_{2,62} + \gamma_{\{65,66\}}b_{2,65} + \gamma_{\{65,66\}}b_{2,66} \equiv 0.
\end{align*}
By computing from the above equalities we get $\gamma_j =0$ for all $j$. This implies $h \equiv 0$,  
$ f \equiv \sum_{4\leqslant u \leqslant 6}\bar \gamma_tq_{2,u}$ and $\mathcal U_{t,2}^{\Sigma_4} = \langle \{ [q_{2,u}] : 4 \leqslant u \leqslant 6\}\rangle$.

By a direct computation we can see that $[q_{2,u}] \in \mathcal V_{t,2}^{\Sigma_4}$ for $7 \leqslant u \leqslant 9$. 

For $t = 2$, the leading monomials of
$q_{2,7}$, $q_{2,8}$, $q_{2,9}$ respectively are $b_{2,90}$, $b_{2,82}$, $b_{2,91}$. Hence, if $[g] \in \mathcal V_{2,2}^{\Sigma_4}$ with $g \in P_4((2)|^2|(1)|^2)$, then there are $\tilde \gamma_u \in \mathbb F_2,\, 7 \leqslant u \leqslant 9$, such that 
$$\tilde h = g +\sum_{7 \leqslant u \leqslant 6}\tilde\gamma_tq_{2,u} \equiv \sum_{67 \leqslant j \leqslant 91,\, j \ne 82, 90, 91}\gamma_jb_{2,j},$$
with suitable $\gamma_j \in \mathbb F_2$. Since $[g],\, [q_{2,u}] \in \mathcal V_{2,2}^{\Sigma_4}$, $7 \leqslant u \leqslant 9$, we have $[\tilde h] \in \mathcal V_{2,2}^{\Sigma_4}$. 
By a direct computation we get
\begin{align*}
\rho_1(\tilde h)+\tilde h &\equiv \gamma_{\{69,72\}}b_{2,67} + \gamma_{73}b_{2,68} + \gamma_{79}b_{2,71} + \gamma_{\{72,76,77\}}b_{2,76}\\ 
&\quad + \gamma_{\{72,76,77\}}b_{2,77} + \gamma_{81}b_{2,81} + \gamma_{81}b_{2,82} + \gamma_{\{70,83,84\}}b_{2,83}\\ 
&\quad + \gamma_{\{70,83,84\}}b_{2,84} + \gamma_{\{85,86\}}b_{2,85} + \gamma_{\{85,86\}}b_{2,86} + \gamma_{\{87,88\}}b_{2,87}\\ 
&\quad + \gamma_{\{87,88\}}b_{2,88} + \gamma_{70}b_{2,89} + \gamma_{\{70,72,89\}}b_{2,90} + \gamma_{70}b_{2,91}  \equiv 0,\\
\rho_2(\tilde h)+\tilde h &\equiv \gamma_{\{67,69,75,78,86,87\}}b_{2,67} + \gamma_{\{68,70\}}b_{2,68} + \gamma_{\{67,69,75,78,87\}}b_{2,69}\\ 
&\quad + \gamma_{\{68,70,84\}}b_{2,70} + \gamma_{\{71,72,87\}}b_{2,71} + \gamma_{\{71,72,86,87\}}b_{2,72}\\ 
&\quad + \gamma_{\{73,83\}}b_{2,73} + \gamma_{\{74,76,83,85,87\}}b_{2,74} + \gamma_{\{73,74,76,79,86,87\}}b_{2,76}\\ 
&\quad + \gamma_{86}b_{2,77} + \gamma_{\{79,85\}}b_{2,79} + \gamma_{\{80,81\}}b_{2,80} + \gamma_{\{80,81\}}b_{2,81}\\ 
&\quad + \gamma_{\{73,83,84\}}b_{2,83} + \gamma_{84}b_{2,84} + \gamma_{\{79,85\}}b_{2,85} + \gamma_{\{88,89\}}b_{2,88}\\ 
&\quad + \gamma_{\{84,88,89\}}b_{2,89} + \gamma_{\{84,86\}}b_{2,90} + \gamma_{84}b_{2,91} \equiv 0,\\
\rho_3(\tilde h)+\tilde h &\equiv \gamma_{\{67,68,72,73\}}b_{2,67} + \gamma_{\{67,68,69\}}b_{2,68} + \gamma_{\{69,72,73\}}b_{2,69}\\ 
&\quad + \gamma_{\{72,79\}}b_{2,71} + \gamma_{\{69,72,73\}}b_{2,73} + \gamma_{\{74,75\}}b_{2,74} + \gamma_{\{74,75\}}b_{2,75}\\ 
&\quad + \gamma_{\{70,72,76,81,83\}}b_{2,76} + \gamma_{\{77,84\}}b_{2,77} + \gamma_{\{78,80\}}b_{2,78}\\ 
&\quad + \gamma_{\{78,80\}}b_{2,80} + \gamma_{\{70,72,76,81,83\}}b_{2,83} + \gamma_{\{77,84\}}b_{2,84}\\ 
&\quad + \gamma_{\{85,87\}}b_{2,85} + \gamma_{\{86,88\}}b_{2,86} + \gamma_{\{85,87\}}b_{2,87} + \gamma_{\{86,88\}}b_{2,88}\\ 
&\quad + \gamma_{89}b_{2,89} + \gamma_{89}b_{2,91} \equiv 0.
\end{align*}
By computing from the above equalities we get $\gamma_j =0$ for all $j$. This implies $\tilde h \equiv 0$,  
$ g \equiv \sum_{7\leqslant u \leqslant 9}\tilde\gamma_uq_{2,u}$ and $\mathcal V_{2,2}^{\Sigma_4} = \langle \{ [q_{2,u}] : 7 \leqslant u \leqslant 9\}\rangle$.

For $t > 2$, the leading monomials of
$q_{2,7}$, $q_{2,8}$, $q_{2,9}$ respectively are $b_{2,90}$, $b_{2,82}$, $b_{2,89}$. Hence, if $[g] \in \mathcal V_{t,2}^{\Sigma_4}$ with $g \in P_4((2)|^2|(1)|^t)$, then there are $\tilde \gamma_u \in \mathbb F_2,\, 7 \leqslant u \leqslant 9$, such that 
$$\tilde h = g +\sum_{7 \leqslant u \leqslant 9}\tilde\gamma_uq_{2,u} \equiv \sum_{67 \leqslant j \leqslant 90,\, j \ne 82, 89,\, 90}\gamma_jb_{2,j},$$
with suitable $\gamma_j \in \mathbb F_2$. Since $[g],\, [q_{2,u}] \in \mathcal V_{t,2}^{\Sigma_4}$, $7 \leqslant u \leqslant 9$, we have $[\tilde h] \in \mathcal V_{t,2}^{\Sigma_4}$. 
By a direct computation we get
\begin{align*}
\rho_1(\tilde h)+\tilde h &\equiv \gamma_{\{69,70,72\}}b_{2,67} + \gamma_{\{70,73\}}b_{2,68} + \gamma_{\{70,84\}}b_{2,71} + \gamma_{\{70,72,76,83\}}b_{2,76}\\ 
& + \gamma_{81}b_{2,83} + \gamma_{\{70,72,76,83\}}b_{2,77} + \gamma_{\{79,80\}}b_{2,81} + \gamma_{\{79,80\}}b_{2,82}\\ 
& + \gamma_{81}b_{2,84} + \gamma_{\{85,86\}}b_{2,85} + \gamma_{\{85,86\}}b_{2,86} + \gamma_{\{87,88\}}b_{2,87} + \gamma_{\{87,88\}}b_{2,88}\\ 
& + \gamma_{\{70,72\}}b_{2,90} \equiv 0,\\
\rho_2(\tilde h)+\tilde h &\equiv \gamma_{\{67,69,75,78,84,86,87\}}b_{2,67} + \gamma_{\{68,70,84\}}b_{2,68} + \gamma_{\{67,69,75,78,87\}}b_{2,69}\\ 
& + \gamma_{\{68,70,84\}}b_{2,70} + \gamma_{\{71,72,84,87\}}b_{2,71} + \gamma_{\{71,72,86,87\}}b_{2,72}\\ 
& + \gamma_{\{73,83\}}b_{2,73} + \gamma_{\{74,76,83,85,87\}}b_{2,74} + \gamma_{\{73,74,76,79,84,86,87\}}b_{2,76}\\ 
& + \gamma_{\{84,86\}}b_{2,77} + \gamma_{\{79,85\}}b_{2,79} + \gamma_{\{80,81\}}b_{2,80} + \gamma_{\{80,81\}}b_{2,81}\\ 
& + \gamma_{\{73,83\}}b_{2,83} + \gamma_{\{79,85\}}b_{2,85} + \gamma_{88}b_{2,88} + \gamma_{88}b_{2,89} + \gamma_{\{84,86\}}b_{2,90} \equiv 0,\\
\rho_3(\tilde h)+\tilde h &\equiv \gamma_{\{67,68,72,73\}}b_{2,67} + \gamma_{\{67,68,69\}}b_{2,68} + \gamma_{\{69,72,73\}}b_{2,69}\\ 
& + \gamma_{\{72,79\}}b_{2,71} + \gamma_{\{69,72,73\}}b_{2,73} + \gamma_{\{74,75\}}b_{2,74} + \gamma_{\{74,75\}}b_{2,75}\\ 
& + \gamma_{\{70,72,76,81,83\}}b_{2,76} + \gamma_{\{77,84\}}b_{2,77} + \gamma_{\{78,80\}}b_{2,78} + \gamma_{\{78,80\}}b_{2,80}\\ 
& + \gamma_{\{70,72,76,81,83\}}b_{2,83} + \gamma_{\{77,84\}}b_{2,84} + \gamma_{\{85,87\}}b_{2,85} + \gamma_{\{86,88\}}b_{2,86}\\ 
& + \gamma_{\{85,87\}}b_{2,87} + \gamma_{\{86,88\}}b_{2,88} \equiv 0.
\end{align*}
By computing from the above equalities we get $\gamma_j =0$ for all $j$. This implies $\tilde h \equiv 0$,  
$ g \equiv \sum_{7\leqslant u \leqslant 9}\tilde\gamma_uq_{2,u}$ and $\mathcal V_{t,2}^{\Sigma_4} = \langle \{ [q_{2,u}] : 7 \leqslant u \leqslant 9\}\rangle$.
The lemma is proved.
\end{proof}

\begin{lems}\label{bdn22} For $s \geqslant 3$, $t \geqslant 2$, we have
$QP_4((2)|^{s}|(1)|^t)^{\Sigma_4} = \langle \{[q_{s,u}] : 1 \leqslant u \leqslant 13\} \rangle$
where $q_{s,u} := q_{t,s,u}$ is defined by \eqref{ctn211} for $1 \leqslant u \leqslant 3$ and 
\begin{align*}
q_{s,4} &= \mbox{$\sum_{j\in \mathbb K_{s,4} =  \{31,\, 32,\, 33,\, 34,\, 35,\, 36,\, 37,\, 38,\, 39,\, 40,\, 41,\, 42,\, 55,\, 56,\, 57,\, 58\}}$}b_{s,j},\\
q_{s,5} &= \mbox{$\sum_{j\in \mathbb K_{s,5} =  \{43,\, 44,\, 45,\, 46,\, 47,\, 48,\, 49,\, 50,\, 51,\, 52,\, 53,\, 54,\, 55,\, 56,\, 57,\, 58\}}$}b_{s,j},\\
q_{s,6} &= \mbox{$\sum_{j\in \mathbb K_{s,6} =  \{63,\, 59,\, 55,\, 64,\, 60,\, 56,\, 65,\, 66,\, 61,\, 62,\, 57,\, 58\}}$}b_{s,j},\\
q_{s,7} &= \mbox{$\sum_{j\in  \mathbb K_{s,7} =  \{90,\, 91,\, 92,\, 93\}}$}b_{s,j},\\
q_{s,8} &= \begin{cases}\sum_{j\in \mathbb K_{3,8} = \{71,\, 72,\, 79,\, 85,\, 86,\, 87,\, 88,\, 89,\, 94,\, 102,\, 103\}}b_{3,j}, & \mbox{if } s = 3,\\ 
\sum_{j\in \mathbb K_{s,8} = \{71,\, 72,\, 79,\, 85,\, 86,\, 87,\, 88,\, 94,\, 102,\, 103,\, 105\}}b_{3,j}, & \mbox{if } s > 3, \end{cases}\\ 
q_{s,9} &= \mbox{$\sum_{j\in \mathbb K_{s,9} = \{69,\, 70,\, 94,\, 95,\, 96,\, 97,\, 98,\, 101\}}$}b_{s,j}, \\ 
q_{s,10} &= \mbox{$\sum_{j\in \mathbb K_{s,10} = \{74,\, 75,\, 76,\, 77,\, 78,\, 80,\, 81,\, 82\}}$}b_{s,j},\\
q_{s,11} &= \begin{cases}\sum_{j\in \mathbb K_{3,11} = \{68,\, 73,\, 83,\, 84,\, 89,\, 94,\, 104,\, 105\}}b_{3,j},  & \mbox{if } s = 3,\\
\sum_{j\in \mathbb K_{s,11} = \{67,\, 68,\, 69,\, 70,\, 71,\, 72,\, 103,\, 104,\, 105\}}b_{s,j}, & \mbox{if } s > 3,\end{cases}\\ 
q_{s,12} &= \begin{cases}\sum_{j\in \mathbb K_{3,12} = \{69,\, 70,\, 71,\, 72,\, 73,\, 83,\, 84,\, 94,\, 99,\, 100,\, 103\}}b_{3,j},& \mbox{if } s = 3,\\ 
\sum_{j\in \mathbb K_{s,12} = \{67,\, 68,\, 73,\, 83,\, 84,\, 89,\, 94,\, 99,\, 100,\, 104\}}b_{s,j}, & \mbox{if } s > 3,\end{cases}\\
q_{s,13} &= \begin{cases}\sum_{j\in \mathbb K_{3,13} = \{67,\, 74,\, 75,\, 76,\, 77,\, 95,\, 97,\, 105\}}b_{3,j},& \mbox{if } s = 3,\\
\sum_{j\in \mathbb K_{s,12} = \{74,\, 75,\, 76,\, 77,\, 89,\, 95,\, 97,\, 99,\, 100,\, 105\}}b_{s,j}, & \mbox{if } s > 3.\end{cases} 
\end{align*}	
Consequently, $\dim QP_4((2)|^{s}|(1)|^t)^{\Sigma_4} = 13$.	
\end{lems}
\begin{proof}
For $s \geqslant 3$, we set
\begin{align*}
\mathcal U_{t,s} &= \langle \{[b_{s,j}] : 31 \leqslant j \leqslant 66\} \subset QP_4^-((2)|^s|(1)|^t),\\
\mathcal V_{t,s} &= \langle \{[b_{s,j}] : 67 \leqslant j \leqslant 89 \mbox{ or } 94 \leqslant j \leqslant 105 \} \subset QP_4^+((2)|^s|(1)|^t).
\end{align*}
Then, $\mathcal U_{t,s}$ and $\mathcal V_{t,s}$ are the $\Sigma_4$-submodules of $QP_4((2)|^s|(1)|^t)$, and we have a direct summand decomposition of $\Sigma_4$-modules:
\begin{align*}
QP_4^+((2)|^s|(1)|^t) \cong	 [\Sigma_4(b_{s,1})]&\bigoplus[\Sigma_4(b_{s,13})]\bigoplus[\Sigma_4(b_{s,19})]\\ & \qquad \bigoplus \mathcal U_{t,s} \bigoplus [\Sigma_4(b_{s,90})]\bigoplus \mathcal V_{t,s}.
\end{align*}
It is easy to see that $[\Sigma_4(b_{s,90})] = \langle \{[b_{s,j}] : 90 \leqslant j \leqslant 93\}\rangle$ and $[\Sigma_4(b_{s,90})]^{\Sigma_4} = \langle [q_{s,7}]\rangle$. We prove that
\begin{align*}
\mathcal U_{t,s}^{\Sigma_4} = \langle \{ [q_{s,u}] : 4 \leqslant u \leqslant 6\}\rangle,\ \
\mathcal V_{t,s}^{\Sigma_4}	= \langle \{[q_{s,u}] : 8 \leqslant u \leqslant 13\}\rangle.
\end{align*} 

By a direct computation we can verify that $[q_{s,u}] \in \mathcal U_{t,s}^{\Sigma_4}$ for $4 \leqslant u \leqslant 6$. The leading monomials of
$q_{s,4}$, $q_{s,5}$, $q_{s,6}$ respectively are $b_{s,42}$, $b_{s,54}$, $b_{s,58}$. Hence, if $[f] \in \mathcal U_{t,s}^{\Sigma_4}$ with $f \in P_4((2)|^s|(1)|^2)$, then there are $\bar \gamma_u \in \mathbb F_2,\, 4 \leqslant u \leqslant 6$ such that 
$$\bar h = f +\sum_{4 \leqslant u \leqslant 6}\bar\gamma_uq_{s,u} \equiv \sum_{31 \leqslant j \leqslant 66,\, j \ne 42, 54, 58}\gamma_jb_{s,j},$$
with suitable $\gamma_j \in \mathbb F_2$. Note that $[\bar h] \in \mathcal U_{t,s}^{\Sigma_4}$.
By a direct computation we get
\begin{align*}
\rho_1(\bar h)+ \bar h &\equiv \gamma_{\{31,34\}}b_{3,31} + \gamma_{\{32,35\}}b_{3,32} + \gamma_{\{33,40\}}b_{3,33} + \gamma_{\{31,34\}}b_{3,34}\\ 
&\quad + \gamma_{\{32,35\}}b_{3,35} + \gamma_{\{36,41\}}b_{3,36} + \gamma_{37}b_{3,37} + \gamma_{\{33,40\}}b_{3,40}\\ 
&\quad + \gamma_{\{36,41\}}b_{3,41} + \gamma_{37}b_{3,42} + \gamma_{\{43,44\}}b_{3,43} + \gamma_{\{43,44\}}b_{3,44}\\ 
&\quad + \gamma_{\{47,49\}}b_{3,47} + \gamma_{\{48,52\}}b_{3,48} + \gamma_{\{47,49\}}b_{3,49} + \gamma_{\{50,53\}}b_{3,50}\\ 
&\quad + \gamma_{51}b_{3,51} + \gamma_{\{48,52\}}b_{3,52} + \gamma_{\{50,53\}}b_{3,53} + \gamma_{51}b_{3,54} + \gamma_{\{55,56\}}b_{3,55}\\ 
&\quad + \gamma_{\{55,56\}}b_{3,56} + \gamma_{\{59,60\}}b_{3,59} + \gamma_{\{59,60\}}b_{3,60} + \gamma_{\{63,64\}}b_{3,63}\\ 
&\quad + \gamma_{\{63,64\}}b_{3,64} + \gamma_{\{38,45,57,61\}}b_{3,65} + \gamma_{\{39,46,62\}}b_{3,66} \equiv 0,\\
\rho_2(\bar h)+ \bar h &\equiv \gamma_{\{31,33\}}b_{3,31} + \gamma_{\{31,33\}}b_{3,33} + \gamma_{\{34,36\}}b_{3,34} + \gamma_{\{35,38\}}b_{3,35}\\ 
&\quad + \gamma_{\{34,36\}}b_{3,36} + \gamma_{\{37,39\}}b_{3,37} + \gamma_{\{35,38\}}b_{3,38} + \gamma_{\{37,39\}}b_{3,39}\\ 
&\quad + \gamma_{\{40,41\}}b_{3,40} + \gamma_{\{40,41\}}b_{3,41} + \gamma_{\{44,45\}}b_{3,44} + \gamma_{\{44,45\}}b_{3,45}\\ 
&\quad + \gamma_{\{46,51\}}b_{3,46} + \gamma_{\{47,48\}}b_{3,47} + \gamma_{\{47,48\}}b_{3,48} + \gamma_{\{49,50\}}b_{3,49}\\ 
&\quad + \gamma_{\{49,50\}}b_{3,50} + \gamma_{\{46,51\}}b_{3,51} + \gamma_{\{52,53\}}b_{3,52} + \gamma_{\{52,53\}}b_{3,53}\\ 
&\quad + \gamma_{\{56,57\}}b_{3,56} + \gamma_{\{56,57\}}b_{3,57} + \gamma_{66}b_{3,58} + \gamma_{\{60,61\}}b_{3,60}\\ 
&\quad + \gamma_{\{60,61\}}b_{3,61} + \gamma_{\{32,43,55,59\}}b_{3,63} + \gamma_{\{64,65\}}b_{3,64}\\ 
&\quad + \gamma_{\{64,65\}}b_{3,65} + \gamma_{66}b_{3,66} \equiv 0,\\
\rho_3(\bar h)+ \bar h &\equiv \gamma_{\{31,32\}}b_{3,31} + \gamma_{\{31,32\}}b_{3,32} + \gamma_{\{34,35\}}b_{3,34} + \gamma_{\{34,35\}}b_{3,35}\\ 
&\quad + \gamma_{\{36,37\}}b_{3,36} + \gamma_{\{36,37\}}b_{3,37} + \gamma_{\{38,39\}}b_{3,38} + \gamma_{\{38,39\}}b_{3,39}\\ 
&\quad + \gamma_{41}b_{3,41} + \gamma_{41}b_{3,42} + \gamma_{\{43,47\}}b_{3,43} + \gamma_{\{44,49\}}b_{3,44} + \gamma_{\{45,46\}}b_{3,45}\\ 
&\quad + \gamma_{\{45,46\}}b_{3,46} + \gamma_{\{43,47\}}b_{3,47} + \gamma_{\{44,49\}}b_{3,49} + \gamma_{\{50,51\}}b_{3,50}\\ 
&\quad + \gamma_{\{50,51\}}b_{3,51} + \gamma_{53}b_{3,53} + \gamma_{53}b_{3,54} + \gamma_{\{33,48,55,63\}}b_{3,55}\\ 
&\quad + \gamma_{\{40,52,56,64\}}b_{3,56} + \gamma_{57}b_{3,57} + \gamma_{57}b_{3,58} + \gamma_{\{61,62\}}b_{3,61}\\ 
&\quad + \gamma_{\{61,62\}}b_{3,62} + \gamma_{\{33,48,55,63\}}b_{3,63} + \gamma_{\{40,52,56,64\}}b_{3,64}\\ 
&\quad + \gamma_{\{65,66\}}b_{3,65} + \gamma_{\{65,66\}}b_{3,66} \equiv 0.
\end{align*}
By computing from the above equalities we get $\gamma_j =0$ for all $j$. This implies $\bar h \equiv 0$,  
$f \equiv \sum_{4\leqslant u \leqslant 6}\bar \gamma_uq_{s,u}$ and $\mathcal U_{t,s}^{\Sigma_4} = \langle \{ [q_{s,u}] : 4 \leqslant u \leqslant 6\}\rangle$.

Now we compute $\mathcal V_{t,s}^{\Sigma_4}$ for $s > 3$. The case $s = 3$ is carried out by a same argument. We can verify that $[q_{s,u}] \in \mathcal V_{t,s}^{\Sigma_4}$ for $8 \leqslant u \leqslant 13$. The leading monomials of $q_{s,u}$, $8\leqslant u \leqslant 13$, respectively are $b_{s,102}$, $b_{s,101}$, $b_{s,82}$, $b_{s,104}$, $b_{s,89}$, $b_{s,77}$. Hence, if $[f] \in \mathcal V_{t,s}^{\Sigma_4}$ with $f \in P_4((2)|^s|(1)|^t)$, then there are $\gamma_u^* \in \mathbb F_2,\, 8 \leqslant u \leqslant 13$, such that 
$$h^* = f +\sum_{8 \leqslant u \leqslant 13}\gamma_u^*q_{s,u} \equiv \sum_{67 \leqslant j \leqslant 105,\, j \notin \mathbb J}\gamma_jb_{2,j},$$
with $\mathbb J = \{77,\, 82,\, 89,\, 90,\, 91,\, 92,\, 93,\, 101,\, 102,\, 104 \}$ and $\gamma_j \in \mathbb F_2$. Since $[f],\, [q_{s,u}] \in \mathcal V_{t,s}^{\Sigma_4}$, we have $[h^*] \in \mathcal V_{t,s}^{\Sigma_4}$.
By a direct computation we obtain
\begin{align*}
\rho_1(h^*)+ h^* &\equiv \gamma_{\{69,70,72,94,105\}}b_{s,67} + \gamma_{\{70,73\}}b_{s,68} + \gamma_{\{72,79,103\}}b_{s,71} + \gamma_{76}b_{s,76}\\ 
&\quad + \gamma_{76}b_{s,77} + \gamma_{81}b_{s,81} + \gamma_{81}b_{s,82} + \gamma_{\{83,84\}}b_{s,83} + \gamma_{\{83,84\}}b_{s,84}\\ 
&\quad + \gamma_{\{85,86\}}b_{s,85} + \gamma_{\{85,86\}}b_{s,86} + \gamma_{\{87,88\}}b_{s,87} + \gamma_{\{87,88\}}b_{s,88}\\ 
&\quad + \gamma_{\{70,72\}}b_{s,89} + \gamma_{\{95,97\}}b_{s,95} + \gamma_{\{96,98\}}b_{s,96} + \gamma_{\{95,97\}}b_{s,97}\\ 
&\quad + \gamma_{\{96,98\}}b_{s,98} + \gamma_{\{70,72,99,100\}}b_{s,99} + \gamma_{\{70,72,99,100\}}b_{s,100}\\ 
&\quad + \gamma_{\{70,72\}}b_{s,105} \equiv 0,\\
\rho_2(h^*)+ h^* &\equiv \gamma_{\{67,69,73,79,83,84,85,86,97,99,103\}}b_{s,67} + \gamma_{\{68,70,84,96\}}b_{s,68}\\ 
&\quad + \gamma_{\{67,69,97,99\}}b_{s,69} + \gamma_{\{68,70,84,96\}}b_{s,70} + \gamma_{\{71,72,86,87\}}b_{s,71}\\ 
&\quad + \gamma_{\{71,72,86,87\}}b_{s,72} + \gamma_{\{73,83\}}b_{s,73} + \gamma_{\{74,76\}}b_{s,74} + \gamma_{\{74,76\}}b_{s,76}\\ 
&\quad + \gamma_{\{79,85\}}b_{s,79} + \gamma_{\{80,81\}}b_{s,80} + \gamma_{\{80,81\}}b_{s,81} + \gamma_{\{73,83\}}b_{s,83}\\ 
&\quad + \gamma_{\{79,85\}}b_{s,85} + \gamma_{88}b_{s,88} + \gamma_{100}b_{s,89} + \gamma_{\{75,78,83,85,94,95\}}b_{s,94}\\ 
&\quad + \gamma_{\{73,75,78,79,94,95\}}b_{s,95} + \gamma_{98}b_{s,98} + \gamma_{\{73,79,83,84,85,86,103\}}b_{s,99}\\ 
&\quad + \gamma_{100}b_{s,100} + \gamma_{98}b_{s,101} + \gamma_{88}b_{s,102} + \gamma_{\{84,86,103\}}b_{s,105} \equiv 0,\\
\rho_3(h^*)+ h^* &\equiv \gamma_{\{67,68,72,73,76,81\}}b_{s,67} + \gamma_{\{67,68,69,76,81,94,105\}}b_{s,68}\\ 
&\quad + \gamma_{\{72,79,103\}}b_{s,71} + \gamma_{\{69,72,73,94\}}b_{s,73} + \gamma_{\{74,75\}}b_{s,74}\\ 
&\quad + \gamma_{\{74,75\}}b_{s,75} + \gamma_{\{78,80\}}b_{s,78} + \gamma_{\{78,80\}}b_{s,80}\\ 
&\quad + \gamma_{\{70,72,76,81,83,99\}}b_{s,83} + \gamma_{\{84,100\}}b_{s,84} + \gamma_{\{85,87\}}b_{s,85}\\ 
&\quad + \gamma_{\{86,88\}}b_{s,86} + \gamma_{\{85,87\}}b_{s,87} + \gamma_{\{86,88\}}b_{s,88}\\ 
&\quad + \gamma_{\{69,72,73,94\}}b_{s,94} + \gamma_{\{76,81,95,96\}}b_{s,95} + \gamma_{\{76,81,95,96\}}b_{s,96}\\ 
&\quad + \gamma_{\{97,98\}}b_{s,97} + \gamma_{\{97,98\}}b_{s,98} + \gamma_{\{70,72,76,81,83,99\}}b_{s,99}\\ 
&\quad + \gamma_{\{84,100\}}b_{s,100} + \gamma_{105}b_{s,104} + \gamma_{105}b_{s,105} \equiv 0.
\end{align*}
The above equalities imply $\gamma_j =0$ for all $j$. Hence, $h^* \equiv 0$ and $f \equiv \sum_{8\leqslant u \leqslant 13}\gamma_u^*q_{s,u}$. This implies $\mathcal V_{t,s}^{\Sigma_4} = \langle \{ [q_{s,u}] : 8 \leqslant u \leqslant 13\}\rangle$. The lemma is completely proved.
\end{proof}

\begin{proof}[Proof of Theorem \ref{dln2}]
Suppose $[g] \in QP_4((2)|^{s}|(1)|^t)^{GL_4}$ with $g \in QP_4((2)|^{s}|(1)|^t)$. Since $\Sigma_4 \subset GL_4$, we have $[g] \in QP_4((2)|^{s}|(1)|^t)^{\Sigma_4}$. 
 
For $s = 2$, by Lemma \ref{bdn21}, we have 
$g \equiv \sum_{1 \leqslant u \leqslant 9}\gamma_u q_{2,u},$
where $\gamma_u \in \mathbb F_2$. 

For $s = t = 2$, by computing $\rho_4(g)+g$ in terms of the admissible monomials, we get
\begin{align*}
\rho_4(g)+g &\equiv \gamma_{\{1,3\}}b_{2,3} + \gamma_{\{1,3\}}b_{2,4} + \gamma_{\{1,4,5,6\}}b_{2,5} + \gamma_{\{1,4,5,6\}}b_{2,6} + \gamma_{\{1,2\}}b_{2,9}\\
& + \gamma_{\{2,6\}}b_{2,14} + \gamma_{\{2,6\}}b_{2,15} + \gamma_{3}b_{2,19} + \gamma_{3}b_{2,20} + \gamma_{3}b_{2,21} + \gamma_{\{3,5\}}b_{2,26}\\
& + \gamma_{\{3,5\}}b_{2,27} + \gamma_{\{4,6,7,8\}}b_{2,31} + \gamma_{\{4,5,6,7,8\}}b_{2,32} + \gamma_{\{5,6\}}b_{2,33} + \gamma_{\{5,6\}}b_{2,36}\\
& + \gamma_{\{5,6\}}b_{2,37} + \gamma_{\{6,8\}}b_{2,43} + \gamma_{\{6,8\}}b_{2,47} + \gamma_{6}b_{2,48} + \gamma_{\{4,6\}}b_{2,50}\\
& + \gamma_{\{4,6\}}b_{2,51} + \gamma_{5}b_{2,55} + \gamma_{4}b_{2,59} + \gamma_{5}b_{2,63} + \gamma_{\{1,2,3,6\}}b_{2,65}\\
& + \gamma_{\{1,2,3,6\}}b_{2,66} + \gamma_{5}b_{2,67} + \gamma_{5}b_{2,68} + \gamma_{4}b_{2,71} + \gamma_{\{3,8\}}b_{2,76} + \gamma_{3}b_{2,77}\\
& + \gamma_{8}b_{2,81} + \gamma_{3}b_{2,83} + \gamma_{3}b_{2,84} + \gamma_{3}b_{2,89} + \gamma_{\{3,6,8\}}b_{2,90} + \gamma_{3}b_{2,91} \equiv 0.
\end{align*}
This equality implies $\gamma_u = 0$ for $1 \leqslant u \leqslant 8$. Hence, $g \equiv \gamma_9 q_{2,9} = \gamma_9\zeta_{2,2}$. The theorem is proved for $s = t = 2$.

For $s= 2,\, t > 2$, we have
\begin{align*}
\rho_4(g)+g &\equiv \gamma_{\{1,3\}}b_{2,3} + \gamma_{\{1,3\}}b_{2,4} + \gamma_{\{1,4,5,6\}}b_{2,5} + \gamma_{\{1,4,5,6\}}b_{2,6} + \gamma_{\{1,2\}}b_{2,9}\\
& + \gamma_{\{2,6\}}b_{2,14} + \gamma_{\{2,6\}}b_{2,15} + \gamma_{3}b_{2,19} + \gamma_{3}b_{2,20} + \gamma_{\{3,9\}}b_{2,21}\\
& + \gamma_{\{3,5\}}b_{2,26} + \gamma_{\{3,5\}}b_{2,27} + \gamma_{\{4,6,7,8\}}b_{2,31} + \gamma_{\{4,5,6,7,8\}}b_{2,32}\\
& + \gamma_{\{5,6,9\}}b_{2,33} + \gamma_{\{5,6\}}b_{2,36} + \gamma_{\{5,6\}}b_{2,37} + \gamma_{\{6,8\}}b_{2,43} + \gamma_{\{6,8\}}b_{2,47}\\
& + \gamma_{6}b_{2,48} + \gamma_{\{4,6\}}b_{2,50} + \gamma_{\{4,6\}}b_{2,51} + \gamma_{\{5,9\}}b_{2,55} + \gamma_{5}b_{2,59}\\
& + \gamma_{\{1,2,3,6\}}b_{2,61} + \gamma_{\{1,2,3,6\}}b_{2,62} + \gamma_{\{4,9\}}b_{2,63} + \gamma_{\{3,5\}}b_{2,67}\\
& + \gamma_{\{3,5\}}b_{2,68} + \gamma_{\{3,4\}}b_{2,71} + \gamma_{8}b_{2,76} + \gamma_{9}b_{2,77} + \gamma_{\{8,9\}}b_{2,83}\\
& + \gamma_{9}b_{2,85} + \gamma_{9}b_{2,87} + \gamma_{\{3,6,8\}}b_{2,90} \equiv 0.
\end{align*}
This equality implies $\gamma_u = 0$ for $1 \leqslant u \leqslant 9$. Hence, $g \equiv 0$. The theorem is proved for $s = 2$ and $t > 2$.

For $s = 3$, we denote $b_{3,j} = b_{t,3,j}$ and $q_{3,u} = q_{t,3,u}$. By using Lemma \ref{bdn22}, we have 
$ f \equiv \sum_{1 \leqslant u \leqslant 13}\gamma_u q_{3,u},$
where $\gamma_u \in \mathbb F_2$. By computing $\rho_4(g)+g$ in terms of the admissible monomials we get
\begin{align*}
\rho_4(g)+g &\equiv \gamma_{\{1,3\}}b_{3,3} + \gamma_{\{1,3\}}b_{3,4} + \gamma_{\{1,4\}}b_{3,5} + \gamma_{\{1,4\}}b_{3,6} + \gamma_{\{1,2\}}b_{3,9}\\ 
&\quad + \gamma_{\{2,5\}}b_{3,14} + \gamma_{\{2,5\}}b_{3,15} + \gamma_{\{3,7\}}b_{3,19} + \gamma_{\{3,7\}}b_{3,20} + \gamma_{\{3,9,12\}}b_{3,21}\\ 
&\quad + \gamma_{\{3,4,5,6\}}b_{3,26} + \gamma_{\{3,4,5,6\}}b_{3,27} + \gamma_{\{4,10,13\}}b_{3,31} + \gamma_{\{4,10,13\}}b_{3,32}\\ 
&\quad + \gamma_{\{4,9,12\}}b_{3,33} + \gamma_{\{1,2,3,4,5,6\}}b_{3,36} + \gamma_{\{1,2,3,4,5,6\}}b_{3,37}\\ 
&\quad + \gamma_{\{1,2,3,5\}}b_{3,41} + \gamma_{\{1,2,3,5\}}b_{3,42} + \gamma_{\{5,10\}}b_{3,43} + \gamma_{\{5,10\}}b_{3,47}\\ 
&\quad + \gamma_{\{5,8,12\}}b_{3,48} + \gamma_{\{5,6\}}b_{3,50} + \gamma_{\{5,6\}}b_{3,51} + \gamma_{\{4,5,6,11,12\}}b_{3,55}\\ 
&\quad + \gamma_{\{6,8\}}b_{3,59} + \gamma_{\{6,8,9,11,12\}}b_{3,63} + \gamma_{\{1,2,3,5\}}b_{3,65} + \gamma_{\{1,2,3,5\}}b_{3,66}\\ 
&\quad + \gamma_{\{4,6,8,10,12\}}b_{3,67} + \gamma_{\{3,4,5,6,7,11\}}b_{3,68} + \gamma_{\{5,6,8,10,12\}}b_{3,71}\\ 
&\quad + \gamma_{\{3,5,7,11,12\}}b_{3,76} + \gamma_{\{3,5,7,10\}}b_{3,77} + \gamma_{\{8,10\}}b_{3,81} + \gamma_{12}b_{3,83}\\ 
&\quad + \gamma_{12}b_{3,85} + \gamma_{12}b_{3,87} + \gamma_{\{7,9\}}b_{3,92} + \gamma_{\{3,5,7,8,9,10,12\}}b_{3,95}\\ 
&\quad + \gamma_{\{9,11\}}b_{3,96} + \gamma_{\{3,5,7,10\}}b_{3,97} + \gamma_{\{8,11\}}b_{3,99} + \gamma_{\{3,5,7,10\}}b_{3,105} \equiv 0.
\end{align*}
From this equality we get $\gamma_{12} = \gamma_{13} = 0$ and $\gamma_{u} = \gamma_{1}$ for $1 \leqslant u \leqslant 11$. Hence, we obtain
$$g \equiv \gamma_1\Big(\sum_{1 \leqslant u \leqslant 11} q_{3,u}\Big) = \gamma_1\Big(\sum_{1 \leqslant j \leqslant 105,\, j \ne 67,\, 89,\, 99,\, 100}b_{3,j}\Big) =\gamma_1 \zeta_{t,3}.$$
This completes the proof of the theorem for $s = 3$.

By using Lemma \ref{bdn22} for $s\geqslant 4$, we have 
$g \equiv \sum_{1 \leqslant u \leqslant 13}\gamma_u q_{s,u},$
where $\gamma_u \in \mathbb F_2$. By computing $\rho_4(g)+g$ in terms of the admissible monomials, we obtain
\begin{align*}
\rho_4(g)+g &\equiv \gamma_{\{1,3\}}b_{s,3} + \gamma_{\{1,3\}}b_{s,4} + \gamma_{\{1,4\}}b_{s,5} + \gamma_{\{1,4\}}b_{s,6} + \gamma_{\{1,2\}}b_{s,9}\\ 
&\quad + \gamma_{\{2,5\}}b_{s,14} + \gamma_{\{2,5\}}b_{s,15} + \gamma_{\{3,7\}}b_{s,19} + \gamma_{\{3,7\}}b_{s,20} + \gamma_{\{3,9,11\}}b_{s,21}\\ 
&\quad + \gamma_{\{3,4,5,6\}}b_{s,26} + \gamma_{\{3,4,5,6\}}b_{s,27} + \gamma_{\{4,10,13\}}b_{s,31} + \gamma_{\{4,10,13\}}b_{s,32}\\ 
&\quad + \gamma_{\{4,9,11\}}b_{s,33} + \gamma_{\{1,2,3,4,5,6\}}b_{s,36} + \gamma_{\{1,2,3,4,5,6\}}b_{s,37} + \gamma_{\{1,2,3,5\}}b_{s,41}\\ 
&\quad + \gamma_{\{1,2,3,5\}}b_{s,42} + \gamma_{\{5,10\}}b_{s,43} + \gamma_{\{5,10\}}b_{s,47} + \gamma_{\{5,8,11\}}b_{s,48}\\ 
&\quad + \gamma_{\{5,6\}}b_{s,50} + \gamma_{\{5,6\}}b_{s,51} + \gamma_{\{4,5,6,12\}}b_{s,55} + \gamma_{\{6,8\}}b_{s,59}\\ 
&\quad + \gamma_{\{6,8,9,12\}}b_{s,63} + \gamma_{\{1,2,3,5\}}b_{s,65} + \gamma_{\{1,2,3,5\}}b_{s,66} + \gamma_{\{3,4,5,6,7,8,11\}}b_{s,67}\\ 
&\quad + \gamma_{\{3,4,5,6,7,11,12\}}b_{s,68} + \gamma_{\{5,6,8,10,11\}}b_{s,71} + \gamma_{\{3,5,7,12\}}b_{s,76}\\ 
&\quad + \gamma_{\{3,5,7,10\}}b_{s,77} + \gamma_{\{8,10\}}b_{s,81} + \gamma_{11}b_{s,83} + \gamma_{11}b_{s,85} + \gamma_{11}b_{s,87}\\ 
&\quad + \gamma_{\{3,5,7,10\}}b_{s,89} + \gamma_{\{7,9\}}b_{s,92} + \gamma_{\{3,5,7,8,9,10,11\}}b_{s,95} + \gamma_{\{9,11,12\}}b_{s,96}\\ 
&\quad + \gamma_{\{3,5,7,10\}}b_{s,97} + \gamma_{\{3,5,7,8,10,11,12\}}b_{s,99} + \gamma_{\{3,5,7,10\}}b_{s,100}\\ 
&\quad + \gamma_{\{3,5,7,10\}}b_{s,105} \equiv 0.
\end{align*}
By computing from this equality we get $\gamma_{11} = \gamma_{13} = 0$ and $\gamma_u = \gamma_1$ for $1 \leqslant u \leqslant 12,\, u \ne 11$. Therefore, we have
$$g \equiv \gamma_1\Bigg(\sum_{1 \leqslant t \leqslant 12,\, t\ne 11} q_{s,t}\Bigg) = \gamma_1\Bigg(\sum_{1 \leqslant j \leqslant 105}b_{s,j}\Bigg) =\gamma_1 \xi_{t,s}.$$
Theorem \ref{dln2} is completely proved.	
\end{proof}

\subsubsection{\textbf{The subcase $s = 2,\, t\geqslant 2$}}\

\medskip
We have $n_{2,t} = 2^{t+2}+2$ and $d_{1,t} = 2^{t+1} - 1$. Consider Kameko's homomorphism
$$(\widetilde{Sq}^0_*)_{(4,n_{2,t})} : (QP_{4})_{2^{t+2}+2}\to (QP_4)_{2^{t+1}-1}.$$
The space $(QP_4)_{2^{t+1}-1}^{GL_4}$ had been determined in our work \cite{su2}.

\begin{props}[\cite{su2}]\label{mdt10} For $t \geqslant 2$, we have $(QP_4)_{2^{t+1}-1}^{GL_4} = \langle [\xi_{t,1}]\rangle$, where
$$\xi_{t,1} = \begin{cases}
\eta_2 + x_1x_2^{2}x_3^{2}x_4^{2}, &\mbox{if } t=2, \\
\eta_t + x_1x_2^{2}x_3^{4}x_4^{2^{t+1}-8}, &\mbox{if } t \geqslant 3.
\end{cases}$$
Here, 
\begin{align*}
\eta_t &= x_4^{2^{t+1}-1} + x_3^{2^{t+1}-1} + x_2^{2^{t+1}-1} + x_1^{2^{t+1}-1} + x_3x_4^{2^{t+1}-2} + x_2x_4^{2^{t+1}-2}\\ &\quad + x_2x_3^{2^{t+1}-2} + x_1x_4^{2^{t+1}-2} + x_1x_3^{2^{t+1}-2} + x_1x_2^{2^{t+1}-2} + x_2x_3^{2}x_4^{2^{t+1}-4}\\ &\quad + x_1x_3^{2}x_4^{2^{t+1}-4} + x_1x_2^{2}x_4^{2^{t+1}-4} + x_1x_2^{2}x_3^{2^{t+1}-4}.
\end{align*}
\end{props}

We prove the following.

\begin{thms} For $s = 2,\, t\geqslant 2$, we have $(QP_4)_{2^{t+2}+2} = \langle [\zeta_{t,2}], [\delta_{t,2}]\rangle$, where $\zeta_{t,2}$ is defined in Theorem \ref{dln2} and 
$$\delta_{t,2} = 
\psi_{2,t}(\xi_{t,1})+ \mbox{$\sum_{1 \leqslant j \leqslant 90,\, j \notin\{31,\, 33,\, 34,\, 36,\, 37,\, 40,\, 41,\, 42,\, 69\}\atop\hskip1.1cm\cup \{71,\, 74,\, 83,\, 84,\, 85,\, 86,\, 87,\, 88,\, 89\}}$}b_{2,j}.$$	 
\end{thms}
\begin{proof}
Let $[f] \in (QP_4)_{2^{t+2}+2}^{GL_4}$ with $f \in (P_4)_{2^{t+2}+2}$, then there is $\lambda_0 \in \mathbb F_2$ and $h \in \mbox{Ker}(\widetilde{Sq}^0_*)_{(4,2^{t+2}+2)}$ such that
$$f \equiv \gamma_0\psi_{2,t}(\xi_{t,1}) + h.$$ 

For $t = 2$, we have $n_{2,2} = 18$. By a direct computation we get
\begin{align*}
\rho_1(\psi_{2,2}(\xi_{2,1}))+\psi_{2,2}(\xi_{2,1}) &\equiv x_1x_2x_3^{2}x_4^{14} + x_1x_2x_3^{6}x_4^{10} + x_1x_2x_3^{14}x_4^{2} + x_1x_2^{3}x_3^{2}x_4^{12}\\ 
& \quad +\, x_1x_2^{3}x_3^{12}x_4^{2} + x_1^{3}x_2x_3^{2}x_4^{12} + x_1^{3}x_2x_3^{12}x_4^{2} + x_1^{3}x_2^{5}x_3^{2}x_4^{8}\\
&= \mbox{$\sum_{j\in\{67,\, 68,\, 71,\, 76,\, 77,\, 83,\, 84,\, 89,\, 91\}}b_{2,j}$},\\
\rho_2(\psi_{2,2}(\xi_{2,1}))+\psi_{2,2}(\xi_{2,1}) &\equiv x_1x_2^{2}x_3^{5}x_4^{10} + x_1x_2^{2}x_3^{13}x_4^{2} + x_1x_2^{3}x_3^{4}x_4^{10} + x_1x_2^{3}x_3^{12}x_4^{2}\\
&=b_{2,73} + b_{2,79} + b_{2,83} + b_{2,85},\\
\rho_3(\psi_{2,2}(\xi_{2,1})) +\psi_{2,2}(\xi_{2,1}) &\equiv x_1x_2^{2}x_3^{3}x_4^{12} + 
x_1x_2^{2}x_3^{12}x_4^{3} =b_{2,74} + b_{2,75}. 
\end{align*} 
By computing from these results and the relations $\rho_i(f) \equiv 0$, $i = 1,\, 2,\, 3$, we obtain
$$f \equiv  \gamma_0(\psi_{2,2}(\xi_{2,1}) + q_{2,0}) + \sum_{1 \leqslant u \leqslant 9}\gamma_uq_{2,u}$$
with $\gamma_u \in \mathbb F_2$ and 
\begin{align*}
q_{2,0} &= x_1x_2x_3^{2}x_4^{14} + x_1x_2x_3^{14}x_4^{2} + x_1x_2^{14}x_3x_4^{2} + x_1x_2x_3^{6}x_4^{10}\\ &\quad + x_1x_2^{6}x_3x_4^{10} + x_1x_2^{2}x_3^{13}x_4^{2} + x_1x_2^{2}x_3^{12}x_4^{3} + x_1x_2^{2}x_3^{5}x_4^{10}\\
&=\mbox{$\sum_{j\in\{67,\, 68,\, 70,\, 71,\, 72,\, 73,\, 75,\, 79\}}b_{2,j}$}.
\end{align*}
	
We have
\begin{align*}
\rho_4(\psi_{2,2}(\xi_{2,1}))&+\psi_{2,2}(\xi_{2,1}) \equiv x_2x_3^{2}x_4^{15} + x_2x_3^{15}x_4^{2} + x_2x_3^{3}x_4^{14} + x_1x_2x_3^{2}x_4^{14}\\ &\quad + x_1x_2x_3^{14}x_4^{2} + x_1x_2x_3^{6}x_4^{10} + x_1^{3}x_2x_3^{2}x_4^{12} + x_1x_2^{7}x_3^{2}x_4^{8}\\ &\quad + x_1x_2^{3}x_3^{12}x_4^{2} + x_1^{3}x_2x_3^{12}x_4^{2} + x_1^{3}x_2^{5}x_3^{2}x_4^{8} + x_1^{3}x_2^{3}x_3^{4}x_4^{8}\\ &\quad + x_1^{3}x_2^{5}x_3^{8}x_4^{2}\\
&=\mbox{$\sum_{j\in\{19,\, 20,\, 31,\, 67,\, 68,\, 71,\, 77,\, 81,\, 83,\, 84,\, 89,\, 90,\, 91\}}b_{2,j}$}.
\end{align*}
	
By using this result and computing $\rho_4(f)+f$ in terms of the admissible monomials, we obtain
\begin{align*}
\rho_4(f)+f &\equiv \gamma_{\{1,3\}}b_{2,3} + \gamma_{\{1,3\}}b_{2,4} + \gamma_{\{1,4,5,6\}}b_{2,5} + \gamma_{\{1,4,5,6\}}b_{2,6} + \gamma_{\{1,2\}}b_{2,9}\\
& + \gamma_{\{2,6\}}b_{2,14} + \gamma_{\{2,6\}}b_{2,15} + \gamma_{\{0,3\}}b_{2,19} + \gamma_{\{0,3\}}b_{2,20} + \gamma_{\{0,3\}}b_{2,21}\\
& + \gamma_{\{3,5\}}b_{2,26} + \gamma_{\{3,5\}}b_{2,27} + \gamma_{\{4,6,7,8\}}b_{2,31} + \gamma_{\{0,4,5,6,7,8\}}b_{2,32}\\
& + \gamma_{\{5,6\}}b_{2,33} + \gamma_{\{5,6\}}b_{2,36} + \gamma_{\{5,6\}}b_{2,37} + \gamma_{\{6,8\}}b_{2,43} + \gamma_{\{6,8\}}b_{2,47}\\
& + \gamma_{\{0,6\}}b_{2,48} + \gamma_{\{4,6\}}b_{2,50} + \gamma_{\{4,6\}}b_{2,51} + \gamma_{\{0,5\}}b_{2,55} + \gamma_{\{0,4\}}b_{2,59}\\
& + \gamma_{\{0,5\}}b_{2,63} + \gamma_{\{1,2,3,6\}}b_{2,65} + \gamma_{\{1,2,3,6\}}b_{2,66} + \gamma_{\{0,5\}}b_{2,67} + \gamma_{\{0,5\}}b_{2,68}\\
& + \gamma_{\{0,4\}}b_{2,71} + \gamma_{\{3,8\}}b_{2,76} + \gamma_{\{0,3\}}b_{2,77} + \gamma_{\{0,8\}}b_{2,81} + \gamma_{\{0,3\}}b_{2,83}\\
& + \gamma_{\{0,3\}}b_{2,84} + \gamma_{\{0,3\}}b_{2,89} + \gamma_{\{0,3,6,8\}}b_{2,90} + \gamma_{\{0,3\}}b_{2,91} \equiv 0.
\end{align*}
The above equality implies $\gamma_t =\gamma_0$ for all $0\leqslant t \leqslant 8$. Hence,  
\begin{align*}
f &\equiv \gamma_0\Big(\psi_{2,2}(\xi_{2,1})+ \sum_{0\leqslant u \leqslant 8}q_{2,u}\Big) +\gamma_9q_{2,9}.
\end{align*}
Since 
$\psi_{2,2}(\xi_{2,1})+\sum_{0\leqslant u \leqslant 8}q_{2,u} = \psi_{2,2}(\xi_{2,1}) + \sum_{1\leqslant j \leqslant 90,\, j \notin \{31,\, 33,\, 34,\, 36,\, 37,\, 40,\, 41,\, 42,\, 69\}\atop\hskip1.1cm\cup\{71,\, 74,\, 83,\, 84,\, 85,\, 86,\, 87,\, 88,\, 89\}}b_{2,j} = \delta_{2,2}$ and $q_{2,9} = \zeta_{2,2}$, we get $f \equiv \gamma_0\delta_{2,2} +\gamma_9\zeta_{2,2}$. 	
The theorem is proved for $t = 2$.

For $t > 2$, we denote $b_{2,j} = b_{t,s,j}$ and $q_{2,j} = q_{t,2,j}$. By a direct computation we obtain
\begin{align*}
\rho_1(\psi_{2,t}(\xi_{t,1}))+\psi_{2,t}(\xi_{t,1}) &\equiv 0,\\ 
\rho_2(\psi_{2,t}(\xi_{t,1}))+\psi_{2,t}(\xi_{t,1}) &\equiv  x_1x_2^{2}x_3^{5}x_4^{2^{t+2}-6} + x_1x_2^{2}x_3^{2^{t+2}-3}x_4^{2} + x_1x_2^{3}x_3^{4}x_4^{2^{t+2}-6}\\ &\quad + x_1x_2^{3}x_3^{2^{t+2}-4}x_4^{2} = b_{2,73} + b_{2,79} + b_{2,83} + b_{2,85},\\
\rho_3(\psi_{2,t}(\xi_{t,1}))+\psi_{2,t}(\xi_{t,1}) &\equiv   x_1x_2^{2}x_3^{3}x_4^{2^{t+2}-4} + x_1x_2^{2}x_3^{2^{t+2}-4}x_4^{3} =   
b_{2,74} + b_{2,75}. 
\end{align*} 
By computing from these results and the relations $\rho_i(f) \equiv 0$, $i = 1,\, 2,\, 3$, we obtain
$$f \equiv  \gamma_0(\psi_{t,2}(\xi_{t,1}) + q_{2,0}) + \sum_{1 \leqslant u \leqslant 9}\gamma_uq_{2,u}$$
with $\gamma_u \in \mathbb F_2$ and 
\begin{align*}
q_{2,0} &= x_1x_2x_3^{2}x_4^{2^{t+2}-2} + x_1x_2x_3^{2^{t+2}-2}x_4^{2} + x_1x_2^{2^{t+2}-2}x_3x_4^{2} + x_1x_2x_3^{6}x_4^{2^{t+2}-6}\\ &\quad + x_1x_2^{6}x_3x_4^{2^{t+2}-6} + x_1x_2^{2}x_3^{2^{t+2}-3}x_4^{2} + x_1x_2^{2}x_3^{2^{t+2}-4}x_4^{3} + x_1x_2^{2}x_3^{5}x_4^{2^{t+2}-6}\\
&=\mbox{$\sum_{j\in\{67,\, 68,\, 70,\, 71,\, 72,\, 73,\, 75,\, 79\}}b_{2,j}$}.
\end{align*}

We have
\begin{align*}
\rho_4(\psi_{2,t}(\xi_{t,1}))&+\psi_{2,t}(\xi_{t,1}) \equiv x_2x_3^{2}x_4^{2^{t+2}-1} + x_2x_3^{3}x_4^{2^{t+2}-2} + x_2x_3^{2^{t+2}-1}x_4^{2}\\ &\quad + x_1x_2^{3}x_3^{2}x_4^{2^{t+2}-4} + x_1x_2^{7}x_3^{2}x_4^{2^{t+2}-8} + x_1^{3}x_2^{3}x_3^{4}x_4^{2^{t+2}-8}\\ &=\mbox{$\sum_{j\in\{19,\, 20,\, 31,\, 76,\, 81,\, 90\}}b_{2,j}$}.
\end{align*}

By using this result and computing $\rho_4(f)+f$ in terms of the admissible monomials, we obtain
\begin{align*}
\rho_4(f)+f &\equiv \gamma_{\{1,3\}}b_{2,3} + \gamma_{\{1,3\}}b_{2,4} + \gamma_{\{1,4,5,6\}}b_{2,5} + \gamma_{\{1,4,5,6\}}b_{2,6} + \gamma_{\{1,2\}}b_{2,9}\\
& + \gamma_{\{2,6\}}b_{2,14} + \gamma_{\{2,6\}}b_{2,15} + \gamma_{\{0,3\}}b_{2,19} + \gamma_{\{0,3\}}b_{2,20} + \gamma_{\{0,3,9\}}b_{2,21}\\
& + \gamma_{\{3,5\}}b_{2,26} + \gamma_{\{3,5\}}b_{2,27} + \gamma_{\{4,6,7,8\}}b_{2,31} + \gamma_{\{0,4,5,6,7,8\}}b_{2,32}\\
& + \gamma_{\{5,6,9\}}b_{2,33} + \gamma_{\{5,6\}}b_{2,36} + \gamma_{\{5,6\}}b_{2,37} + \gamma_{\{6,8\}}b_{2,43} + \gamma_{\{6,8\}}b_{2,47}\\
& + \gamma_{\{0,6\}}b_{2,48} + \gamma_{\{4,6\}}b_{2,50} + \gamma_{\{4,6\}}b_{2,51} + \gamma_{\{0,5,9\}}b_{2,55} + \gamma_{\{0,5\}}b_{2,59}\\
& + \gamma_{\{1,2,3,6\}}b_{2,61} + \gamma_{\{1,2,3,6\}}b_{2,62} + \gamma_{\{0,4,9\}}b_{2,63} + \gamma_{\{3,5\}}b_{2,67}\\
& + \gamma_{\{3,5\}}b_{2,68} + \gamma_{\{3,4\}}b_{2,71} + \gamma_{\{0,8\}}b_{2,76} + \gamma_{9}b_{2,77} + \gamma_{\{0,8,9\}}b_{2,83}\\
& + \gamma_{9}b_{2,85} + \gamma_{9}b_{2,87} + \gamma_{\{0,3,6,8\}}b_{2,90} \equiv 0.
\end{align*}
The above equality implies $\gamma_9 = 0$ and $\gamma_t =\gamma_0$ for all $0\leqslant t \leqslant 8$. Hence,  
\begin{align*}
f &\equiv \gamma_0\Big(\psi_{2,t}(\xi_{t,1})+  \sum_{0\leqslant u \leqslant 8}q_{2,u}\Big) = \gamma_0\delta_{t,2}  
\end{align*}
Since $\sum_{0\leqslant u \leqslant 8}q_{2,u}= \sum_{1\leqslant j \leqslant 90,\, j \notin \{31,\, 33,\, 34,\, 36,\, 37,\, 40,\, 41,\, 42,\, 69\}\atop\hskip1.1cm\cup\{71,\, 74,\, 83,\, 84,\, 85,\, 86,\, 87,\, 88,\, 89\}}b_{2,j}$, we get $f \equiv \gamma_0\delta_{t,2}$. 
Thus, the theorem is completely proved.
\end{proof}

\begin{rems}
In \cite[Page 471, Line 15$\uparrow$]{pp25}, the author stated that $(QP_4)_{n_{2,t}}^{GL_4}$ is generated by $[\psi_{2,t}(\xi_{1,t})]$ but this is false because the class	$[\psi_{2,t}(\xi_{1,t})]$ is not an $GL_4$-invariant in $(QP_4)_{n_{2,t}}$.
\end{rems}

\subsubsection{\textbf{The subcase $s = 3,\, t \geqslant 2$}}\ 

\medskip
We prove the following.
\begin{thms} For $s = 3,\, t\geqslant 2$, we have $(QP_4)_{2^{t+3}+6}^{GL_4} = \langle [\zeta_{t,3}], [\delta_{t,3}]\rangle$, where $\zeta_{t,3}$ is defined in Theorem \ref{dln2} and 
$$\delta_{t,3} = \begin{cases} \psi_{3,2}(\xi_{2,2})+ \sum_{j\in \{67,\, 71,\, 73,\, 83,\, 86,\, 88,\, 94,\, 99\}}b_{3,j}, &\mbox{if } t = 2,\\ 
\psi_{3,3}(\xi_{3,2})+ \mbox{$\sum_{j \in\{71,\, 74,\, 77,\, 81,\, 82,\, 86,\, 89,\, 94,\, 95,\, 100,\, 102,\, 105\}}$}b_{3,j}, &\mbox{if } t = 3,\\
\psi_{3,t}(\theta_{t,2}) + \sum_{j \in\{68,\, 71,\, 74,\, 75,\, 77,\, 81,\, 82,\, 84,\, 86,\, 89,\, 95,\, 102,\, 104\}}b_{3,j}, &\mbox{if } t \geqslant 4.\end{cases}$$
Here, $\xi_{2,2},\, \xi_{3,2}$, $\theta_{t,2}$ are defined in Theorem \ref{dlt21}, Proposition \ref{mdt32}, Theorem \ref{d1t4}. 	 
\end{thms}
\begin{proof} Let $[f] \in (QP_4)_{2^{t+3}+6}^{GL_4}$ with $f \in (P_4)_{2^{t+3}+6}$, then there is $\lambda_0 \in \mathbb F_2$ and $h \in \mbox{Ker}(\widetilde{Sq}^0_*)_{(4,2^{t+3}+6)}$ such that
	$$f \equiv \gamma_0\psi_{3,t}(\xi_{t,2}) + h.$$
	 
For $t = 2$, we have $n_{3,2} = 38$. From Theorem \ref{dlt21} we have
\begin{align*}
\xi_{2,2} &= x_1x_2x_3x_4^{14} + x_1x_2x_3^{14}x_4 + x_1x_2^{3}x_3x_4^{12} + x_1x_2^{3}x_3^{12}x_4\\ &\quad + x_1^{3}x_2x_3x_4^{12} + x_1^{3}x_2x_3^{12}x_4 + x_1^{3}x_2^{5}x_3x_4^{8} + x_1^{3}x_2^{5}x_3^{8}x_4.
\end{align*} 
By a direct computation we get
\begin{align*}
\rho_1(\psi_{3,2}(\xi_{2,2}))&+\psi_{3,2}(\xi_{2,2}) \equiv x_1x_2x_3^{6}x_4^{30} + x_1x_2x_3^{30}x_4^{6} + x_1x_2^{3}x_3^{6}x_4^{28} + x_1x_2^{3}x_3^{28}x_4^{6}\\ & \quad + x_1^{3}x_2x_3^{6}x_4^{28} + x_1^{3}x_2x_3^{28}x_4^{6} + x_1x_2^{3}x_3^{12}x_4^{22} + x_1^{3}x_2x_3^{12}x_4^{22}\\ & \quad + x_1x_2^{3}x_3^{14}x_4^{20} + x_1^{3}x_2x_3^{14}x_4^{20}\\
&= \mbox{$\sum_{j\in\{67,\, 68,\, 83,\, 84,\, 85,\, 86,\, 87,\, 88,\, 99,\, 100\}}b_{3,j}$},\\
\rho_2(\psi_{3,2}(\xi_{2,2}))&+\psi_{3,2}(\xi_{2,2}) \equiv x_1^{3}x_2x_3^{6}x_4^{28} + x_1^{3}x_2x_3^{14}x_4^{20} + x_1^{3}x_2^{5}x_3^{2}x_4^{28} + x_1^{3}x_2^{5}x_3^{6}x_4^{24}\\ & \quad + x_1^{3}x_2^{13}x_3^{2}x_4^{20} = b_{3,88} + b_{3,89} + b_{3,100} + b_{3,102} + b_{3,105},\\
\rho_3(\psi_{3,2}(\xi_{2,2}))&+\psi_{3,2}(\xi_{2,2}) \equiv 0. 
\end{align*} 
By computing based on these results, Lemma \ref{bdn22} and the relations $\rho_i(f)+ f \equiv 0$, $i = 1,\, 2,\, 3$, we obtain
$$f \equiv  \gamma_0(\psi_{3,2}(\xi_{2,2}) + q_{3,0}) + \sum_{1 \leqslant u \leqslant 13}\gamma_uq_{3,u}$$
with $\gamma_u \in \mathbb F_2$ and 
\begin{align*}
q_{3,0} &= x_1x_2x_3^{6}x_4^{30} + x_1x_2x_3^{14}x_4^{22} + x_1x_2^{2}x_3^{29}x_4^{6} + x_1x_2^{3}x_3^{28}x_4^{6}\\ &\quad+ x_1^{3}x_2x_3^{12}x_4^{22} + x_1^{3}x_2x_3^{14}x_4^{20} + x_1x_2^{2}x_3^{5}x_4^{30} + x_1x_2^{3}x_3^{6}x_4^{28}\\
&=b_{3,67} + b_{3,71} + b_{3,73} + b_{3,83} + b_{3,86} + b_{3,88} + b_{3,94} + b_{3,99}.
\end{align*}

We have
\begin{align*}
\rho_4(\psi_{3,2}(\xi_{2,2}))+\psi_{3,2}(\xi_{2,2}) &\equiv x_2^{3}x_3^{5}x_4^{30} + x_2^{3}x_3^{29}x_4^{6} + x_1x_2^{3}x_3^{12}x_4^{22} + x_1x_2^{3}x_3^{14}x_4^{20}\\
&= b_{3,55} + b_{3,63} + b_{3,85} + b_{3,87}.
\end{align*}

By using this result and computing $\rho_4(f)+f$ in terms of the admissible monomials, we obtain
\begin{align*}
\rho_4(f)&+f \equiv \gamma_{\{1,3\}}b_{3,3} + \gamma_{\{1,3\}}b_{3,4} + \gamma_{\{1,4\}}b_{3,5} + \gamma_{\{1,4\}}b_{3,6} + \gamma_{\{1,2\}}b_{3,9}\\ 
&\quad + \gamma_{\{2,5\}}b_{3,14} + \gamma_{\{2,5\}}b_{3,15} + \gamma_{\{3,7\}}b_{3,19} + \gamma_{\{3,7\}}b_{3,20} + \gamma_{\{3,9,12\}}b_{3,21}\\ 
&\quad + \gamma_{\{3,4,5,6\}}b_{3,26} + \gamma_{\{3,4,5,6\}}b_{3,27} + \gamma_{\{4,10,13\}}b_{3,31} + \gamma_{\{4,10,13\}}b_{3,32}\\ 
&\quad + \gamma_{\{4,9,12\}}b_{3,33} + \gamma_{\{1,2,3,4,5,6\}}b_{3,36} + \gamma_{\{1,2,3,4,5,6\}}b_{3,37} + \gamma_{\{1,2,3,5\}}b_{3,41}\\ 
&\quad + \gamma_{\{1,2,3,5\}}b_{3,42} + \gamma_{\{5,10\}}b_{3,43} + \gamma_{\{5,10\}}b_{3,47} + \gamma_{\{5,8,12\}}b_{3,48}\\ 
&\quad + \gamma_{\{5,6\}}b_{3,50} + \gamma_{\{5,6\}}b_{3,51} + \gamma_{\{4,5,6,11,12\}}b_{3,55} + \gamma_{\{6,8\}}b_{3,59}\\ 
&\quad + \gamma_{\{6,8,9,11,12\}}b_{3,63} + \gamma_{\{1,2,3,5\}}b_{3,65} + \gamma_{\{1,2,3,5\}}b_{3,66} + \gamma_{\{4,6,8,10,12\}}b_{3,67}\\ 
&\quad + \gamma_{\{3,4,5,6,7,11\}}b_{3,68} + \gamma_{\{5,6,8,10,12\}}b_{3,71} + \gamma_{\{3,5,7,11,12\}}b_{3,76}\\ 
&\quad + \gamma_{\{3,5,7,10\}}b_{3,77} + \gamma_{\{8,10\}}b_{3,81} + \gamma_{12}b_{3,83} + \gamma_{12}b_{3,85} + \gamma_{12}b_{3,87}\\ 
&\quad + \gamma_{\{7,9\}}b_{3,92} + \gamma_{\{3,5,7,8,9,10,12\}}b_{3,95} + \gamma_{\{9,11\}}b_{3,96} + \gamma_{\{3,5,7,10\}}b_{3,97}\\ 
&\quad + \gamma_{\{8,11\}}b_{3,99} + \gamma_{\{3,5,7,10\}}b_{3,105} \equiv 0.
\end{align*}
Computing from this equality gives $\gamma_{12} =\gamma_{13} = 0$ and $\gamma_u =\gamma_1$ for all $0\leqslant u \leqslant 11$. Hence, we obtain 
\begin{align*}
f &\equiv \gamma_0\Big(\psi_{3,2}(\xi_{2,2}) + q_{3,0}\Big) +\gamma_1 \Big(\sum_{1\leqslant u \leqslant 11}q_{3,u}\Big).
\end{align*}
Since $\psi_{3,2}(\xi_{2,2}) + q_{3,0} = \delta_{3,2}$ and $\sum_{1\leqslant u \leqslant 11}q_{3,u} = \zeta_{3,2}$ we obtain $f \equiv \gamma_0\delta_{3,2} +\gamma_1\zeta_{3,2}$. Hence, the theorem is proved for $t = 2$.

For $t = 3$, we have $n_{3,3} = 70$. From Theorem \ref{mdt32} we have
\begin{align*}
\xi_{3,2} = \bar\xi_{3,2} + \sum_{j\in \{29,\, 36,\, 37,\, 39,\, 41,\, 42,\, 45\}}\widetilde a_{3,j}.
\end{align*} 
By a direct computation we get
\begin{align*}
\rho_1(\psi_{3,3}(\xi_{3,2}))&+\psi_{3,3}(\xi_{3,2}) \equiv  \mbox{$\sum_{j\in\{71,\, 76,\, 77,\, 85,\, 86,\, 95,\, 97,\, 99,\, 100,\, 105\}}b_{3,j}$},\\
\rho_2(\psi_{3,3}(\xi_{3,2}))&+\psi_{3,3}(\xi_{3,2}) \equiv \mbox{$\sum_{j\in\{74,\, 76,\, 80,\, 81,\, 88,\, 89,\, 100,\, 102,\, 105\}}b_{3,j}$},\\
\rho_3(\psi_{3,3}(\xi_{3,2}))&+\psi_{3,3}(\xi_{3,2}) \equiv \mbox{$\sum_{j\in\{67,\, 71,\, 73,\, 74,\, 75,\, 83,\, 86,\, 88,\, 94,\, 99\}}b_{3,j}$}. 
\end{align*} 
By computing based on these results, Lemma \ref{bdn22} and the relations $\rho_i(f)+ f \equiv 0$, $i = 1,\, 2,\, 3$, we obtain
$$f \equiv  \gamma_0(\psi_{3,3}(\xi_{3,2}) + q_{3,0}) + \sum_{1 \leqslant u \leqslant 13}\gamma_uq_{3,u}$$
with $\gamma_u \in \mathbb F_2$ and 
\begin{align*}
q_{3,0} &= \mbox{$\sum_{j\in\{69,\, 70,\, 71,\, 73,\, 75,\, 76,\, 78,\, 79,\, 80,\, 83,\, 84,\, 85,\, 87,\, 88,\, 94,\, 95,\, 99,\, 105\}}b_{3,j}$}.
\end{align*}
We have
\begin{align*}
\rho_4(\psi_{3,3}(\xi_{3,2}))&+\psi_{3,2}(\xi_{3,2}) \equiv \mbox{$\sum_{j\in\{31,\, 55,\, 71,\, 77,\, 85,\, 97,\, 99,\, 105\}}b_{3,j}$}.
\end{align*}

By using this result and computing $\rho_4(f)+f$ in terms of the admissible monomials, we obtain
\begin{align*}
\rho_4(f)&+f \equiv \gamma_{\{1,3\}}b_{3,3} + \gamma_{\{1,3\}}b_{3,4} + \gamma_{\{1,4\}}b_{3,5} + \gamma_{\{1,4\}}b_{3,6} + \gamma_{\{1,2\}}b_{3,9}\\
& + \gamma_{\{2,5\}}b_{3,14} + \gamma_{\{2,5\}}b_{3,15} + \gamma_{\{3,7\}}b_{3,19} + \gamma_{\{3,7\}}b_{3,20} + \gamma_{\{0,3,9,12\}}b_{3,21}\\
& + \gamma_{\{3,4,5,6\}}b_{3,26} + \gamma_{\{3,4,5,6\}}b_{3,27} + \gamma_{\{0,4,10,13\}}b_{3,31} + \gamma_{\{0,4,10,13\}}b_{3,32}\\
& + \gamma_{\{0,4,9,12\}}b_{3,33} + \gamma_{\{1,2,3,4,5,6\}}b_{3,36} + \gamma_{\{1,2,3,4,5,6\}}b_{3,37} + \gamma_{\{1,2,3,5\}}b_{3,41}\\
& + \gamma_{\{1,2,3,5\}}b_{3,42} + \gamma_{\{0,5,10\}}b_{3,43} + \gamma_{\{0,5,10\}}b_{3,47} + \gamma_{\{5,8,12\}}b_{3,48}\\
& + \gamma_{\{5,6\}}b_{3,50} + \gamma_{\{5,6\}}b_{3,51} + \gamma_{\{6,8,9,11,12\}}b_{3,55} + \gamma_{\{1,2,3,5\}}b_{3,57}\\
& + \gamma_{\{1,2,3,5\}}b_{3,58} + \gamma_{\{0,4,5,6,11,12\}}b_{3,59} + \gamma_{\{0,6,8\}}b_{3,63} + \gamma_{\{0,4,6,8,10,12\}}b_{3,67}\\
& + \gamma_{\{3,4,5,6,7,11\}}b_{3,68} + \gamma_{\{0,5,6,8,10,12\}}b_{3,71} + \gamma_{\{0,3,5,7,11,12\}}b_{3,76}\\
& + \gamma_{\{0,3,5,7,10\}}b_{3,77} + \gamma_{\{8,10\}}b_{3,81} + \gamma_{\{0,12\}}b_{3,83} + \gamma_{\{0,12\}}b_{3,85}\\
& + \gamma_{\{0,12\}}b_{3,87} + \gamma_{\{7,9\}}b_{3,92} + \gamma_{\{0,3,5,7,8,9,10,12\}}b_{3,95} + \gamma_{\{9,11\}}b_{3,96}\\
& + \gamma_{\{0,3,5,7,10\}}b_{3,97} + \gamma_{\{0,8,11\}}b_{3,99} + \gamma_{\{0,3,5,7,10\}}b_{3,105} \equiv 0.
\end{align*}
Computing from this equality gives $\gamma_{13} = 0$ and $\gamma_u =\gamma_1$ for all $0\leqslant u \leqslant 11$, $u \ne 8,\, 10$, $\gamma_{12} = \gamma_0$ and $\gamma_8 = \gamma_{10} = \gamma_0 + \gamma_1$. Hence, we obtain 
\begin{align*}
f &\equiv \gamma_0\Big(\psi_{3,3}(\xi_{3,2}) + q_{3,0} + q_{3,8} + q_{3,10} + q_{3,12}\Big) +\gamma_1 \Big(\sum_{1\leqslant u \leqslant 11}q_{3,u}\Big). 
\end{align*}
We have $q_{3,0} + q_{3,8} + q_{3,10} + q_{3,12} = \sum_{j\in \{25,\, 26,\, 27,\, 28,\, 35,\, 56,\, 62,\, 65,\, 66\}\atop \hskip0.5cm\cup\{67,\, 68,\, 69,\, 78,\, 81,\, 82,\, 83,\, 84\}}\widetilde a_{3,j} = \delta_{3,3}$ and $\sum_{1\leqslant u \leqslant 11}q_{3,u} = \zeta_{3,3}$. Hence $f\equiv \gamma_0\delta_{3,3} +\gamma_1\zeta_{3,3}$. The theorem is proved for $t = 3$.	

For $t \geqslant 4$, from Theorem \ref{d1t4} we have
\begin{align*}
\theta_{t,2} = 
\bar \xi_{t,2} + \mbox{$\sum_{j\in \{25,\, 26,\, 27,\, 28,\, 29,\, 32,\, 35,\, 37,\, 39,\, 43,\, 44,\, 45\}}$}\widetilde a_{2,j}.
\end{align*} 
By a direct computation we get
\begin{align*}
\rho_1(\psi_{3,t}(\theta_{t,2}))&+\psi_{3,t}(\theta_{t,2}) \equiv  \mbox{$\sum_{j\in\{67,\, 68,\, 71,\, 76,\, 77,\, 83,\, 84,\, 85,\, 86,\, 95,\, 97,\, 105\}}b_{3,j}$},\\
\rho_2(\psi_{3,t}(\theta_{t,2}))&+\psi_{3,t}(\theta_{t,2}) \equiv \mbox{$\sum_{j\in\{74,\, 76,\, 80,\, 81,\, 88,\, 89,\, 100,\, 102,\, 105\}}b_{3,j}$},\\
\rho_3(\psi_{3,t}(\theta_{t,2}))&+\psi_{3,t}(\theta_{t,2}) \equiv \mbox{$\sum_{j\in\{68,\, 71,\, 84,\, 86,\, 88,\, 89,\, 100,\, 104\}}b_{3,j}$}. 
\end{align*} 
By computing based on these results, Lemma \ref{bdn22} and the relations $\rho_i(f)+ f \equiv 0$, $i = 1,\, 2,\, 3$, we obtain
$$f \equiv  \gamma_0(\psi_{3,t}(\theta_{t,2}) + q_{3,0}) + \sum_{1 \leqslant u \leqslant 13}\gamma_uq_{3,u}$$
with $\gamma_u \in \mathbb F_2$ and 
\begin{align*}
q_{3,0} &= \mbox{$\sum_{j\in\{69,\, 70,\, 71,\, 76,\, 78,\, 79,\, 80,\, 84,\, 85,\, 87,\, 88,\, 89,\, 94,\, 95,\, 99,\, 100,\, 105\}}b_{3,j}$}.
\end{align*}
We have
\begin{align*}
\rho_4(\psi_{3,t}(\theta_{t,2}))&+\psi_{3,t}(\theta_{t,2}) \equiv \mbox{$\sum_{j\in\{31,\, 32,\, 68,\, 71,\, 77,\, 85,\, 96,\, 97,\, 99,\, 105\}}b_{3,j}$}.
\end{align*}

By using this result and computing $\rho_4(f)+f$ in terms of the admissible monomials, we obtain
\begin{align*}
\rho_4(f)&+f \equiv \gamma_{\{1,3\}}b_{3,3} + \gamma_{\{1,3\}}b_{3,4} + \gamma_{\{1,4\}}b_{3,5} + \gamma_{\{1,4\}}b_{3,6} + \gamma_{\{1,2\}}b_{3,9}\\
& + \gamma_{\{2,5\}}b_{3,14} + \gamma_{\{2,5\}}b_{3,15} + \gamma_{\{3,7\}}b_{3,19} + \gamma_{\{3,7\}}b_{3,20} + \gamma_{\{0,3,9,12\}}b_{3,21}\\
& + \gamma_{\{3,4,5,6\}}b_{3,26} + \gamma_{\{3,4,5,6\}}b_{3,27} + \gamma_{\{0,4,10,13\}}b_{3,31} + \gamma_{\{0,4,10,13\}}b_{3,32}\\
& + \gamma_{\{0,4,9,12\}}b_{3,33} + \gamma_{\{1,2,3,4,5,6\}}b_{3,36} + \gamma_{\{1,2,3,4,5,6\}}b_{3,37} + \gamma_{\{1,2,3,5\}}b_{3,41}\\
& + \gamma_{\{1,2,3,5\}}b_{3,42} + \gamma_{\{0,5,10\}}b_{3,43} + \gamma_{\{0,5,10\}}b_{3,47} + \gamma_{\{5,8,12\}}b_{3,48}\\
& + \gamma_{\{5,6\}}b_{3,50} + \gamma_{\{5,6\}}b_{3,51} + \gamma_{\{0,6,8,9,11,12\}}b_{3,55} + \gamma_{\{1,2,3,5\}}b_{3,57}\\
& + \gamma_{\{1,2,3,5\}}b_{3,58} + \gamma_{\{4,5,6,11,12\}}b_{3,59} + \gamma_{\{0,6,8\}}b_{3,63} + \gamma_{\{0,4,6,8,10,12\}}b_{3,67}\\
& + \gamma_{\{0,3,4,5,6,7,11\}}b_{3,68} + \gamma_{\{0,5,6,8,10,12\}}b_{3,71} + \gamma_{\{3,5,7,11,12\}}b_{3,76}\\
& + \gamma_{\{0,3,5,7,10\}}b_{3,77} + \gamma_{\{8,10\}}b_{3,81} + \gamma_{\{0,12\}}b_{3,83} + \gamma_{\{0,12\}}b_{3,85}\\
& + \gamma_{\{0,12\}}b_{3,87} + \gamma_{\{7,9\}}b_{3,92} + \gamma_{\{0,3,5,7,8,9,10,12\}}b_{3,95} + \gamma_{\{0,9,11\}}b_{3,96}\\
& + \gamma_{\{0,3,5,7,10\}}b_{3,97} + \gamma_{\{8,11\}}b_{3,99} + \gamma_{\{0,3,5,7,10\}}b_{3,105} +  \equiv 0.
\end{align*}
Computing from this equality gives $\gamma_{13} = 0$ and $\gamma_u =\gamma_1$ for all $0\leqslant u \leqslant 9,\, u \ne 8$, $\gamma_{8} = \gamma_{10} = \gamma_{11} = \gamma_0 + \gamma_1$ and $\gamma_{12} = \gamma_{0}$. Hence, we obtain 
\begin{align*}
f &\equiv \gamma_0\Big(\psi_{3,t}(\xi_{t,2}) + q_{3,0} + q_{3,8} + q_{3,10} + q_{3,11} + q_{3,12}\Big) +\gamma_1 \Big(\sum_{1\leqslant u \leqslant 11}q_{3,u}\Big). 
\end{align*}
We have $\psi_{3,t}(\xi_{t,2}) + q_{3,0} + q_{3,8} + q_{3,10} + q_{3,11} + q_{3,12} = \delta_{t,3}$ and $\sum_{1\leqslant u \leqslant 11}q_{3,u} = \zeta_{t,3}$. Hence $f\equiv \gamma_0\delta_{t,3} +\gamma_1\zeta_{t,3}$. This completes the proof of the theorem.	
\end{proof}

\subsubsection{\textbf{The subcase $s \geqslant 4$, $t \geqslant 2$}}\
	
\medskip	
We prove the following.
\begin{thms} For $s \geqslant 4,\, t\geqslant 2$, we have $$(QP_4)_{n_{s,t}}^{GL_4} = \begin{cases}
\langle [\zeta_{t,4}]\rangle, &\mbox{if } s = 4,\, 2 \leqslant t \leqslant 3,\\
\langle [\zeta_{t,4}],[\delta_{t,4}]\rangle, &\mbox{if } s = 4,\, t \geqslant 4,\\
\langle [\zeta_{t,s}], [\chi_{t,s}]\rangle, &\mbox{if } s\geqslant 5,\, 2 \leqslant t \leqslant 3,\\
\langle [\zeta_{t,s}], [\delta_{t,s}], [\chi_{t,s}]\rangle, &\mbox{if } s \geqslant 5, t \geqslant 4,\\
\end{cases}$$ 	
where $\zeta_{t,s}$ is defined in Theorem \ref{dln2}, and 
$$\begin{cases}
\delta_{t,s} = \psi_{s,t}(\theta_{t,s-1}) +\sum_{j \in\{68,\, 74,\, 75,\, 84,\, 85,\, 86\}\atop\hskip0.1cm\cup\{87,\, 88,\, 100,\, 104,\, 105\}}b_{s,j}, &\mbox{for }s \geqslant 4,\, t \geqslant 4,\\
\chi_{t,s} = \psi_{s,t}(\xi_{t,s-1}) +\sum_{j \in\{67,\, 76,\, 77,\, 83,\, 84,\, 88,\, 90,\, 91,\, 92\}\atop\hskip0.1cm\cup\{93,\, 94,\, 96,\, 98,\, 99,\, 101,\, 102,\, 104\}}b_{s,j}, &\mbox{for }s \geqslant 5.
\end{cases}$$
Here, $\xi_{t,s-1}$ and $\theta_{t,s-1}$ are respectively defined in Theorems \ref{dlt21} and \ref{d1t4}. 	 
\end{thms}
\begin{proof} Recall that Kameko's homomorphism
$$(\widetilde{Sq}^0_*)_{(4,n_{s,t})} : (QP_{4})_{n_{s,t}} \longrightarrow (QP_4)_{d_{s-1,t}}$$
is a homomorphism of $GL_4$-modules, hence it induces a one also denote by  
$$(\widetilde{Sq}^0_*)_{(4,n_{s,t})} : (QP_{4})_{n_{s,t}}^{GL_4} \longrightarrow (QP_4)_{d_{s-1,t}}^{GL_4}.$$

By Theorem \ref{dln2}, for $s \geqslant 4,\, t\geqslant 2$, $$\mbox{Ker}\big((\widetilde{Sq}^0_*)_{(4,n_{s,t})}\big)^{GL_4} = \langle [\zeta_{t,s}] \rangle \subset (QP_{4})_{n_{s,t}}^{GL_4}.$$  
	
For $s = 4$ and $2 \leqslant t \leqslant 3$, by Theorem \ref{thm1}, $(QP_4)_{d_{3,t}}^{GL_4} = 0$. Hence, 
$$(QP_4)_{n_{4,t}}^{GL_4} = \mbox{Ker}\big((\widetilde{Sq}^0_*)_{(4,n_{4,t})}\big)^{GL_4} = \langle [\zeta_{t,4}] \rangle.$$ 
Thus, the theorem holds for $s = 4$ and $2 \leqslant t \leqslant 3$.

For $s \geqslant 4,\, t \geqslant 4$, by Theorem \ref{d1t4}, $[\theta_{t,s-1}] \in (QP_4)_{d_{s-1,t}}^{GL_4}$. We find an element $h  \in \mbox{Ker}\big((\widetilde{Sq}^0_*)_{(4,n_{s,t})}\big)$ such that $f = \psi_{s,t}(\theta_{t,s-1}) + h$ and $[f] \in (QP_4)_{n_{s,t}}^{GL_4}$.

By a direct computation we get
\begin{align*}
\rho_1(\psi_{s,t}&(\theta_{t,s-1}))+\psi_{s,t}(\theta_{t,s-1}) \equiv x_1x_2x_3^{2^{s}-2}x_4^{2^{s+t}-2} + x_1x_2x_3^{2^{s+t}-2}x_4^{2^{s}-2}\\ 
& + x_1x_2^{3}x_3^{2^{s+t}-4}x_4^{2^{s}-2} + x_1^{3}x_2x_3^{2^{s+t}-4}x_4^{2^{s}-2} + x_1x_2^{3}x_3^{2^{s}-2}x_4^{2^{s+t}-4}\\ 
& + x_1^{3}x_2x_3^{2^{s}-2}x_4^{2^{s+t}-4} = \mbox{$\sum_{j\in\{67,\, 68,\, 83,\, 84,\, 99,\, 100\}}b_{s,j}$},\\
\rho_2(\psi_{s,t}&(\theta_{t,s-1}))+\psi_{s,t}(\theta_{t,s-1}) \equiv x_1x_2^{2}x_3^{2^{s}-1}x_4^{2^{s+t}-4} + x_1x_2^{2}x_3^{2^{s+1}-3}x_4^{2^{s+t}-2^{s}-2}\\ 
& + x_1x_2^{3}x_3^{2^{s}-4}x_4^{2^{s+t}-2} + x_1x_2^{3}x_3^{2^{s+1}-4}x_4^{2^{s+t}-2^{s}-2} + x_1x_2^{2^{s}-1}x_3^{2}x_4^{2^{s+t}-4}\\ 
& + x_1^{3}x_2x_3^{2^{s}-2}x_4^{2^{s+t}-4} + x_1^{3}x_2x_3^{2^{s+1}-2}x_4^{2^{s+t}-2^{s}-4} + x_1^{3}x_2^{5}x_3^{2^{s}-6}x_4^{2^{s+t}-4}\\ 
& + x_1^{3}x_2^{2^{s}-3}x_3^{2}x_4^{2^{s+t}-4} + x_1^{3}x_2^{2^{s+1}-3}x_3^{2}x_4^{2^{s+t}-2^{s}-4}\\ 
&= \mbox{$\sum_{j\in\{74,\, 76,\, 79,\, 85,\, 88,\, 89,\, 95,\, 100,\, 102,\, 105\}}b_{s,j}$},\\
\rho_3(\psi_{s,t}&(\theta_{t,s-1}))+\psi_{s,t}(\theta_{t,s-1}) \equiv 0. \end{align*} 
By computing based on these results, Lemma \ref{bdn22} and the relations $\rho_i(f)+ f \equiv 0$, $i = 1,\, 2,\, 3$, we obtain
$$f \equiv  \psi_{s,t}(\theta_{t,s-1})) + q_{s,0} + \sum_{1 \leqslant u \leqslant 13}\gamma_uq_{s,u}$$
with $\gamma_u \in \mathbb F_2$ and 
$q_{s,0} = \sum_{j\in\{67,\, 69,\, 70,\, 71,\, 72,\, 74,\, 75,\, 84,\, 85,\, 86,\, 87,\, 88,\, 100,\, 103\}}b_{s,j}$.
	
We have
\begin{align*}
\rho_4(\psi_{s,t}&(\theta_{t,s-1}))+\psi_{s,t}(\theta_{t,s-1}) \equiv x_2x_3^{2^{s}-1}x_4^{2^{s+t}-2} + x_2x_3^{2^{s+t}-2}x_4^{2^{s}-1}\\
& + x_1x_2x_3^{2^{s}-2}x_4^{2^{s+t}-2} + x_1x_2x_3^{2^{s+t}-2}x_4^{2^{s}-2} + x_1x_2^{3}x_3^{2^{s}-4}x_4^{2^{s+t}-2}\\
& + x_1x_2^{3}x_3^{2^{s+1}-4}x_4^{2^{s+t}-2^{s}-2} + x_1x_2^{3}x_3^{2^{s+1}-2}x_4^{2^{s+t}-2^{s}-4} + x_1x_2^{3}x_3^{2^{s+t}-2}x_4^{2^{s}-4} \\ 
&= \mbox{$\sum_{j\in\{31,\, 32,\, 67,\, 68,\, 95,\, 85,\, 87,\, 96\}}b_{s,j}$}
\end{align*}
	
By using this result and computing $\rho_4(f)+f$ in terms of the admissible monomials, we obtain
\begin{align*}
\rho_4(f)&+f \equiv \gamma_{\{1,3\}}b_{s,3} + \gamma_{\{1,3\}}b_{s,4} + \gamma_{\{1,4\}}b_{s,5} + \gamma_{\{1,4\}}b_{s,6} + \gamma_{\{1,2\}}b_{s,9}\\
& + \gamma_{\{2,5\}}b_{s,14} + \gamma_{\{2,5\}}b_{s,15} + \gamma_{\{3,7\}}b_{s,19} + \gamma_{\{3,7\}}b_{s,20} + (\gamma_{\{3,9,11\}}+1)b_{s,21}\\
& + \gamma_{\{3,4,5,6\}}b_{s,26} + \gamma_{\{3,4,5,6\}}b_{s,27} + \gamma_{\{4,10,13\}}b_{s,31} + \gamma_{\{4,10,13\}}b_{s,32}\\
& + (\gamma_{\{4,9,11\}}+1)b_{s,33} + \gamma_{\{1,2,3,4,5,6\}}b_{s,36} + \gamma_{\{1,2,3,4,5,6\}}b_{s,37} + \gamma_{\{1,2,3,5\}}b_{s,41}\\
& + \gamma_{\{1,2,3,5\}}b_{s,42} + \gamma_{\{5,10\}}b_{s,43} + \gamma_{\{5,10\}}b_{s,47} + (\gamma_{\{5,8,11\}}+1)b_{s,48}\\
& + \gamma_{\{5,6\}}b_{s,50} + \gamma_{\{5,6\}}b_{s,51} + \gamma_{\{4,5,6,12\}}b_{s,55} + \gamma_{\{6,8\}}b_{s,59} + \gamma_{\{6,8,9,12\}}b_{s,63}\\
& + \gamma_{\{1,2,3,5\}}b_{s,65} + \gamma_{\{1,2,3,5\}}b_{s,66} + (\gamma_{\{3,4,5,6,7,8,11\}}+1)b_{s,67}\\
& + (\gamma_{\{3,4,5,6,7,11,12\}}+1)b_{s,68} + (\gamma_{\{5,6,8,10,11\}}+1)b_{s,71} + \gamma_{\{3,5,7,12\}}b_{s,76}\\
& + \gamma_{\{3,5,7,10\}}b_{s,77} + \gamma_{\{8,10\}}b_{s,81} + (\gamma_{11}+1)b_{s,83} + (\gamma_{11}+1)b_{s,85}\\
& + (\gamma_{11}+1)b_{s,87} + \gamma_{\{3,5,7,10\}}b_{s,89} + \gamma_{\{7,9\}}b_{s,92} + (\gamma_{\{3,5,7,8,9,10,11\}}+1)b_{s,95}\\
& + (\gamma_{\{9,11,12\}}+1)b_{s,96} + \gamma_{\{3,5,7,10\}}b_{s,97} + (\gamma_{\{3,5,7,8,10,11,12\}}+1)b_{s,99}\\
& + \gamma_{\{3,5,7,10\}}b_{s,100} + \gamma_{\{3,5,7,10\}}b_{s,105} \equiv 0.
\end{align*}
Computing from this equality gives and $\gamma_u =\gamma_{12} = \gamma_1$ for $1\leqslant u \leqslant 10$,  $\gamma_{11} = 1$ and $\gamma_{13} = 0$. Hence, we obtain 
\begin{align*}
f &\equiv \psi_{s,t}(\theta_{t,s-1}) + q_{s,0} + q_{s,11} +\gamma_1 \Big(\sum_{1\leqslant u \leqslant 12,\, u \ne 11}q_{s,u}\Big).
\end{align*}
Since $q_{s,0} + q_{s,11} = \sum_{j \in\{68,\, 74,\, 75,\, 84,\, 85,\, 86\}\atop\hskip0.1cm\cup\{87,\, 88,\, 100,\, 104,\, 105\}}b_{s,j}$ and $\sum_{1\leqslant u \leqslant 12,\, u \ne 11}q_{s,u} = \zeta_{t,s}$ we obtain $f \equiv \delta_{t,s} +\gamma_1\zeta_{t,s}$.

For $s = 4$ and $t \geqslant 4$, by Theorem \ref{d1t4}, $(QP_4)_{d_{3,t}}^{GL_4} = \langle [\theta_{t,3}]\rangle$, hence we get $(QP_4)_{n_{4,t}}^{GL_4} = \langle [\zeta_{t,4}], [\delta_{t,4}]\rangle$. The theorem is proved for the cases $s = 4$ and $t \geqslant 4$. 
	
For $s \geqslant 5$, by Theorem \ref{dlt21} we have $[\xi_{t,s-1}] \in (QP_4)_{d_{s-1,t}}^{GL_4}$. We find an element $g  \in \mbox{Ker}\big((\widetilde{Sq}^0_*)_{(4,n_{s,t})}\big)$ such that $f = \psi_{s,t}(\xi_{t,s-1}) + g$ and $[f] \in (QP_4)_{n_{s,t}}^{GL_4}$.

By a direct computation we get
\begin{align*}
\rho_1(\psi_{s,t}&(\xi_{t,s-1}))+\psi_s(\xi_{t,s-1}) \equiv x_1x_2x_3^{2^{s}-2}x_4^{2^{s+t}-2} + x_1x_2x_3^{2^{s+1}-2}x_4^{2^{s+t}-2^{s}-2}\\ 
& + x_1x_2^{3}x_3^{2^{s+1}-2}x_4^{2^{s+t}-2^{s}-4} + x_1^{3}x_2x_3^{2^{s+1}-2}x_4^{2^{s+t}-2^{s}-4} + x_1x_2^{3}x_3^{2^{s}-2}x_4^{2^{s+t}-4}\\ 
& + x_1^{3}x_2x_3^{2^{s}-2}x_4^{2^{s+t}-4} =  b_{s,67} + b_{s,71} + b_{s,87} + b_{s,88} + b_{s,99} + b_{s,100},\\
\rho_2(\psi_{s,t}&(\xi_{t,s-1}))+\psi_s(\xi_{t,s-1}) \equiv x_1x_2x_3^{2^{s}-2}x_4^{2^{s+t}-2} + x_1x_2^{2}x_3^{2^{s+t}-3}x_4^{2^{s}-2}\\ 
& + x_1x_2^{2}x_3^{2^{s}-1}x_4^{2^{s+t}-4} + x_1x_2^{2^{s}-1}x_3^{2}x_4^{2^{s+t}-4} + x_1x_2^{3}x_3^{2^{s+t}-4}x_4^{2^{s}-2}\\ 
& + x_1x_2^{3}x_3^{2^{s}-4}x_4^{2^{s+t}-2} + x_1x_2^{3}x_3^{2^{s}-2}x_4^{2^{s+t}-4}\\ 
&= b_{s,67} + b_{s,73} + b_{s,74} + b_{s,76} + b_{s,83} + b_{s,95} + b_{s,99},\\
\rho_3(\psi_{s,t}&(\xi_{t,s-1}))+\psi_s(\xi_{t,s-1}) \equiv  x_1x_2x_3^{2^{s}-2}x_4^{2^{s+t}-2} + x_1x_2x_3^{2^{s+1}-2}x_4^{2^{s+t}-2^{s}-2}\\ 
& + x_1x_2^{2}x_3^{2^{s}-3}x_4^{2^{s+t}-2} + x_1x_2^{2}x_3^{2^{s+t}-3}x_4^{2^{s}-2} + x_1x_2^{3}x_3^{2^{s}-2}x_4^{2^{s+t}-4}\\ 
& + x_1x_2^{3}x_3^{2^{s+t}-4}x_4^{2^{s}-2} + x_1^{3}x_2x_3^{2^{s+1}-4}x_4^{2^{s+t}-2^{s}-2} + x_1^{3}x_2x_3^{2^{s+1}-2}x_4^{2^{s+t}-2^{s}-4}\\ 
& = \mbox{$\sum_{j\in\{67,\, 71,\, 73,\, 83,\, 86,\, 88,\, 94,\, 99\}}b_{s,j}$}
\end{align*} 
By computing from these results and the relations $\rho_i(f) + f \equiv 0$, $i = 1,\, 2,\, 3$, we obtain
$$f \equiv  \psi_{s,t}(\xi_{t,s-1}) + q_{s,0} + \sum_{1 \leqslant u \leqslant 13}\gamma_uq_{s,u}$$
with $\gamma_u \in \mathbb F_2$ and 
$q_{s,0} = \sum_{j \in \{68,\, 69,\, 70,\, 71,\, 72,\, 73,\, 74,\, 75,\, 79,\, 85,\, 86,\, 87,\, 99,\, 103\}}b_{s,j}.$

A direct computation gives
\begin{align*}
\rho_4(\psi_{s,t}&(\xi_{t,s-1}))+\psi_s(\xi_{t,s-1}) \equiv x_2x_3^{2^{s}-2}x_4^{2^{s+t}-1} + x_2x_3^{2^{s+t}-1}x_4^{2^{s}-2}\\
& + x_2^{3}x_3^{2^{s}-3}x_4^{2^{s+t}-2} + x_1x_2^{3}x_3^{2^{s}-4}x_4^{2^{s+t}-2} + x_1x_2^{3}x_3^{2^{s}-2}x_4^{2^{s+t}-4}\\
& + x_1x_2^{3}x_3^{2^{s+1}-2}x_4^{2^{s+t}-2^{s}-4} + x_1x_2^{2^{s+1}-1}x_3^{2}x_4^{2^{s+t}-2^{s}-4} + x_1^{3}x_2x_3^{2^{s}-4}x_4^{2^{s+t}-2}\\
& + x_1^{3}x_2x_3^{2^{s}-2}x_4^{2^{s+t}-4} + x_1^{3}x_2^{5}x_3^{2^{s}-6}x_4^{2^{s+t}-4} + x_1^{3}x_2^{2^{s}-3}x_3^{2}x_4^{2^{s+t}-4}\\
& + x_1^{2^{s}-1}x_2x_3^{2}x_4^{2^{s+t}-4} = \mbox{$\sum_{j \in \{19,\, 20,\, 63,\, 77,\, 81,\, 87,\, 89,\, 95,\, 97,\, 99,\, 100,\, 105\}}b_{s,j}$}.
\end{align*}

By using this result and computing $\rho_4(f)+f$ in terms of the admissible monomials, we obtain
\begin{align*}
\rho_4(f)+f &\equiv \gamma_{\{1,3\}}b_{s,3} + \gamma_{\{1,3\}}b_{s,4} + \gamma_{\{1,4\}}b_{s,5} + \gamma_{\{1,4\}}b_{s,6} + \gamma_{\{1,2\}}b_{s,9}\\ 
&\quad + \gamma_{\{2,5\}}b_{s,14} + \gamma_{\{2,5\}}b_{s,15} + (\gamma_{\{3,7\}}+1)b_{s,19} + (\gamma_{\{3,7\}}+1)b_{s,20}\\ 
&\quad + (\gamma_{\{3,9,11\}}+1)b_{s,21} + \gamma_{\{3,4,5,6\}}b_{s,26} + \gamma_{\{3,4,5,6\}}b_{s,27}\\ 
&\quad + (\gamma_{\{4,10,13\}}+1)b_{s,31} + (\gamma_{\{4,10,13\}}+1)b_{s,32} + (\gamma_{\{4,9,11\}}+1)b_{s,33}\\ 
&\quad + \gamma_{\{1,2,3,4,5,6\}}b_{s,36} + \gamma_{\{1,2,3,4,5,6\}}b_{s,37} + \gamma_{\{1,2,3,5\}}b_{s,41} + \gamma_{\{1,2,3,5\}}b_{s,42}\\ 
&\quad + \gamma_{\{5,10\}}b_{s,43} + \gamma_{\{5,10\}}b_{s,47} + (\gamma_{\{5,8,11\}}+1)b_{s,48} + \gamma_{\{5,6\}}b_{s,50}\\ 
&\quad + \gamma_{\{5,6\}}b_{s,51} + (\gamma_{\{4,5,6,12\}}+1)b_{s,55} + (\gamma_{\{6,8\}}+1)b_{s,59}\\ 
&\quad + (\gamma_{\{6,8,9,12\}}+1)b_{s,63} + \gamma_{\{1,2,3,5\}}b_{s,65} + \gamma_{\{1,2,3,5\}}b_{s,66}\\ 
&\quad + \gamma_{\{3,4,5,6,7,8,11\}}b_{s,67} + \gamma_{\{3,4,5,6,7,11,12\}}b_{s,68} +( \gamma_{\{5,6,8,10,11\}}+1)b_{s,71}\\ 
&\quad + \gamma_{\{3,5,7,12\}}b_{s,76} + (\gamma_{\{3,5,7,10\}}+1)b_{s,77} + (\gamma_{\{8,10\}}+1)b_{s,81} + \gamma_{11}b_{s,83}\\ 
&\quad + \gamma_{11}b_{s,85} + \gamma_{11}b_{s,87} + (\gamma_{\{3,5,7,10\}}+1)b_{s,89} + \gamma_{\{7,9\}}b_{s,92}\\ 
&\quad + (\gamma_{\{3,5,7,8,9,10,11\}}+1)b_{s,95} + \gamma_{\{9,11,12\}}b_{s,96} + (\gamma_{\{3,5,7,10\}}+1)b_{s,97}\\ 
&\quad + (\gamma_{\{3,5,7,8,10,11,12\}}+1)b_{s,99} + (\gamma_{\{3,5,7,10\}}+1)b_{s,100}\\ 
&\quad +( \gamma_{\{3,5,7,10\}}+1)b_{s,105} \equiv 0.
\end{align*}
Computing from this equality gives $\gamma_{11} = 0$, $\gamma_u =\gamma_1 = \gamma_{10}$ for $1 \leqslant u \leqslant 6$,  $\gamma_{13} =1$ and $\gamma_{u} = \gamma_1 + 1$ for $u = 7,\, 8,\, 9,\, 12$. Hence, we obtain 
\begin{align*}
f \equiv \psi_{s,t}(\xi_{t,s-1})) + q_{s,0} + q_{s,7} &+ q_{s,8} + q_{s,9} + q_{s,12} + q_{s,13}\\ &\quad
 +\gamma_1 \big(\mbox{$\sum_{1\leqslant t \leqslant 12,\, u \ne 11}$}q_{s,u}\big).  
\end{align*}
We have  
$$q_{s,0} + q_{s,7} + q_{s,8} + q_{s,9} + q_{s,12} + q_{s,13} = \sum_{j \in \{67,\, 76,\, 77,\, 83,\, 84,\, 88,\, 90,\, 91,\, 92\}\atop\hskip0.1cm\cup\{93,\, 94,\, 96,\, 98,\, 99,\, 101,\, 102,\, 104\}}b_{s,j}$$ and $\sum_{1\leqslant t \leqslant 12,\, u \ne 11}q_{s,u} = \zeta_{t,s}$. Hence, we get $f \equiv \chi_{t,s} + \gamma_1\zeta_{t,s}$.

For $s \geqslant 5$ and $2 \leqslant t \leqslant 3$, by Theorem \ref{dlt21}, 
$$(QP_4)_{d_{s-1},t}^{GL_4} = QP_4((3)|^{s-1}|(1)|^{t+1}) = \langle [\xi_{t,s-1}]\rangle.$$
Hence, we get $(QP_4)_{n_{s,t}}^{GL_4}	= \langle [\zeta_{t,s}],[\chi_{t,s}]\rangle$ and $\dim(QP_4)_{n_{s,t}}^{GL_4} = 2$.

For $s \geqslant 5$ and $t \geqslant 4$, by Theorems \ref{dln2} and \ref{d1t4}, we have
$$\dim(QP_4)_{n_{s,t}}^{GL_4} \leqslant \dim \mbox{Ker}\big((\widetilde{Sq}^0_*)_{(4,n_{s,t})}\big)^{GL_4} + \dim(QP_4)_{d_{s-1,t}}^{GL_4} = 3.$$
Since the set $\{[\zeta_{t,s}],[\delta_{t,s}],[\chi_{t,s}]\}$ is linearly independent in $(QP_4)_{n_{s,t}}^{GL_4}$, we obtain $(QP_4)_{n_{s,t}}^{GL_4} = \langle[\zeta_{t,s}],[\delta_{t,s}],[\chi_{t,s}] \rangle$ and $\dim(QP_4)_{n_{s,t}}^{GL_4} = 3$. The theorem is completely proved.
\end{proof}


\end{document}